\documentclass[12pt]{article}

\usepackage{amsmath}
\usepackage{amsfonts}
\usepackage{amsthm}
\usepackage{amssymb}
\usepackage{url}
\usepackage{color,graphics}
\usepackage{psfrag,epsfig}

\usepackage[medium]{titlesec}

\usepackage{todonotes}

\presetkeys{todonotes}{color=green!30, size = \tiny}{}
\usepackage{longtable}
\usepackage{microtype}

\usepackage[text={6in,8in},centering]{geometry}

\usepackage[breaklinks = true]{hyperref}

\newtheorem{thm}{Theorem}[section]
\newtheorem{kthm}{Key Theorem}
\newtheorem{lem}{Lemma}[section]

\newtheorem{cor}[lem]{Corollary}
\newtheorem{prop}[lem]{Proposition}

\theoremstyle{remark}
\newtheorem{rmk}{Remark}[section]

\theoremstyle{definition}
\newtheorem{defn}[lem]{Definition}

\newcommand{\supp}{\operatorname{supp}}
\newcommand{\diag}{\operatorname{diag}}
\newcommand{\dist}{\operatorname{dist}}
\newcommand{\Z}{\mathbb{Z}}

\newcommand{\R}{\mathbb{R}}

\newcommand{\NN}{\mathbb{N}}

\newcommand{\mC}{\ensuremath{\mathcal{C}}}

\newcommand{\cP}{\ensuremath{\mathcal{P}}}

\newcommand{\mM}{\ensuremath{\mathcal{M}}}

\newcommand{\mO}{\ensuremath{\mathcal{O}}}

\newcommand{\Tm}{\ensuremath{\mathbb{T}}}
\newcommand{\T}{\ensuremath{\mathbb{T}}}

\newcommand{\e}{\ensuremath{\epsilon}}
\newcommand{\bH}{\mathbb{H}}
\newcommand{\Te}{\sqrt{\epsilon}\mathbb{T}}
\newcommand{\Hse}{H^s_{\epsilon}}
\newcommand{\bA}{\mathbb{A}}

\newcommand{\tcM}{\widetilde{\mathcal M}}
\newcommand{\tcA}{\widetilde{\mathcal A}}
\newcommand{\tcN}{\widetilde{\mathcal N}}

\newcommand{\cC}{\mathcal{C}}
\newcommand{\cG}{\mathcal{G}}
\newcommand{\cV}{\mathcal{V}}
\newcommand{\cU}{\mathcal{U}}
\newcommand{\cA}{\mathcal{A}}
\newcommand{\cN}{\mathcal{N}}
\newcommand{\cM}{\mathcal{M}}

\newcommand{\barc}{\bar{c}}

\newcommand{\bV}{\mathbb{V}}

\def\textb{\textcolor{blue}}
\def\textr{\textcolor{red}}

\newcommand{\cK}{\ensuremath{\mathcal{K}}}

\newcommand{\cI}{\ensuremath{\mathcal{I}}}

\newcommand{\Phg}{\Phi_{\mathrm{glob}}}
\newcommand{\Phl}{\Phi_{\mathrm{loc}}}

\newcommand{\bmat}[1]{\begin{bmatrix} #1 \end{bmatrix}}
\newcommand{\Span}{\mathrm{Span}}

\newcommand{\tdt}{\tilde{t}}

\newcommand{\tdC}{\tilde{C}}
\newcommand{\LF}{\mathcal{LF}}
\newcommand{\Id}{\mathrm{Id}}
\newcommand{\Mane}{Ma\~ne\ }

\newcommand{\tcS}{\widetilde{\mathcal{S}}}
\newcommand{\tcH}{\widetilde{\mathcal{H}}}
\newcommand{\tcI}{\widetilde{\mathcal{I}}}
\newcommand{\Inv}{\mathrm{Inv}}
\newcommand{\loc}{\mathrm{loc}}
\newcommand{\cH}{\mathcal{H}}
\newcommand{\cF}{\mathcal{F}}

\newcommand{\id}{\mathrm{id}}

\newcommand{\pre}{\mathrm{pre}}

\newcommand{\cW}{\mathcal{W}}
\newcommand{\st}{\mathrm{st}}

\newcommand{\vp}{\varpi\,}

\newcommand{\N}{\mathbb{N}}
\newcommand{\chG}{\check{\mathcal{G}}}

\newcommand{\chT}{\check{T}}

\newcommand{\bL}{\mathbb{L}}

\def\to{\longrightarrow}

\def\leq{\leqslant}
\def\geq{\geqslant}

\def\bdef{\begin{definition}}
\def\endef{\end{definition}}
\def\bthm{\begin{thm}}
\def\ethm{\end{thm}}
\def\bkthm{\begin{kthm}}
\def\ekthm{\end{kthm}}
\def\blm{\begin{lemma}}
\def\elm{\end{lemma}}
\def\brm{\begin{rmk}}
\def\erm{\end{rmk}}
\def\bprop{\begin{prop}}
\def\eprop{\end{prop}}
\def\bcor{\begin{cor}}
\def\ecor{\end{cor}}
\def\be{\begin{eqnarray}}
\def\ee{\end{eqnarray}}
\def\beal{\begin{aligned}}
\def\enal{\end{aligned}}

\def\om{\omega}
\def\al{\alpha}

\def\bt{\beta}
\def\eps{\epsilon}

\def\M{\mathcal M}
\def\R{\mathbb R}

\def\T{\mathbb T}
\def\Q{\mathbb Q}
\def\Z{\mathbb Z}

\def\cR{\mathcal R}

\def\cU{\mathcal U}

\def\cC{\mathcal C}
\def\cL{\mathcal L}

\def\cG{\mathcal G}
\def\gm{\gamma}

\def\Lb{\Lambda}
\def\th{\theta}
\def\dt{\delta}
\def\lb{\lambda}
\def\cS{\mathcal S}

\def\be {\begin{equation}}
\def\ee {\end{equation}}
\def\bdef{\begin{definition}}
\def\endef{\end{definition}}
\def\blm{\begin{lem}}
\def\elm{\end{lem}}
\def\beal{\begin{aligned}}
\def\enal{\end{aligned}}

\newtheorem{definition}{Definition}

\newcommand{\argmin}{\mathrm{argmin}}
\newcommand{\con}{\mathrm{con}}
\numberwithin{equation}{section}
\numberwithin{figure}{section}

\begin{document}
\title{Arnold diffusion for smooth systems of two and a half degrees of freedom}

\author{
 V. Kaloshin\footnote{University of Maryland at College Park
 (\texttt{vadim.kaloshin\@ gmail.com})},\ \ \
 K. Zhang\footnote{University of Toronto (\texttt{kzhang\@ math.utoronto.edu})}}

\maketitle

\begin{center}
{\it \quad \quad \quad 
Dedicated to the memory of John Mather: 
\newline a great  mathematician  and a remarkable person }
\end{center}

\vskip 0.3in

\begin{abstract}
In the present paper we prove a strong form of Arnold diffusion.
Let $\T^2$ be the two torus  and $B^2$ be the unit ball around the origin
in $\R^2$. Fix $\rho>0$.  Our main result says that for a ``generic''
time-periodic perturbation of an integrable system of two degrees
of freedom
\[
H_0(p)+\eps H_1(\th,p,t),\quad \ \th\in \T^2,\ p\in B^2,\ t\in \T=\R/\Z,
\]
with a strictly convex $H_0$, there exists a $\rho$-dense
orbit $(\th_{\e},p_{\e},t)(t)$ in $\T^2 \times B^2 \times \T$, namely,
a $\rho$-neighborhood of the orbit contains $\T^2 \times B^2 \times \T$.

Our proof is a combination of geometric and variational methods.
The fundamental elements of the construction are usage of
crumpled normally hyperbolic invariant cylinders from \cite{BKZ},
flower and simple normally hyperbolic invariant manifolds
from  as well as their kissing property at a strong 
double resonance. This allows us to build a ``connected'' net 
of $3$-dimensional normally hyperbolic invariant manifolds.
To construct diffusing orbits along this net we employ 
a version of Mather variational method \cite{Ma2} 
proposed by Bernard in \cite{Be}. This version is equipped 
with weak KAM theory \cite{Fa}.
\end{abstract}

\newpage

\tableofcontents


\section{Introduction}
The famous question called the ergodic hypothesis, formulated by
Maxwell and Boltzmann, suggests that for a typical Hamiltonian on
a typical energy surface all, but a set of zero measure of initial
conditions, have trajectories dense in this energy surface. However,
KAM theory showed that for an open set of (nearly integrable) Hamiltonian 
systems there is a set of initial conditions of positive measure with almost
periodic trajectories. This disproved the ergodic
hypothesis and forced to reconsider the problem.

A quasi-ergodic hypothesis, proposed by Ehrenfest \cite{E} and Birkhoff \cite{Bi},
asks if a typical Hamiltonian on a typical energy surface has a dense orbit.
A definite answer whether this statement is true or not is still far out
of reach of modern dynamics. There was an attempt to prove this statement
by E. Fermi \cite{Fe}, which failed (see \cite{Ga} for a more detailed account).
To simplify the problem, Arnold \cite{Ar3}
asks:

{\it Does there exist a real instability in many-dimensional problems of
perturbation theory when the invariant tori do not divide the phase space?}

For nearly integrable systems of two degrees (resp. of one and a half) of freedom
the invariant tori do divide the phase space and an energy surface respectively.
This implies that instability do not occur.
We solve a weaker version of this question for systems with two and a half
and $3$ degrees of freedom. This corresponds autonomous perturbations
of integrable systems with three degrees of freedom (resp. time-periodic 
perturbations of integrable systems with two degrees of freedom).

\subsection{Statement of the result}

Let $(\th,p)\in \T^2 \times B^2$ be the phase space of an integrable
Hamiltonian system $H_0(p)$ with $\T^2$ being $2$-dimensional torus $\T^2=\R^2/\Z^2\ni \th=(\th_1,\th_2)$ 
and $B^2$ being the unit ball around $0$ in $\R^2$, 
$p=(p_1,p_2)\in B^2$.  $H_0$ is assumed to be strictly convex with the following uniform estimate: there exists $D>1$ such that 
\[
	D^{-1} I \le \partial^2_{pp} H_0 \le D I, \quad |H_0(0)|, \quad \|\partial_p H_0(0)\| \le D. 
\]
where $I$ is the $2\times 2$ identity matrix.

Consider a smooth time periodic perturbation
\[
	H_\eps(\th,p,t)=H_0(p)+\eps H_1(\th,p,t),\quad t\in \T=\R/\T.
\]
We study  Arnold diffusion for this system, namely,
\begin{center}
	{\it topological instability in the $p$ variable. }
\end{center}
Arnold \cite{Ar1} proved existence of such orbits 
for an example and conjectured that they exist for 
a typical perturbation (see e.g. \cite{Ar3,Ar2,Ar5}).

Denote $\Z^3_* = \Z^3 \setminus (0, 0, 1)\Z$, then integer relations $k \cdot (\partial_p H_0,1)=0$ with $k=(\vec k_1,k_0) \in \Z^3_*$ and $\ \cdot\ $ being the inner product define {\it a resonant submanifold}. The strict convexity of $H_0$ implies that $\partial_p H_0:B^2 \to \R^2$
is a diffeomorphism and each resonant line defines a smooth curve embedded into action space
$$
S_{k}=\{p\in \R^2:\ k \cdot (\partial_p H_0,1)=0\}.
$$
If curves $S_{k}$ and $S_{k'}$ are given by
two linearly independent resonances vectors $\{k,k'\}$,
they either have no intersection or intersect at a single point in $B^2$.
We call a  resonance $S_{k}$ {\it space irreducible} if 
the greatest common divisor of components of $\vec k_1$ is one. 
Notice that space irreducible resonances are dense.

\medskip 

\begin{figure}[t]
	\centering
	\includegraphics[width=2in]{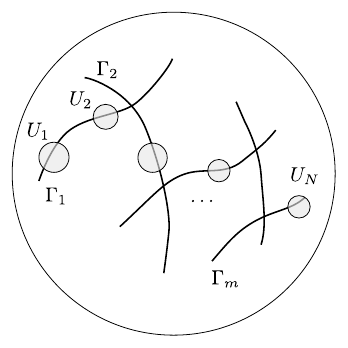}
  \caption{Resonant net}
  \label{fig:res-net}
\end{figure}

	Consider a finite collection of tuples:
	\[
		\cK = \left\{  (k, \Gamma_k): \quad k \in \Z^3_*, \quad \Gamma_k \subset S_k\cap B^2\right\}, 
	\]
	where $k$ is \emph{space irreducible}, 
	\footnote{This condition is not really necessary, we assume it as it helps simplify the presentation for single resonances.}
	and $\Gamma_{k} \subset S_{k}$ is a closed segment. 
We say $\cK$ defines a \emph{diffusion path} if 
\[
	\cP = \bigcup_{(k, \Gamma_k) \in \cK} \Gamma_k
\]
is a connected set. We would like to construct diffusion orbits along the path $\cP$ (see Figure \ref{fig:res-net}).


\begin{thm} 
\label{main} Let $\mathcal P\subset B^2$ be a diffusion path,
$5\le r < +\infty$, and $U_1, \dots,  U_N$ be open sets such that $U_i \cap \cP \ne \emptyset$, $i = 1, \dots, N$.
Then there exist:
\begin{itemize}
	\item a $C^r$ open and dense set 
	$\mathcal U = \mathcal U(\mathcal P)\subset \cS^r$ 
	depending on $\mathcal P$, 
	\item a nonnegative  lower semi-continuous
	function $\eps_0=\eps_0(H_1)$ with $\eps_0|_{\mathcal U}>0$ 
	and 
	\item a ``cusp'' set 
	$$
	{\mathcal V}:=\mathcal V(\mathcal U,\eps_0):=
	\{\epsilon H_1:\ H_1\in {\mathcal U} ,\ \ 0<\eps<\eps_0(H_1)\},
	$$
	\item 
	a $C^r$ open and dense subset of $\eps H_1 \in \mathcal W
	\subsetneq \mathcal V$ 
\end{itemize}
such that for each $\epsilon H_1 \in \cW$  there is an orbit $(\theta,p)(t)$ of $H_\epsilon$
and times $0<T_1< \cdots <T_N$ with the property  
$$
p(T_i) \in U_i, \qquad i=1,\dots,N. 
$$
\end{thm}
\begin{rmk}
The condition that $\epsilon_0$ is lower semi-continuous implies the set $\cV$ is open. 
\end{rmk}

\begin{figure}[t]
	\centering
	\includegraphics[width=3.5in]{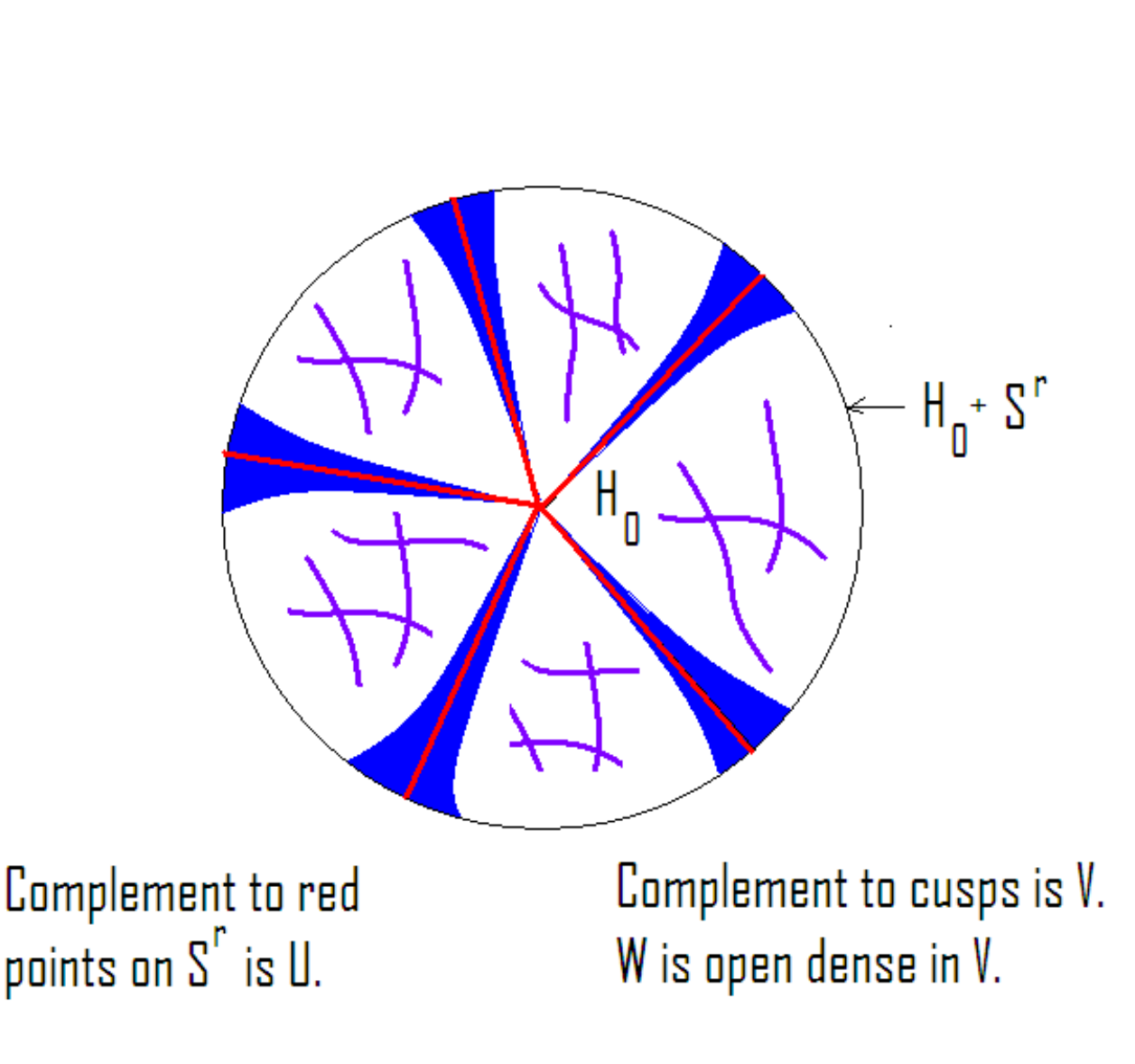}
  \caption{Description of generic perturbations}
  \label{fig:cusp-resid}
\end{figure}

\begin{rmk}
Note that the notion of genericity we use is not standard.
We show that in a neighborhood of perturbations of $H_0$ 
the set of good directions $\mathcal U$ is open dense in $\cS^r$. 
Around each exceptional (nowhere dense) direction we remove 
a cusp and call the complement $\mathcal V$. For this set 
of perturbations we establish connected collection of invariant 
manifolds. Then in the complement to some exceptional 
perturbations $\mathcal W$ in $\mathcal V$ we show that there 
are diffusing orbits ``shadowing'' these cylinders. Mather calls it {\it cusp residual}. See 
Figure~\ref{fig:cusp-resid}.
\end{rmk}

Consider $(k_1, \Gamma_{k_1}) \in \cK$, for $\lambda >0$, in section \ref{sec:intro-SR} we define a quantitative 
nondegeneracy hypothesis relative to the resonant segment $\Gamma_{k_1}$. 
\[
	SR (k_1, \Gamma_{k_1}, \lb) = [SR1_\lambda] - [SR3_\lambda]
\]
Let us denote by $\cU_{SR}^\lambda(k_1, \Gamma_{k_1})$ the set of $H_1$ 
which satisfy $SR(k_1, \Gamma_{k_1}, \lambda)$.  

Suppose $H_1$ satisfies the conditions $SR(k_1,\Gamma_{k_1},  \lambda)$, we define a finite subset of $\Z^3_*$ called the \emph{strong additional resonances}. In Theorem~\ref{thm:c-equiv-sr} we define a large constant $K = K(k_1, \Gamma_{k_1}, \lambda)$, and call $k_2 \in \Z^3_*$ \emph{strong} if:
\begin{itemize}
 \item either there is $(k_2, \Gamma_{k_2}) \in \cK$ such that $\Gamma_{k_1} \cap \Gamma_{k_2} \ne \emptyset$; 
 \item or $|k_2| \le K(k_1, \Gamma_{k_1}, \lambda)$ and $S_{k_2}\cap \Gamma_{k_1} \ne \emptyset$. 
\end{itemize}
We emphasize that strong additional resonances are taken from the set $\Z^3_*$, not just the space irreducible ones. Denote the set of strong additional resonances $\cK^\st(k_1, \Gamma_{k_1}, \lambda)$,  If $k_2$ is strong, then it defines a unique double resonance $\Gamma_{k_1, k_2} = \Gamma_{k_1} \cap S_{k_2}$. 

For each double resonance $\Gamma_{k_1, k_2}$, we associate 
non-resonance conditions of two types: 
\begin{itemize}
	\item high energy $[DR1^h] -[DR3^h]$ 
	(section \ref{sec:intro-DR}), 
	\item low energy $[DR1^c] - [DR4^c]$ 
	(section \ref{sec:intro-DR}).
\end{itemize} 
 For each pair $k_2 \in \mathcal K^{st}(k_1, \Gamma_{k_1}, \lambda)$ consider the set of 
$H_1$ which satisfy the above conditions and denote it by 
$\cU_{DR}(k_1, \Gamma_{k_1}, k_2)$. 

\begin{rmk}
All our non-degenerate conditions at a double resonance is stated relative to a single resonance. Namely, the condition $DR(k_1, \Gamma_{k_1}, k_2)$ may differ from the condition $DR(k_2, \Gamma_{k_2}, k_1)$.  
\end{rmk}

The following theorem is immediate given that:
\begin{itemize}
 \item (Proposition~\ref{prop:sr-non-deg}) Each $\cU_{SR}^\lambda(k_1, \Gamma_{k_1})$ is open and the union $\bigcup_{\lambda > 0} \cU_{SR}^\lambda(k_1)$ is dense;
 \item (Proposition~\ref{prop:Dr-non-deg}) The set $\cU_{DR}(k_1,\Gamma_{k_1}, k_2)$ is open and dense. 
\end{itemize}
\begin{thm}
	The set 
	\[
		\mathcal U = \cU(\cP) := \bigcup_{\lambda > 0} \bigcap_{(k_1, \Gamma_{k_1}) \in \mathcal K} 
		\left( \cU_{SR}^\lambda(k_1, \Gamma_{k_1}) \cap \bigcap_{k_2 \in \mathcal K^{st}(k_1, \Gamma_{k_1}, \lambda)} 
		\cU_{DR}(k_1,\Gamma_{k_1}, k_2) \right) 
	\]
	is open and dense in $\mathcal S^r$. 
\end{thm}

As a corollary of Theorem~\ref{main}, we obtain:
\begin{thm}[Almost Density Theorem]\label{thm:alomst-dense}
For any $\rho >0$ there are 
\begin{itemize}
	\item an open dense set $\mathcal U = \cU(\rho) \subset \mathcal S^r$, 
	\item a nonnegative 
	lower semi--continuous function $\eps_0:\mathcal S^r\to \R_+$ with 
	$\eps_0|_{\mathcal U}>0$, 
	\item a cusp set 
	$
	{\mathcal V}:=\mathcal V(\mathcal U,\eps_0):=
	\{\epsilon H_1:\ H_1\in {\mathcal U} ,\ \ 0<\eps<\eps_0(H_1)\},
	$
	\item an open dense subset $\mathcal W\subsetneq \mathcal V$ 
\end{itemize}
all depending on $\rho$ such that for any $\eps H_1 \in \mathcal W$   
there is a  $\rho$--dense orbit on $\mathbb T^2 \times B^2 \times \mathbb T.$ 
\end{thm}
\begin{proof}
	Given a vector $(\omega, 1) \in \R^3$, let us call $\omega$ being $\rho$-irrational if there exists $T>0$ such that  $\{t(\omega, 1): \, t\in [-T, T]\} \subset \T^3$ is $\rho$-dense, and let $T(\omega)$ be the smallest such $T$. Using the fact that $\dot{p} = O(\epsilon)$, $\dot{\theta} = \nabla H_0(p) + O(\epsilon)$, there is $\epsilon_0 >0$ depending on $\rho$ and $T(\omega)$ such that if $0 < \epsilon < \epsilon_0$
	\begin{equation}
		\label{eq:rho-dense}
		\beal
		B_{2\rho}\left( \bigcup_{t \in \R}\phi^t_{H_\epsilon}(\theta_0, p_0, t_0) \right) & \supset \ \T^2 \times B_\rho(p_*) \times \T , \\
\text{ for all } &p_0 \in B_\rho(p_*), \, (\theta_0, t_0) \in \T^2 \times \T.
\enal  
\end{equation}

Any vector that is not $\rho$-irrational (called $\rho$-rational) must be resonant: namely $k \cdot (\omega, 1) = 0$ for some 
$k \in \Z^3_*$. Moreover, there are only finitely many resonances 
that corresponds to $\rho$-rational vectors. Since there are 
infinitely many space irreducible resonances, there is a diffusion 
path $\cP$ such  consisting only of $\rho/2$-irrational 
resonances. Moreover, we may choose the path $\cP$ to be 
$\rho/2$-dense in $B^2$, since space irreducible resonances 
are dense. 

We now apply Theorem~\ref{main} to the path $\cP$, and pick $p_i$, $i = 1, \dots, N \in \cP$ such that $(\nabla H_0(p_i), 1)$ is $(\rho/2)-$irrational, and such that $\bigcup_{i = 1}^n B_{\rho/2}(p_i) \supset B^2$. According to our theorem, there is an orbit whose $p$ component  visit every $B_\rho(p_i)$. Then the orbit must be $\rho-$dense in view of \eqref{eq:rho-dense}. 
\end{proof}

\subsection{Discussions of the result}

\subsubsection*{Relation with Mather's approach}

Theorem \ref{main} was announced by Mather in \cite{Ma4}, where 
he proposed a plan to prove it. Some parts of the proof are written 
in \cite{Ma7}.
Our work realizes Mather's general plan using weak KAM theory and 
Hamiltonian point of view. New techniques and tools are introduced, 
below we summarize them.  

\begin{itemize}
	\item We utilize Bernard's forcing relation to simplify the construction of diffusion orbit. This allows a more Hamiltonian treatment of the variational concepts, and allows us to reduce the main theorem to local forcing equivalence of cohomology classes. 

	\item We use Hamiltonian normal forms to construct a {collection} of normally hyperbolic invariant cylinders along the chosen diffusion path. We obtain precise control of the normal forms (via an anisotropic $C^2$ norm) at both single and double resonances. Mather's method uses mostly the Lagrangian point of view. 

	\item We introduce the concept of \emph{Aubry-Mather type}, which generalizes the work done in 
Bernard--Kaloshin--Zhang \cite{BKZ} to a more abstract setting, applicable to both single and double resonances. Heuristically, this means the Aubry sets behaves like Aubry-Mather sets in twist maps. Our approach can be seen as a generalization of the variational technique for a priori unstable systems from \cite{Be, BKZ, CY2}. 

	\item One important obstacle is the problem of regularity 
	of barrier functions (see section \ref{sec:reg-barrier}), 
	which outside of the realm of twist maps is difficult to overcome. 
	Our definition of Aubry-Mather type allows proving this statement in a general setting. 
	It is our understanding that Mather \cite{Ma9} handles 
	this problem without proving existence of invariant cylinders.

	\item In a double resonance we also construct normally hyperbolic
	invariant cylinders. This leads to a fairly simple and explicit
	structure of minimal orbits near a double resonance. In particular,
	in order to switch from one resonance to another we need {\it only one jump}
	(see section \ref{sec:intro-forcing} for the formulation of the statement).

	\item It is our understanding that Mather's approach \cite{Ma9} requires
	an implicitly defined number of jumps. His approach resembles
	his proof of existence of diffusing orbits for twist maps inside
	a Birkhoff region of instability \cite{Ma}.
\end{itemize}

\subsubsection*{Other results on apriori stable systems}

In \cite{Cheng2017}, Theorem 5.1 a weaker result to Theorem  \ref{main}
is stated. The set of admissible perturbations $\mathfrak R_a$ in 
Theorem 5.1 is residual, while our set of admissible perturbations 
$\mathcal U$ is open and dense. The size of admissible perturbation 
$a_P$ from Theorem 5.1 is analog of $\eps_0$ in Theorem  \ref{main}. 
Regularity of dependence of the size of admissible perturbation $a_P$ 
on $P$ is not discussed. Therefore, the genericity of perturbations from 
Theorem 5.1 is up to the reader's interpretation. 
Notice that the proof in \cite{Cheng2017} is variational and, as well 
as our proof, fundamentally relies on Mather's ideas. 

In \cite{Marco2016a} (see also \cite{GM2017, Marco2016}),  Theorem 1 
is nearly identical to Theorem  \ref{main}. A slight difference is a higher 
regularity requirement. This proof is geometrical and does not  use 
variational methods. 

An earlier version of the current paper was available (\cite{KZ2012}) since 2012, the current version is a thorough revision. We introduce a more general concept 
of Aubry-Mather type, propose a different way to perform normal forms, and also 
use a different method to handle the transition from single to double resonances. 


In \cite{KZb} we propose a way to prove 
Arnold diffusion in the same setting as in the present paper, namely, 
for generic time-periodic perturbations of integrable systems of 
three degrees of freedom with strictly convex unperturbed $H_0$.
The key element of the construction is to find a diffusion path such that 
at every strong resonance the associated averaged mechanical system 
is dominant. For dominant systems we introduce dimension reduction 
and prove existence of $3$-dimensional normally hyperbolic invariant 
cylinders. Moreover, we show existence of families of Aubry sets of 
Aubry-Mather type (see section \ref{sec:AM-type} for a definition).
 Finally, we use these sets to construct diffusing orbits. 

In \cite{KZ2017} for dominant mechanical systems in any dimension 
we prove analogous statement on existence of $3$-dimensional cylinders 
carying a family of Aubry-Mather type sets.

%
%

\subsubsection*{Autonomous version}

Let $n=3,\ 
\widetilde p=(\widetilde p_1,\widetilde p_2,\widetilde p_3)\in B^3,$ 
and $\widetilde H_0(p)$ be a strictly convex Hamiltonian. 
{Fix $i\in\{1,2,3\}$ and a regular value $a$ of $\widetilde H_0$.
Since $\widetilde H_0$ is strictly convex, there is only 
one critical value of $\widetilde H_0$}. Consider a convex 
connected open set $W$ in the region $\partial_{p_i}H_0>\rho>0$ 
for some $\rho>0$ and two open sets $U$ and $U'$ in $W$ {intersecting 
the same energy surface 
$S_a=\{\widetilde H_0^{-1}(a)\}$}. Then for a cusp generic 
(autonomous) perturbation
$\widetilde H_0(\widetilde p)+
\eps \widetilde H_1(\widetilde \th,\widetilde p)$ 
there is an orbit $(\widetilde \th_\eps,\widetilde p_\eps)(t)$ 
{on the energy surface $S_a$ connecting} $U$ with $U'$, 
namely, $\widetilde p_\eps(0)\in U$ and 
$\widetilde p_\eps(t)\in U'$ for some $t=t_\eps$. 

This can be shown using energy reduction to a time periodic system 
of two and a half degrees of freedom (see e.g. \cite[Section 45]{Ar4}). 
%
%

\subsubsection*{Generic instability of resonant totally elliptic points}

In \cite{KMV} stability of resonant totally elliptic fixed points
of symplectic maps in dimension $4$ is studied. It is shown that
generically a convex, resonant, totally elliptic
point of a symplectic map is Lyapunov unstable.

\subsubsection*{ Non-convex Hamiltonians}

In the case the Hamiltonian $H_0$ is non-convex 
{or non-strictly convex} for all $p\in B^2$, for example, 
$H_0(p)=p_1^2+p_2^3$, the problem of global Arnold diffusion 
is wide open. Some results for the Hamiltonian 
$H_0(p)=p_1^2-p_2^2$ are in \cite{BCT,BK}. 

To apply variational approach one faces another deep open problem of extending Mather theory and weak KAM theory 
beyond convex Hamiltonians or developping a new technique to construct diffusing orbits. \\

\subsubsection*{Other diffusion mechanisms}

Here we would like to give a short review of other diffusion mechanisms.
In the case $n = 2$ Arnold proposed the following example
\[
	H(q, p, \phi,I, t) =\frac{I^2}{2}+
	\frac{p^2}{2}+ \epsilon (1 - \cos q)(1 + \mu(\sin \phi + \sin t)).
\]
This example is a perturbation of the product of a a one-dimensional
pendulum and a one-dimensional rotator. The main feature of this example
is that it has a $3$-dimensional normally hyperbolic invariant cylinder.
There is a rich literature on Arnold example and we do not intend to
give extensive list of references; we mention \cite{AKN, Be4, BB, Bs1, Zha}, and 
references therein. This example gave rise to a family of examples of
systems of $n+1/2$ degrees of freedom of the form
\[
	H_\eps(q, p,\phi,I, t) = H_0(I) + K_0(p, q) + \eps H_1( q, p, \phi, I, t),
\]
where $(q, p) \in \T^{n-1} \times \R^{n-1}, I \in \R, \phi, t \in T$.
Moreover, the Hamiltonian $K_0(p, q)$ has a saddle fixed point at the origin
and $K_0(0, q)$ attains its strict maximum at $q = 0$. For small $\eps$
a 3-dimensional NHIC $\cC$ persists. Several geometric mechanisms of diffusion have evolved:

--– In \cite{DLS, DGLS, DH, GL} the authors carefully analyze two types of dynamics 
induced on the cylinder $\cC$. These two dynamics are given by so-called inner 
and outer maps. In \cite{DS2017a}, \cite{DS2017}, these techniques are applied to
 a general perturbation of Arnold's example. 

—-- In \cite{T1,T2, T3, DT2016} a return (separatrix) map along invariant manifolds of $\cC$
is constructed. A detailed analysis of this separatrix map gives diffusing orbits.

--- In \cite{CK2015,GKZ2016,KZZ2015} for an open set of perturbations of Arnold's example, one 
constructs an probability measure $\mu$ in the phase space such that the pushforward of 
$\mu$ projected onto the $I$ component in the proper time scale weakly converges to the stochastic 
diffusion process. This, in particular, implies existence of diffusing orbits. 

-- In \cite{GT2017}, the authors treats the a priori chaotic setting, but prove diffusion in 
the real analytic category, which is much more difficult. A different mechanism related 
to the slow-fast system is given by the same authors in \cite{GT}. 

As we mentioned on several other occasions the other two groups \cite{Be, CY1, CY2} are
inspired and influenced by Mather variation method \cite{ Ma, Ma1, Ma2} and build
diffusing orbits variationally. Recently a priori unstable structure was established
for the restricted planar three body problem \cite{FGKR}. It turns out that for
this problem there are no large gaps.

A multidimensional diffusion mechanism of different nature, but also based on
existence and persistence of a $3$-dimensional NHIC $\cC$ is proposed in \cite{BK}.\\


We start with an outline of our proof with a sufficient condition 
for Arnold diffusion.

\subsection{Scheme of diffusion}
\label{sec:intro-scheme}

For all $\epsilon \le 1$, the Hamiltonian $H_\epsilon$ satisfy the Tonelli property of superlinearity, strict convexity, and completeness (see Section~\ref{sec:forcing-def}). Mather theory (\cite{Fa}, \cite{Ma2}) implies that for each cohomology class  $c \in \R^2 \simeq H^1(\T^2, \R)$, the Hamiltonian $H_\epsilon$ admits families of invariant sets of the Hamiltonian flow on $\T^2\times \R^2\times \T$, called the Mather, Aubry, and \Mane sets, satisfying
\[
	\tcM_{H_\epsilon}(c) \subset \tcA_{H_\epsilon}(c) \subset \tcN_{H_\epsilon}(c) \subset  \T^2 \times B_{C\sqrt{\epsilon}}(c) \times \T, 
\]
where $C$ is a constant depending only on $D$ (see Corollary~\ref{cor:min-near-int}). We use $\tcM^0_{H_\epsilon}(c)$, $\tcA^0_{H_\epsilon}(c)$ and $\tcN^0_{H_\epsilon}(c)$ to denote their intersection with the section $\{t=0\}$, which are invariant under the time-$1$-map. Throughout the paper, we may switch between the two equivalent settings: either consider continuous invariant sets of the flow, or discrete invariant sets under the time-$1$-map. 

Our main strategy is then to pick a subset $\Gamma_* \subset \R^2$ of cohomologies very close to the diffusion path $\cP$, then find an orbit that shadows a sequence of Aubry sets $\tcA_{H_\epsilon}(c_i)$, $i = 1, \dots, N$, $c_i \in \Gamma_*$. This requires the existence of non-degenerate \emph{heteroclinic connections} between the Aubry sets, the family of invariant sets with heteroclinic connections is called a \emph{transition chain} by Arnold \cite{Ar1}. To do this, 
we show that the cohomologies satisfy one of the four \emph{diffusion mechanisms}. 

We give a general introduction to  these mechanisms below, and refer to Section \ref{sec:suff-diffusion} for precise definitions.

\subsubsection*{Mather mechanism} For a twist map, it is known since Birkhoff that a region free of essential invariant curves is unstable, namely there exists orbits that drifts from one boundary of the region to another. Mather (\cite{Ma2}) gave a conceptual description of this phenomenon, and generalized it into higher dimension. 

We say that the pair $(H_\epsilon, c)$ satisfy the \emph{Mather mechanism} if 
\[
	\pi_\theta \tcN^0_{H_\epsilon}(c) \subset \T^2
\]
is contractible. (Note  in the twist map case this means $\tcN^0(c)$ is not a rotational 
invariant curve.) 
Mather proved that in this case, $\tcA^0_{H_\epsilon}(c)$ admits a heteroclinic connecting orbit to  
$\tcA^0_{H_\epsilon}(c')$ if $c, c'$ are close. 

\subsubsection*{Arnold mechanism} In Arnold's original paper \cite{Ar1}, Arnold showed the existence of 
a family of invariant tori, whose own stable and unstable manifolds intersects transversally. In our setting, 
the tori are the Aubry sets $\tcA_{H_\epsilon}(c)$ contained in a \emph{normally hyperbolic invariant cylinder}
\footnote{More precisely, these cylinders are weakly invariant, i.e. the associated vector field is tangent to 
them, but orbits may escape through the boundary}. 
To consider homoclinic connections, we lift the system to a double covering map, then the homoclinic 
connections becomes \emph{heteroclinic} connections between the two copies. 

We say that the pair $(H_\epsilon, c)$ satisfy the \emph{Arnold mechanism} if $\tcA^0_{H_\epsilon}(c)$ is 
an invariant curve, and there exists a symplectic double covering map $\Xi: \T^2 \times \R^2 \to \T^2 \times \R^2$, 
such that the set 
\[
	\tcN_{H_\epsilon \circ \Xi}^0(\Xi^*c) \setminus \tcA^0_{H_\epsilon \circ \Xi}(\Xi^*c)
\]
is totally disconnected. If $\tcA^0_{H_\epsilon}(c)$ is a smooth invariant curve with transversal 
intersection of stable and unstable manifolds, then the above set is discrete. 

\subsubsection*{Bifurcation mechanism} This is technically similar to the Arnold mechanism, but 
happens when the Aubry set $\tcA_{H_\epsilon}(c)$ is contained in two disjoint normally hyperbolic 
invariant cylinders. 

We say that the pair $(H_\epsilon, c)$ satisfy the \emph{Bifurcation mechanism} if the set 
\[
	\tcN^0_{H_\epsilon}(c) \setminus \tcA^0_{H_\epsilon}(c)
\]
is totally disconnected.

\begin{figure}[t]
 \centering 
 \includegraphics[width=2in]{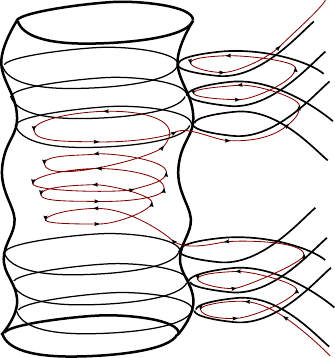} 
 \caption{Diffusion along a cylinder} \label{fig:arnold-mather}
\end{figure}

\subsubsection*{Normally hyperbolic invariant cylinders} The three mechanisms given do not apply to all cases, even after a generic perturbation. The main observation is that they do apply to cohomologies that satisfy:
\begin{enumerate}
	\item (2D NHIC) The Aubry set $\tcA^0_{H_\epsilon}(c)$ is contained in a two-dimensional normally hyperbolic invariant cylinder. 
	\item (1D Graph Theorem) The Aubry set $\tcA^0_{H_\epsilon}(c)$ is contained in a Lipschitz graph over the circle $\T$. 
\end{enumerate}
 In this case, we say that the pair $(H_\epsilon, c)$ is of \emph{Aubry-Mather type}.

Under these assumptions, the Aubry set resembles the \emph{Aubry-Mather sets} for twist maps, and in particular, generically we have the following dichotomy: either $\pi \tcN_{H_\epsilon}^0(c)$ is contractible, or $\tcA_{H_\epsilon}^0(c)$ is a Lipschitz invariant curve. In the latter case, we can show that Arnold mechanism applies after an additional perturbation. Since either Mather or Arnold mechanism applies, we conclude that $\tcA(c)$ is connected to $\tcA(c')$ for $c, c'$ close. Moreover, this argument can be continued if $c'$ is also of Aubry-Mather type. Dynamically, the orbit is either diffusing along the heteroclinic orbits of invariant curves, or diffusing in a Birkhoff region of instability within the cylinder. See Figure~\ref{fig:arnold-mather}. We now briefly describe the non-degeneracy conditions.

While a cohomology of Aubry-Mather type is robust, namely it can be extended along a continuous curve, in 
a one-parameter family one may encounter a bifurcation where the Aubry set jumps from one cylinder to another 
one. At the bifurcation, the Aubry set is contained in both cylinders. We say that the pair $(H_\epsilon, c)$ 
if of \emph{Bifurcation Aubry-Mather type} if the Aubry set is possibly contained in two cylinders. 

Technically we have to involve a different bifurcation type, called \emph{asymmetric bifurcation type}. This is very similar to the bifurcation Aubry-Mather type, the main difference is on one side of the bifurcation, the Aubry set is a Aubry-Mather type set contained in a invariant cylinder, while the other side we have a hyperbolic periodic orbit. This happens when we cross double resonance, see Definition~\ref{def:asym-bif}.

\subsubsection*{Forcing relation}

The rigorous formulation of the three diffusion mechanisms will be given using the concept of \emph{forcing equivalence} defined by Bernard in \cite{Be} (which is  generalization of a equivalence relation defined by Mather, see \cite{Ma2}). If $c,c'$ are forcing equivalent (denoted $c \dashv \vdash c'$), then there is a heteroclinic orbit connecting the associated Aubry sets. Moreover, there exists orbit shadowing an arbitrary sequence of cohomologies, as long as they are all equivalent. See Section \ref{sec:suff-diffusion} for more details. 

The main theorem reduces to Theorem~\ref{thm:c-equiv-path}, which proves forcing equivalence of a net of cohomologies, called $\Gamma_*$. The set $\Gamma_*$ consists of finitely many smooth curves. On each of the smooth curves, we prove the cohomologies are of Aubry-Mather type, and therefore one of the three mechanisms apply. 

We prove forcing equivalence of different connected components directly, using the definition of the forcing relation. We call this the \textbf{Jump mechanism}.

\subsection{Three regimes of diffusion}

Recall that our plan is to choose a net $\Gamma_*$ of cohomology classes, and prove their forcing equivalence, by first proving they are of Aubry-Mather type or Bifurcation Aubry-Mather type. This is done in three distinct regimes. 

\subsubsection*{Single resonance}

Let $(k_1, \Gamma_{k_1}) \in \cK$ be one of the single resonant component, and let $\cK^\st(k_1, \Gamma_{k_1}, \lambda)$ be the collection of the strong additional resonances. Then for $p$ in a $O(\sqrt \varepsilon)$-neighborhood of the set
\[
	\Gamma_{k_1}^{SR}{(M,\lb)} 
:= \Gamma_{k_1} \setminus \left( \bigcup_{k_2 \in \cK^\st(k_1, \Gamma_{k_1}, \lambda)} B_{M\sqrt{\epsilon}}(\Gamma_{k_1, k_2}) \right),
\]
where $M$ is a large parameter, the system admits the normal form
\[
	N_\epsilon^{SR} = H_0 + \epsilon Z(\theta^s, p) + O(\epsilon \delta), \quad
	(\theta^s, \theta^f, t)\in \T^3. 
\]
where how small $\delta$ is depends on how many double resonances we exclude. 
Under the non-degeneracy conditions $SR(k_1, \Gamma_{k_1},\lambda)$, the above system 
admits three dimensional (for the flow) normally hyperbolic invariant cylinders, and one can prove each 
$c \in \Gamma_{k_1}^{SR}{(M,\lb)}$ is of AM or of Bifurcation AM type. 

\subsubsection*{Double resonance, high energy}

Let $p_0 = \Gamma_{k_1, k_2}$ be a double resonance. On the set 
\[
	B_{M\sqrt{\epsilon}}(p_0), \quad p_0 = \Gamma_{k_1} \cap \Gamma_{k_2}, 
\]
we perform a normal form transformation, and then the $p$ variable via $I = (p-p_0)/\sqrt{\epsilon}$. 
One can show that the system is conjugate to:
\[
	\frac{1}{\beta}\left( K(I) - U(\varphi) + O(\sqrt{\epsilon}) \right), 
\]
where $K: \R^2 \to \R$ is a positive quadratic form and $U: \T^2 \to \R$, and $\beta>0$ is
 a constant depending only on $k_1, k_2$. The system $H^s = K(I) - U(\varphi)$ is a two degrees of 
freedom mechanical system. {Below we use  the shifted energy $E:= H^s(\varphi, I)+\min U(\varphi)$
as a parameter.}

When the shifted energy $E$ is not too close to $0$, we are in the 
{\it the high energy regime}. By imposing the conditions $[DR1^h] - [DR3^h]$, 
one shows existence of two-dimensional normally hyperbolic invariant 
cylinders associated to the shortest loops for the associated Jacobi metric, 
with the shifted energy as a parameter. This cylinder persists under perturbation, 
and one can show that the associated cohomologies are of Aubry-Mather type. 

\subsubsection*{Double resonance, low energy}

As the energy decreases, the cylinder constructed in the high energy may 
not persist. Under the non-degeneracy conditions $DR(k_1, k_2)$, 
we distinguish two separate cases:
\begin{enumerate}
	\item Simple cylinder: in this case the cylinder extends across zero {shifted} 
energy to negative {shifted} energy. 	In this case one can still show the associated 
cohomologies are of Aubry-Mather type. 
	\item Non-simple cylinder: In this case the cylinder may be destroyed 
before {the shifted} energy becomes zero. However, we show 
the existence of two simple cylinders near the non-simple one, and 
one can ``jump'' from one cylinder to another one. 

	This is the only case where the \textbf{Jump mechanism} is used. 
\end{enumerate}

\section{Forcing relation}

\subsection{Sufficient condition for Arnold diffusion}
\label{sec:suff-diffusion}

Recall that we will utilize the concept of forcing equivalence, denoted $c \dashv \vdash c'$. The actual definition will not be important for the current discussions, instead, we state its main application to Arnold diffusion. 

\begin{prop}[\cite{Be}, Proposition~0.10]
	\label{prop:forcing-implications}
Let $\{c_i\}_{i = 1}^N$ be a sequence of cohomology classes which are forcing equivalent. For each $i$, let $U_i$ be neighborhoods of the discrete Mather sets $\tcM_H^0(c_i)$, then there is a trajectory of the Hamiltonian flowvisiting all the sets $U_i$. 
\end{prop}

Let $\sigma>0$ and $\cV_\sigma(H)$ denote the $\sigma$ neighborhood of $H$ in the space $C^r(\T^2\times B^2 \times \T)$ with respect to the natural $C^r$ topology. The following statement is a ``local'' version of our main theorem, where we state that given $H_1 \in \cU$, we can: 

(1) Choose $\epsilon_0$ to be locally constant on a neighborhood of $H_1$; 

(2) Prove forcing equivalence on a residual subset of a neighborhood of $H_\epsilon$. \\

\begin{thm}\label{thm:c-equiv-path}
Let $\cP$ be a diffusion path and $U_1, \dots, U_N$ be open sets intersecting 
$\cP$.  Then there is and open and dense subset $\cU(\cP) \subset \cS^r$, and 
for each $H_1 \in \cU(\cP)$, there are $\delta = \delta(H_0, H_1) > 0$,  
$\epsilon_1 = \epsilon_1(H_0, H_1) > 0$ such that for each 
\[
	H_1' \in \cV_\delta(H_1), \quad 0 < \epsilon < \epsilon_1, 
\]
there is a subset $\Gamma_*(\epsilon, H_0, H_1') \subset \R^2$ satisfying  
\[
	\Gamma_* = \Gamma_*(\epsilon, H_0, H_1') \cap U_i \ne \emptyset, 
\quad i = 1, \dots, N,  
\]
with the property that 
there is  $\sigma = \sigma(\epsilon, H_0, H_1') > 0$, and a residual subset 
$\cR_\sigma(H_0 + \epsilon H_1') \subset \cV_\sigma(H_0 + \epsilon H_1')$, 
such that for each $H' \in \cR$, with respect to the Hamiltonian $H'$, all 
the $c \in \Gamma_*(\epsilon, H_0, H_1')$ are forcing equivalent. 
\end{thm}

Proposition~\ref{prop:forcing-implications} and Theorem~\ref{thm:c-equiv-path} imply our main theorem. 
\begin{proof}[Proof of Theorem~\ref{main}]
	First of all, let us define the lower semi-continuous function $\epsilon_0$. For each $H_1 \in \cU$, define
	\[
		\epsilon_2^{H_1}(\cdot) = \epsilon_1(H_0, H_1) \mathbf{1}_{\cV_\delta(H_1)}(\cdot)  
	\]
	where $\mathbf{1}_{\cV}$ denote the indicator function of $\cV$. An indicator function of an open set is lower semi-continuous by definition. For $H_1 \in \cS^r \setminus \cU$, let $\epsilon_2^{H_1} \equiv 0$.  We then define 
	\[
		\epsilon_0(\cdot) =  \sup_{H_1 \in \cS^r} \epsilon_2^{H_1}(\cdot), 
	\]
	which is lower semi-continuous, being an (uncountable) supremum of lower semi-continuous function. Note 
that $\epsilon_0$ is positive on each $H_1 \in \cU$ since $\epsilon_2^{H_1}(H_1) = \epsilon_1(H_0, H_1) > 0$. 

	Consider now $H_1 \in \cU$ and $0 < \epsilon < \epsilon_0(H_1)$ as defined above. Let $\Gamma_*(\epsilon_0, H_0, H_1)$ be as in Theorem~\ref{thm:c-equiv-path}.  Let $c_i \in U_i \cap  \Gamma_*(\epsilon, H_0, H_1)$. For any $\|H - H_0\|_{C^r} \le \epsilon$, there is $C>0$ depending only on $D$ such that  $\tcM^0_H(c) \subset \T^2 \times  B_{C\sqrt{\epsilon}}(c)$. As a result, reducing $\epsilon_0$ if necessary (note that minimum of an lower semi-continuous function and a constant is still lower semi-continuous) , we have $\tcM^0_H(c_i) \subset  \T^2 \times U_i$. Since $c_i$ are all forcing equivalent by Theorem~\ref{thm:c-equiv-path}, Proposition~\ref{prop:forcing-implications} implies the existence of an orbit visiting each neighborhood 
	$ \T^2 \times U_i$. 

	Since the above discussion applies to all $H \in \cR_\sigma(H_0 + \epsilon H_1)$ where  $H_1 \in \cU$, $0 < \epsilon < \epsilon_0(H_1)$, we conclude that for a dense subset of $\cV(\cU, \epsilon_0)$ (as defined in Theorem~\ref{main}), there is an orbit visiting each $\T^2 \times U_i$. Since this property is open due to the smoothness of the flow, it holds on an open and dense subset $\cW$ of $\cV$. 
\end{proof}

The set $\Gamma_*(\epsilon, H_0, H_1)$ will be chosen to be the union of finitely many smooth curves, and will coincide with $\cP$ except on finitely many neighborhoods of size $O(\sqrt{\epsilon})$ of strong double 
resonances.

\subsection{Diffusion mechanisms via forcing equivalence}

We reformulate the diffusion mechanisms introduced in Section~\ref{sec:intro-scheme} using forcing equivalence. We start with Mather mechanism. 
\begin{prop}[\cite{Be}, Theorem 0.11]\label{prop:mather-mech}
Suppose 
\begin{equation}
	\label{eq:mather-mech}\tag{$\mathbf{Ma}$}
	\cN_H^0(c) \text{ is contractible.}
\end{equation}
as a subset of $\T^2$, then there is $\sigma>0$ such that $c$ is forcing equivalent to all $c' \in B_\sigma(c)$. 
\end{prop}

To define the Arnold mechanism, we consider a finite covering of our space. Let 
\[
	\xi: \T^n \to \T^n
\]
be a linear double covering map, for example: $(\theta_1, \theta_2) \mapsto (2\theta_1, \theta_2)$. Then $\xi$ 
lifts to a symplectic map 
\[
	\Xi: \T^n \times \R^n  \to \T^n \times \R^n , \quad
	\Xi(\theta, p) = (\xi \theta, \xi^*p), 
\]
where $\xi^*(p)$ is defined by the relation $\xi^*(p) \cdot v = p \cdot d\xi(v)$ for all $v \in \R^2$. For example, if $n = 2$ and $\xi(\theta_1, \theta_2) = (2\theta_1, \theta_2)$ we have $\xi^*(p_1, p_2) = (p_1/2, p_2)$. This allows us to consider the lifted Hamiltonian $H\circ \Xi$. 

\begin{lem}(\cite{Be}, Section 7)\label{lem:covering-aub}
We have 
\[
\tcA^0_{H \circ \Xi}(\xi^*c) = \Xi^{-1} \tcA^0_H(c), \quad \tcN^0_{H \circ \Xi}(\xi^*c) \supset \Xi^{-1} \tcN^0_H(c). 
\]
Moreover, $\xi^*c \vdash \xi^* c'$ relative to $H\circ \Xi$ implies $c \vdash c'$ relative to $H$. 
\end{lem}

The Aubry set can be decomposed into disjoint invariant sets called {\it static classes}, which gives important insight into 
the structure of the Aubry set. In particular, when there is only one static class, then $\tcA_H^0(c) = \tcN_H^0(c)$. 
In the case $\tcA_H^0(c) \ne \tcN_H^0(c)$, the difference $\tcN_H^0(c)\setminus \tcA_H^0(c) $ 
consists of heteroclinic orbits from one static class to another (\cite{Be}). Using Lemma~\ref{lem:covering-aub}, when $\tcA_H^0(c) = \tcN_H^0(c)$, it may happen that $\tcA^0_{H \circ \Xi}(\xi^*c) \subsetneq \tcN^0_{H \circ \Xi}(\xi^*c)$, and the difference provides additional heterclinic orbits to the Aubry set that is not contained in the Ma\~ne set before the lifting. This can be exploited to create diffusion orbits.

\begin{prop}[\cite{Be}, Theorem 9.2, Proposition 7.3]\label{prop:arnold-mech}
Suppose, either:
\begin{equation}
	\label{eq:bifurcation-mech} \tag{$\mathbf{Bif}$}
	\tcA_H^0(c) \text{ has two static classes, and } \tcN_H^0(c) \setminus \tcA_H^0(c) 
	\text{ is totally disconnected, }
\end{equation}
or:
\begin{equation}
	\label{eq:arnold-mech} \tag{$\mathbf{Ar}$}
	\tcA_H^0(c) \text{ has one static class, and } \tcN^0_{H \circ \Xi}(\xi^* c)  \setminus \tcA^0_{H \circ \Xi}(\xi^* c) \text{ is totally disconnected.}
\end{equation}
Then  there is $\sigma>0$ such that $c$ is forcing equivalent to all $c' \in B_\sigma(c)$. 
\end{prop}

As implied by the labeling, the first item is called the bifurcation mechanism, and the second the Arnold mechanism. 
We obtain the following immediate corollary:
\begin{cor}[Mather-Arnold mechanism]
	Suppose $\Gamma \subset B^2$ is a continuous curve, and 
	for each $c \in \Gamma$ one of the diffusion mechanisms \eqref{eq:mather-mech}, \eqref{eq:bifurcation-mech}, or \eqref{eq:arnold-mech} holds. Then all $c \in \Gamma$ are forcing equivalent. 
\end{cor}

Recall that $\Gamma_*(\epsilon, H_0, H_1)$ can be chosen as a union of finitely many smooth curves. 
\begin{itemize}
	\item We will later show that for $c \in \Gamma_*(\epsilon, H_0, H_1)$ in Theorem~\ref{thm:c-equiv-path}, one of the two applies: Proposition \ref{prop:mather-mech}
	or  Proposition \ref{prop:arnold-mech}. 

	As a result, each connected component of $\Gamma_*(\epsilon, H_0, H_1)$ 
	consists of equivalent $c$'s.

	\item We prove the forcing equivalence between different connected components 
	using directly the definition of forcing relation. We call this {\it the ``jump'' mechanism.} 
\end{itemize}

\subsection{Invariance under the symplectic coordinate changes}

A diffeomorphism $\Psi = \Psi(\theta, p): \T^n \times \R^n \to \T^n \times \R^n$ is called exact symplectic if $\Psi^* \lambda - \lambda$ is an exact one-form, where $\lambda = \sum_{i=1}^n p_i d\theta_i$ is the canonical form. We say $\Phi: \T^n \times \R^n \times \T \to \T^n \times \R^n \times \T$ is {\it exact symplectic} if  
\[
\Phi(\theta, p, t) = (\Phi_1(\theta, p, t), t), 
\]
and there is $\widetilde{E} = \widetilde{E}(\theta, p, t)$, such that 
\[
\Psi(\theta, p, t, E) = (\Phi(\theta, p, t), E + \widetilde{E}(\theta, p, t))
\]
(called the autonomous extension of $\Phi$) is exact symplectic. The new term $\widetilde{E}(\theta, p, t)$ is defined up to adding a function $f'(t)$, where $f(t)$ is periodic in $t$. Let us assume $\widetilde{E}(0, 0, t) \equiv 0$, therefore the choice of $\widetilde{E}$ is unique. 

Let $H = H(\theta, p, t)$, and $\Phi$ is exact symplectic with extension $\Psi$. Then for $G(\theta, p, t, E) = H (\theta, p, t) + E$, we define
\begin{equation}
  \label{eq:Phi-star-H}
  \Phi^* H = \Psi^* H = G \circ \Psi (\theta, p, t, E) - E = H \circ \Phi + \widetilde{E}(\theta, p, t). 
\end{equation}

The Aubry, Mather, Ma\~ne sets are invariant under exact symplectic coordinate change in the following sense. 
\begin{prop}[\cite{Be2}, \cite{MS2016}] \label{prop:symp-inv}
Suppose $H$ and $\Phi^*H$  are Tonelli, and let $\Psi$ be the extension of $\Phi$. Let $(c, \alpha) \in \R^n \times \R \simeq H^1(\T^n \times \T, \R)$, and let $\Psi^*(c, \alpha) = (c^*, \alpha^*)$ be the push forward of the  cohomology class via the identification $H^1(\T^2 \times \T, \R)\simeq H^1(\T^n \times \R^n \times \T \times \R, \R)$. Then 
\[
\alpha = \alpha_H(c) \quad \Longleftrightarrow \quad \alpha^* = \alpha_{\Phi^*H}(c^*). 
\]
Let us denote 
\[
(\Phi^*_H c, \alpha^*) = \Psi^*(c, \alpha_H(c)), 
\]
then
\[
	\tcM_H(c) = \Phi\left(  \tcM_{\Phi^*H}(\Phi_H^*c) \right), \ 
	\tcA_H(c) = \Phi\left(  \tcA_{\Phi^*H}(\Phi_H^*c) \right), \ 
	\tcN_H(c) = \Phi\left(  \tcN_{\Phi^*H}(\Phi_H^*c) \right).
\]
Note in the particular case when $\Phi$ is homotopic to identity, $\Phi_H^*c = c$ and $\alpha_{\Phi^*H}(c) = \alpha_H(c)$. 
\end{prop}

\begin{lem}\label{lem:symp-inv-normal-form}
Let $\Phi$ be an exact symplectic coordinate change. The tuple $(H, c)$ satisfies \eqref{eq:bifurcation-mech}  or \eqref{eq:arnold-mech} if and only if $(\Phi^*H, \Phi_H^*c)$ satisfies the same conditions. The property \eqref{eq:mather-mech} with the additional condition  that $\tcA_H^0(c) = \tcN_H^0(c)$ is also invariant under exact symplectic coordinate changes. 
\end{lem}
\begin{proof}
	Since our symplectic coordinate changes are always identity in the $t$ components, the invariance of Aubry and Ma\~ne sets imply the invariance of their zero section under the map $\Phi(\cdot, \cdot, 0)$. The invariance of  \eqref{eq:bifurcation-mech} follows. For the invariance of \eqref{eq:mather-mech}, note that due to the graph property, $\tcA_H^0(c)$ is contractible in $\T^2 \times \R^2$ if and only if $\cA_H^0(c)$ is contractible in $\T^2$. Therefore the contractibility of Aubry set is invariant, and since the Aubry set coincide with the Ma\~ne set by assumption, \eqref{eq:mather-mech} is invariant. 

	For \eqref{eq:arnold-mech}, let $\Psi$ be the extension of $\Phi$, and let us extend $\Xi$ trivially to $\T^n \times \R^n \times \T$ or $\T^n \times \R^n \times \R$ without changing its name. Let $\Phi_1$ be an exact symplectic change homotopic to $\Phi$, with extension $\Psi_1$, such that 
	\[
	  	\Phi \circ \Xi = \Xi \circ \Phi_1, \quad \Psi \circ \Xi = \Xi \circ \Psi_1. 
	\]
	Let us note $\Xi^*c$, defined as the push forward of $H^1(\T^n \times \R^n \times \T \times \R, \R)$ under the identification with $H^1(\T^n \times \T, \R)$, is identical to $\xi^*c$. We have:
	\[
	\Psi_1^*(\Xi^*c, \alpha_{H \circ \Xi}(\Xi^*c)) = \Psi_1^* \Xi^* (c, \alpha_H(c)) = \Xi^* \Psi^* (c, \alpha_H(c)). 
	\]
	Since $\Xi$ is independent of $t$,  $\Xi^*$ is identity in the last component. We conclude that 
	$(\Phi_1)_{H \circ \Xi}^* \,\Xi^*c   = \Xi^*(\Phi_H^* c)$. Moreover,
	\[ 	   
 	   (\Phi \circ \Xi)^* H = (\Phi^*H )\circ \Xi, \quad
 	   (\Xi \circ \Phi_1)^* H = \Phi_1^* (H \circ \Xi). 
	\]

	Then 
	\[
		\Phi_1 \left( \tcN_{(\Phi \circ \Xi)^* H}((\Phi \circ \Xi)_H^* c) \right) = \Phi_1 \left( \tcN_{\Phi_1^* (H \circ \Xi)}((\Phi_1)_{H \circ \Xi}^* \ \Xi^*c) \right) =  \tcN_{H \circ \Xi}(\Xi^* c)
	\]
	and 
	\[
		\Phi_1\left( \Xi^{-1} \tcN_{\Phi^*H}(\Phi_H^* c) \right) = \Xi^{-1} \Phi \left( \Xi^{-1} \tcN_{\Phi^*H}(\Phi_H^* c) \right) = \Xi^{-1}\left(  \tcN_H (c) \right), 
	\]
	therefore 
	\[
	\Phi_1\left( \tcN_{(\Phi^*H) \circ \Xi}(\Xi^*(\Phi_H^* c)) \setminus  \Xi^{-1} \tcN_{\Phi^*H}(\Phi_H^* c)   \right) = \tcN_{H \circ \Xi}(\Xi^* c) \setminus \Xi^{-1} \tcN_H (c). 
	\]
	This implies invariance of $\eqref{eq:arnold-mech}$ after considering the zero section of the above equality. 
\end{proof}

Our definition of exact symplectic coordinate change for time-periodic system is somewhat restrictive, and in particular, it does not apply directly to the linear coordinate change performed at the double resonance. In that setting, we will prove invariance of Mather, Aubry and Ma\~ne set directly. 

\subsection{Normal hyperbolicity and Aubry-Mather type}
\label{sec:AM-type}

Call a two dymensional normally hyperbolic invariant cylinder {\it symplectic} if the restriction 
of the canonical form for this cylinder is non-degenerate on the domain of definition.  Loosely speaking, a pair 
$(H_*, c_*)$ is called of {\it Aubry-Mather type} (AM type for short, refer to Definition~\ref{defn:AM} for details) if: 
\begin{enumerate}
	\item The discrete Aubry set $\tcA_{H_*}^0(c_*)$ is contained in \emph{two dimensional} 
	normally hyperbolic invariant cylinder, the restriction of the symplectic form is non-degenerate on the cylinder. 
	\item There is $\sigma > 0$ such that the following holds for $c \in B_\sigma(c_*)$ and 
$H \in \cV_\sigma(H_*)$\ 
:
	\begin{enumerate}
		\item  The discrete Aubry set satisfies the graph property under the local coordinates of the cylinder.
		\item When the Aubry set is an invariant graph, then locally the unstable manifold of the Aubry set is a graph over the configuration space $\T^2$. 
	\end{enumerate}
\end{enumerate}
This definition gives an abstract version of the setting seen in the \emph{a priori unstable}  systems.

\begin{thm}[See Theorem~\ref{thm:AM-non-deg}]
	\label{thm:intro-AM-non-deg}
Suppose $H_* \in C^r$, $r \ge 2$ and $(H_*, c_*)$ is of Aubry-Mather type, $\Gamma \subset \R^2$ is a smooth curve containing $c_*$ in the relative interior. Then there is $\sigma > 0$ such that for all 
$c \in \overline{B_\sigma(c_*)}\cap \Gamma$, the following dichotomy holds 
for a $C^r$-residual subset of $H \in \cV_\sigma(H_*)$: 
\begin{enumerate}
	\item  Either the projected Ma\~ne set $\cN^0_{H}(c)$ is contractible as a subset of $\T^2$ (Mather mechanism \eqref{eq:mather-mech});
	\item  Or there is a finite covering map $\Xi$ such that the set 
	\[
		\tcN^0_{H \circ \Xi}(\xi^*c) \setminus  \Xi^{-1} \tcN^0_{H}(c)
	\]
	is totally disconnected (Arnold mechanism \eqref{eq:arnold-mech}). 
\end{enumerate}
\end{thm}

We say $(H_*, c_*)$ is of \emph{bifurcation Aubry-Mather type} if there exists two normally hyperbolic invariant cylinders, such that the \emph{local Aubry set} restricted to each cylinder satisfy the conditions of Aubry-Mather type. The precise definition is given in Definition~\ref{def:bif}. 

We will also consider a particular (and simpler) bifurcation. We say $(H_*, c_*)$ is of \emph{asymmetric bifurcation type} if there exists one normally hyperbolic invariant cylinder, and a hyperbolic periodic orbit, such that the Aubry set is either contained in the union of the cylinder (and of Aubry-Mather type), and the periodic orbit, see Definition~\ref{def:asym-bif}. 

We state the consequence of these definitions in terms of diffusion. 

\begin{thm}[See Theorem~\ref{thm:bifur-non-deg}]
	\label{thm:intro-bifur-non-deg}
	Suppose $H_* \in C^r$, $r \ge 2$ and $(H_*, c_*)$ is of bifurcation Aubry-Mather type or asymmetric bifurcation type, and $\Gamma \subset \R^2$ is 
a smooth curve containing $c_*$ in the relative interior. Then there is 
$\sigma>0$ and an open and dense subset $\cR\subset \cV_\sigma(H_*)$ 
such that for each $H \in \cR$ and each $c \in \Gamma \cap B_\sigma(c_*)$,  \eqref{eq:bifurcation-mech} holds on at most finitely many $c$'s, and for all 
other $c$'s either \eqref{eq:mather-mech} or \eqref{eq:arnold-mech} holds.
\end{thm}

The following Proposition is a direct consequence of Theorem~\ref{thm:intro-AM-non-deg} and Theorem~\ref{thm:intro-bifur-non-deg}. 
\begin{prop}\label{prop:forcing-connected}
Suppose $\Gamma \subset B^2$ is a piecewise smooth curve of cohomologies such that for each $c \in \Gamma$, such that the pair $(H_*, c)$ is of  Aubry-Mather type, bifurcation AM type or asymmetric bifurcation type. Then there is $\sigma>0$ 
and a residual subset $\cR_\sigma(H_*) \subset \cV_\sigma(H_*)$, such that either \eqref{eq:mather-mech}, \eqref{eq:bifurcation-mech}, or \eqref{eq:arnold-mech} holds for each $H \in \cR_\sigma(H_*)$ and each $c \in \Gamma$. 
\end{prop}
\begin{proof}
For a piecewise smooth $\Gamma = \bigcup_{i = 1}^m \Gamma_i$, we can extend each $\Gamma_i$ to $\Gamma_i'$ smoothly, such that $\Gamma_i$ is contained in the relative interior of $\Gamma_i'$. We then apply Theorem~\ref{thm:intro-AM-non-deg} and Theorem~\ref{thm:intro-bifur-non-deg} to each $c \in \Gamma_i$ relative to the smooth curve $\Gamma_i'$, to get the conclusion of our proposition for $c' \in B_{\sigma(c)}(c)$, and $H \in \cR_{\sigma(c)}(H_*) \subset \cV_{\sigma(c)}(H_*)$. The proposition then follows by considering a finite covering of $\Gamma$ by $B_{\sigma(c_j)}(c_j)$, and taking finite intersection of residual subsets $\cR_{\sigma(c_j)}(H_*)$. 
\end{proof}

We now describe the selection of cohomologies and prove AM type in each of the two regimes. Single resonance is covered in Section~\ref{sec:SR-overview}, and double resonance is split into two sections, Section~\ref{sec:DR-geometrical} covers the geometrical part, while Section~\ref{sec:coh-DR} covers the variational part.

\section{Normal forms and cohomology classes at single resonances}
\label{sec:SR-overview}

\subsection{Resonant component and non-degeneracy conditions}
\label{sec:intro-SR}
Let $(k_1, \Gamma_{k_1}) \in \cK$ be a resonant segment in the diffusion path. Define the resonant component of $H_1$ relative to the single resonance $k_1$ as follows:
\[
	[H_1]_{k_1}(\theta, p, t) = \sum_{k \in k_1 \Z} h_k(p) e^{2\pi i k \cdot (\theta, t)}, 
\]
where $h_k$ are the Fourier coefficients of $H_1(\theta, p, t)$. Since $[H_1]_{k_1}$ only depends on the variables $k_1 \cdot (\theta, t)$ and $p$, we define 
\[
	Z_{k_1}: \T \times \R^2 \to \R, \quad 
	Z_{k_1}(k_1 \cdot (\theta, t), p) = [H_1]_{k_1}(\theta, p, t). 
\]

For  $p_0 \in \Gamma_k$, define the following conditions: 
\begin{itemize}
	\item[{[SR$1_\lambda$]}] For all $p \in B_\lambda(p_0)$, the function $Z_{k_1}(\cdot, p)$ achieves a global maximum at $\theta_*^s(p) \in \T$, and 
	\[
		Z_{k_1}(\theta^s, p) - Z_{k_1}(\theta^s_*(p), p) < \lambda d(\theta^s, \theta^s_*(p))^2. 
	\]
	\item[{[SR$2_\lambda$]}] For all  $p\in B_\lambda(p_0)$,
	there exists two local maxima
	$\theta^s_1(p)$ and 
	$\theta^s_2(p)$ of the function $Z_{k_1}(.,p)$ 
	in $ \Tm^{n-1}$ satisfying  
	\begin{align*}
		\partial^2_{\theta^s} Z_{k_1}(\theta^s_1(p),p)< \lambda I
		\quad,\quad 
		\partial^2_{\theta^s} Z_{k_1}(\theta^s_2(p),p)< \lambda I,\qquad \qquad \qquad \\
		Z_{k_1}(\theta^s,p) < \max \{Z_{k_1}(\theta^f_1(p),p),Z_{k_1}(\theta^f_2(p),p)\}- \lambda
		\big(\min\{d(\theta^s-\theta^s_1), d(\theta^s-\theta^s_2)\}\big)^2. 
	\end{align*}
\end{itemize} 

\begin{defn}
	We say that $H_1$ satisfy the condition $SR(k_1, \Gamma_{k_1},  \lambda)$ if 
	for each $p_0 \in \Gamma_{k_1}$, at least one of 
	[SR$1_\lambda$] and [SR$2_\lambda$] holds for 
	the function $Z_{k_1}(\theta^s, p)$. 
\end{defn}

\begin{prop}\label{prop:sr-non-deg}
	The set of $H_1 \in \cS^r$ such that $SR(k_1, \Gamma_{k_1}, \lambda)$ holds for some $\lambda > 0$ is open and dense. 
\end{prop}
\begin{proof}
It suffices to show that a generic one-parameter family $f(x, a)$, $x \in \T$, $a \in [a_1, a_2]$ the function $f(\cdot, a)$ has  unique non-degenerate maximum, with the exception of up to finitely many $a$'s for which there are two non-degenerate maxima. We can first show that the following property is open and dense: any local maxima in $x$ is non-degenerate ($\partial^2_{xx} f(x) < 0$). The main observation is that it's implied by a co-dimension two condition: we require whenever $\partial_x f  =0$, and $\partial_{xx} f =0$, we have $\partial_{xxx} f \ne 0$. Then any degenerate critical point cannot be a maxima. 

We obtain a finite family of local minima. We then can ``slide'' them against each other so that they intersect transversally. 
\end{proof}

Let $K$ be a large parameter, recall the strong additional resonances are defined by
\[
	\cK^\st(k_1, \Gamma_{k_1}, K) = \{ k_2 \in \Z^3_*: \,  
	|k_2| \le K, \, \Gamma_{k_1} \cap S_{k_2} \ne \emptyset \}. 
\]
We show generic forcing equivalence on each connected components of $\Gamma_{k_1}$ minus $O(\sqrt{\epsilon})-$neigh\-borhoods 
of the strong double resonances, called {\it punctures}. The following theorem is the main result of this section, the proof is given in Section~\ref{sec:res-comp} assuming propositions proved in the later sections. 
For $M,K>0$ denote 
\begin{equation}
	\label{eq:punctured-sr}
	\Gamma_{k_1}^{SR}(M,K) := 
\Gamma_{k_1} \setminus \left(  \bigcup_{k_2 \in \cK^\st(k_1, \Gamma_{k_1}, K)} \quad
	B_{2M\sqrt{\epsilon}}(\Gamma_{k_1, k_2})
	\right).
\end{equation}

\begin{thm}\label{thm:c-equiv-sr}
Suppose $H_1$ satisfy the condition $SR(k_1,\Gamma_{k_1}, \lambda)$ on 
$\Gamma_{k_1}$. Then there is $K = K(D, k_1, \lambda)$,  
$M = M(D, k_1, \lambda)$,  $\epsilon_1 = \epsilon_1(D, k_1, \lambda)>0$,  
$\sigma = \sigma(k_1, H_0,\epsilon, H_1)>0$, and for every 
$0 < \epsilon < \epsilon_1$, a residual subset 
$\cR \subset \cV_\sigma(H_0 + \epsilon H_1)$, such that the following 
hold for all $H \in \cR$: for each $c \in \Gamma^{SR}_{k_1}(M,K)$ 
{the associated Aubry or Ma\~ne sets satisfy} either 
\eqref{eq:mather-mech}, \eqref{eq:bifurcation-mech} or 
\eqref{eq:arnold-mech}. As a result, each  connected components of 
$\Gamma_{k_1}^{SR}(M,K)$ is contained in one forcing equivalent class. 
\end{thm}

\begin{figure}[t]
	\centering 
	\includegraphics[width=2.35in]{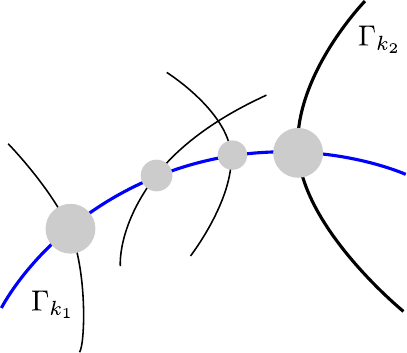} 
	\caption{Single resonance after removing punctures}
	\label{fig:sr-puncture}
\end{figure}

Refer to Figure~\ref{fig:sr-puncture} for an illustration.

\subsection{Normal form}

Then the classical partial averaging theory indicates that after a coordinate change, the system has 
the normal form $H_0 + \epsilon Z_{k_1} + h.o.t$ away from punctures. In order to state the normal form, 
we need an anisotropic norm adapted to the perturbative nature of the system. Define
\begin{equation}
	\label{eq:rescaled-norm}
	\|H_1(\theta, p, t)\|_{C^r_I} = \sup_{|\alpha| + |\beta| \le r} 
	\epsilon^{\frac{|\beta|}{2}}  \sup \left| \partial^\alpha_{(\theta, t)} \partial^\beta_p H_1(\theta, p, t) \right| ,
\end{equation}
where $\alpha \in (\Z_+)^3, \beta \in (\Z_+)^2$ are multi-indices and $|\cdot|$ denote the sum of the indices. The rescaled norm is similar to  $C^r$ norm, but replace the $p$ derivatives by the derivatives in $I = p/\sqrt{\epsilon}$, hence the name.

\begin{thm}[See end of this section]
	\label{thm:normal-form-sr} With the notations above there is $C = C(D, k_1) > 1$ such that 
the following hold. Let $\delta > 0$ be a small parameter and set $K = 1/\delta^2$ and $M = K^2$. 
Then for any $c \in \Gamma_{k_1}^{SR}(M,K)$ there exists $p_c \in \Gamma_{k_1}$ such that 
\begin{equation}
	\label{eq:rat-cover}
	c \in B_{CK\sqrt{\epsilon}/2}(p_c), 
\end{equation}
and an $C^\infty$ exact symplectic coordinate change homotopic to the identity
\[
	\Phi_\epsilon: \T^2 \times B_{CK \sqrt{\epsilon}/2}(p_c) \times \T \to 
	\T^2 \times B_{CK \sqrt{\epsilon}}(p_c) \times \T
\]
such that:
\begin{enumerate}
	\item 
	\begin{equation}
		\label{eq:sr-normal}
		(\Phi_\epsilon)^*H = H_0 + \epsilon [H_1]_{k_1} + \epsilon R, \qquad
		\|R\|_{C^2_I} \le C \delta. \qquad \qquad \qquad
	\end{equation}
	\item \quad \quad $\|\Pi_\theta(\Phi_\epsilon - Id)\|_{C^2_I} \le C \delta^4$ \  and \ $\|\Pi_p(\Phi_\epsilon - Id)\|_{C^2_I} \le C \delta^4 \sqrt{\epsilon}$.
\end{enumerate}
Here the $C^2_I$ norm is evaluated on the set $\T^2 \times B_{CK\sqrt{\epsilon}/2}(p_c) \times \T$. 
\end{thm}
\begin{rmk}
The $C^\infty$ coordinate change is obtained by approximating a coordinate change that is only $C^{r-1}$, see Appendix~\ref{sec:norm-form-double}. The reason this can be done is that we only need $C^2$ estimates of the coordinate change. 
\end{rmk}

We use the idea of Lochak (see for example \cite{Loc93}) to cover the action space with double resonances. 
A double resonance $p_0 = S_{k_1} \cap S_{k_2}$ corresponds to a periodic orbit of 
the unperturbed system $H_0$. More precisely, we have $\omega_0 = \Omega_0(p_0): =  \nabla H_0(p_0)$ which satisfies $\R (\omega, 1) \cap \Z^3 \ne \emptyset$. Denote by 
$T_{\omega_0} = \min\{t>0:\, t(\omega_0,1) \in \Z^3\}$ the minimal period. 

The resonant lattice for $p_0$ is $\Lambda = \Span_\R \{k_1, k_2\} \cap \Z^3$, and the resonant component is 
\[
	[H_1]_{k_1, k_2} = \sum_{k \in \Lambda} h_k(p) e^{2\pi i k \cdot (\theta, t)}. 
\]

We have the following general normal form theorem, the proof is given at the end of Section~\ref{sec:norm-form-double}. 
\begin{prop}
	\label{prop:normal-form-dr-gen}
	Let $p_0 = \Gamma_{k_1, k_2},$ $T = T_{\omega_0(p_0)}$. Then for a parameter $C_1 >1$,  there exists  $C = C(r, C_1) > 1$, $\epsilon_0 = \epsilon_0(r) >0$ such that if $K_1 > C$ satisfies
	\[
		T < \frac{C_1}{K_1^2 \sqrt{\epsilon}}, 
	\]
    then for each $0 < \epsilon < \epsilon_0$ there exists a $C^\infty$ exact symplectic map 
	\[
		\Phi: \T^2 \times B_{K_1\sqrt{\epsilon}/2} \times \T \to \T^2 \times B_{K_1\sqrt{\epsilon}} \times \T 
	\]
	such that 
	\[
		(\Phi)^*H_\epsilon = H_0 + \epsilon [H_1]_{k_1, k_2} + \epsilon R_1, 
	\]
	where 
	\[
		\|R_1\|_{C^2_I} \le C K_1^{-1}, 
	\]
	and 
	\[
		\left\| \Pi_\theta (\Phi - Id) \right\|_{C^2_I} \le C K_1^{-2}, \quad 
		\left\| \Pi_p (\Phi - Id) \right\|_{C^2_I} \le C K_1^{-2} \sqrt{\epsilon}. 
	\]
\end{prop}

Let us denote $\Lambda_1 = \Span_\R \{k_1\} \cap \Z^3$ and $\Lambda_2 = \Span_\R \{k_1, k_2\} \cap \Z^3$. 
\begin{lem}\label{lem:add-res}
There is an absolute constant $C>0$  such that if 
$$
\min\{|k|: \, k \in \Lambda_2 \setminus \Lambda_1\} \ge K > 0,
$$ we have 
\[
	\left\|  [H_1]_{k_1, k_2} - [H_1]_{k_1} \right\|_{C^2} \le C K^{-\frac12}. 
\]
\end{lem}
\begin{proof}
	First let us note that there is an absolute constant $C >0$ such that for each two dimensional lattice $\Lambda \subset \Z^3$, we have 
	\[
		\sum_{k \in \Lambda \setminus \{0\}} |k|^{-2 -\frac12} < C. 
	\]
	To see this, we can bound the sum above using the integral $\int_{|z| \ge 1, \, z \in \Span_\R \Lambda} |z|^{-2-\frac12}dA$. Then using the fact that $\|h_k\|_{C^2} \le |k|^{2-r}\|H_1\|_{C^r}$, we have 
	\[
		\left\|  [H_1]_{k_1, k_2} - [H_1] \right\|_{C^2} \le  \sum_{k \in \Lambda_2 \setminus \Lambda_1} |k|^{2-r} \le K^{-\frac12} \sum_{k \in \Lambda_2 \setminus\{0\}} |k|^{2 + \frac12 - r} < C K^{-\frac12}. 
	\]
\end{proof}

The following lemma is an easy consequence of the Dirichlet theorem (see \cite{Loc93}). 
\begin{lem}\label{lem:dirichlet} There is $C = C(D, k_1) > 0$ such that for each $Q_1>1$ and each $c \in S_{k_1}$, there is a double resonance $p_0$ with  $T_{\omega_0} < CQ_1$, and $\|c - p_0\| < C (T_{\omega_0}Q_1)^{-1}$, where $\omega_0 = \nabla H_0(p_0)$. 
\end{lem}

\begin{proof}
	[Proof of Theorem~\ref{thm:normal-form-sr}]
	Denote $\tau=K^2 \sqrt \varepsilon$. 
	First we apply Lemma~\ref{lem:dirichlet} using the parameter $Q_1 = C{\tau}^{-1}$, 
then each $c \in S_{k_1}$ is contained in the $\frac{\tau}{T(p_c)} \le K^2 \sqrt{\epsilon}$ neighborhood of a double 
resonance $p_c$,  whose period $T(\omega_c)$ is at most $C Q_1$, $\omega_c = \nabla H_0(p_c)$. Note that for $c$ in 
the set \eqref{eq:punctured-sr} (with $M = K^2$), we have 
	\[
		p_c \notin \cK^\st(k_1, \Gamma_{k_1}, K). 
	\]
	Let $p_c = \Gamma_{k_1, k_2}$, necessarily $\Lambda_2 \setminus \Lambda_1$ (see Lemma~\ref{lem:add-res}) contains only vectors larger than $K$, in this case we have $T_{\omega_c} \ge C^{-1} K$ where $C$ may depend on $k_1$. This lead to a better estimate
	\[
 	\|p_c - c\| \le \frac{K^2 \sqrt{\epsilon}}{T_{\omega_c}} \le  C K\sqrt{\epsilon}. 
	\]
	Moreover, from Lemma~\ref{lem:add-res}, 
	\[
		\|[H_1]_{k_1, k_2} - [H_1]_{k_1}\|_{C^2} < C K^{- \frac12}. 
	\]

	Let $K_1 = C K$, we have $T_{\omega_c} \le \frac{C}{K^2 \sqrt{\epsilon}} = \frac{C^3}{K_1^2 \sqrt{\epsilon}}$, therefore Proposition~\ref{prop:normal-form-dr-gen} applies with the parameter $C_1 = C^3$ and $K_1 = C K$, we obtain, for a different constant $C_2$
	\[
		\|R_1\|_{C^2_I} \le C_2 K_1^{-1} = C_2 C^{-1} K^{-1}, 
	\]
	therefore 
	\[
		H_\epsilon \circ \Phi = H_0 + \epsilon [H_1]_{k_1} + \epsilon R, 
	\]
	where 
	\[
	\|R\|_{C^2_I} = \|\epsilon([H_1]_{k_1, k_2} - [H_1]_{k_1}) + \epsilon R_1\|_{C^2_I} \le C_2 C^{-1} K^{-1}  +  C K^{-\frac12} \le 2 C K^{-\frac12}
	\]
	if $K$ is large enough. Moreover,
	\[
		\|\Pi_\theta(\Phi - Id)\|_{C^2_I} \le C_2 K_1^{-2} \le K^{-2} \le C_2 \delta^4,\] 
\[
		\|\Pi_p(\Phi -\Id)\|_{C^2_I} \le C_2 K_1^{-2} \sqrt{\epsilon} \le C_2 \delta^4 \sqrt{\epsilon}. 
	\]
\end{proof}

\subsection{The resonant component}
\label{sec:res-comp}

Using the fact that $k_1 = (k_1^1, k_1^0) \in \Z^2 \times \Z$ is space irreducible, there is $k_2 = (k_2^1, k_2^0)$ such that $B_0^T := \bmat{k_1^1 \\ k_2^1} \in SL(2, \Z)$.
\footnote{This is the only part where the space irreducibility of resonance is used. }
Define: 
\[
	B= \bmat{ (k_1)^T \\ (k_2)^T \\ \bmat{0 & 0 & 1} } \in SL(3, \Z), 
\]
and 
\[
	\Phi_L(\theta,  p, t) =  (\theta^s, \theta^f, p^s, p^f, t), \quad
	\bmat{\theta^s \\ \theta^f \\ t} = B \bmat{\theta \\ t},\quad
	\bmat{p^s \\ p^f} = (B_0^T)^{-1} p. 
\]
One verifies that $\Phi_L$ is a linear exact symplectic coordinate change. Note that $\theta^s = k_1 \cdot (\theta, t)$ and $[H_1]_{k_1}\circ \Phi_L$ depends only on $\theta^s, p^s, p^f$. Let us write 
\begin{equation}
	\label{eq:N-sr}
	N_\epsilon = (\Phi_\epsilon^* H_\epsilon) \circ \Phi_L = 
	H_0(p^s, p^f) + \epsilon Z(\theta^s, p^s, p^f) + \epsilon R(\theta^s, \theta^f, p^s, p^f, t),
\end{equation}
where we abused notation by keeping the name of $H_0$ and $R$ after the coordinate change. Let us also abuse notation by writing  $\theta = (\theta^s, \theta^f)$ and $p = (p^s, p^f)$. Let us note that $N_\epsilon$ is defined on the set 
\[
	\T^2 \times B_{K_1\sqrt{\epsilon}}(p_1) \times \T, \quad \text{ where }\quad 
	K_1 = \frac{2K}{\|B^{-1}\|}, \quad p_1 = (B_0^T)^{-1}\, p_0. 
\] 
and the resonant segment $\Gamma_{k_1}$ is represented by $\Gamma^s = \{p: \partial_{p^s}H_0 = 0\}$ in the new coordinates. 

Let us consider the following set:
\[
	\cR(\epsilon, \delta, p_1) = 
	\left\{ N_\epsilon = H_0 + \epsilon Z(\theta^s, p) + \epsilon R(\theta, p, t), 
	\quad 
	\|R\|_{C^2_I(\T^2 \times B_{K_1 \sqrt{\epsilon}}(p_1)\times \T )} < \delta
	\right\}.
\]

We show that the system $N_\epsilon$ admits a three dimensional normally hyperbolic invariant cylinder 
of the type 
\[
(\theta^s, p^f) = (\Theta^s, P^s)(\theta^f, p^f, t).
\]
and if we consider the discrete Aubry set, it is a graph over $\theta^f$ component. 
The details will be given 
in Section~\ref{sec:sr-AM}, here we state the consequences of those results:

\begin{prop}[See Theorem~\ref{thm:sr-AM-type-single}]\label{prop:single-max}
Assume that $Z(\theta^s, p)$ satisfies condition $[SR1_\lambda]$ at $p_1\in \Gamma^s$. 
Then there is $\delta_0, \epsilon_0>0$ depending on $D, \lambda$ such that if $0< \epsilon < \epsilon_0$ and $0 < \delta < \delta_0$, for each $N \in \cR(\epsilon, \delta, p_1)$, each $c \in B_{K_1 \sqrt{\epsilon}/2}(p_1) \cap \Gamma^s$ 
the pair $(N, c)$ is of Aubry-Mather type. 
\end{prop}

\begin{prop}[See Theorem~\ref{thm:sr-double-loc}]\label{prop:double-max}
Consider $N_\epsilon$ as in Proposition~\ref{prop:single-max}, and assume that $Z(\theta^s, p)$ satisfies condition $[SR2_\lambda]$ at $p_2 \in \Gamma^s$. Then there is $\delta_0, \epsilon_0 > 0$ depending on $D, \lambda$ such that if $0< \epsilon < \epsilon_0$ and $0 < \delta < \delta_0$, there is an open and dense subset $\cR_1 \subset \cR(\epsilon, \delta, p_1)$, such that each $c \in B_{K_1 \sqrt{\epsilon}/2}(p_1) \cap \Gamma^s$ 
the pair $(N, c)$ is of bifurcation Aubry-Mather type. 
\end{prop}

\begin{proof}
	[Proof of Theorem~\ref{thm:c-equiv-sr}]
	Choose $K_0 = 1/\delta_0$. Let $\Gamma$ be a connected component of \eqref{eq:punctured-sr}, and consider $c_0 \in \Gamma$. Then there is $p_0 \in \Gamma_k$ such that $c_0 \in B_{K_2\sqrt{\epsilon}/2}(p_0)$ (this is possible by choosing $Q$ large in \eqref{eq:rat-cover}), where $K_2 = K/\|B^{-1}\|^2$. After the coordinate change, $c$ is mapped to $c_1 = (M_0^T)^{-1}(c)$ which is contained in $B_{K_1/2}(p_1)$ with $K_1 = K/\|B^{-1}\|$. We now apply either Proposition~\ref{prop:single-max} or \ref{prop:double-max} depending on the condition, and conclude that on the curve $B_{K_1 \sqrt{\epsilon}/2}(p_1) \cap \Gamma^s$, each $c$ is of Aubry-Mather or bifurcation AM type, relative to $N_\epsilon$. We now apply Proposition~\ref{prop:forcing-connected}, to conclude that to conclude that  either \eqref{eq:mather-mech}, \eqref{eq:bifurcation-mech}, or \eqref{eq:arnold-mech} holds for each $c \in B_{K_1 \sqrt{\epsilon}/2}(p_1) \cap \Gamma^s$, on a $C^r$-residual subset $\cR_\sigma(N_\epsilon)$ of $N \in \cV_\sigma(N_\epsilon)$, 
for some $\sigma>0$. 

	We now revert the coordinate change. The fact that the coordinate change is $C^\infty$ implies the mapping $H_\epsilon \mapsto (\Phi_L \circ \Phi_\epsilon)^*H_\epsilon =: \Phi^* H_\epsilon$ is a homeomorphism between $C^r$ spaces, and in particular, open neighborhoods and residual subsets are preserved between coordinate changes. As a result, there is $\sigma' >0$ and a residual subset $\cR_{\sigma'}(H_\epsilon)$ of $\cV_{\sigma'}(H_\epsilon)$,  such that $H \in \cR_{\sigma'}(H_\epsilon)$ implies $\Phi^* H_\epsilon \in \cR_\sigma(N_\epsilon)$. Then Lemma~\ref{lem:symp-inv-normal-form} (invariance of diffusion mechanism under symplectic coordinate changes) implies for each $c \in B_{K_2\sqrt{\epsilon}/2}(p_0)\cap \Gamma$, and for each $H \in \cR_{\sigma'}(H_\epsilon)$, one of \eqref{eq:mather-mech}, \eqref{eq:bifurcation-mech}, or \eqref{eq:arnold-mech} hold. 

	We now apply the above argument to each $c \in \Gamma$, and establish \eqref{eq:mather-mech}, \eqref{eq:bifurcation-mech}, or \eqref{eq:arnold-mech} for an neighborhood $B_{\sigma(c)}(c)$ of $c$, on a $C^r$ residual subset $\cR_{c}$ of $\cV_{\sigma(c)}(H)$. By compactness, $\Gamma$ 
can be covered by finitely many $B_{\sigma(c_i)}(c_i)$'s, then 
by taking intersections over $\cR_{c_i}$, we conclude that our conditions hold on all $c \in \Gamma$, over a residual subset 
of $\cV_{\sigma_0}(H_\epsilon)$, where 
$\sigma_0 = \min \sigma(c_i)$. The theorem follows. 
\end{proof}

\section{Double resonance: geometric description}
\label{sec:DR-geometrical}

In this section we describe the non-degeneracy condition at the double resonance. We then describe the normally hyperbolic cylinders in this regime. In next section, we will return to variational setting, define the cohomology classes and prove their forcing equivalence. 

\subsection{The slow system}

We now consider one of the strong double resonance. Let $k \in \cK^\st(k_1,\Gamma, K)$, and denote $p_0 = \Gamma_{k_1, k}$, and $\omega_0 = \nabla H_0(p)$. Define 
\[
	\Lambda = \Span_\R \{k_1, k\} \cap \Z^3, 
\]
and choose $k_2 \in \Z^3_*$ such that $\Lambda = \Span_\Z\{k_1, k_2\}$. It is always possible to choose $|k_2| \le |k_1| + K$. 

Given $H_1 = \sum_{k \in \Z^3} h_k(p) e^{2\pi i k \cdot (\theta, t)}$, we define
\[
	[H_1]_\Lambda = \sum_{k \in \Lambda} h_k(p) e^{2\pi i k \cdot (\theta, t)}. 
\]
Then after a symplectic coordinate change defined on the set $K \sqrt{\epsilon}$ (see Theorem~\ref{double-norm-form}), the system has the normal form:
\[
	N_\epsilon = \Phi_\epsilon^* H_\epsilon = H_0 + \epsilon [H_1]_\Lambda + O(\epsilon^{\frac32}). 
\]

The system is conjugate to a two degrees of freedom mechanical system 
after a coordinate change and an energy reduction. The details are given 
in Appendix~\ref{slow-fast-section}, here we give a brief description. 
Let $k_3 \in \Z^3$ be such that 
\[
	B^T = \bmat{k_1 & k_2 & k_3} \in SL(3, \Z). 
\]
To define a symplectic coordinate change, we consider the corresponding autonomous system $N_\epsilon(\theta, p, t) + E$, and consider the coordinate change 
\begin{equation}
	\beal   \label{eq:Phi-linear}
	(\theta, p, t, E)\  &=& \Phi_L(\varphi, I, \tau, F), 
\qquad \qquad \qquad \\
	\bmat{\theta \\ t} = B^{-1} \bmat{\varphi \\ \tau/\sqrt{\epsilon}}, & \qquad & 
	\bmat{p - p_0 \\ E + H_0(p_0)} = B^T \bmat{\sqrt{\epsilon} I \\ \epsilon F}. 
	\enal 
\end{equation}
One checks that 
\begin{multline*}
	\left( 	\T^2 \times \R^2 \times \T \times \R, \quad  d\theta \wedge dp + dt \wedge dE \right)  \\
	\overset{\Phi_L}{\to}  	\left( \T^2 \times \R^2 \times \T \times \R, \quad \frac{1}{\sqrt{\epsilon}}(d\varphi \wedge dI + d\tau \wedge dF)  \right)
\end{multline*}
is an exact symplectic coordinate change. The transformed Hamiltonian 
$(N_\epsilon + E) \circ \Phi_L$ is no longer Tonelli in the standard sense, 
however by using a standard energy reduction on the energy level $0$, 
with $\tau$ as the new time takes the system to 
\[
	\frac{1}{\beta}\left( K(I) - U_0(\varphi) + \sqrt{\epsilon} P(\varphi, I, \tau) \right), \quad \varphi \in \T^2, \, I \in \R^2, \, \tau \in \sqrt{\epsilon}\T, 
\]
where 
\[
	\beta = k_3 \cdot (\omega_0, 1), \quad K(I) = \frac12 \left( B_0 \partial^2_{pp} H_0(p_0) B_0^T \right), \quad B_0 = \bmat{k_1^T \\ k_2^T}, 
\]
and 
\[
	U(k_1 \cdot( \theta, t), k_2\cdot (\theta, t)) = - [H_1]_\Lambda(\theta, p_0, t). 
\]
The system 
\be \label{eq:slow-mechanical}
	H^s(\varphi, I) = K(I) - U(\varphi) = K - U
\ee
is called the slow mechanical system, and the non-degeneracy conditions at 
the double resonance $p_0$ is stated for this system. 

\subsection{Non-degeneracy conditions for the slow system}
\label{sec:intro-DR}

We consider the (shifted) energy as a parameter. 
For each $E>0$, by the Maupertuis principle,  the Hamiltonian dynamics on 
the energy surface $\cS_E := \{(\varphi, I): \, H^s(\varphi, I) +\min U(\varphi)= E \}$, 
is the time change of the geodesic flow for the Jacobi metric
\[
	g_E(\varphi)(v) = \textb{2(E + U(\varphi))}\,K^{-1}(v), 
\]
where $K^{-1}(v) = \frac12 \left( \partial^2_{II}K \right)^{-1} v \cdot v$ is the Lagrangian associated to the Hamiltonian $K(I)$. 

We will be interested in a special homology class $h = (0, 1) \in \Z^2 \simeq H_1(\T^2, \Z)$. They represent classes of the original system satisfying $k_1 \cdot(\dot\theta, 1) = 0$, i.e. orbits that travel close to the resonance $\Gamma_{k_1}$. We assume the following non-degeneracy conditions: 
\begin{enumerate}
	\item[{$[DR1^h]$}] For each $E \in (0, \infty)$, each shortest closed geodesic (called a loop) of $g_E$ in the homology class $h$ is a hyperbolic orbit of the geodesic flow. 
	\item[{$[DR2^h]$}] At all but finitely many bifurcation values, there is only one $g_E$-shortest loop. At each bifurcation value $E$, there are exactly two shortest $g_E$ loops denoted 
$\gamma_h^E$ and $\bar\gamma_h^E$. 
	\item[{$[DR3^h]$}] At bifurcation value $E_*$,
	\[
		\dfrac{d(\ell_E(\gm_h^E))}{dE}|_{E=E^*} \ne
		\dfrac{d(\ell_E(\overline \gm_h^E))}{dE} |_{E=E^*},
	\]
	where $l_E$ denote the $g_E$ length of a loop. 
\end{enumerate} 

We now discuss the conditions at the critical shifted energy $E=0$. 
The Jacobi metric $g_0$ becomes degenerate at one point 
$\varphi_* = \argmin_\varphi U(\varphi)$. By performing a translation, 
we may assume $\varphi_* = 0$. Let $\gamma_h^0$ be a shortest loop 
of $g_0$ in the homology $h$. Consider the following cases:
\begin{enumerate}
	\item $0 \in \gm^0_h$ and $\gm^0_h$ is not self-intersecting. 
	Call such homology class $h$ {\it simple critical} and
	the corresponding geodesic $\gm^0_h$ {\it simple loop}.

	\item $0 \in \gm^0_h$ and $\gm^0_h$ is self-intersecting.
	Call such homology class $h$ {\it non-simple} and
	the corresponding geodesic $\gm_h^0$ {\it non-simple}.

	\item  $0 \not\in \gm_h^0$, then  $\gm_h^0$ is a regular geodesic.
	Call such homology class $h$ {\it simple non-critical}.
\end{enumerate}
Mather \cite{Ma8} proved that  generically only these three cases occur (see below for the precise claim). 

\begin{lem}
	\label{lem:non-simple-loop-decomposition}
	Let $h$ be a non-simple homology class. Then for a generic potential $U$ 
	the curve $\gm^0_h$  is the concatenation of of two simple loops, possibly 
	with multiplicities. More precisely, given $h\in H_1(\T^2,\Z)$ generically 
	there are simple homology classes $h_1,\ h_2 \in H_1(\T^s,\Z)$ and integers 
	$n_1,\ n_2\in\Z_+$ such that the corresponding  minimal geodesics 
	$\gm^0_{h_1}$ and $\gm^0_{h_2}$ are simple and $h=n_1h_1+n_2h_2$.
\end{lem}

We call $\gamma_h^0$ extensible 
 if there exists a family of shortest curves $\gamma_h^E$ converging to it in the Hausdorff topology. 
Consider the lift of $\gamma_h^E$ to the universal cover $\R^2$, then as $E \to 0$ it converges to 
a periodic curve in $\R^2$ consists of concatenation of $\gamma_{h_1}^0$ and $\gamma_{h_2}^0$. 
Let $(\sigma_1, \dots, \sigma_n) \in \{0, 1\}^N$ be the order that $\gamma_{h_i}^0$ are traced. 

\begin{lem}[See Section~\ref{sec:homology}]
	\label{lem:homotopy}
Assume that $\gamma_h^0$ is extensible, then the sequence 
$(\sigma_n, \dots, \sigma_n)$, as described, 
is uniquely determined up to cyclic permutation. We write
$$
\gamma_h^E \to  \gamma_{h_{\sigma_1}}^0 * \cdots
* \gamma_{h_{\sigma_n}}^0, \quad \text{as} \quad E \to 0. 
$$
\end{lem}

We note that $(0,0)$ is a fixed point of the Hamiltonian flow $H^s(\varphi, I)$ which is hyperbolic if $\partial^2_{\theta\theta}U(0) > 0$. Any simple $g_0$--shortest loop corresponds to a homoclinic orbit of 
the fixed point $(0,0)$. We impose the following non-degeneracy conditions:
\begin{itemize}
	\item[{[DR1$^c$]}] $(0,0)$ is a hyperbolic fixed point with distinct eigenvalues 
	$-\lambda_2 < - \lambda_1 <0 < \lambda_1 < \lambda_2$. Let $v_1^\pm, v_2^\pm$ be 
	the eigendirections for $\pm \lambda_1$, $\pm \lambda_2$. 
	\item[{[DR2$^c$]}] There is a unique $g_0$--shortest loop in the homology $h$. If it is non-simple, 
	then it is the concatenation of two simple loops $\gamma_{h_1}^0$ and $\gamma_{h_2}^0$. 
	\item[{[DR3$^c$]}] If $\gamma_h^0$ is simple critical, then it is \emph{not} tangent to 
	the $\langle v_2^+, v_2^-\rangle$ plane. If $\gamma_h^0$ is non-simple, then each of 
	$\gamma_{h_1}^0$ and $\gamma_{h_2}^0$ are not tangent to the $\langle v_2^+, v_2^-\rangle$ plane. 
	\item[{[DR4$^c$]}] 
	\begin{itemize}
		\item If $\gamma_h^0$ is simple non-critical, then $\gamma_h^0$ is hyperbolic. 
		\item If $\gamma_h^0$ is simple critical, then $\gamma_h^0$ is non-degenerate 
		in the sense that it is the transversal intersection of the stable and unstable manifolds of $(0,0)$. 
		\item If $\gamma_h^0$ is non-simple, then each of $\gamma_{h_1}^0$ and 
		$\gamma_{h_2}^0$ is non-degenerate. 
	\end{itemize}
\end{itemize}

The genericity of these conditions are summarized in the following statement. 
\begin{prop}
	\label{prop:Dr-non-deg}
The conditions $[DR1^h-DR3^h]$ and $[DR1^c-DR4^c]$ 
hold on an open and sense set of potentials $U \in C^r(\T^2)$,
for $r\ge 2$. 
\end{prop}

\subsection{Normally hyperbolic cylinders}

Conditions $[DR1^h] - [DR3^h]$ ensures that for each $E_0 > 0$, the set 
$$
\bigcup_{E \in (E_0- \delta, E_0 + \delta)} \eta_h^E
$$ 
is a normally hyperbolic invariant cylinder. This cylinder does not 
necessarily extend to the shifted energy $E=0$. The following 
statement ensures existence of cylinders near critical energy, 
using the conditions $[DR1^c] - [DR4^c]$. 

Recall that $\gamma_h^0$ is a shortest curve in the critical energy. 
The corresponding set in the phase space is called $\eta_h^0$. Due to 
the symmetry of the system, we also have the shortest curve $\gamma_{-h}^0$ 
which coincide with $\gamma_{h}^0$ but has a different orientation. 

\begin{thm}[See Section~\ref{sec:NHIC-DR-proof}]
	\label{thm:NHIC-DR}
Suppose that $H^s$ satisfies conditions 
$[DR1^c] - [DR4^c]$\,\footnote{Note that we do not assume 
$[DR1^h] - [DR3^h]$}. 
\begin{enumerate}
	\item If $\gamma_h^0$ is simple, then there is $e>0$ depending on $H^s$ such that:
	\begin{enumerate}
		\item  For each $0 < E < e$, there exist periodic orbits $\eta_h^E$ and $\eta_{-h}^E$, such that the projections $\gamma_h^E \to \gamma_h^0$ and $\gamma_{-h}^E \to \gamma_{-h}^0$ in the Hausdorff topology. 
		\item For each $- e < E < 0$, there exists a periodic orbit $\eta_c^E$ which shadows the concatenation of $\eta_h^0$ and $\eta_{-h}^0$. 
	\end{enumerate}  
	{Then the union }
	\[
		\bigcup_{0 < E < e} \left( \eta_h^E \cup  \eta_{-h}^E \right) \cup \eta_h^0 \cup \eta_{-h}^0 \cup \bigcup_{-e < E < 0} \eta_c^E
	\]
	is a $C^1$ normally hyperbolic invariant manifold containing the homoclinics $\eta_{\pm h}^0$. 
	\item If $\gamma_h^0$ is non-simple: Let $\sigma_1, \dots, \sigma_n$ be the sequence determined in Lemma~\ref{lem:homotopy}, $[DR1^c] - [DR4^c]$ ensures $\gamma_h^0$ is extensible. More precisely, there is $e>0$ such that for each $0 < E < e$, there is a periodic orbit $\gamma_h^E$ such that $\gamma_h^E \to \gamma_{h_{\sigma_1}}^0 * \cdots
	* \gamma_{h_{\sigma_n}}^0$. Moreover, each $\gamma_h^E$ is hyperbolic. 
\end{enumerate}
\end{thm}

\begin{figure}[t]
	\centering 
	\includegraphics[height=2.65in]{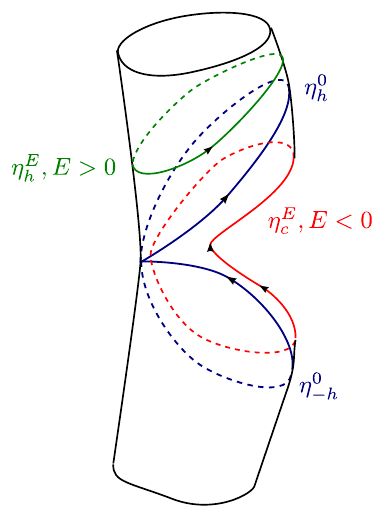} 
	\hskip .15in
	\includegraphics[height=2.65in]{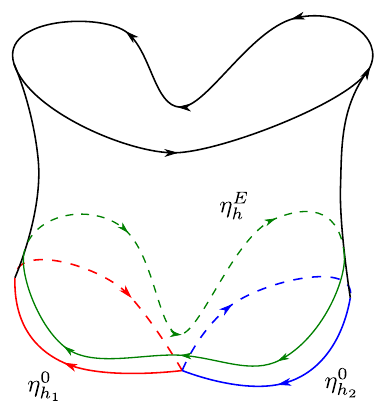} 
	\caption{Extension of homoclinics to periodic orbits, simple and non-simple.}
\end{figure}


Theorem~\ref{thm:NHIC-DR} allows us to prove Proposition~\ref{prop:Dr-non-deg} together with the following statements. 

\begin{thm}[See Appendix~\ref{generic-mechanical-system}]
	\label{thm:DR-non-deg-HE}
	Fix two parameters $0 < e_0 < \bar{E}$, then the set of potentials $U$ such 
	that $[DR1^h] - [DR3^h]$ hold on the smaller interval 
	$E \in [e_0, E]$ is open and dense in $C^r(\T^2)$, 
	for $r\ge 2$. 
\end{thm}

The following statement is proved in Section~\ref{sec:very-high-energy}. 
\begin{prop}
\label{prop:DR-very-high-energy}
The set of potentials $U\in C^r(\T^2)$ such that all $\gamma_h^E$  for $E \ge \bar{E}$ is unique and hyperbolic is open and dense for $r\ge 2$. 
\end{prop}

The following statement is proved in Section~\ref{sec:non-deg-crit}.
\begin{prop}
	\label{prop:DR-non-deg-Crit}
The set of potentials $U$ such that $[DR1^c] - [DR4^c]$ hold is open and dense in $C^r(\T^2)$, for $r\ge 2$.
\end{prop}

\begin{proof}
	[Proof of Proposition~\ref{prop:Dr-non-deg}]
	Proposition~\ref{prop:DR-non-deg-Crit} implies the set of potentials which satisfy $[DR1^c]-[DR4^c]$ is open and dense. By Theorem~\ref{thm:NHIC-DR} the set of potentials $\cU_{crit}$ such that there is $e_0 > 0$ such that all $\gamma_h^E$ are unique and hyperbolic for all $0 < E < e_0$ is open and dense.

	The set $\cU_{High}$ of $U$'s such that that there exists $\bar{E} > 0$ 
such that all $\gamma_h^E$ are unique and hyperbolic is open and dense 
by Proposition~\ref{prop:DR-very-high-energy}. 
By Theorem~\ref{thm:DR-non-deg-HE} the set of potentials
$\cU_{med}^{e_0, \bar{E}}$ such that for given $0<e_0 < \bar{E}$, 
$[DR1^h] - [DR3^h]$ holds on $E \in [e_0, E]$ is also open and dense. 
As a result the set of potentials, where $[DR1^h] - [DR3^h]$ holds on 
$E \in (0, \infty)$ is 
	\[
		\cU_{crit} \cap 
\bigcup_{0 < e_0 < \bar{E}} \cU_{med}^{e_0, \bar{E}} \,\cap\, \cU_{high} 
	\]
	which is open and dense. 
\end{proof}

\subsubsection*{Diffusion across a double resonance: a geometric description}

\begin{figure}[t]
	\centering 
	\includegraphics[width=2.5in]{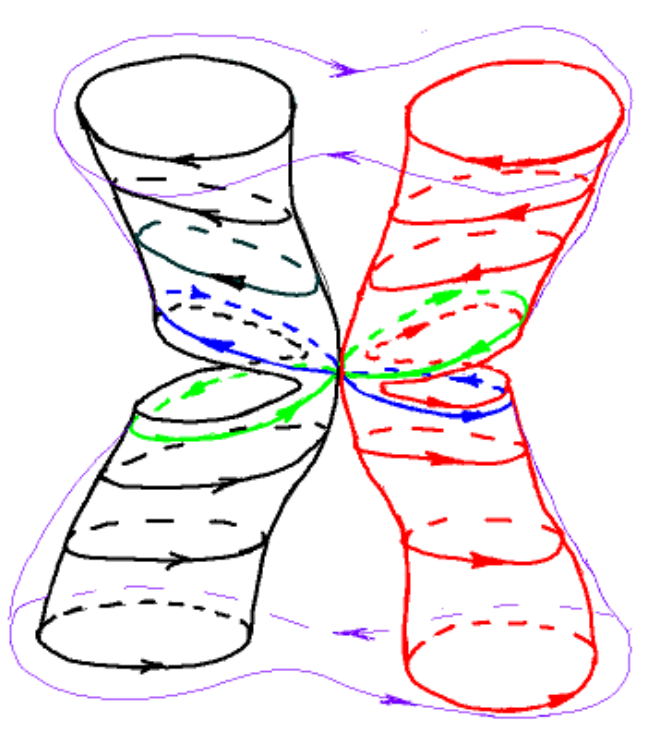}
	\caption{Hyperbolic invariant manifolds with kissing property} 
	\label{fig:kissing}
\end{figure}

The diffusion across a double resonance may be described heuristically as follows:
\begin{itemize}
	\item If $h$  is a simple homology, then the cylinder extends to 
the shifted energy $E < 0$ and connecting with the homology $-h$. 
As a result, the family of periodic orbits $\gamma_h^E$, $E \ge 0$ and $\gamma_{-h}^E$, $E  \ge 0$ are all contained in a normally hyperbolic invariant manifold. This corresponds to a continuous curve of cohomologies that are of Aubry-Mather type. Moreover, this picture survives small perturbations of $H^s$. 
	\item If $h$ is non-simple, then the cylinder ``pinches'' 
at $E = 0$. In particular, after considering the perturbation 
$H^s + \sqrt{\epsilon} P$ of the slow system, the cylinder may 
not survive the perturbation for $E$ sufficiently close to $0$. However, for each simple homology $h_1, h_2$, there exists 
a simple cylinder due to Theorem~\ref{thm:NHIC-DR}, item 1. 
The two simple cylinders are tangent to 
the weak stable/unstable directions plane 
at the fixed point $(0, 0)$. See Figure~\ref{fig:kissing}. 

	To diffuse across a double resonance, we ``jump'' from 
the cylinder for homology $h$ to the cylinder with homology $h_1$, then diffuse across to homology $-h_1$ since $h_1$ 
is now simple, then jump back to homology $-h$. All of these 
are realized by choosing the appropriate cohomology curves 
the lie on these cylinders. This construction is detailed in 
Section~\ref{sec:coh-DR}.  See also Figure \ref{fig:sim-non-simp-coh}.
\end{itemize}

\subsection{Local maps and global maps}
\label{sec:intro-local-map}

In this section we outline the basic approach to proving Theorem~\ref{thm:NHIC-DR}, based on ideas of Shil'nikov and others (\cite{BR}, \cite{Shil67}, \cite{ST89}). The full proofs are given in Section~\ref{sec:NHIC-DR-proof}. 

Let us describe the simple homology case first. Let $\eta^+ = \eta_h^0$ be the  homoclinic orbit to the hyperbolic fixed point $O = (0, 0)$, and $\eta^- = \eta_{-h}^0$ its time-reversal. Condition $[DR2^c]$ ensures that $\eta^\pm$ are not tangent to the strong stable/unstable directions, which implies they must be tangent to the weak stable/unstable directions.

\begin{figure}[t]
	\centering 
	\includegraphics[width=3.5in]{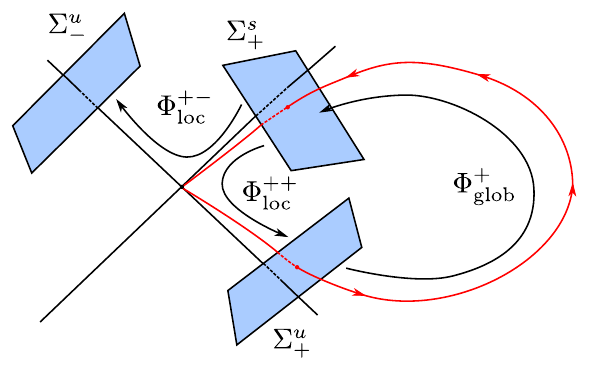} 
	\caption{A local map and a global map}
	\label{fig:loc-glob}
\end{figure}

Consider four (three dimensional) sections 
\[
	\Sigma^u_\pm, \quad \Sigma^s_\pm
\]
transverse to the weak stable and unstable eigen-directions,  
sufficiently close to the origin $O$, on each side of the equilibrium (see \eqref{eq:Sigma-su}). Up to renaming, we may 
assume that $\Sigma^{s/u}_+$ is transverse to $\eta^+$, and $\Sigma^{s/u}_-$ is transverse to $\eta^-$. 
We then define four local maps:
\[
	\Phl^{ij}: U_i\,(\subset \Sigma^s_i) \to \Sigma^u_j, \quad  i,j \in \{+, -\}
\]
as the Poincar\'e map between the corresponding sections. Note that 
these maps are {\it not defined on the whole section,} in particular, 
they are undefined along the (full) stable/unstable manifolds $W^{s/u}(O)$. 
However, they are defined on open sets. Moreover, we have a pair of 
global maps
\[
	\Phg^+: \Sigma^u_+ \to \Sigma^s_+, \quad \Phg^-: \Sigma^u_- \to \Sigma^s_-
\]
which are the Poincar\'e maps along the orbits $\eta^+$ and $\eta^-$. These maps are well defined from a neighborhood of $\Sigma^u_\pm \cap \eta^\pm$ to 
a neighborhood of $\Sigma^s_\pm \cap \eta^\pm$. See Figure~\ref{fig:loc-glob}. 

The periodic orbits obtained in Theorem~\ref{thm:NHIC-DR} corresponds to 
the fixed points of compositions of local and global maps, when restricted to 
the suitable energy surfaces. More precisely:
\begin{itemize}
	\item The orbits $\eta_{h}^E$, $E>0$ correspond to the fixed point of 
	$\Phg^+ \circ \Phl^{++}|_{S_E}$, where $S_E$ denote the energy surface 
$\{H^s +\min U(\varphi) = E\}$, 
	and similarly for $\eta_{-h}^E$. 
	\item The orbits $\eta_c^E$, $E < 0$ correspond to the fixed point of 
	$\Phg^+ \circ \Phl^{-+}\circ \Phg^-\circ \Phl^{+-}|_{S_E}$. 
\end{itemize}

In the non-simple case, we similarly consider two homoclinics $\eta_1 = \eta_{h_1}$ and $\eta_2 = \eta_{h_2}$, and assume that both crosses the sections $\Sigma^{s/u}_+$. Then we have the same local maps, and different global maps $\Phg^{1/2}$. The periodic orbit $\eta^E_h$ corresponds to the fixed point of 
\[
	\prod_{i=n}^1 \left(\Phg^{\sigma_i} \circ  \Phl^{++} \right)|_{S_E},
\]
where $\sigma_i$ is the sequence in Lemma~\ref{lem:homotopy}.

\section{Double resonance: choice of cohomology and forcing equivalence}
\label{sec:coh-DR}

\subsection{Choice of cohomologies for the slow system}

As in the case of single resonance, our strategy is to choose a continuous curve in the cohomology space, and prove forcing equivalence up to a residual perturbation. To do this, we need to use the duality between homology and cohomology. 

Let $L$ be the Lagrangian associated to the Hamiltonian $H$, and let $\mu$ denote an invariant measure of 
the Euler--Lagrange flow. The rotation vector of $\mu$ is 
given by $\rho(\mu) = \int v \,d\mu(\theta, t, t)$. Then Mather's 
alpha and beta functions are defined as:
\[
	\alpha_H(c) = - \inf_\mu \int \left( L(\theta, v, t) - c \cdot v \right)d\mu, \quad
	\beta_H(\rho) = \inf_{\rho(\mu) = \rho} \int L(\theta, v, t) d\mu. 
\]
A measure reaching the minimum in the definition of $\alpha_H(c)$ is called $c-$minimal. 
Then $\alpha$ and $\beta$ are both convex and Fenchel dual of each other:
\[
	\beta_H(\rho) = \sup\{\rho \cdot c - \alpha_H(c)\}. 
\]
The Legendre-Fenchel transform of $\beta$ is defined as 
\[
	\LF_{\beta_H}(\rho) = \{ c \in \R^2: \quad \beta_H(\rho) = \rho \cdot c - \alpha_H(c)\}. 
\]
Geometrically, 
\[
	\LF_{\beta_H}(\rho) = conv \left\{c:\quad  \text{there is a c-minimal } \mu \text{ such that  }
	\rho(\mu) = \rho  \right\}
\]
where $conv$ denotes the convex hull.

Let $\gamma_h^E$ be a shortest loop for the Jacobi metric $g_E$. Let $T(\gamma_h^E)$ denotes its period under the Hamiltonian flow, and if $\gamma_h^E$ is unique, we define
$$ \lambda_h^E =1/T(\gamma_h^E). $$
If $H^s$ satisfies $[DR1^h] - [DR3^h]$, then there are at most finitely many $E$'s such that there are two shortest loops
$\gamma_h^E$ and $\bar{\gamma}_h^E$. 
We will show that the set 
$\LF_{\beta_{H^s}}(\lambda_h^E h) = \LF_{\beta_{H^s}}(\bar{\lambda}_h^E h)$, and, therefore, 
the set $\LF_{\beta_H^s}(\lambda_h^E h)$ is independent of the choice of $\gamma_h^E$.

Each $\LF_{\beta_H^s}(\lambda_h^E h)$ is a segment of nonzero length parallel to $h^\perp$, 
and depends continuously on $E$. We call the union
\begin{equation}
	\label{channel}
	\bigcup_{E>0}\LF_{\beta_{H^s}}(\lambda_h^E h)
\end{equation}
the \emph{channel} associated to the homology $h$, and we will choose
a curve of cohomologies in the \emph{interior} of this channel. The channel is connected at 
the bottom to the set $\LF_{\beta_H}(0)$, which is a convex set with non-empty interior. 
The following proposition summarizes the channel picture and the relation to the Aubry sets.

\begin{prop}[See Section~\ref{sec:var-slow}]
	\label{slow-localization}
Assume that $H^s$ satisfies the conditions $[DR1^h-DR3^h]$ and 
$[DR1^c- DR4^c]$. Then each $\LF_{\beta_{H^s}}(\lambda_h^E h)$ is a segment of non-zero length orthogonal to $h$, which varies continuously with respect to $E$. 

 For $\bar{E} > 0$,  let
$\bar{c}_h: (0, \bar{E}]\to H^1(\T^s, \R)$ be 
a $C^1$ function such that $\bar{c}_h(E)$ is 
in the relative interior of $\LF_{\beta_H}(\lambda_h^E h)$.
The following hold. 
\begin{enumerate}
	\item If $E$ is not a bifurcation energy, then 	$ \cA_{H^s}(\bar{c}_h(E)) = \gamma_h^E $. 
	\item If $E$ is a bifurcation energy, 
	$ \cA_{H^s}(\bar{c}_h(E)) = \gamma_h^E \cup \bar\gamma_h^E $. 
	\item If $h$ is simple, then the  limit $\lim_{E \to 0}\LF_{\beta_{H^s}}(\lambda_h^E h)$ contains 
	a segment of non-zero width. We assume, in addition, that $\bar{c}_h(0)$ is in the relative interior of 
	this segment. 
	\begin{enumerate}
		\item If $h$  is simple critical, then $ \cA_{H^s}(\bar{c}_h(0)) = \gamma_h^0 $;
		for each $0 \le \lambda <1$,
		$ \cA_{H^s}(\lambda\bar{c}_h(0)) = \{\varphi =0\} $. 
		\item If $h$ is simple non-critical, then 
		$  \cA_{H^s}(\bar{c}_h(0)) = \gamma_h^0 \cup \{\varphi =0\} $;
		for each $0 \le \lambda <1$,
		$ \cA_{H^s}(\lambda\bar{c}_h(0)) = \{\varphi =0\}$. 
	\end{enumerate}
	\item If $h$ is non-simple, then the limit $\lim_{E \to 0}\LF_{\beta_{H^s}}(\lambda_h^E h)$ is a single point. 
\end{enumerate}
\end{prop}

Let us also note that due to symmetry of the system, $\LF_{\beta_H}(- \lambda h) = - \LF_{\beta_H}(\lambda h)$. Denote $\bar{c}_{-h}(E) = - \bar{c}_h(E)$. We now choose the cohomology classes for $H^s$ as follows: 
if $h$ is simple (either critical or non-critical),  we choose 
\begin{equation}
	\label{eq:bargm-simple}
\begin{aligned}
 & 	\bar\Gamma_h = \bar\Gamma_h(\bar{E}) =\left(  \bigcup_{0 \le E \le \bar{E}} \bar{c}_h(E) \right) \cup 
	\left( \bigcup_{0 \le  \lambda \le 1} \lambda \bar{c}_h(0) \right), \\
	&	\bar\Gamma^{DR}_h = \bar\Gamma^{DR}_h(\bar{E}) = \bar\Gamma_h \cup \bar\Gamma_{-h} = \bar\Gamma_h \cup \left( - \bar\Gamma_h \right). 
\end{aligned}
\end{equation}
The curve $\bar\Gamma_h^{DR}$ is a continuous curve connecting $\bar{c}_h(\bar{E})$, $0,$ and $\bar{c}_{-h}(\bar{E})$. 

\begin{figure}[t]
	\centering 
	\includegraphics[width=1.95in]{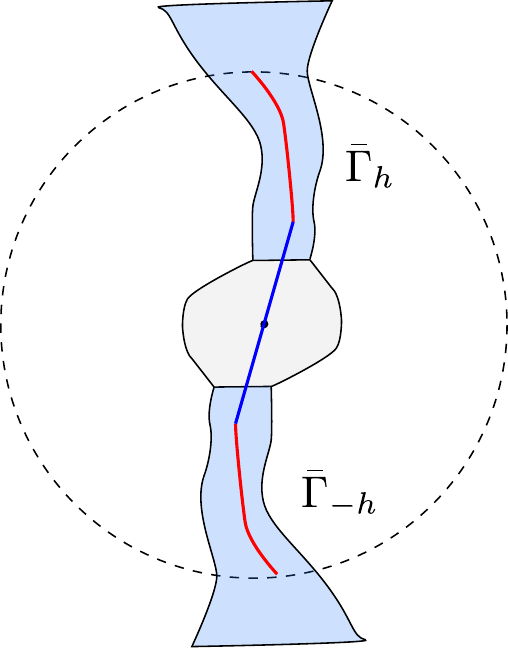} 
	\hspace{.5in}
	\includegraphics[width=1.95in]{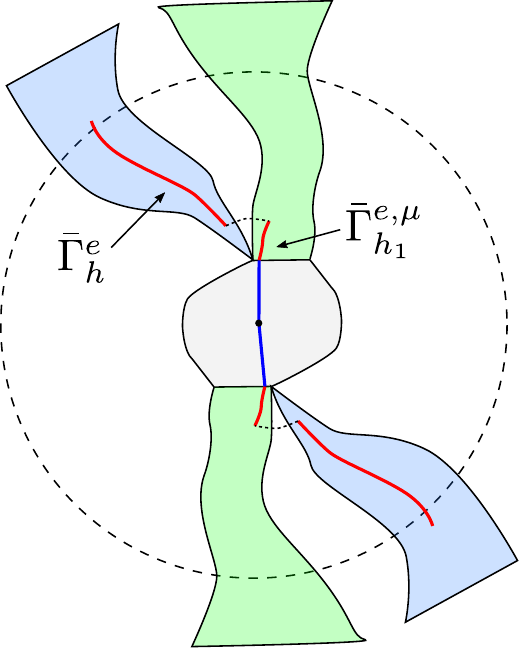}
	\caption{Choice of cohomology for $H^s$. Left: simple case; right: non-simple case}
	\label{fig:sim-non-simp-coh}
\end{figure}

If $h$ is non-simple, then $h = n_1 h_1 + n_2 h_2$ is a combination of simple homologies $h_1$ and $h_2$. Let $0 < \mu < e$ be parameters. Define 
\begin{equation}
	\label{eq:bargm-non-simp}
	\bar\Gamma_h^e = \bigcup_{e  \le E \le \bar{E}} \bar{c}_h(E). 
\end{equation}
Since $h_1$ is simple critical, we choose a continuous curve $\bar{c}_{h_1}^\mu(E)$ such that 
\[
	\|\bar{c}_{h_1}^\mu(0) - \bar{c}_h(0)\| < \mu.
\]
We then define 
\begin{equation}
	\label{eq:bargm-aux}
	\bar\Gamma_{h_1}^{e, \mu} = \left(  \bigcup_{0 \le E \le e + \mu} \bar{c}^\mu_{h_1}(E) \right) \cup 
	\left( \bigcup_{0 \le  \lambda \le 1} \lambda \bar{c}_{h_1}(0) \right). 
\end{equation}
$\bar\Gamma_{h_1}^{e, \mu}$ is a continuous curve connecting 
$\bar{c}_{h_1}^\mu(e + \mu)$ with $0$ (see Fig. 5.1, right). We then define 
\begin{equation}
	\label{eq:bargm-full}
	\bar\Gamma_h = \bar\Gamma_h^e \cup \bar\Gamma_{h_1}^{e, \mu}, \quad
	\bar\Gamma_h^{DR} = \bar\Gamma_h \cup \bar\Gamma_{-h}. 
\end{equation}
Let us note that the set $\bar\Gamma_h^{DR}$ consists of 
three connected component, $\bar\Gamma_h^e$, 
$\bar\Gamma_{-h}^e$ and 
$\bar\Gamma_{h_1}^{e, \mu} \cup 
\bar\Gamma_{-h_1}^{e, \mu}$. Since the mechanisms 
we described so far can only prove forcing equivalence 
along a connected set, we will use a different mechanism 
called ``jump'' to prove forcing-equivalence of different connected components.

\subsection{Aubry-Mather type at a double resonance}

\begin{prop}\label{prop:AM-DR}
Suppose $H^s = K(I) - U(\varphi)$ satisfies the non-degeneracy assumptions. Given $C>0$, there is $\epsilon_0, \delta > 0$ depending on $H^s$ and $C$ such that for each $0 < \epsilon < \epsilon_0$, if $U' \in \cV_\delta(U)$ and $\|P\|_{C^2} < C$, then for 
\[
	H^s_\epsilon = \frac{1}{\beta}\left(K(I) - U'(\varphi) + \sqrt{\epsilon} P(\varphi, I, \tau) \right). 
\]
each $\bar{c} \in \bar\Gamma_h^{DR}$ is of one of fours types: Aubry-Mather type, bifurcation Aubry-Mather type,  asymmetric bifurcation type, or $\tcA_{H^s_\epsilon}$ is a hyperbolic periodic orbit. Note that this applies to all types of homologies: simple critical, simple non-critical, and non-simple. 
\end{prop}
\begin{proof} 
We prove our proposition by referring to technical statement proved in later sections. 

(1) \emph{Simple, critical homology}. Denote $c_0 = \bar{c}_h(0)$.  Theorem~\ref{thm:dr-am-crit} states that there exists $\epsilon_0, \epsilon, \delta >0$ such that for all $0 < \epsilon < \epsilon_0$, $U' \in \cV_\delta(U)$, and $c \in B_e(c_0)$, the pair $H_\epsilon^s, c$ is of Aubry-Mather type. The same theorem also states $H_\epsilon^s, \lambda c_0$ for $0 \le \lambda \le 1$ is of Aubry-Mather type. This covers the cohomologies:
\[
c \in \bigcup_{0 \le E \le e} \bar{c}_h(E) \, \cup \, \bigcup_{0 \le \lambda \le 1} \lambda \bar{c}_h(0). 
\]
The cohomologies 
$\bigcup_{e \le E \le \bar{E}}\bar{c}_h(E)$, is covered in Theorem~\ref{thm:dr-high-AM}, which states that for each $e>0$, for there is $\epsilon, \delta$ as in our proposition, such that with respect to $H^s_\epsilon$, the cohomology $\bar{c}_h(E)$, $E \ge e$ is of Aubry-Mather type if $\gamma_h^E$ is the unique shortest curve, and of bifurcation Aubry-Mather type if there are two shortest curves. As a result, all cohomologies in $\bar{\Gamma}_h$ are of AM or bifurcation AM type, and by symmetry, so does $\bar{\Gamma}_{-h}$. This proves our proposition in the simple homology case, see \eqref{eq:bargm-simple}. 

(2) \emph{Non-simple homology}. In the non-simple case. the cohomology curve $\bar{\Gamma}_h$ is the disjoint union of two parts, namely 
\[
\bar{\Gamma}_h^e = \bigcup_{e \le E \le \bar{E}}\bar{c}_h(E), \quad \bar{\Gamma}_{h_1}^{e, \mu} = \left(  \bigcup_{0 \le E \le e + \mu} \bar{c}^\mu_{h_1}(E) \right) \cup 
	\left( \bigcup_{0 \le  \lambda \le 1} \lambda \bar{c}_{h_1}(0) \right).
\]
We note that $h_1$ is a simple homology and  therefore each $c \in \bar{\Gamma}_{h_1}^{e, \mu}$ are of AM type in the same way as in case (1). On the other hand, each homology in $\bar{\Gamma}_h^e$ is of AM or bifurcation AM type since Theorem~\ref{thm:dr-high-AM} applies the same way to simple and non-simple homology. We conclude that our Proposition holds in the non-simple case. 

(3) \emph{Simple, non-critical homology}. The high energy case follows from 
Theorem~\ref{thm:dr-high-AM} in the same way as before. The critical energy 
follows from Theorem~\ref{thm:dr-asym-bif}, the main difference is that the 
critical energy is an asymmetric bifurcation (see Definition~\ref{def:asym-bif}). 
For $c = \lambda \bar{c}_h(0)$ where $0 \le \lambda < 1$, the Aubry set is 
a single periodic orbit as a perturbation of the hyperbolic fixed point $(0, 0)$. 
\end{proof}

\begin{cor}\label{cor:DR-equiv}
Suppose $H_\epsilon^s$ is from Proposition~\ref{prop:AM-DR}, then there is $\sigma > 0$ and a residual subset $\cR$ of $\cV_\sigma(H_\epsilon^s)$, such that if $H' \in \cR$, then each $\bar{c} \in \bar\Gamma_h^{DR}$ satisfies one of the diffusion mechanisms 
\eqref{eq:mather-mech}, 
\eqref{eq:bifurcation-mech} and \eqref{eq:arnold-mech}. 
\end{cor}
\begin{proof}
Note that Proposition~\ref{prop:forcing-connected} applies when $c$ is either Aubry-Mather type, bifurcation Aubry-Mather type,  asymmetric bifurcation type. Then Using Proposition~\ref{prop:AM-DR}, we only need to show the same conclusion hold if we add a fourth case, when the Aubry set 
is a single hyperbolic periodic orbit. However, in this case, 
$\tcN^0 = \tcA^0$ is discrete, so \eqref{eq:mather-mech} applies. 
\end{proof}

We now revert to the original coordinate system. Denote $p_0 = \Gamma_{k_1, k_2}$ the double resonance point. For a given cohomology class $\bar{c} \in \R^2$, we consider the pair $c, \alpha$ defined by the equality
\begin{equation}
	\label{eq:c-barc}
	\bmat{c - p_0 \\ - \alpha + H_0(p_0)} = M^T \bmat{ \bar{c}\sqrt{\epsilon} \\  \alpha_{H_\epsilon^s}(\bar{c}) \, \epsilon}
\end{equation}
then $\alpha = \alpha_{N_\epsilon}(c)$ (compare to \eqref{eq:Phi-linear}). If \eqref{eq:c-barc} is satisfied, we denote
\[
(c, -\alpha) = \Phi_L^*(\bar{c}, \alpha_{H_\epsilon^s}(\bar{c})), \quad 
c = \Phi_{L, H_\epsilon^s}^*(\bar{c}). 
\]

Moreover, let us consider the autonomous version of the Aubry set:  
\[
	\tcA_{H^s_\epsilon + F}(\bar{c}) = 
	\{ (\varphi, I, \tau, - H_s^\epsilon(\varphi, I, \tau)):\quad (\varphi, I, \tau) \in \tcA_{H_\epsilon^s}(\bar{c}) \},
\]
and similarly define $\tcM$ and $\tcN$. Then 
by Proposition~\ref{prop:He-He-slow-var}
\begin{equation}
	\label{eq:symp-inv-linear}
	\tcA_{H_\epsilon + E}(c) = \Phi_L \left( \tcA_{H^s_\epsilon + F}(\bar{c}) \right), \quad
	\tcN_{H_\epsilon + E}(c) = \Phi_L \left( \tcN_{H^s_\epsilon + F}(\bar{c}) \right). 
\end{equation}
\begin{prop}\label{prop:DR-symp-inv}
Suppose $c, \bar{c}$ are related by \eqref{eq:c-barc}. Then with respect to $N_\epsilon$,  
$c$ satisfies  one of the diffusion mechanisms \eqref{eq:mather-mech}, \eqref{eq:bifurcation-mech} 
and \eqref{eq:arnold-mech} if and only if $\bar{c}$ does the same with respect to $H_\epsilon^s$. 
\end{prop}
\begin{proof}
	The main issue here is that \eqref{eq:mather-mech}, \eqref{eq:bifurcation-mech} and \eqref{eq:arnold-mech} are defined for the zero section $\tcA^0$ and $\tcN^0$, but the symplectic coordinate change $\Phi_L$ does not preserve the zero section. However, the property that there \emph{exists} a global section of the Hamiltonian flow such that one of \eqref{eq:mather-mech}, \eqref{eq:bifurcation-mech} or \eqref{eq:arnold-mech} applies to the intersection of $\tcA$ or $\tcN$  with the section, is invariant under the coordinate change $\Phi_L$, due to \eqref{eq:symp-inv-linear}. Since any two global sections are related by the Poincare map, any topological properties of invariant sets are equivalent across different sections. 
\end{proof}

\subsection{Connecting to $\Gamma_{k_1, k_2}$ and $\Gamma^{SR}_{k_1}$}

At this point it is natural to consider the cohomology class 
\begin{equation}
	\label{eq:Phi-L-star}
	\Phi_{L, H_\epsilon^s}^* \left( \bar\Gamma_h^{DR} \right) := \pi_c \left( \{ \Phi_L^*(\bar{c}, \alpha_{H^s_\epsilon}(\bar{c})):\quad  \bar{c} \in \Gamma_h^{DR}\} \right)
\end{equation}
(see \eqref{eq:c-barc})
where $\pi_c$ denote the projection $(c, -\alpha) \mapsto c \in \R^2$. We note that $\alpha$ is automatically determined by $c$ via the relation $\alpha = \alpha_{N_\epsilon}(c) = \alpha_{H_\epsilon}(c)$. 

\begin{figure}[t]
	\centering 
	\includegraphics[width=3in]{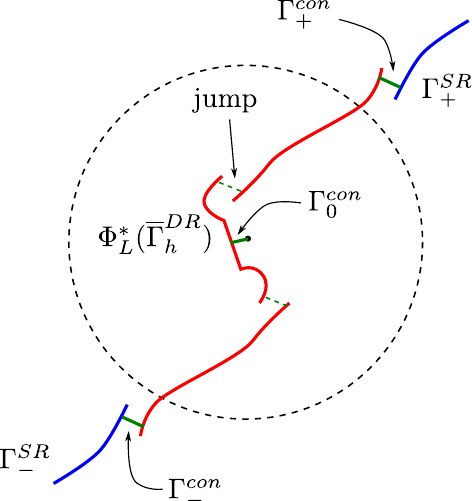} 
	\caption{Cohomology curve at double resonance, with connectors}
	\label{fig:gamma-DR}
\end{figure}

We would like to choose the cohomology $\Phi_{L, H_\epsilon^s}^* \left( \bar\Gamma_h^{DR} \right)$ for the original system, but due to the $\epsilon$ dependence of the map, the new set does not necessarily contain the double resonance point $\Gamma_{k_1, k_2}$, nor does it connect to the cohomology $\Gamma^{SR}_{k_1}$ already chosen at single resonance. To solve this problem, we will add three pieces of ``connectors'' to the set: $\Gamma^{\con}_0$ is used to connect to $\Gamma_{k_1, k_2}$, and $\Gamma^{\con}_\pm$ to connect to $\Gamma^{SR}_{k_1, k_2}$. Then 
\begin{equation}
	\label{eq:c-DR}
	\Gamma_{k_1, k_2}^{DR}(\epsilon, H_1) = \Phi_{L, H_\epsilon^s}^* \left( \bar\Gamma_h^{DR} \right) \cup \Gamma_0^{\con} \cup \Gamma^{\con}_- \cup \Gamma^{\con}_+. 
\end{equation}
See Figure~\ref{fig:gamma-DR}. 

\subsubsection{Connecting to the double resonance point}
We define
\[
	c_0 = \Gamma_{k_1, k_2}, \quad c_0^\epsilon = \Phi_{L, H_\epsilon^s}^*(0, \alpha_{H^s_\epsilon}(0)),  
	\quad \Gamma^\con_0 = \bigcup_{s \in [0, 1]}\{s c_0 + (1-s)c_0^\epsilon\}
\]
and define $h$ to be the homology class corresponding to $k_1$.

Notice that we add a small segment to connect directly to 
$\Gamma_{k_1, k_2}$. 

\begin{prop}\label{prop:center-connector}
Suppose $H_1$ satisfies the satisfies the conditions 
$[DR1^h-DR3^h]$ and $[DR1^c-DR4^c]$ at the double resonance 	$\Gamma_{k_1, k_2}$ 
relative to $\Gamma_{k_1}$. Then there is $\epsilon_0, \delta >0$ depending only on $H^s$ such that for 
$0 < \epsilon < \epsilon_0$, $H_1 ' \in \cV_\delta(H_1)$ and $0 \le s \le 1$, for $H_\epsilon' = H_0 + \epsilon H_1'$: 
\[
	\cN_{H_\epsilon'}^0(s c_0 + (1-s)c_0^\epsilon) = 
	\cA_{H_\epsilon'}^0(s c_0 + (1-s)c_0^\epsilon) = 
	\cA_{H_\epsilon'}^0(c_0) = \cA_{H_\epsilon'}^0(c_0^\epsilon)
\]
is contractible in $\T^2$. As a result the cohomologies in $\Gamma^\con_0$ are forcing equivalent. 
\end{prop}
\begin{proof}
Consider $\bar{c}_0^\epsilon \in \R^2$, $\bar{\alpha}_0^\epsilon \in \R$ defined by the formula
\[
\bmat{ \bar{c}_0^\epsilon\sqrt{\epsilon} \\ \bar{\alpha}_0^\epsilon\epsilon} = 
M^{-T} \bmat{ 0 \\ - \alpha_{H_\epsilon}(c_0) + H_0(c_0)}, 
\]
then according to \eqref{eq:c-barc}, $\Phi_{L, H_\epsilon^s}^*(\bar{c}_0, \bar{\alpha}_0) = c_0$. According to Lemma~\ref{lem:pert-alpha}, we have $\|\alpha_{H_\epsilon}(c_0) - H_0(c_0)\| \le C\epsilon$ for some $C$ depending only on $H_0$, therefore $\|\bar{c}_0^\epsilon\| \le C \sqrt{\epsilon}$. 

 Recall that $O = (0, 0)$ is a hyperbolic fixed point of the system $H^s$, and 
$\tcA_{H^s}(0) = O$.  By standard perturbation theory of hyperbolic sets, 
there is a neighborhood $V \ni O$, such that the system 
$H^s_\epsilon = K - U' + \sqrt{\epsilon}P$ with $U' \in \cV_\delta(U)$ and 
$0 < \epsilon < \epsilon_0$, $H^s_\epsilon$ admits a unique hyperbolic 
periodic orbit $O_\epsilon$ contained in $V$. Moreover, using the upper 
semi-continuity of the Aubry set Corollary~\ref{cor:aubry-unique-static}, 
by possibly choosing $\delta$ and $\epsilon_0$ smaller, we ensure for all 
$0 \le \lambda \le 1$,  
$\tcA_{H^s_\epsilon}(\lambda \bar{c}_0^{\,\epsilon}) \subset V$, and therefore 
$\tcA_{H^s_\epsilon}(\lambda\bar{c}_0^{\,\epsilon}) = 
\tcA_{H^s_\epsilon}(0) = O_\epsilon$. We note that in this case the Aubry set 
has a unique static class, hence Aubry set coincides with the Ma\~ne set, also 
the discrete Aubry set is finite and therefore contractible. We now apply 
symplectic invariance \eqref{eq:symp-inv-linear} to get the same for 
the original system. 
\end{proof}

\subsubsection{Connecting single and double resonance}

Recall that the single resonance cohomologies 
$\Gamma^{SR}_{k_1}{(M,K)}$ is defined in \eqref{eq:punctured-sr}.
For each double resonance $\Gamma_{k_1, k_2}$, 
$\Gamma^{SR}_{k_1,k_2,\pm}$ be the two connected component 
of \eqref{eq:punctured-sr} adjacent to $\Gamma_{k_1, k_2}$. 
We define connectors $\Gamma^\con_\pm$ which connect 
the double resonance cohomology curve 
$\Phi_L^* \left( \bar\Gamma_h^{DR} \right)$ to 
$\Gamma^{SR}_{k_1,k_2,\pm}$, respectively. 

Let $c \in \R^2 \simeq H^1(\T^2, \R)$ be a cohomology, and let $\rho_H(c)$ 
denote the convex hull of rotation vectors of all $c$-minimal measures. 
This coincides with the Legendre-Fenchel transform relative to the alpha 
function 
\[
	\rho_H(c) = \LF_{\alpha_H}(c) \subset \R^2. 
\]
Note that $\rho_H(c)$ is a set valued function taking values in convex sets. 
A common feature of the cohomology classes we've chosen is that they stay on the rational line
\[
	\Omega_{k_1} = \{\omega: \quad k_1 \cdot (\omega, 1) = 0\}. 
\]
Let us also denote $\Omega_{k_1, k_2} = \Omega_{k_1} \cap \Omega_{k_2}, \ c_0 = \Gamma_{k_1, k_2},$ and 
$\omega_0 = \Omega_{k_1, k_2}.$

\begin{figure}[t]
	\centering
	\includegraphics[width=2in]{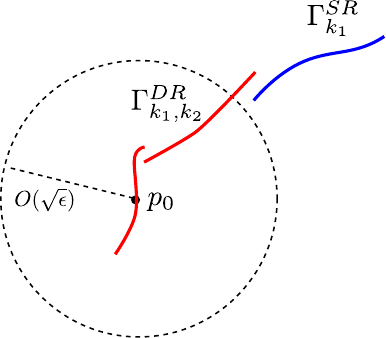}
	\includegraphics[width=2.5in]{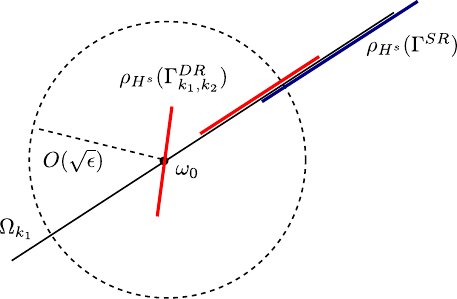}
	\caption{Left: cohomology curves, right: rotation vectors.} \label{fig:transition}
\end{figure}

The main observation is that the rotation vector of the curve 
$\Gamma^{SR}_{k_1}{(M,K)}$ and $\Gamma^{DR}_{k_1, k_2}$ overlap 
on the line $\Omega_{k_1}$, see Figure~\ref{fig:transition}. 
To prove this statement, we show that the rotation vector of 
$c \in \Gamma^{SR}_{k_1}{(M,K)}$ is $O(\sqrt{\epsilon})$ close 
to $\omega_0$ at its nearest point, while the rotation vector of 
$c \in \Gamma^{DR}_{k_1, k_2}$ is $C^{-1} E\sqrt{\epsilon}$ away from 
$\omega_0$ when $E$ is sufficiently large. Since both sets lie on 
the same line, they have to overlap. 

\begin{prop}\label{prop:overlap}
	Let $\Gamma^{SR}_{k_1,k_2,\pm}$ and 
	$\Gamma_{k_1, k_2}$ be as before. Suppose $H_1$ the conditions 
	$[DR1^h-DR3^h]$ and $[DR1^c-DR4^c]$ relative to 
	$\Gamma_{k_1}$. Then there is $\epsilon_0, \delta >0$ depending 
	only on $H^s$ such that for $0 < \epsilon < \epsilon_0$ such that 
	if $0< \epsilon < \epsilon_0$, and $H_1' \in \cV_\delta(H_1)$, for 
	$H_\epsilon' = H_0 + \epsilon H_1'$:
	\begin{enumerate}
		\item For each $c \in \Gamma^{SR}_{k_1}{(M,K)}$,  
$\rho_{H_\epsilon'}(c)$ is a single point contained in $\Omega_{k_1}$. 
The function $\rho_{H_\epsilon}(c)$ is continuous on $\Gamma^{SR}$, 
and there is $C> 0$ such that
		\[
			\|\rho_{H_\epsilon}(c) - \omega_0\| \le C \|c - c_0\|, \quad \forall \ c \in \Gamma^{SR}_{k_1}{(M,K)}. 
		\]
		\item Let $c_h(E) = \Phi_L^*\left( \bar{c}_h(E)\right) \in \Gamma_{k_1, k_2}^{DR}$. Then each $\rho_{H_\epsilon'}(c_h(E))$  is a single point contained in $\Omega_{k_1}$, and there are $C, E_0>0$ such that 
		\[
			\|\rho_{H_\epsilon'}(c_h(E)) - \omega_0\| \ge C^{-1} E \sqrt{\epsilon}, \quad E \ge E_0. 
		\]
	\end{enumerate}
	As a result, there exists $\bar{E}>0$ such that  for $E > \bar{E}$, the rotation vectors 
	$\rho_{H_\epsilon'}(c_h(E))$ coincide with one of 
	$\rho_{H_\epsilon'}(\Gamma^{SR}_{k_1,k_2,\pm})$, 
	see Figure~\ref{fig:transition}. Similarly, for $E$ sufficiently 
	large, the rotation vectors $\rho_{H_\epsilon'}(c_{-h}(E))$ 
	coincide with  one of 
	$\rho_{H_\epsilon'}(\Gamma^{SR}_{k_1,k_2,\pm})$ 
	(which are not covered in the first case).
\end{prop}
\begin{proof}
In the single resonance regime, after a linear coordinate change (Section~\ref{sec:res-comp}) the system is converted to $H_0 + \epsilon Z(\theta^s, p) + \epsilon R(\theta^s, \theta^f, p^s, p^f, t)$, where the invariant cylinder is given by $(\theta^s, p^s) = (\Theta^s, P^s)(\theta^f, p^f, t)$, and the Aubry set for any $c \in \Gamma^s = \{(\partial_p H_0(0), 1) \cdot (1, 0, 0) = 0\}$ is a graph over $(\theta^f, t)$ (see Theorem~\ref{thm:nhic-sr}, Theorem~\ref{thm:var-local-sr}) This implies that any $c$ admits a unique rotation vector, as different rotation vectors will resulting in intersecting minimizing orbits, which violates the graph theorem. By reverting the coordinate change, we also obtain $\rho_{H_\epsilon}(c) \in \Omega_{k_1}$. We then apply Corollary~\ref{cor:near-int-rho}, to get (in general) 
\[
\|\rho_{H_\epsilon}(c) - \nabla H_0(c)\| \le C \sqrt{\epsilon}, 
\]
while 
\[
\|\nabla H_0(c) - \omega_0\| = \|\nabla H_0(c) - \nabla H_0(c_0)\| \le C\|c - c_0\|
\]
and recall that $\|c - c_0\| \ge K \sqrt{\epsilon}$ for any $c \in \Gamma_{k_1}^{SR}$, item (1) follows. 

In the double resonance case, the cohomology curve $\bar{c}_h(E)$ at 
high energy are of Aubry-Mather type corresponding to the homology class 
$h = (1, 0)$. The same arguments as in the single resonance case implies 
$\rho_{H_\epsilon^s}(c)$ is unique, and contained in the line $\omega_1 = 0$. 
Again by reverting the coordinate change, we obtain the desired property for 
$\rho_{H-\epsilon}(c)$. For the in equality, note that if any point $(\varphi, I, t)$ 
in the Aubry set $\tcA_{H^s_\epsilon}(c)$ satisfies $H^s(\varphi, I) \ge E_0$, 
we have $K(I) \ge H^s(\varphi, I) - \|U\|_{C^0} \ge \frac12 H^s(\varphi, I)$ if 
$E_0 \ge 2\|U\|_{C^0}$. As a result
\[
\min_{(\varphi, I, \tau) \in \tcA_{H_\epsilon^s}(\bar{c}_h(E)) }\|I\| \ge C^{-1} \min_{(\varphi, I, \tau) \in \tcA_{H_\epsilon^s}(\bar{c}_h(E)) } \sqrt{K(I)} \ge C^{-1} \sqrt{E}
\]
for a suitable $C>1$. Reverting the coordinate change $\Phi_L$ implies
\[
\min_{(\theta, p, t) \in \tcA_{H_\epsilon}(c_h(E))} \|p - p_0\| \ge C^{-1} \|I\| \sqrt{\epsilon} \ge C^{-1} \sqrt{E}\sqrt{\epsilon}. 
\]
Finally, we apply Corollary~\ref{cor:near-int-rho} again to get our inequality. 
\end{proof}

We have the following lemma:
\begin{lem}[\cite{Mas2003}, Proposition 6]\label{lem:int-aubry}
	Let $c_1, c_2 \in H^1(\T^2, \R)$ and $\rho \in H_1(\T^2, \R)$ satisfy 
$\rho_H(c_1) = \rho_H(c_2) = \rho$, and both $c_1, c_2$ lie in 
the relative interior of $\LF_{\beta_H}(\rho)$. Then $\tcA_H(c_1) = \tcA_H(c_2)$. 
\end{lem}

\begin{prop}\label{prop:SR-connector}
Suppose $H_1$ satisfies the conditions $[DR1^h-DR3^h]$ and $[DR1^c-DR4^c]$ 
relative to $\Gamma_{k_1}$. Then there is $\epsilon_0, \delta >0$ depending only 
on $H^s$ {such that for $0 < \epsilon < \epsilon_0$ 
and a residual subset} $\cR_\delta(H_1)$ of  $\cV_\delta(H_1)$,  for each 
$H_1' \in \cR_\delta(H_1)$, there is $c_1^\pm \in \Gamma^{SR}_{k_1}(M,K)$ 
and $c_2^\pm \in \Gamma^{DR}_{k_1, k_2}$, such that for all 
$0 \le \lambda \le 1$ {and $H_\epsilon' = H_0 + \epsilon H_1'$}
\[
	\cN_{H_\epsilon'}^0(\lambda c_1^\pm + (1-\lambda)c_2^\pm) = 
\cA_{H_\epsilon'}^0(\lambda c_1^\pm + (1-\lambda)c_2^\pm) = 
\cA_{H_\epsilon'}^0(c_1^\pm) = \cA_{H_\epsilon'}^0(c_2^\pm)
\]
is contractible in $\T^2$. As the result both cohomology curves
\[
	\Gamma^\con_+ = \bigcup_{s \in [0, 1]} \{sc_1^+ + (1-s)c_2^+\}, \quad
	\Gamma^\con_- = \bigcup_{s \in [0, 1]} \{sc_1^- + (1-s)c_2^-\}
\]
are contained in a single forcing equivalent class. 
\end{prop}

\begin{proof}
Proposition~\ref{prop:overlap} implies the curves $\rho_{H_\epsilon}(\Gamma^{SR}_{k_1})$ and $\rho_{H_\epsilon}(\Gamma^{DR}_{k_1, k_2})$ overlap on an interval contained in $\Omega_{k_1}$, for all $H_1' \in \cV_\delta(H_1)$.  In particular, there must be $c_1 \in \Gamma^{SR}_{k_1}$ and $c_2 \in \Gamma^{DR}_{k_1, k_2}$ where they both have rational rotation vectors.  We now assume that $H_\epsilon'=H_0 + \epsilon H_1'$ satisfies the residual condition that all Aubry sets with rational rotation vector is supported on a hyperbolic periodic orbit, in this setting $\cA = \cN$, $c_1, c_2$ are contained in the relative interior of $\LF_\beta(\rho)$. Lemma~\ref{lem:int-aubry} now implies all the Aubry sets $\cA_{H_\epsilon}(\lambda c_1 + (1-\lambda)c_2)$, $0 \le \lambda \le 1$ coincide, which implies \eqref{eq:mather-mech} applies to the whole segment. The proposition follows. 
\end{proof}

\subsection{Jump from non-simple homology to simple homology}
\label{sec:intro-jump}

As described, when $h$ is not a simple homology, the cohomology class $\Gamma_{k_1, k_2}^{DR}$ as chosen is not connected. More precisely, $\Gamma_{k_1,k_2}^{DR}$ consists of three connected components (see also \eqref{eq:bargm-non-simp} and \eqref{eq:c-DR}):
\[
	\Phi_L^*\bar\Gamma_h^e \cup \Gamma^\con_+,\quad  \Phi_L^*\bar\Gamma_{-h}^e \cup \Gamma^\con_-, \quad \Phi_L^*\left( \bar\Gamma_{h_1}^{e, \mu} \cup \bar\Gamma_{-h_1}^{e, \mu} \right) \cup \Gamma^\con_0.
\]
We will show forcing equivalence of the components by the following:
\begin{thm}[Section~\ref{sec:forcing-for-kissing}]\label{thm:forcing-jump}
Suppose the slow system $H^s$ satisfies the conditions $[DR1^c] - [DR4^c]$, and that the associated homology $h = n_1 h_1 + n_2 h_2$ is non-simple. Then there exists $e, \mu, \epsilon_0, \delta >0$ depending on $H^s$, such that there is a residual subset $\cR_\delta(H_1)$ of $\cV_\delta(H_1)$, and for each $H_1' \in \cR_\delta(H_1)$, there is $E_1, E_1 \in (e, e + \mu)$ such that 
\[
	\Phi_L^*\left( \bar{c}_h(E_1) \right),\quad \Phi_L^*\left(  \bar{c}_{h_1}^{\,\mu}(E_2) \right)
\]
are forcing equivalent, with respect to $H_0 + \epsilon H_1'$. The same conclusions apply when $h, h_1$ are replaced with $-h, -h_1$. 
\end{thm}

See the dashed line in Figure~\ref{fig:gamma-DR}. 

\subsection{Forcing equivalence at the double resonance}
\label{sec:intro-forcing}

We summarize all of our constructions in the following theorem.
\begin{thm}
	\label{thm:c-equiv-dr}
	Suppose $H_1$ satisfies all non-degeneracy conditions of $\Gamma_{k_1}$ and $\Gamma_{k_1, k_2}$. 
	Then there is $\epsilon_1 = \epsilon_1(H_0, H_1) > 0$, $\delta = \delta(H_0, H_1) > 0$, such that for all 
	$H_1' \in \cV_\delta(H_1)$, and $0 < \epsilon < \epsilon_1$, there is a subset $\Gamma^{DR}_{k_1, k_2} = \Gamma^{DR}_{k_1, k_2}(H_0, H_1', \epsilon)$ satisfying:
	\begin{enumerate}
		\item $\Gamma_{k_1, k_2} \subset \Gamma^{DR}_{k_1, k_2}$;
		\item $\Gamma^{DR}_{k_1, k_2}$ intersects both $\Gamma^{SR}_{k_1, k_2,+}$ and $\Gamma^{SR}_ {k_1, k_2,-}$;
	\end{enumerate}
	there is $\sigma = \sigma(H_0, H_1, \epsilon) > 0$ and a residual subset $\cR_\sigma(H_0 + \epsilon H_1')$ of $\cV_\sigma(H_0 + \epsilon H_1')$, such that for all $H \in \cR(H_0 + \epsilon H_1')$, all of $\Gamma^{DR}_{k_1, k_2}$ are forcing equivalent. 
\end{thm}
\begin{proof}
	The two properties of $\Gamma_{k_1, k_2}^{DR}$ hold by construction. During the proof, we say a property hold ``after a residual perturbation'' if it holds on an residual subset of a neighborhood of the corresponding Hamiltonian. Note that if several properties each hold after a residual perturbation, then they hold simultaneously after a residual perturbation. Moreover, if a property holds after a residual perturbation, and is invariant under coordinate changes, then it also holds in the new coordinate after a residual perturbation, provided the coordinate change is smooth enough. This is the case for our system since all coordinate changes are $C^\infty$. 

	Let $\delta$ be the smallest parameter depending on $H^s$ such that all of Corollary~\ref{cor:DR-equiv}, Proposition~\ref{prop:center-connector} and Proposition~\ref{prop:SR-connector} holds. 

	Let $\delta_1 = \delta_1(H_0, H_1)$ be such that for any $H_1' \in \cV_{\delta_1}(H_1)$, the associated slow system is in $\cV_{\delta}(H^s)$. Then for $0 < \epsilon < \epsilon_0$, Lemma~\ref{lem:energy-reduction} implies the system $H_0 + \epsilon H_1'$ is reduced to a system $(G^s + \sqrt{\epsilon} P)/\beta$ satisfying the conclusions of Corollary~\ref{cor:DR-equiv}. Moreover, by Proposition~\ref{prop:DR-symp-inv}, relative to the the normal form system $N_\epsilon$, one of the three diffusion mechanisms hold for $\Phi_L^*(\bar{\Gamma}_h^{DR})$, after taking a residual perturbation. By Lemma~\ref{lem:symp-inv-normal-form}, the same holds for relative to the original system $H_0 + \epsilon H_1'$.  

	Recall that (see \ref{eq:c-DR}) $\Gamma_{k_1, k_2}^{DR}$ consists of the set $\Phi_L^*(\bar{\Gamma}^{DR}_h)$ and the connector sets $\Gamma_0^\con$ and $\Gamma_\pm^{\con}$. Proposition~\ref{prop:center-connector} and \ref{prop:SR-connector} implies forcing equivalence of the connectors after a residual perturbation. This implies each connected component of $\Gamma_{k_1, k_2}^{DR}$ are forcing equivalent.

	Finally, Theorem~\ref{thm:forcing-jump} implies all three components of $\Gamma_{k_1, k_2}^{DR}$ are equivalent to each other after a residual perturbation.
\end{proof}

Assuming all the propositions and theorems formulated thus far, we prove Theorem~\ref{thm:c-equiv-path} which implies our main theorem. 
\begin{proof}[Proof of Theorem~\ref{thm:c-equiv-path}]
	Let $H_1 \in \cU$, which means that $H_1$ satisfies $SR(k_1, \lambda)$ for some $\lambda>0$ and all $k_1 \in \cK$, and that for all the strong double resonances $k_2 \in \bigcup_{k_1 \in \cK} \bigcup \cK^\text{st}(k_1, \Gamma_{k_1},  \lambda)$, $H_1$ satisfies the non-degeneracy conditions $DR(k_1, \Gamma_{k_1},  k_2)$. 

	For each $k_1$, Theorem~\ref{thm:c-equiv-sr} applies. Therefore, there exists $\epsilon_1^{k_1}(H_0, \lambda) > 0$ such that the theorem applies for each $H_0 + \epsilon H_1$ with $0 < \epsilon < \epsilon_1$. Since $SR(k_1, \Gamma_{k_1}, \lambda)$ is an open condition, there exists $\delta^{k_1} = \delta^{k_1}(H_0, H_1)$ such that the conclusion of the theorem to all $H_1' \in \cV_{\delta^{k_1}}(H_1)$. 

	For each $k_1, k_2$, the conclusion of Theorem~\ref{thm:c-equiv-dr} holds to all $H_1' \in \cV_{\delta^{k_1, k_2}}(H_1)$ and $0 < \epsilon < \epsilon_1^{k_1, k_2}(H_0, H_1)$. Define 
	\[
		\epsilon_1(H_0, H_1) = \min\left\{\min_{k_1, k_2} \epsilon_1^{k_1, k_2}, \min_{k_1} \epsilon_1^{k_1}  \right\}, \quad
		\delta(H_0, H_1) = \min\left\{\min_{k_1, k_2} \delta^{k_1, k_2}, \min_{k_1} \delta^{k_1}  \right\}. 
	\]
	Then the conclusion of both theorems apply to $H_1' \in \cV_\delta(H_1)$ and $0 < \epsilon < \epsilon_1$.

	Define 
	\[
		\Gamma_*(H_0, H_1', \epsilon) = \bigcup_{k_1 \in \cK} \left( \Gamma^{SR}_{k_1} \cup \bigcup_{k_2 \in \cK^{st}(k_1, \lambda)} \Gamma^{DR}_{k_1, k_2} \right). 
	\]
	For each single resonance $k_1$, the union 
	$$
	\Gamma_*^{SR}(k_1) = 	\Gamma^{SR}_{k_1} \cup 
	\bigcup_{k_2 \in \cK^{st}(k_1, \lambda)} 
	\Gamma^{DR}_{k_1, k_2}
	$$ are contained in a single equivalent class, since each $\Gamma^{DR}_{k_1, k_2}$ are forcing equivalent, and they connect all the disconnected pieces from $\Gamma^{SR}_{k_1}$. If two single resonances $\Gamma_{k_1}$ and $\Gamma_{k_1'}$ intersect at $\Gamma_{k_1, k_1'}$, then $\Gamma_*^{SR}(k_1)$ and $\Gamma_*^{SR}(k_1')$ also intersect at $\Gamma_{k_1, k_1'}$, since $\Gamma_{k_1, k_1'}$ is contained in both $\Gamma^{DR}_{k_1, k_1'}$ and $\Gamma^{DR}_{k_1', k_1}$. As a result, the entire $\Gamma_*(H_0, H_1', \epsilon)$ is contained in a single forcing equivalent class. 

	Finally, if $U_1, \dots, U_N$ are open sets which intersect $\cP$, by setting $\epsilon_0$ small enough, they also intersect $\bigcup_{k_1 \in \cK}\Gamma^{SR}_{k_1}$, since the said union is obtained from $\cP$ by removing finitely many neighborhoods of size $O(\sqrt{\epsilon})$. Therefore, $U_1, \dots, U_N$ also intersect $\Gamma_*(H_0, H_1', \epsilon)$. 
\end{proof}

\section{Weak KAM theory and forcing equivalence}
\label{sec:forcing-def}

In this section we give an introduction to weak KAM theory and forcing relation. Most of the presentation follow \cite{Be} and \cite{Fa}. 

The forcing relation is introduced by Bernard (\cite{Be}). It generalizes a similar equivalence relation defined by Mather (\cite{Ma2}), using the point of view of Fathi (\cite{Fa}), see also \cite{KO}. This approach frees us from needing to construct a variational principle for the global diffusion orbit, as the approaches in \cite{Ma1,Ma,Ma3,Ma4,CY1,CY2} requires. Another advantage is that Mather's approach works only on continuous curves of cohomologies, as the his definition is local. The forcing relation allows disconnected pieces of cohomology as long as we prove the forcing relation by definition, enabling the \emph{jump mechanism}. 

\subsection{Periodic Tonelli Hamiltonians}
\label{sec:intro-Tonelli}

A $C^2$ Hamiltonian 
\[
	H: \T^n \times \R^n \times \R \to \R
\]
is called (time-periodic) Tonelli if it satisfies:
\begin{enumerate}
	\item (Periodicity) There is $0 < \vp = \vp_H$ such that $H(\theta, p, t + \vp) = H(\theta, p, t)$. 
	\item (Convexity) $\partial^2_{pp}H(x, p, t)$  is strictly positive definite as a quadratic form. 
	\item (Superlinearity) $\lim_{\|p\| \to \infty} H(x, p, t)/\|p\| \to \infty$. 
	\item (Completeness) The Hamiltonian vector field generates a complete flow on $\T^n \times \R^n$. We denote by $\phi_H^{s,t}$ the flow from time $s$ to time $t$, and by $\phi_H$ the flow $\phi_H^{0,1}$. 
\end{enumerate}
We will denote by 
\[
	L = L_H(\theta, v, t) = \sup\{p \cdot v - H(\theta, p, t)\}
\]
its Legendre transform. For the most part, we will restrict to Hamiltonians with $\vp = 1$, namely, defined on $\T^n \times \R^n \times \T$, but near double resonances we need to consider Hamiltonians that are $\sqrt{\epsilon}-$periodic. 

It is helpful to consider a family of Hamiltonians which satisfy these properties uniformly. For $D>0$, consider 
\[
	\begin{aligned}
		\bH(D) = 
		\Bigl\{ 
		& H \in C^2(\T^n \times \R^n \times \R): \\
		& \vp_H \le 1, \quad D^{-1} I \le \partial^2_{pp}H \le D I, \quad 
		\|H(\cdot, 0)\|_{C^0}, \|\partial_p H(\cdot, 0)\|_{C^0} \le D. 
		\Bigr\}
	\end{aligned}
\]
We then check that each $H \in \bH(D)$ is Tonelli, and it satisfy a list of uniform estimates called \emph{uniform family} in \cite{Be}. In particular, if $H_0 \in \bH(D)$, then $H_\epsilon \in \bH(2D)$ for all $\epsilon \le \epsilon_0 = \epsilon_0(D)$.

Given $C >0$, we say a function $u: \R^n \to \R$ is $C$ semi-concave if for every $x \in \R^n$, there is a linear function $l_x: \R^n \to \R$, such that 
\[
	u(y) - u(x) \le l_x(y-x) + C \|y-x\|^2, \quad y \in \R^n. 
\]
The linear form $l_x$ is called a super-differential at $x$. The set of all super-differentials at $x$ is denoted $\partial^+ u(x)$. It is easy to see if $u$ is differentiable at $x$, then $\partial^+ u(x) = \{ du(x)\}$. 
A function $u: \T^n \to \R$ is semi-concave if it's semi-concave as a function on $\R^n$. 

\begin{lem}(\cite{Be})
If $u: \T^n \to \R$ is $C$ semi-concave, then it is 
$C\sqrt{n}$--Lipschitz. The super-differential set 
$\partial^+ u(x) \subset \{\|p\| \le C \sqrt{n}\}$.
\end{lem}

Given $s < t \in \R$, $x, y \in \T^n$, we define the Lagrangian action
\begin{equation}
	\label{eq:def-action}
	A(x, s, y, t) = \inf_{\gamma(s) = x, \, \gamma(t) = y} \int_s^t L_H(\gamma(\tau), \dot{\gamma}(\tau), \tau) d\tau, 
\end{equation}
where the infimum is taken over all absolutely continuous $\gamma$. 
We outline a series of useful results:
\begin{prop}
	\label{prop:action-concave}
	Let $H \in \bH(D)$, then:
	\begin{enumerate}
		\item  (Tonelli Theorem) (\cite{Ma1}, Appendix 2) The infimum in \eqref{eq:def-action} is always reached. It is then $C^2$, which solves the Euler-Lagrange equation. Such a $\gamma$ is called a minimizer.
		\item (A priori compactness) (\cite{Be}, Section B.2) Let $\delta >0$, then there is a constant $C_\delta>0$ depending only on $\delta$ and $D$, such that for $t -s > \delta$, any minimizer of $A(s, x, t, y)$ satisfies $\|\dot{\gamma}\| \le C_\delta$. 
		\item (Uniform semi-concavity) (\cite{Be}, Theorem B.7, see also Proposition~\ref{prop:near-auto-semi-concave}) For $t - s > \delta$, the function $A(s, x; t, y)$ is $C_\delta$ semi-concave in $(x,s)$ and $(y, t)$. Moreover, if $\gamma:[s,t] \to \T^n$ is a minimizer, then 
		\begin{equation}
			\label{eq:sup-grad-action}
			\begin{aligned}
				& 			\left( p(s), - H(\gamma(s), p(s), s) \right)  \in - \partial^+_{(x, s)}A(x, s, y, t), \\
				& \left( p(t), - H(\gamma(t), p(t), t) \right)  \in - \partial^+_{(y, t)}A(x, s, y, t),
			\end{aligned}
		\end{equation}
		where $p(\tau) = \partial_v L_H(\gamma(\tau), \dot\gamma(\tau), \tau)$. 
	\end{enumerate}
\end{prop}
\begin{rmk}
For item 3 of Proposition~\ref{prop:action-concave}, \cite{Be} only treated the semi-concavity in $x$. A proof that also include time-dependence can be obtained from Proposition~\ref{prop:near-auto-semi-concave} by choosing a non-small $\epsilon$. 
\end{rmk} 

For each $c \in \R^n$, define 
\[
	L_{H, c}(\theta, v, t) = L_H(\theta, v, t) - c \cdot v, 
\]
and denote $A_{H, c} = A_{L_{H, c}}$. Then \eqref{eq:sup-grad-action} becomes
\[
	\left( p(t) - c,  - H(\gamma(t), p(t), t) \right) \in \partial_{(y, t)}A_{H,c}(x, s, y, t)
\]
and similarly for $s$.

\subsection{Weak KAM solution}
\label{sec:intro-weak-kam-solution}

We now define the (continuous) Lax-Oleinik semi-group $T_c^{s, t}: C(\T^n) \to C(\T^n)$ via the formula 
\[
	T_c^{s, t} u(x) = \min_{z \in \T^n} \left\{  u(z) + A_{H,c}(z, s; x, t)  \right\} . 
\]
For a $\vp$-periodic Hamiltonian, the associated discrete semi-group is generated by the operator: $T_c u = T_c^{0, \vp} u$. 

\begin{lem}\label{lem:inf-concave}
Let $\{u_\zeta\}_{\zeta \in Z}$ be a (possibly uncountable) family of $C$-semi-concave functions $\T^n \to \R$, and $v = \inf u_\zeta$ is bounded, then $\inf_{\zeta \in Z} u_\zeta$ is also $C$ semi-concave.

Moreover, suppose for $x_0 \in \Z$, the infimum $v(x_0)$ is reached at $u_{\zeta_0}(x_0)$, then $\partial^+ u_{\zeta_0}(x_0) = \partial^+ v(x_0)$.
\end{lem}
Using Proposition~\ref{prop:action-concave}, the functions $T^n_\eta u$ are $C$ semi-concave with the constant depending only on uniform family. 

\begin{prop}[Ma\~ne's critical value, see \cite{Be}, Proposition~3.1]
	There is a unique $\alpha \in \R$ such that $T_c^n u(x) + n \alpha$ stays bounded for all $n \in \N$. This value coincides with Mather's alpha function: 
	\begin{equation}
		\label{eq:alpha}
		\alpha_H(c) = - \inf_{\mu}\left\{   \int L_{H, c}(\theta, v, t) d\mu(\theta, v, t) \right\}
	\end{equation}
	where the infimum is taken over all invariant probability measure of the Euler-Lagrange flow on $\T^n \times \R^n \times \T_\vp$. 
\end{prop}
Let us also point out an alternative definition of the alpha function (see \cite{Sor2015}), namely, we can replace the class of minimal measures with the class of \emph{closed measures}, which satisfies 
\[
	\int df(\theta, t) \cdot (v, 1) d\mu(\theta, v, t) = 0 
\]
for every $C^1$ function $f: \T^n \times \R \to \R$ satisfying $f(\theta, t + \varpi) = f(\theta, t)$. 

A function $w: \T^n \times \T_\varpi \to \R$ is called 
a {\it weak KAM solution} if 
\[
	T_c^{s, t} w(\cdot, s) + \alpha_{H}(c) (t-s)= w(\cdot, t), \quad s< t \in \R,
\]
i.e. the family $w(\cdot, t)$ is invariant under the semi-group up to a linear drift. The function $u(\theta) = w(\theta, 0)$ is then a fixed point of the operator $T_c + \vp \alpha_H(c)$.

\begin{prop}[Existence of weak KAM solution, \cite{Be}, Proposition~3.2]
	For $c \in \R^n$ and $u_0 \in C(\T^n)$, the function 
	\[
		w(\theta, t) = \liminf_{N \to - \infty} (T_c^{t - N \vp ,t} u_0  + N \vp \alpha(c)), \quad N \in \N
	\]
	is  a $\vp$-periodic weak KAM solution. 
\end{prop}

The following is a easy consequence of Lemma~\ref{lem:inf-concave}. 
\begin{lem}\label{lem:weak-KAM-semi-concave}
Weak KAM solutions are uniformly semi-concave over all  $H \in \bH(D)$. 
\end{lem}

A function $w:\T^n \times \R \to \R$ is called {\it dominated} 
by $L_{H, c}, \alpha$ if 
\[
	w(y, t) - w(x, s) \le A_{H, c}(x, s, y, t) + (t-s)\alpha, 
	\quad x, y \in \T^n, \quad t < s. 
\]
$\gamma: I \to \T^n$, where $I$ is an interval in $\R$ is  
called {\it calibrated} by $w$ if 
\[
	w(\gamma(t), t) - w(\gamma(s), s) = A_{H, c}(x, s, y, t) + (t-s)\alpha. 
\]
\begin{prop}[\cite{Fa}, Proposition~4.1.8]
	$w: \T^n \times \R \to \R$ is a weak KAM solution if and only if it is dominated by $L_c, \alpha_H(c)$ and for every $(y, t)$ there is a calibrated curve $\gamma: (-\infty, t] \to \T^n$ such that $\gamma(t) = y$. 
\end{prop}

\subsection{Pseudographs, Aubry, Ma\~ne and Mather sets}

Let $u: \T^n \to \R$ be semi-concave. By the Radamacher theorem, $u$ is differentiable almost everywhere. For $c \in \R^n$,  we define the (overlapping) psudograph 
\[
	\cG_{c, u} = \cG_{c, H, u} = \left\{  (x, c +  \nabla u(x)), \quad  \nabla u(x) \text{ exists}  \right\}.
\]
In the one-dimensional case, at every discontinuity of the function $du(x)$, the left limit is larger than the right limit. In the time-dependent setting if $w: \T^n \times \R \to \R$ is semi-concave, we write 
\[
	\cG_{c, w} = \{(x, c + \partial_x u(x, t), t): \quad \partial_x u \text{ exists}\}. 
\]

The evolution by the Lax-Oleinik semi-group generates an evolution operator on the psudograph. The following statement outline its relation with the Hamiltonian dynamics.
\begin{prop}[\cite{Be}] \label{prop:graph-inv}
For each $s < t$, we have 
\[
	\overline{\cG_{c, T_c^{s,t} u}} \subset  \phi_H^{s, t}  \left(  \cG_{c, u} \right),
\]
here $\phi_H^{s,t}$ denotes the Hamiltonian flow. 
\end{prop}

\begin{cor}\label{cor:back-inv}
Suppose $w(\theta, t)$ is a (time-periodic) weak KAM solution of $L_{H, c}$, then 
\[
	\left( \phi_H^{s, t} \right)^{-1} \, \overline{ \cG_{c, w(\cdot, t)}}\subset \cG_{c, w(\cdot, s)}. 
\]
In particular, $(\phi_H^{-1})\,  \overline{\cG_{c, u}} \subset \cG_{c, u}$, where $u = w(\cdot, 0)$ and $\phi_H = \phi_H^{0, \vp}$ is the associated discrete dynamics.
\end{cor}

Let $w = w(\theta, t)$ be a continuous weak KAM solution for $L_c$. Then Corollary~\ref{cor:back-inv} implies the set 
\[
	\cG_{c, w} = \left\{ (\theta, p): \quad (\theta, p) \in \cG_{c, w(\cdot, t)} \right\}
\]
is backward invariant under the flow $\phi_H^{s, t}$. We then define
\[
	\tcI(c, w) = 
	\left\{ (\theta, p, t): \quad
	(\theta, p) \in \bigcap_{s < t} \left( \phi_H^{s, t} \right)^{-1} \, \overline{ \cG_{c, w(\cdot, t)}}
	\right\},  
\]
in other words, $\tcI(c, w)$ is the invariant set generated by the family of psudographs $\cG_{c, w(\cdot, t)}$, in the extended phase space $\T^n \times \R^n \times \T_\varpi$. 

The Aubry and Ma\~ne sets admit the following equivalent definitions: 
\[
	\tcA(c) = \bigcap_{w} \tcI(c, w), \quad \tcN(c) = \bigcup_{w} \tcI(c, w), 
\]
where $w: \T^n \times \R \to \T_\varpi$ is taken over all $L_c$ continuous-time weak KAM solutions. The Mather set $\tcM(c)$ is then the support of all $\phi_H^{s, t}$ invariant measures contained in $\tcA(c)$. Note that if we consider a discrete weak KAM solution $u: \T^n \to \R$, then the analogous definitions give us $\tcA^0$ and $\tcN^0$.

\subsection{The dual setting, forward solutions}

There is a dual setting which corresponds to forward dynamics (as apposed to the backward invariant sets obtained before). 
Define
\[
	\chT_c^{s, t} u (x) = \max_z \left\{ u(z) - A_c(x, s; z, t)  \right\},
\]
and note the following: 
\begin{enumerate}
	\item $\alpha = \alpha_H(c)$ is the unique number such that a weak KAM solution may exist. 
	\item $-\chT_c^{s, t} u$ are uniformly semi-concave (if $t-s > \tau$), fixed points of $\chT^{s,t}_c - \alpha$ exists, and are called forward weak KAM solution.
	\item For a semi-concave function $u$, we define $\chG_{c, u} = \{(\theta, c + \nabla_x u(\theta))\}$, and call it an anti-overlapping psudographs.
	\item Analogs of the previous section apply with appropriate changes. 
\end{enumerate}

Let $w(\theta, t)$ be a weak KAM solution for $L_{H, c}$, and $w^+$ a forward weak KAM solution. We say that $w, w^+$ are conjugate if they coincide on the set $\cM_H(c)$. 
\begin{prop}[\cite{Fa}, Theorem 5.1.2]
	For each weak KAM solution $w$, there exists a forward solution $w^+$ conjugate to $w$ satisfying $w^+ \ge w$, and 
	\[
		\tcI(c, w) = \tcI(c, w^+) = \tcI(c, w, w^+) := \{(x, c+ \partial_x w(x, t), t): \quad w(x, t) = w^+(x, t)\}. 
	\]
\end{prop}

\subsection{Peierls barrier, static classes, elementary solutions}

We define 
\[
	h_{H, c}(x, s, y, t) = \liminf_{N \to \infty} A_{H, c}(x, s, y, t + N \vp) + N \vp \alpha_H(c)
\]
called the Peirels barrier by Mather (\cite{Ma2}). The projected Aubry set $\cA_H(c)$ has the following alternative characterization: 
\[
	\cA_H(c) = \{(x, s) \in \T^n \times \T_\vp:  h_{H,c}(x, s, x, s) = 0\}. 
\]
We will also consider the discrete barrier:
\begin{equation}
	\label{eq:dis-barrier}
	h_{H, c}(x, y) = h_{H, c}(x, 0, y, 0), 
\end{equation}
then $\cA^0_{H}(c) = \{x:\,  h_{H, c}(x, x) = 0\}$.

On the set $\cA_H(c)$ we define the Mather semi-distance: 
\[
	\bar{d}(x, s, y, t) = h_{H, c}(x, s, y, t) + h_{H, c}(y, t, x, s), 
\]
then the condition $\bar{d}(x, s, y, t) = 0$ defines an equivalence relation $(x, s) \sim (y, t)$ on $\cA_H(c)$. The equivalence classes of this relation are called the \emph{static classes}. Let $\cS \subset \cA_H(c)$ be a static class, it corresponds uniquely to an invariant set in the phase space:
\[
	\tcS = \pi_{(x, t)}^{-1}\Bigr|_{\cA_H(c)} \cS. 
\]

For $(\zeta, \tau) \in \cS$, the function 
\[
	h_{H, c}(\zeta, \tau, \cdot, \cdot) 
\]
is a weak KAM solution for $L_{H, c}$, called the \emph{elementary solution}. The elementary solution is independent of the choice of $(\zeta, \tau) \in \cS$, up to an additive constant. 

We have the following useful statements concerning elementary solutions. 
\begin{prop} \label{prop:rep-formula}
\begin{enumerate}
	\item (Representation formula, see \cite{Fa}, Theorem 8.6.1 and \cite{CIS2013}, Theorem 7) Let $w(x, t)$ be an $L_{H, c}$ weak KAM solution. Then 
	\[
		w(x, t) = \min_{(\zeta, \tau) \in \cA_H(c)}\{ w(\zeta, \tau) + h_{H, c}(\zeta, \tau, x, t) \}. 
	\]
	\item (See \cite{Be}, Proposition 4.3) Every orbit in the \Mane set $\tcN_H(c)$ is a heteroclinic orbit between two static classes $\tcS_1, \tcS_2$. We have $\tcA_H(c) = \tcN_H(c)$ if and only if $\tcA_H(c)$ has only one static class. 
\end{enumerate}
\end{prop}

\subsection{The forcing relation}
\label{sec:def-forcing}

\begin{defn}\label{defn:forcing}
Let $u: \T^n \to \R$ be semi-concave, and $N \in \N$. We say that $\cG_{c, u} \vdash_N c'$, if there exists a semi-concave function $v: \T^n \to \R$ such that 
\[
	\overline{ \cG_{c', v} } \subset \bigcup_{k = 0}^N  \phi_H^{k} \left(  \cG_{c, u} \right). 
\]
We say the $c \vdash c'$ if there exists $N \in \N$ such that 
\[
	\cG_{c, u} \vdash_N c'
\]
for every psudograph $\cG_{c, u}$. We say $c \dashv \vdash c'$ if $c \vdash c'$ and $c' \vdash c$. 
\end{defn}

In view of Proposition~\ref{prop:graph-inv}, we always have $c \vdash c$. The relation is transitive  by definition. Therefore $\dashv \vdash$ defines an equivalence relation. 

We summarize the property of the forcing relation below. 
\begin{itemize}

	\item (Proposition~\ref{prop:forcing-implications})

	Let $\{c_i\}_{i = 1}^N$ be a sequence of cohomology classes which are forcing equivalent. For each $i$, let $U_i$ be neighborhoods of the discrete Mather sets $\tcM_H^0(c_i)$, then there is a trajectory of the Hamiltonian flow $\phi^t$ of $H$ visiting all the sets $U_i$. 

	\item (Mather mechanism, Proposition~\ref{prop:mather-mech}) 

	Suppose $\cN_H^0(c)$ is contractible
	as a subset of $\T^2$, then there is $\sigma>0$ such that $c$ is forcing equivalent to all $c' \in B_\sigma(c)$. 

	\item (Arnold and Bifurcation mechanism, Proposition~\ref{prop:arnold-mech})

	Suppose, either: $\tcA_H^0(c)$ has only two static classes and $\tcN_H^0(c) \setminus \tcA_H^0(c)$ is totally disconnected, or $\tcA_H^0(c)$ has only one static class, and there is a symplectic double covering map $\xi$ such that $\tcN^0_{H \circ \Xi}(\xi^* c)  \setminus  \Xi^{-1} \tcN^0_{H}$ is totally disconnected, then there is $\sigma>0$ such that $c$ is forcing equivalent to all $c' \in B_\sigma(c)$. 
\end{itemize}

\section{Perturbative Weak KAM theory}
\label{sec:pert-weak-kam}

By perturbative weak KAM theory, we mean two things: 
\begin{itemize}
	\item How  do the weak KAM solutions and the Mather, Aubry, \Mane sets respond to limits of the Hamiltonian;
	\item  How do the weak KAM solutions change when we perturb a system, in particular, what happens when we perturb (1) completely integrable systems, and (2) autonomous systems. 
\end{itemize}
In this section, we state and prove results in both aspects, as a technical tool for proving forcing equivalence.

\subsection{Semi-continuity}

Let $\bH(D)$ be the uniform family defined before, note that they are periodic of period $0 < \vp_H \le 1$, but not necessarily of the same period. It is known that in the case that all periods are fixed at $1$, the weak KAM solutions are upper semi-continuous (see precise statements below) under $C^2$ convergence over compact sets (\cite{Be3}). The results  generalize to the case when the periods are not the same, as we now show. 

Let us remark that if $H_n \in \bH(D)$ is a family of Hamiltonians, and $H_n \to H$ uniformly over compact sets on $\T^n \times \R^n \times \R$, then $H$ is necessarily periodic of \emph{some} period, and therefore $H \in \bH(D)$. 

We define the upper limit $\limsup$ for a sequence of sets $A_n$ to be the set of all accumulation points of all sequences $x_n \in A_n$. 
\begin{lem}[\cite{Be3}, Lemma 7]\label{lem:conv-weak-KAM}
Suppose Hamiltonians $H_k \in \bH(D)$, $k \in \N$ are  a family of periodic Tonelli Hamiltonians. Suppose $H_n \to H$ in $C^2$ over compact sets,   
$c_k \to c \in \R^n$, then $\alpha_{H_k}(c_k) \to \alpha_H(c)$. 
If $w_k: \T^n \times \R \to \R$ is a sequence of weak KAM solution of 
$H_k, c_k$ which converges uniformly to $w: \R^n \times \R \to \R$, then 
$w$ is a weak KAM solution of $H, c$.

Moreover, we have 
\[
	\limsup_{k\to \infty} \cG_{c_k, H_k, w_k}  \subset \cG_{c, H, w}, \quad
	\limsup_{k\to \infty} \tcI_{H_k}(c_k, w_k) \subset \tcI_H(c, w), \quad
	\limsup_{k \to \infty} \tilde{\cN}_{H_k}(c_k) \subset \tilde{\cN}_H(c). 
\]
\end{lem}
\begin{proof}
	The proof is an elaboration of \cite{Be3}, Lemma 7. 
	First, note that if $H_k \in \bH$ and $c_k$ uniformly bounded, 
	then $\alpha_{H_k}(c_k)$ is uniformly bounded (\eqref{eq:alpha}). By restricting to a subsequence, 
	we may assume $\alpha_{H_k}(c_k) \to \alpha \in \R$. Then 
	taking limit in 
	\[
		w_k(\gamma(t),t) - w_k(\gamma(s), s) \le \int_s^t L_{H_k}(\gamma(\tau), \dot{\gamma}(\tau), \tau) - c_k \cdot \dot{\gamma}(\tau) + \alpha_{H_k}(c_k) \, d\tau
	\]
	we obtain
	\[
		w(\gamma(t),t) - w(\gamma(s), s) \le \int_s^t L_{H}(\gamma(\tau), \dot{\gamma}(\tau), \tau) - c \cdot \dot{\gamma}(\tau) + \alpha \, d\tau,
	\]
	implying $\alpha \ge \alpha_{H}(c)$, using the definition of the $\al$-function (\ref{eq:alpha}). 

	Now given $x \in \T^n$, $t \in \R$, let $\gamma_k: (-\infty, t]\to \T^d$ be a $L_{H_k, c_k}$,  $w_k$ calibrated curve,  then Proposition~\ref{prop:action-concave} implies $\gamma_k$ are uniformly Lipschitz. Then $\gamma_k$ has a subsequence that converges in $C^1_{loc}$ to a limit $\gamma(t)$, which is $L_{H, c} + \alpha$, $u$ calibrated. This implies both $\alpha_H(c) = \alpha$ and that $u$ is a $H, c$ weak KAM solution. 

	We now prove the ``moreover'' part. Denote
	\[
		\cG_k = \cG_{c_k, H_k, w_k}, \quad \cG = \cG_{c, H, w}
		\quad \tcI_k = \tcI_H(c_kw_k), \quad \cI = \cI_H(c, w). 
	\]
	then  if $(x_k, p_k, t_k) \in \cG_k$, then there exists $L_{H_k, c_k}$, $w_k$ calibrated curves $\gamma_k: (-\infty, t_k]\to \T^d$ with $\gamma_k(t_k) =x_k$ and $\partial_p H_k(x_k, p_k, t_k) = \dot{\gamma}_n$. Then by the same argument, after restricting to a subsequence, $\gamma_n$ converges in $C^1_{loc}$ to $\gamma: (-\infty, t] \to \T^d$, and $\gamma$ is a $L_{H,c}$, $w$ calibrated curve. This implies $(t, x, p) \in \cG$. 

	For the set $\cI$, let us prove for each fixed $T$, we have 
	\[
		\limsup_{n\to \infty} \phi_{H_k}^{-T}\, \cG_k \subset \phi_H^{-T} \, \cG. 
	\]
	Indeed, for any $(x_k, p_k,t_k) \in 
	\phi_{H_k}^{-T} \cG_{c_k, H_k, w_k}$, there exists 
	$L_{H_k, c_k}$, $w_k$ calibrated curves 
	$\gamma:(-\infty, t_k +T] \to \T^d$ such that 
	$\gamma(t_k) = x_k$, $\partial_p H_k(x_k, p_k, t_k) = \dot{\gamma}_k(t_k)$. Then exactly the same argument 
	as before implies $\gamma_n$ accumulates to a $L_{H,c}$, 
	$u$ calibrated curve $\gamma:(-\infty, t+T)\to \T^d$, and implying $(t, x, p) \in \Phi_H^{-T}\cG$. We obtain 
	\[
		\limsup_{n \to \infty}\ \tcI_n \subset \tcI.
	\]

	Finally, since 
	\[
		\tilde{\cN}_{H_k}(c_k) = \bigcup_{w} \tcI_H(c_k, w), 
	\]
	where the union is over all $L_{H_k, c_k}$ weak KAM 
	solutions. For $(x_k, p_k, t_k) \in \tilde{\cN}_{H_k}(c_k)$, 
	there exists $u_k$ such that $(x_k, p_k, t_k) \in 
	\cG_{c_k, H_k, u_k}$. Since all $w_k$ are equi-continuous 
	and equi-bounded (see e.g. (A.3), \cite{Be}), 
	there exists a subsequence that converges to $u$ uniformly 
	on the interval $\T^d \times [0,K]$. Since all $w_k$'s are 
	periodic with period bounded by $1$, 
	this implies $w_k \to w$ on $\T^d \times \R$ as well. We 
	then apply the semi-continuity of $\cG$ sets to get 
	semi-continuity of $\tilde{\cN}$.
\end{proof}

The theory built in \cite{Be3} allows one to pass from semi-continuity of pseudographs to semi-continuity of Aubry set under a condition called the \emph{coincidence hypothesis}. A sufficient condition for this hypothesis is when the Aubry set has finitely many static classes. 
\begin{cor}\label{cor:aubry-unique-static}
Suppose $\tilde{\cA}_H(c)$ has at most finitely many static classes. Then if $H_n \in \bH(D)$ $C^2$ converges to $H$ over compact sets, and $c_n \to c$, we have 
\[
	\limsup_{n\to \infty}\ \tilde{\cA}_{H_n}(c_n) \subset \tilde{\cA}_{H}(c)
\]
as subsets of $\T^n \times \R^n \times \R$. 
\end{cor}
\begin{proof}
The proof  follows in the same way as the proof of Theorem~1 in \cite{Be3}. We also refer to \cite{KZ2017}, Section~6.2 where this is carried out in detail. 
\end{proof}

\subsection{Continuity of the barrier function}

In general,  the barrier function $h_{H,c}$ may be discontinuous with respect to $H$ and $c$. However,
the continuity properties hold in the particular case when the limiting
Aubry set contains only one static class.

\begin{prop}\label{cont-barrier}
Assume that a sequence $H_k \in \bH(D)$ converges to $H$ in $C^2$ over compact sets,
and $c_k \to c \in \R^n$. Assume that the projected Aubry set
$\cA_{H}(c)$ contains at a unique static class. 
Let $(x_k,s_k) \in \cA_{H_k}(c_k)$ with $(x_k,s_k) \to (x,s)$, then the barrier
functions $h_{H_k,c_k}(x_k, s_k; \cdot, \cdot)$ converges to
$h_{H,c}(x,s; \cdot, \cdot)$ uniformly.

Similarly, for $(y_k, t_k)\in \cA_{H_k}(c_k)$ and $(y_k, t_k)\to (y,t)$,
the barrier functions $h_{H_k,c_k}(\cdot, \cdot; y_k, t_k)$ converges to
$h_{H,c}(\cdot, \cdot; y,t)$ uniformly.
\end{prop}

\begin{prop}\label{unif-barrier-conv}
Assume that a sequence $H_n$ converges to $H$ in $C^2$ 
over compact sets,
$c_k \to c \in \R^n$ and the Aubry set $\cA_{H}(c)$
contains a unique static class.
\begin{enumerate}
	\item For any $(x,s)\in \cA_H(c)$ we have
	$$ \lim_{n\to \infty}\ 
	\sup_{(y,t)\in M\times \T}
	|h_{H_k,c_k}(x_k,s_k; y,t) - h_{H,c}(x,s;y,t)| =0 $$
	uniformly over  $(x_k,s_k)\in \cA_{H_k}(c_k)$ and $(y,t)\in M\times \R$.

	\item For any  $l_k\in \partial_y^+ h_{H_k,c_k}(x_k,s_k;y,t)$ and $l_k \to l$,
	we have $l\in \partial^+_y h_{H,c}(x,s;y,t)$. Moreover, the convergence
	is uniform in the sense that
	$$ \lim_{n\to \infty} \inf_{l_k \in \partial^+_y h_{H_k,c_k}(x_k,s_k;y,t)}
	d(l_k, \partial_y^+ h_{H,c}(x,s;y,t)) =0  $$
	uniformly in $(x,s)\in \cA_{H}(c)$, $(x_k,s_k)\in \cA_{H_k}(c_k)$.
\end{enumerate}
\end{prop}

The following statement follows easily from the representation formula (Proposition~\ref{prop:rep-formula}). 
\begin{lem}\label{unique-weak-kam}
Assume that $\cA_{H}(c)$ has a unique static class.
Let $(x_1, t_1)\in \cA_{H}(c)$, then any weak KAM solution differs from
$h_{H,c}(x_1, t_1; \cdot, \cdot)$ by a constant.
\end{lem}

\begin{proof}[Proof of Proposition~\ref{cont-barrier}]

	We prove the second statement. By Proposition~\ref{prop:action-concave},
	all functions $h_{H_k,c_k}(x_k, s_k; \cdot, \cdot)$ are uniformly
	semi-concave, and hence equi-continuous. By Arzela-Ascoli,
	any subsequence contains a uniformly convergent subsequence, whose limit is
	$$ h_{H,c}(x,s; \cdot, \cdot) + C$$
	due to Lemma~\ref{lem:conv-weak-KAM} and Lemma~\ref{unique-weak-kam}.
	Moreover,
	$$
	h_{H,c}(x_k,s_k; x,s) \to h_{H,c}(x,s; x,s)=0,
	$$
	so $C=0$.
	It follows that  $h_{H_k,c_k}(x_k, s_k; \cdot, \cdot)$ converges to
	$h_{H,c}(x,s; \cdot, \cdot)$ uniformly.

	Statement 1 follows from the definition of the projected Aubry set
	$$ \cA_{H}(c)=\{(x,s)\in M\times \T: h_{H,c}(x,s;x,s)=0\}$$
	and statement 2.
\end{proof}

\begin{proof}[Proof of Proposition~\ref{unif-barrier-conv}]
	\textit{Part 1.}  We argue by contradiction. Assume that there exist
	$\delta>0$, and by restricting to a subsequence,
	$$ \inf_{C\in \R}\ \sup_{(y,t)}\
	|h_{H_k,c_k}(x_k,s_k, y,t) - h_{H,c}(x,s, y,t)-C| > \delta. $$
	By compactness, and by restricting to a subsequence again, we may assume
	that $(x_k,s_k)\to (x^*,s^*)$, $(y_k, t_k) \to (y,t)$.
Using Proposition~\ref{cont-barrier}, take limit as $n\to \infty$, we have
$$  \sup_{(y,t)}
|h_{H,c}(x^*,s^*, y,t) - h_{H,c}(x,s, y,t)| > \delta. $$
By Lemma~\ref{unique-weak-kam}, the left hand side is $0$, which is
a contradiction.

\textit{Part 2.} 
$h_{H_k,c_k}(x_k,s_k, \cdot,t)$ converges to $h_{H,c}(x,s, \cdot,t)$
uniformly. Convergence of super-differentials follows directly from
Proposition~\ref{prop:action-concave}. It suffices to prove uniformity.
Assume, by contradiction, that by restricting to a subsequence,
we have $(x_k,s_k)\to (x,s) \in \cA_{H}(c)$,
$l_k \in \partial^+_y h_{H_k,c_k}(x_k,s_k, y,t)$
and $(x,s) \in \cA_{H}(c)$ such that
$$\lim_{n\to \infty} l_k \notin \partial^+_y h_{H,c}(x,s,  y, t). $$
By Proposition~\ref{prop:action-concave}, $l_n \to l \in
\partial^+_y h_{H,c}(x^*,s^*,  y, t)$, but we also have
$\partial^+_y h_{H,c}(x,s,  y, t) = \partial^+_y h_{H,c}(x^*,s^*,  y, t)$
since the functions differ by a constant using Lemma~\ref{unique-weak-kam}.
This is a contradiction.
\end{proof}

\subsection{Lipschitz estimates for nearly integrable systems}

In this section we record uniform estimates for the nearly integrable system 
\[
	H_\epsilon = H_0(p) + \epsilon H_1(\theta, p, t), \quad (\theta, p, t)\in \T^n \times \R^n \times \T,
\]
with the assumptions $H_0 \in \bH(D)$, $\|H_1\|_{C^2} \le 1$. 

\begin{prop}[Proposition 4.3, \cite{BKZ}]\label{prop:semi-concave-near-int}
For $H_\epsilon$ as given,  for any $c \in \R^n$, any $L_{H_\epsilon, c}$ weak KAM solution $u(x,t)$ is $6 D \sqrt{\epsilon}-$semi-concave and $6D \sqrt{d \epsilon}-$Lipschitz in $x$. 
\end{prop}

\begin{cor}\label{cor:min-near-int}
For any $L_{H_\epsilon, c}$ weak KAM solution $u: \T^n\times \T \to \R$, let $\gamma(-\infty, t_0] \to \T^n$ be a calibrated curve. Then for any $t \in (-\infty, t_0)$, 
\[
	p(t) : = \partial_v L_{H_\epsilon}(\gamma(t), \dot{\gamma}(t), t) 
\]
satisfies $\|p(t) - c\| \le 6D\sqrt{n\epsilon}$. 

Moreover, any $(x,p,t)\in \tcN_{H_\epsilon}(c)$ satisfies the same estimate $\|p - c\| \le 6D \sqrt{n\epsilon}$. 
\end{cor}
\begin{proof}
	$p(t)-c$ is a super gradient of $u(\cdot, t)$ at $\gamma(t)$. The conclusion follows from the fact that $u(\cdot, t)$ is $6D\sqrt{n\epsilon}$--Lipschitz. The statement for $\tcN$ follows from the fact that any $(x,p,t) \in \cG \subset \cN$ is the end point of a calibrated curve. 
\end{proof}

\begin{cor}
\label{cor:near-int-rho}
Let $\mu$ be any $c$-minimal measure of $H_\epsilon$, then there is $C>0$ depending only on $H_0$ such that 
\[
\|\rho(\mu) - \nabla H_0(c)\| \le C \sqrt{\epsilon}. 
\]
\end{cor}
\begin{proof}
Recall that $\mu$ is a measure on $\T^n \times \R^n\times \R$, invariant under the Euler-Lagrange flow. By considering the Legendre transform $\bL$, we obtain 
\[
\rho(\mu) = \int v d\mu(\theta, v, t) = \int \partial_p H_\epsilon(\theta, p, t) d\bL_*\mu(\theta, p, t). 
\]
By Corollary~\ref{cor:min-near-int}, we have $\partial_p H_\epsilon(\theta, p, t) = \nabla H_0(c) + O(\sqrt{\epsilon})$ for all 
$$(\theta,p, t) \in \supp \bL_*\mu \subset \tcA_{H_\epsilon(}c),
$$ 
therefore, $\rho(\mu) = \nabla H_0(c) + O(\sqrt{\epsilon})$. 
\end{proof}

\subsection{Estimates for nearly autonomous systems}

 The goal of this section is to derive a special Lipshitz estimate of weak KAM solutions for  perturbations of autonomous systems. More precisely, consider
\[
	H_\epsilon(x, p, t) = H_1(x, p) + \epsilon H_2(x, p, t). 
\]
Assume that $H_1 \in \bH(D/2)$ and $\|H_2\|_{C^2} = 1$. Then for $\epsilon$ small enough all $H_\epsilon \in \bH(D)$. Let $L_\epsilon$ denote the associated Lagrangian. 

We first state an estimate for the alpha function. 
\begin{lem}
	\label{lem:pert-alpha}
	There is $C>0$ depending only on $\|H_1\|_{C^2}$ and $D$, such that 
	\[
		\|\alpha_{H_\epsilon}(c) - \alpha_{H_1}(c)\| < C \epsilon. 
	\]
\end{lem}
\begin{proof}
	Let $C$ always denote a generic constant depending on $\|H_1\|_{C^2}$. Let $L_1$ be the Lagrangian for $H_1$, then  $\|L_\epsilon- L_1\|_{C^0} \le C \epsilon $. As a result, the functionals
	\[
		\int (L_\epsilon - c \cdot v) d\mu, \quad \int (L - c \cdot v) d\mu,
	\]
	defined on the space of closed probability measures, differ by at most $C \epsilon$. The lemma follows immediately using the closed measure version for the definition of the alpha function, see Section \ref{sec:intro-weak-kam-solution}. 
\end{proof}

Let $u(x, t)$ be a weak KAM solution to $L_\epsilon - c\cdot v$, and $w(x, t)$ a forward weak KAM solution conjugate to $u$.
Recall that 
\[
	\cI_{u, w}(c) =  \argmin_{(x, t)} (u - w), \quad \tcI_{u, w}(c) = \{(x, t, c +  d_x u(x, t)):\ (x, t) \in \cI_{u, w}\}. 
\]

\begin{thm}\label{thm:energy-lip}
Let $u, w$ be conjugate pair of weak KAM solutions to $L_\epsilon - c \cdot v$, then there is $C$ depending only $D$ and $\|c\|$  such that 
\[
	|H_\epsilon (x_1, p_1, t_1) - H_\epsilon (x_2, p_2, t_2)| \le 
	C \sqrt{\epsilon}\, \|(x_2 - x_1, t_2 - t_1)\|, \quad (x_i, p_i, t_i) \in \tcI_{u, w}. 
\]
In particular, the above estimates holds on the Aubry set $\tcA_{L_\epsilon}(c)$. 
\end{thm}

The proof relies on semi-concavity of weak KAM solution, following Fathi (\cite{Fa}). However to get an improved estimate we need a notion of semi-concavity that is ``stronger'' in the $t$ direction.  

Let $\Omega \subset \R^n$ be an open convex set. Let $A$ be a symmetric $n\times n$ matrix. We say that $f: \Omega \to \R$ 
is $A$-semi-concave if for each $x \in \R^n$, there is $l_x \in \R^n$ such that 
\[
	f(y) - f(x) -  l_x \cdot (y - x) \le \frac12 A (y-x)^2, \quad x, y \in \Omega
\]
where $Ax^2$ denotes $Ax \cdot x$. These definitions generalizes the standard semi-concavity, as $A$-semi-concave functions are $\frac12 \|A\|$-semi-concave.  We say $f$ is $A$-semi-convex if $-f$ is $A$-semi-concave. The following lemma follows from a direct computation.
\begin{lem}
	$f$ is $A$-semi-concave if and only if $f_A(x) = f(x) - \frac12 Ax^2$ is concave. 
\end{lem}

The following lemma is proved in \cite{Zha2017}.
\begin{lem}[See \cite{Zha2017}, Lemma 3.2]
	\label{lem:aniso-lipshitz}
Suppose $f: \R^n \to \R$ is $B$-semi-concave and $g: \R^n \to \R$ is $(-A)$-semi-convex, and $S = B - A$ is positive definite. Suppose $f(x)\ge g(x)$ and $M$ is the set on which $f - g$ reaches its minimum. 

Then for all $x_1, x_2 \in M$, we have 
\[
	\|df(x_2) - df(x_1) - \frac12(A + B)(x_2 - x_1)\|_{S^{-1}} \le \frac12 \|x_2 - x_1\|_S, 
\]
where $\|x\|_S = \sqrt{S x^2}$. 
\end{lem}

\subsection{Semi-concavity of viscosity solutions}

Recall that the action function $A_L(x, t, y, s)$ is the minimal Lagrangian action of curves with end points $\gamma(t) = x$, $\gamma(s) = y$. The goal of this section is to prove:
\begin{prop}\label{prop:near-auto-semi-concave}
For the Hamiltonian $H_\epsilon$, there is a constant $C>0$ depending only on $D$ such that if $1/(2\sqrt{\epsilon}) \le s - t \le 1/\sqrt{\epsilon}$, the action function $A(x, t, y, s)$ is $S_\epsilon$-semi-concave in $(x, t)$ and $(y,s)$, where 
\[
	S_\epsilon = C\bmat{ \Id_{n \times n} & 0 \\ 0 & \sqrt{\epsilon}}. 
\]
\end{prop}

Our proposition follows directly from the following technical lemma by choosing $1/(2\sqrt{\epsilon}) \le  T  \le 1/\sqrt{\epsilon}$. 
\begin{lem}\label{lem:lin-drift}
There is a constant $C$ depending $D$ such that if $\gamma: [t_0, t_1] \to \T^n$ is an extremal curve with $T = t_1 - t_0>1$, then for $h \in \R^n$, $\Delta T \in [-T/2, T/2]$,  $p(t) = \partial_v L_\epsilon(\gamma(t), \dot{\gamma}(t), t)$,  
\[
	\begin{aligned}
		& 	A(\gamma(t_0), t_0, \gamma(t_1) + h, t_1 + \Delta T) - A(\gamma(t_0), t_0, \gamma(t_1), t_1) \\
		&  \le p(t_1) \cdot h + H_\epsilon(\gamma(t_1), p(t_1), t_1) \cdot (\Delta T) +  C  \|h\|^2 + C \left( \frac{1}{T} + \epsilon T \right) (\Delta T)^2.   
	\end{aligned}
\]
\end{lem}
\begin{proof}
	We omit the subscript in $L_\epsilon$ within the proof. Moreover, by considering the Lagrangian $L(\cdot, \cdot, \cdot + t)$, it suffices to consider $t = 0$, $s = T$. 

	Let $\gamma:[0, T] \to \T^n$ be a extremal curve, write $\lambda = \Delta T/ T$ and define
	\[
		\gamma_h^\lambda: [0, T + \Delta T] \to \T^n,
	\]
	\[
		\gamma_h^\lambda(t) = 
		\begin{cases}
			\gamma\left( \frac{t}{1+ \lambda} \right), & t \in [0, (T-1)(1+\lambda)]; \\
			\gamma\left( \frac{t}{1+\lambda} \right) + \left( \frac{t}{1+\lambda} - T + 1 \right) h , & 
			t \in [(T-1)(1 + \lambda), T(1+\lambda)].
		\end{cases}
	\]
	The curve is obtained by adding a linear drift in $x$ on the time interval $[T-1, T]$, then reparametrize time to the interval $[0, T + \Delta T]$, see Figure~\ref{fig:lin-drift}. 
	\begin{figure}[t]
		\centering
		\includegraphics[width=3.5in]{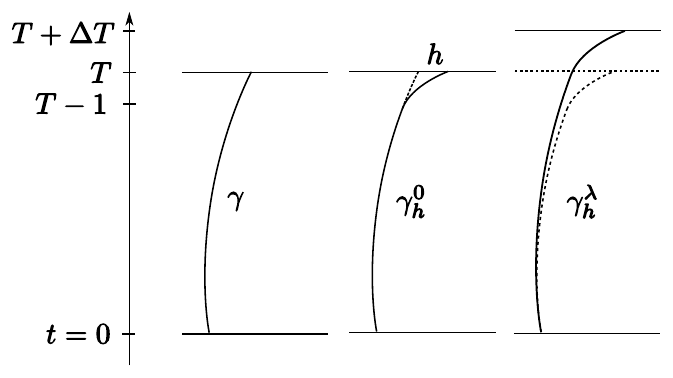}
		\caption{Defining the curve $\gamma_h^\lambda$} \label{fig:lin-drift}
	\end{figure}

	Let $C$ denote an unspecified constant. Note that $\|\dot{\gamma}\|\le C$. 
	\[
		\begin{aligned}
			& 	\bA_L(\gamma_h^\lambda)	= \int_0^{(1+\lambda)T} L(\gamma_h^\lambda, \dot\gamma_h^\lambda, t) dt  \\
			& = (1+ \lambda) \int_0^{T-1} L\left(  \gamma(s), \frac{1}{1+\lambda}\dot\gamma(s) , (1+ \lambda )s\right) ds \\
			& \qquad  +  (1+ \lambda)\int_{T-1}^T  L\left( \gamma(s) + (s - T + 1)h , \frac{1}{1+\lambda} \left( \dot{\gamma}(s) + h\right) , (1+ \lambda)s \right) ds, 
		\end{aligned}
	\]
	then
	\[
		\begin{aligned}
			& \bA_{L}(\gamma_h^\lambda)	 \le  (1+ \lambda)\int_0^T L(\gamma, \dot\gamma, s)ds + 
			\int_0^{T-1} \left( - \lambda\partial_v L \cdot \dot\gamma + (1 + \lambda)  \partial_t L \cdot \lambda s \right) ds  + \\
			& \quad (1+ \lambda) \int_{T-1}^T \left( \partial_x L \cdot (s-T +1) h + \frac{1}{1+\lambda}\partial_vL \cdot ( - \lambda \dot\gamma + h ) + \partial_t L \cdot \lambda s
			\right) ds \\
			& \quad + C T \left(  \|\partial^2_{vv}L\|  + \|\partial^2_{vt}L\|  T + \|\partial^2_{tt}L\| T^2\right) \lambda^2  + C \|\partial^2_{(x, v)^2} L\| \left( \|h\|^2  + |\lambda| \|h\| + \lambda^2\right)   \\
			& \quad + C \left( \|\partial_{xt}L\| |\lambda| T \|h\|  + \|\partial_{vt} L \| (|\lambda| + \|h\|)\lambda T + \|\partial_{tt}L\| \lambda^2 T^2\right). 
		\end{aligned}
	\]
	Using $\|\partial^2 L\| \le C$, $\|\partial^2_{(x, v, t) t}L\| \le C \epsilon$, and plug in $\Delta T = \lambda T$, we have
	\[
		\begin{aligned}
			& \bA_L(\gamma_h^\lambda) - \bA_L(\gamma) \le 
			\lambda \int_0^T \left(  L - \partial_v L \cdot \dot{\gamma} +  \partial_t L \right) ds + 
			\int_0^1\left(  \partial_x L \cdot sh + \partial_v L \cdot h \right)ds \\
			& \quad + C\left( \frac{1}{T} + \epsilon T \right) (\Delta T)^2 
			+ C\left( \|h\|^2 + \frac{1}{T} |\Delta T| \|h\| + \frac{1}{T^2} (\Delta T)^2 \right) \\
			& \quad+ C\epsilon \left( |\Delta T| \|h\| + (\Delta T)^2 \right)
		\end{aligned}
	\]
	Using the Euler-Lagrange equation, we have for $p(t) = \partial_v L(\gamma, \dot\gamma, t)$, 
	\[
		\frac{d}{dt}(- tH) =  \frac{d}{dt}(t(L - \partial_v L \cdot \dot{\gamma})) = (L - \partial_v L \cdot \dot{\gamma}) + t \partial_t L, 
	\]
	\[
		\frac{d}{dt}\left( p \cdot th \right) = \frac{d}{dt} \left(  \partial_v L \cdot th \right) = \partial_v L \cdot h + \partial_x L \cdot th , 
	\]
	we get 
	\[
		\begin{aligned}
			& \bA_L(\gamma_h^\lambda) - \bA_L(\gamma) \le \lambda (-tH)\Bigr|_0^T + \left( p(T-1+t) \cdot h t\right)\Bigr|_0^1 \\
			& \quad + C \|h^2\| + C\left(  \frac{1}{T} + \epsilon  \right)|\Delta T| \|h\| + C\left( \frac{1}{T} + \epsilon T \right) (\Delta T)^2 \\
			& \le H(\gamma(T), p(T), T) (\Delta T) + p(T) \cdot h 
			+ 2C \|h\|^2  + 2C\left( \frac{1}{T} + \epsilon T \right) (\Delta T)^2 . 
		\end{aligned}
	\]
\end{proof}

Let $u(x, t)$ be a weak KAM solution to $L_\epsilon - c\cdot v$, and $w(x, t)$ a forward weak KAM solution conjugate to $u$. 
Then Proposition~\ref{prop:near-auto-semi-concave} implies that $u$ is $S_\epsilon$-semi-concave, while $w$ is $(-S_\epsilon)$-semi-convex. 
\begin{proof}[Proof of Theorem~\ref{thm:energy-lip}]
	Note that for each $(x, p, t) \in \tcI_{u, w}$, $- H(x, p, t) + \alpha_{L_\epsilon}(c) = \partial_t u(x, t)$. Apply Proposition~\ref{prop:near-auto-semi-concave} and Lemma~\ref{lem:aniso-lipshitz}, we obtain 
	\[
		\left\| \bmat{p_2 - p_1 \\ H(x_2, p_2, t_2) - H(x_1, p_1, t_1) } - \bmat{C(x_2 - x_1) \\ C\sqrt{\epsilon}(t_2 - t_1)}   \right\|_{S_\epsilon^{-1}}
		\le \left\|  \bmat{x_2 - x_1 \\ t_2 - t_1} \right\|_{S_\epsilon}, 
	\]
	therefore 
	\[
		\begin{aligned}
			& \frac{1}{\sqrt{\epsilon}} |H(x_2, p_2, t_2) - H(x_1, p_1, t_1) - C\sqrt{\epsilon}(t_2 - t_1)|  \\
			& \le 		C \|x_2 - x_1\| + C\|t_2 - t_1\|,
		\end{aligned}
	\]
	the proposition follows. 
\end{proof}

\section{Cohomology of Aubry-Mather type}

\subsection{Aubry-Mather type and diffusion mechanisms}

Let 
\[
	H: \T^n \times \R^n \times \T \to \R
\]
be a $C^r$ Tonelli Hamiltonian with $r \ge 2$ contained in the family $\bH(D)$. 

\begin{defn}\label{defn:AM}
We say that the pair $(H_*, c_*)$ is of \emph{Aubry-Mather type} if it satisfies 
the following conditions:
\begin{enumerate}
	\item  There is an embedding 
	\[
		\chi: \,  \T \times (-1, 1) \to \T^n \times \R^n,
	\]
	such that $\cC(H) = \chi (\T \times (-1, 1))$ is a normally hyperbolic weakly invariant cylinder 
	under the time-$1$-map $\Phi^{H}$. We require $\cC(H)$ is \emph{symplectic}, 
	i.e. the restriction of the symplectic form $\omega$ to $\overline{\cC}$ is non-degenerate. 

	\item There is $\sigma > 0$, such that the following holds for all 
	$$
	\cV_{\sigma}(H_*):=\{\|H - H_*\|_{C^2}\} < \sigma \quad , \quad 
	B_{\sigma}(c_*):=\{\|c - c_*\| < \sigma\}.
	$$ 
	(We then say that the property below holds robustly at $(H_*,c_*)$.)
	\begin{enumerate}
		\item  
		The discrete Aubry set $\tcA^0(c) \subset \cC(H)$. Moreover, the pull back of the Aubry set 
		is contained in a Lipschitz graph, namely, the map
		\[
			\pi_x: \, \chi^{-1} \tcA^0_H(c) \subset \T\times (-1,1) \to \T
		\]
		is bi-Lipschitz. 
		\item If $\chi^{-1} \tcA^0_H(c)$ projects onto $\T$, then there is a neighborhood $V^0$ of 
		$\cA^0_H(c)$ such that the strong unstable manifold $W^u(\tcA^0_H(c))\cap \pi_\theta^{-1} V^0$ 
		is a Lipschitz graph over the $\theta$ component. 
	\end{enumerate}
\end{enumerate}
\end{defn}

If $(H_*,c_*)$ is of Aubry-Mather type, let $h \in \Z^n \simeq H_1(\T^n, \Z)$ be the homology 
class of the curve $\chi(\T \times \{0\})$.

\begin{rmk}\label{rmk:AM-hyp-fixed-pt}
The definition of Aubry-Mather type includes a much simpler case, that is when the Aubry set $\tcA_{H}(c)$ is a hyperbolic periodic orbit, still contained in a NHIC $\cC$. The condition 2(a) is satisfied since $\tcA_{H}^0(c)$ is discrete. Condition (b) is vacuous since the $\chi^{-1} \tcA_{H}^0(c)$ never projects onto $\T$. This holds, in particular, at double resonance when the energy is critical. 
\end{rmk}

\begin{defn} \label{def:bif}
We say that the pair $(H_*, c_*)$ is of \emph{bifurcation Aubry-Mather type} if there exist
$\sigma >0$ and  open sets $V_1, V_2 \subset \T^n$ with 
$\overline{V_1} \cap \overline{V_2} = \emptyset$, and a smooth bump function 
\[
	f: \T^n \to [0, 1], \quad f|_{V_1} = 0, \quad f|_{V_2} = 1,
\]
such that for all $c \in \overline{B_{\sigma}(c_*)}$ and $H \in \cV_{\sigma}(H_*)$: 
\begin{enumerate}
	\item Each of the Hamiltonians $H^1 = H - f$ and $H^2 = H - (1-f)$ satisfies 
	\[
		\cA_{H^1}(c) \subset V_1 \times \T , \quad  \cA_{H^2}(c) \subset V_2 \times \T. 
	\]
	It follows that $\tcA_{H^1}(c)$, $\tcA_{H^2}(c)$ are both invariant sets of $H$, called 
	\emph{the local Aubry sets.}
	\item $(H^1, c_*)$ (resp. $(H^2, c_*)$) are of Aubry-Mather type, with the invariant cylinders 
	$\cC_1$ and $\cC_2$. Moreover, they have the same homology class $h$. 
	\item The Aubry set 
	\[
		\tcA_H(c) \subset \tcA_{H^1}(c) \cup \tcA_{H^2}(c). 
	\]
\end{enumerate}
\end{defn}
According to item (3) above, the Aubry set $\tcA_H(c)$ may be contained in one of the local component, or both. 

There is another type of bifurcation which happens at double resonance, when the homology $h$ is simple non-critical (i.e. the shortest geodesic loop does not contain the saddle fixed point). We call this an asymmetric bifurcation.
\begin{defn}
\label{def:asym-bif}
We say that the pair $(H_*, c_*)$ is of \emph{asymmetric bifurcation type} if there exist $\sigma>0$ and open sets $V_1, V_2 \subset \T^n$ with 
$\overline{V_1} \cap \overline{V_2} = \emptyset$, and a smooth bump function 
\[
	f: \T^n \to [0, 1], \quad f|_{V_1} = 0, \quad f|_{V_2} = 1,
\]
such that for all $c \in \overline{B_{\sigma}(c_*)}$ and $H \in \cV_{\sigma}(H_*)$: 
\begin{enumerate}
	\item Item 1 of Definition~\ref{def:bif} holds. 
	\item $(H^1, c_*)$ is of Aubry-Mather type, with the invariant cylinder 
	$\cC_1$. The Aubry set $\tcA_{H^2}(c_*)$ is a single hyperbolic periodic orbit. 
	\item The Aubry set 
	\[
		\tcA_H(c) \subset \tcA_{H^1}(c) \cup \tcA_{H^2}(c). 
	\]
\end{enumerate}
\end{defn}
\vskip 0.1in 

\begin{thm}\label{thm:AM-non-deg}
Suppose a pair $(H_*, c_*)$ is of Aubry-Mather type, and let $\Gamma \ni c_*$ be a smooth curve in $\R^2$. Then there are $\sigma_1, \sigma_2 > 0$ 
such that for all $c \in \Gamma_1:= \overline{B_{\sigma_1}(c_*)}\cap \Gamma$, we have the following dichotomy for a residual subset 
of $H \in \cV_\sigma(H_*)$: 
\begin{enumerate}
	\item  The projected Ma\~ne set $\cN^0_{H}(c)$ is contractible as a subset of $\T^n$;
	\item  There is a double covering map $\Xi$ such that the set 
	\[
		\tcN^0_{H \circ \Xi}(\xi^*c) \setminus \Xi^{-1} \tcN^0_{H}(c)
	\]
	is totally disconnected. 
\end{enumerate}
\end{thm}

\begin{thm}
	\label{thm:bifur-non-deg}
	Suppose $(H_*, c_*)$ is of either
	\begin{itemize}
	 \item bifurcation  AM type, or 
	 \item  asymmetric bifurcation type,
	\end{itemize}
	and $c_* \in \Gamma$, where $\Gamma$ is a smooth curve in $\R^2$, 
	then there is $\sigma>0$ and an open and dense subset $\cR \subset \cV_\sigma(H_*)$ such that for
	$c \in B_\sigma(c_*) \cap \Gamma$ and $H \in \cR$, either:
	\begin{enumerate}
		\item $\tcA^0_H(c)$ has a unique static class, and one of the following holds.
		\begin{enumerate}
			\item  The projected Ma\~ne set $\cN^0_{H}(c)$ is contractible as a subset of $\T^n$;
			\item  There is a double covering map $\Xi$ such that the set 
			\[
				\tcN^0_{H \circ \Xi}(\xi^*c) \setminus  \Xi^{-1}  \tcN^0_{H}(c)
			\]
			is totally disconnected.  
		\end{enumerate}
		\item Or $\tcA^0_H(c)$ has two static classes, and \[
			\tcN^0_H(c) \setminus \tcA^0_H(c)
		\]
		is a discrete set. 
	\end{enumerate}
\end{thm}

Analogous result in the a priori unstable setting is due to Mather and 
Cheng-Yan (\cite{CY1}, \cite{CY2}). Our definition is more general and applies, 
as will be seen, to both single and double resonant settings. The prove is base on 
the result of  \cite{BKZ}, which applies to our setting with appropriate changes. 
Here we describe the changes needed for the proof in \cite{BKZ} to apply.

If $(c_*, H_*)$ is of Aubry-Mather type, let $h \in \Z^n \simeq H_1(\T^n, \Z)$ be its 
homology class. Then any minimal measure contained in $\tcA^0(c)$ must have 
rotation vector $\lambda h$, $\lambda \in \R$. Moreover, similar to the case of twist map, 
all such measures has the same rotation vector $\lambda h$. We say the rotation vector 
is rational/irrational if $\lambda$ is rational/irrational. 

The following proposition is a consequence of the Hamiltonian Kupka-Smale theorem, 
see for example \cite{RiRu11}. 
\begin{prop} \label{prop:kupta-smale}
There are $\sigma_1, \sigma_2>0$, and a residual subset $\cR_1$ of the ball $\cV_{\sigma_2}(H_*)$,
such that  for all $c \in \Gamma_1 := \overline{B_{\sigma_1}(c^*)} \cap \Gamma$,
$H \in \cR_1$, if $c$ has rational rotation vector, then $\pi \chi^{-1} \tcA^0(c)\ne \T $. 
\end{prop}

Let $1 \le j \le n$ be such that $e_j \nparallel h$, define the covering map $\xi: \T^n \to \T^n$ by
\[
	\xi(\theta_1, \cdots, \theta_n) = (\theta_1, \cdots, 2\theta_j, \cdots, \theta_n).
\]
\begin{itemize}
	\item For $H \in \cR_1$, define
	\begin{equation}
		\label{eq:gamma-star}
		\Gamma_*(H) = \{ c \in \Gamma_1: \,  \pi\chi^{-1}(\cN^0_H(c)) = \T\}. 
	\end{equation}
	Proposition~\ref{prop:kupta-smale} implies that  each $c \in \Gamma_*(H)$ has an irrational rotation vector, and the Aubry set has a unique static class, 
	and $\cA^0_H(c) = \cN^0_H(c)$. Each $\Gamma_*(H)$ is compact due to upper semi-continuity of the Ma\~ne set. 
	\item For $H \in \cR_1$ and $c \in \Gamma_*(H)$, the lifted Aubry set $\tcA^0_{H\circ \Xi}(\xi^*c) =  \Xi^{-1} \tcA_H(c)$ has two static classes, denote them $\tcS_1, \tcS_2$ (and projections $\cS_1, \cS_2$). 

	We have the decomposition 
	\[
		\tcN^0_{H \circ \Xi}(\xi^*c) = \tcS_1 \cup \tcS_2 \cup \tcH_{12} \cup \tcH_{21},
	\]
	where $\tcH_{ij}$ consists of heteroclinic orbits from $\tcS_i$ to $\tcS_j$. We will use 
	the notation $\tcS_i(H,c), \tcH(H, c)$ to show dependence on the pair $(H, c)$, and 
	$\cH_{ij}$ for projection of $\tcH_{ij}$. 

	\item For $H \in \cR_1$ and $c \in \Gamma_*(H)$, consider the (discrete) Peierl's barrier (see \eqref{eq:dis-barrier})
	\[
		h(\zeta_1, \cdot), \, h(\zeta_2, \cdot), \quad
		h(\cdot, \zeta_1), \, h(\cdot, \zeta_2),\quad
		\zeta_i \in \cS_i, \, i = 1, 2, 
	\]
	where $h = h_{H\circ \Xi, \, \xi^*c}$. These functions are independent of the choice of $\zeta_i$ 
	except for an additive constant. We define
	\[
		b^-_{H, c}(\theta) = h(\zeta_1, \theta) + h(\theta, \zeta_2) - h(\zeta_1, \zeta_2)
	\]
	and $b^+_{H,c}$ by switching $\zeta_1, \zeta_2$. The functions $b^\pm_{H,c}$ are non-negative 
	and vanish on $\cH_{12} \cup \cS_1 \cup \cS_2$ and $\cH_{21} \cup \cS_1 \cup \cS_2$, respectively. 

	\item Consider small neighborhoods $V_1, V_2$ of $\cS_1(H_*, c_*), \cS_2(H_*, c_*)$, and define 
	$K = \T^n \setminus (V_1 \cup V_2)$. By semi-continuity of the Aubry set, for sufficiently small 
	$\sigma_1, \sigma_2$, $K$ is disjoint from $\cS_i(H, c)$ for all $c \in \Gamma_1$ and $H \in \cR_1$. 
	Moreover, $\pi^{-1}K$ intersect every orbit of $\tcH_{12}(H,c)$ and $\tcH_{12}(H, c)$. 
\end{itemize}

\begin{lem}[\cite{BKZ}, Lemma 5.2] \label{lem:K}
For each $(H, c) \in \cR_1 \times \Gamma_1$, the set 
$$
\tcN_{H \circ \Xi}(\xi^*c)  \setminus \Xi^{-1} \tcN_H(c)
$$ 
is totally disconnected if and only if 
\[
	\cN_{H \circ \Xi}(\xi^* c) \cap K = (\cH_{12} \cup \cH_{21}) \cap K
\]
is totally disconnected. 
\end{lem}

We will show that the set of $H \in \cR_1$ with the following property contains a dense $G_\delta$ set: for each $c \in \Gamma_*(H)$, $\cN_{H \circ \Xi}(\xi^* c) \cap K$ is 
totally disconnected. The following lemma implies the $G_\delta$ property. 

\begin{lem}[\cite{BKZ}, Lemma 5.3] \label{lem:G-dt}
Let $K \subset \T^n$ be compact, then the set of $H \in \cR_1$ such that for all $c \in \Gamma_*(N)$, 
the set $\cN_{H \circ \Xi}(\xi^*c)$ is totally disconnected, is a $G_\delta$ set. 
\end{lem}

The following proposition allows local perturbations of $b^\pm_{H, c}$ \emph{simultaneously} for 
all $c$'s in a small ball. Let $B_\sigma(x)$ denote the ball of radius $\sigma$ at $x$ in a metric space. 
\begin{prop}[\cite{BKZ}, Proposition~5.2]\label{prop:perb-b}
Let $H_* \in \cR_1$, $c_* \in \Gamma_*(H)$, and $K \cap \cA_{H \circ \Xi}(\xi^*c) = \emptyset$.  
Then there is $\sigma>0$ such that for all
\[
	H \in \cR_1 \cap B_\sigma(H_*), \quad
	\theta_0 \in K \cap \cH_{12}(H_*, c_*), \quad
	\varphi \in C_c^r(B_\sigma(\theta_0)) \text{ with }
	\|\varphi\|_{C^r} < \sigma,
\]
there is a Hamiltonian $H_\varphi$ such that:
\begin{enumerate}
	\item For all $c \in B_\sigma(c_*)$, the Aubry sets $\tcA_{H_\varphi \circ \Xi}(\xi^*c)$ 
	coincides with $\tcA_{H \circ \Xi}(\xi^*c)$ with the same static classes. In particular, 
	$B_\sigma(c_*) \cap \Gamma_*(H) = B_\sigma(c_*) \cap \Gamma_*(H_\varphi)$. 
	\item  For all $c \in B_\sigma(c_0)\cap \Gamma_*(H)$, there exists a constant $e\in \R$ such that 
	\begin{equation}
		\label{eq:b-plus-pert}
		b_{H_\varphi, c}^+ (\theta) = b_{H,c}^+(\theta) + \varphi(\theta)+e, \quad \theta \in B_\sigma(\theta_0).
	\end{equation}	
	The same holds for $\theta_0 \in K \cap \cH_{21}(H_*, c_*)$, with $b^+$ replaced 
	with $b^-$ in \eqref{eq:b-plus-pert}. 
	Moreover, for each $H\in \cR_1\cap B_{\sigma}(H_*)$, $\|H_{\varphi}-H\|_{C^r}\to 0$ 
	when $\|\varphi\|_{C^r} \to 0$.
\end{enumerate}
\end{prop}

We will use Proposition~\ref{prop:perb-b} to locally perturb the functions $b^+_{H, c}$, therefore, 
perturbing 
$$
\cH_{12}(H, c) \cap B_\sigma(\theta_0) = \argmin\ b^+_{H, c} \cap B_\sigma(\theta_0).
$$ 
Similarly for $b^-$. However, as observed by Mather and Cheng-Yan, this requires additional information 
on how $b^\pm_{H, c}$ depends on $c$. 

\begin{prop}[Section~\ref{sec:reg-barrier}] \label{prop:barrier-holder}
There is $\beta> 0$ such that  each $H \in \cR_1$, the maps $c \mapsto b^\pm_{H, c}$
from $\Gamma_*(H)$ to $C^0(\T^n, \R)$  are $\beta$-H\"older. 
\end{prop}

As a result, the set $\{b^\pm_{H,c}: \, c \in \Gamma^*(N)\}$ has Hausdorff dimension at most 
$1/\beta$ in $C^0(\T^n, \R)$. The following lemma allows us to take advantage of this fact.
\begin{lem}[\cite{BKZ}, Lemma 5.6] \label{lem:disc}
Let $\cF\subset C^0([-1,1]^n,\R)$ be a compact set of finite Hausdorff dimension.
The following property is satisfied on a residual set of functions $\varphi\in C^r(\R^n,\R)$ 
(with the uniform $C^r$ norm):

For each $f\in \cF$, the set of minima of the function $f+\varphi$ on $[-1,1]^n$ is 
totally disconnected.

As a consequence, for each open neighborhood $\Omega$ of $[-1,1]^n$ in $\R^n$, 
there exists arbitrarily $C^r$-small compactly supported functions $\varphi : \Omega\to \R$ 
satisfying this property.
\end{lem}
\begin{proof}[Proof of Theorem~\ref{thm:AM-non-deg}]
	Let $(H_*,c_*)$ be of Aubry-Mather type, and let $\sigma > 0$ be as in Definition~\ref{defn:AM}. Let $K$ be as in Lemma~\ref{lem:K} and denote $\Gamma_1 = \overline{B_{\sigma}(c_*)} \cap \Gamma$. Let  $\cR_2 \subset \cR_1 \cap \cV_{\sigma}(H_*)$ be the set of Hamiltonians such that for all 
	$c \in \Gamma_1 \cap \Gamma_*(H)$, the set
	$\tcN^0_{H \circ \Xi}(\xi^*c) \setminus \tcN^0_{H}(c)$ is totally disconnected. 
	According to Lemma~\ref{lem:G-dt}, this set is $G_\delta$, we will  to show 
	that it is dense on $\cV_{\sigma_2}(H_*)$ for some $0 < \sigma_2 < \sigma$. 

	Consider $c \in \Gamma_1 \cap \Gamma_*(H_*)$, let $\sigma_c>0$ be small enough so that Proposition \ref{prop:perb-b} applies to to the pair $(H_*, c)$ on the set $K$. For each $\theta_0 \in \tcN^0_{H_* \circ \Xi}(\xi^* c) \cap K$, define 
	$$
	D_{\sigma_c}(\theta_0) = \{\theta : \max_i|\theta^i - \theta_0^i| \le \sigma_c/{(2\sqrt{n})}\} 
	\subset B_{\sigma_c}(\theta_0).
	$$
	Proposition \ref{prop:barrier-holder} implies that the family of functions
	the family 	of functions 
	\[
		b^{\pm}_{H,c},\ c\in \Gamma_1\cap \Gamma_*(N)
	\]
	has Hausdorff dimension at most $1/\beta$, therefore we can apply Lemma \ref{lem:disc}  
	on the cube $D_{\sigma_c}(\theta_0)$ for each $H\in \cR_1$. We find arbitrarily small functions 
	$\varphi$ compactly supported in $D_{\sigma_c}(\theta_0)$ and such that each of the functions 
	\[
		b^{\pm}_{N,c}+\varphi, c\in  \Gamma_1 \cap \Gamma^*(N)
	\]
	have a totally disconnected set of minima in $D_{\sigma_c}(\theta_0)$. We then apply 
	Proposition \ref{prop:perb-b} to get Hamiltonians $H_{\varphi}$ approximating $H$.  We obtain:
	\begin{itemize}
		\item The set of Hamiltonians $H$ such that $ {\cN}_{H \circ \Xi}(\xi^* c) \cap D_{\sigma_c}(\theta_0)$ 
		is totally disconnected for each $c\in \Gamma_*(H)$	is dense in $\cR_1\cap \cV_{\sigma_c}(H_*)$. 
		By Lemma \ref{lem:G-dt}, it is also $G_{\delta}$, therefore residual.
	\end{itemize}
	
	Since $K$ is compact, there is a finite cover $K \subset \bigcup_{i=1}^k D_{\sigma_i}(\theta_i)$, such hat
	\begin{itemize}
		\item For a residual set $\cR^i(c)$ of $H \in B_{\sigma_c}(H_*)$, the set    $ {\cN}_{H \circ \Xi}(\xi^* c) \cap D_{\sigma_c}(\theta_i)$  is totally disconnected for all $i=1, \ldots, k$ and $c \in \Gamma_*(H)$.
	\end{itemize}
	Take $\cR_c = \bigcap_{i}\cR^i$, then for $H \in \cR_c$, the set  $ {\cN}_{H \circ \Xi}(\xi^* c) \cap K$ is totally disconnected. Finally, we consider a finite covering $\Gamma_1 \subset \bigcup_{j} B_{\sigma_{c_j}}(c_j)$ and repeat the above argument. 
\end{proof}

The next two sections are dedicated to proving Proposition~\ref{prop:barrier-holder}. 

\subsection{Weak KAM solutions are unstable manifolds}

An important consequence of the Aubry-Mather type is that the local unstable manifold coincide with 
an elementary weak KAM solution. 
\begin{prop}\label{prop:local-unst-wkam}
Suppose $(H_*, c_*)$ is of Aubry-Mather type. Then for each $(H, c) \in \cV_\sigma(H_*) \times B_\sigma(c_*)$ 
such that $\chi^{-1} \tcA^0_H(c)$ projects onto $\T$, we have 
\[
	W^u(\tcA_H^0(c)) \cap \pi_\theta^{-1} V^0 = \{(\theta, c + \nabla u(\theta)): \quad \theta \in V^0\}, 
\]
where $u(\theta) = h_{H, c}(\zeta, \theta)$ for some $\zeta \in \cA_H^0(c)$. 
\end{prop}

Proposition~\ref{prop:local-unst-wkam} follows from the following general statement 
(for the continuous Aubry sets). 
\begin{lem}\label{lem:unstable-weak-kam}
Assume that: 
\begin{enumerate}
	\item $\tcA_H(c)$ is a partially hyperbolic set of the dynamics. 
	\item There is a neighborhood $V \supset \cA_H(c) \subset \T^n \times \T$  such that  
	(strong) unstable manifold $W^u(\tcA_H(c)) \cap \pi_{(\theta, t)}^{-1} V$ is a Lipschitz graph 
	over $V \subset \T^n \times \T$. 
\end{enumerate}
Then there is a function $u: V \to \R$ solving the Hamilton-Jacobi equation 
\[
	u_t + H(\theta, c + \nabla u, t) = \alpha_H(c), 
\]
and $(\theta, \nabla u(\theta, t), t) = W^u(\tcA_H(c))$ for all $(\theta, t) \in V$.

Assume, in addition, that $\tcA_H(c)$ has a unique static class. Then there is 
$C \in \R$ and $(\zeta, \tau) \in \cA_H(c)$ such that
\[
	u(\theta, t)  = h_{H, c}(\zeta, \tau, \theta, t) + C, \quad (\theta, t) \in V. 
\]
\end{lem}

The idea of the Lemma is that the unstable manifold is a backward invariant Lagrangian manifold, which can be used to construct a solution of the Hamilton-Jacobi equation. The proof is a standard application of the method of characteristics. 

\begin{proof}
	By considering the Hamiltonian $H(\theta, c + p) - \alpha$ instead of $H$, we can always assume $c=0$ and $\alpha_H(c) =0$. 

	let $z = z_{\vartheta, s} \in \T^n \times \R^n \times \T$ be the unique point such that 
	$z \in W^u(\tcA_H(c))$ and $\pi_{(x, t)}z = (\vartheta, s)$ (Since $W^u(\tcA_H(c))$ is locally a graph over $(\theta, t)$), and let 
	$\gamma_{\vartheta, s}, p_{\vartheta, s}: (-\infty, s] \to \T^n \times \R^n$ be 
	the backward orbit of $z$. Let $z' = z'_{\vartheta, s} \in \tcA_H(c)$ be the unique such point 
	with $z_{\vartheta, s} \in W^u(z'_{\vartheta, s})$, and let 
	$\gamma_{\vartheta, s}',  p_{\vartheta, s}': (-\infty, s] \to \T^n \times \R^n$ be the backward orbit of $z'$. 
	See Figure~\ref{fig:unstable-weak-kam}. 

	\begin{figure}[t]
		\centering 
		\includegraphics[width=3in]{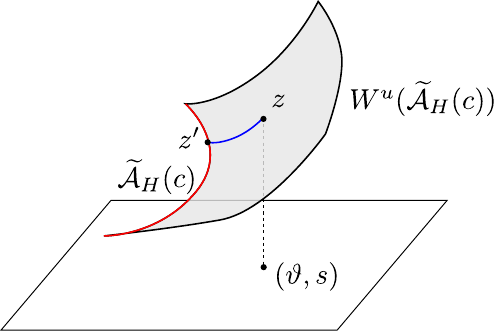} 
		\caption{Proof of Lemma~\ref{lem:unstable-weak-kam}} \label{fig:unstable-weak-kam}
	\end{figure}

	Let $w$ be an arbitrary weak KAM solution for $H$ at cohomology $0$. We then define 
	\[
		u_T(\vartheta, s) =   w(\gamma_{\vartheta, s}'(-T), -T) + \int_{-T}^0 L(d \gamma_{\vartheta, s}(t)) dt,  
	\]
	where $d\gamma(t) = (\gamma, \dot\gamma, t)$. 
	Note that 
	\[
		\begin{aligned}
			u_T(\vartheta, s) & =   \int_{-T}^0 L(d\gamma_{\vartheta, s}(t)) - L(d\gamma_{\vartheta, s}'(t)) dt + w(\gamma_{\vartheta, s}'(-T), -T) +  \int_{-T} L(d\gamma_{\vartheta, s}'(t)) dt \\
			& = \int_{-T}^0 L(d\gamma_{\vartheta, s}(t)) - L(d\gamma_{\vartheta, s}'(t)) dt + w(\gamma_{\vartheta, s}'(s), s),
		\end{aligned}
	\]
	where the second line is due to the fact that $(\gamma_{\vartheta, s}'(t), t) \in \cA_H(c)$.  Since the  first integral converges since $dist(d\gamma_z(t), d\gamma_{\vartheta, s}'(t)) \to 0$ exponentially fast,  the limit
	\[
		u(\vartheta, s) := \lim_{T \to \infty} u_T(\vartheta, s)
	\]
	converges exponentially fast. Consider $-T_1 < s$, then
	\[
		\begin{aligned}
			& 	u(\vartheta, s)  - u(\gamma_{\vartheta, s}(-T_1), -T_1)  = \lim_{T \to \infty} u_T(\vartheta, s) - \lim_{T - T_1 \to \infty} u_{T - T_1}(\gamma_{\vartheta, s}(-T_1), -T_1) \\
			& = \lim_{T \to \infty}  
			\left( 	w(\gamma_{\vartheta, s}'(-T), -T) + \int_{-T}^0 L(d \gamma_{\vartheta, s}(t)) 	dt \right)	\\
			& \quad - \lim_{T \to \infty} \left(  w(\gamma_{\vartheta, s}'(-T), -T) + \int_{-T}^{T_1} L(d \gamma_{\vartheta, s}(t)) dt\right) \\
			& = \int_{T_1}^T  L(d \gamma_{\vartheta, s}(t)) dt. 
		\end{aligned}
	\]
	In other words, each curve $\gamma_{\vartheta, s}:(-\infty, s] \to \T^n$ is $(u, L_c)$-calibrated. Moreover, we have 
	\[
		\begin{aligned}
			& 	d_\vartheta u(\vartheta, s) = d_\vartheta w(\gamma_{\vartheta, s}'(-T), -T)
			+ d_\vartheta  \int_{-T}^0 L(d\gamma_{\vartheta, s}(t))dt  \\
			& = \partial_v L(d\gamma_{\vartheta, s}'(-T)) \cdot \frac{\partial \gamma_{\vartheta, s}'(-T)}{\partial \vartheta}  + \left( \partial_v L(d\gamma_{\vartheta, s}(t)) \cdot 
			\frac{\partial\gamma_{\vartheta, s}(t)}{\partial \vartheta} \right) \Bigr|_{t = -T}^0, 
		\end{aligned}
	\]
	assuming the derivatives exist. The fact that $\gamma'_{\vartheta, s}(-T)$ and $\gamma_{\vartheta, s}(-T)$ are exponentially close to each other uniformly in $\vartheta, s$ implies
	\[
		\nabla u(\vartheta, s) := d_\vartheta u(\vartheta, s) = \partial_v L(d\gamma_{\vartheta, s}(0)), 
	\]
	and as a result
	\[
		\nabla u(\gamma_{\vartheta, s}(t), t) = \partial_v L(\gamma_{\vartheta, s}(t), \dot{\gamma}_{\vartheta, s}(t), t) = p_{\vartheta, s}(t) - c. 
	\]
	We now compute
	\[
		\begin{aligned}
			& 	L(\gamma_{\vartheta, s}(t), \dot{\gamma}_{\vartheta, s}(t), t) =  \frac{d}{dt}u(\gamma_{\vartheta, s}(t), t) = \nabla u(\gamma_{\vartheta, s}(t), t) \cdot \dot{\gamma}_{\vartheta, s}(t) + u_t(\gamma_{\vartheta, s}(t), t) \\
			& = 	u_t(\gamma_{\vartheta, s}(t), t) + p_{\vartheta, s}(t) \cdot \dot{\gamma}_{\vartheta, s}(t),  
		\end{aligned}
	\]
	we get 
	\[
		\begin{aligned}
			& 	\partial_t u(\gamma_{\vartheta, s}(t), t)  = L(\gamma_{\vartheta, s}(t), \dot{\gamma}_{\vartheta, s}(t), t)  - p_{\vartheta, s}(t) \cdot \dot{\gamma}_{\vartheta, s}(t) \\
			& = - H(\gamma_{\vartheta, s}(t), p_{\vartheta, s}(t), t) = - H(\gamma_{\vartheta, s}(t), \nabla u(\gamma_{\vartheta, s}(t), t), t)
		\end{aligned}
	\]
	which is exactly the Hamilton-Jacobi equation.

	We now prove the ``moreover'' part.  Consider a subsequence $T_n \to \infty$ such that $\gamma_{\vartheta, s}(-T_n) \to (\zeta, \tau) = (\zeta_{\vartheta, s}, \tau_{\vartheta, s}) \in \cA_H(0)$.  We then have 
	\[
		\begin{aligned}
			&	u(\vartheta, s)  \ge  \liminf_{n \to \infty} u(\gamma_{\vartheta, s}(-T_n), -T_n) +  \int_{T_n}^s L(d\gamma_{\vartheta, s}) dt \\
			&	\ge u(\zeta, \tau) + \liminf_{n \to \infty} A_{H}(\gamma_{y, s}(-T_n), -T_n, \vartheta, s) = u(\zeta, \tau) + h_{H}(\zeta, \tau, \vartheta, s). 
		\end{aligned}
	\]
	On the other hand, according to Proposition 4.1.8 of \cite{Fa}, when $u$ solves the Hamilton-Jacobi equation, $u$ is dominated by $L_c$, namely for any absolutely continuous curve $(\theta_1, t_1), (\theta_2, t_2) \in V$, we have $u(\theta_2, t_1) - u(\theta_1, t_1) \le A_H(\theta_1, t_1, \theta_2, t_2)$. As a result, 
	\[
		u(\vartheta, s) \le u(\zeta, \tau) + \liminf_{n \to \infty} A_{H}(\gamma_{y, s}(-T_n), -T_n, \vartheta, s)
	\]
	which implies $u(\vartheta, s) = u(\zeta, \tau) + h_{H}(\zeta, \tau, \vartheta, s)$. Notice that we have not proven what's needed since $(\zeta, \tau) = (\zeta_{\vartheta, s}, \tau_{\vartheta, s})$ depends on $(\vartheta, s)$. However when there is only one static class we have 
	\[
		h_{H}(\zeta_1, \tau_1, \vartheta, s) =  h_{H}(\zeta_1, \tau_1, \zeta_2, \tau_2) + h_H(\zeta_2, \tau_2, \vartheta, s), \quad (\zeta_1, \tau_1), (\zeta_2, \tau_2) \in \cA_H(0), 
	\]
	which allows a consistent choice of $(\zeta, \tau)$  for all $(\vartheta, s)$. 
\end{proof}

\begin{proof}
	[Proof of Proposition~\ref{prop:local-unst-wkam}]
	By converting Definition~\ref{defn:AM}, (2)(b), to its continuous counterpart, the condition of Lemma~\ref{lem:unstable-weak-kam} is satisfied for $\tcA_H(c)$. The proposition follows by taking the zero section of the weak KAM solution. 
\end{proof}

\subsection{Regularity of the barrier functions}
\label{sec:reg-barrier}

In this section we prove Proposition~\ref{prop:barrier-holder}. Let $c_* \in \Gamma_*(H_*)$, $H \in \cV_\sigma(H_*)$ and $c \in B_\sigma(c_*) \cap \Gamma_*(H)$,  then according to our assumption, the Aubry set $\tcA_H(c)$ is contained in the cylinder $\cC$, and has a unique static class. 

Let us now consider the lift $\tcA_{H \circ \Xi}^0(\xi^*c)$, which has two components $\tcS_1, \tcS_2$. The cylinder $\cC$ lifts to two disjoint cylinders $\cC_1, \cC_2$. For the rest of the discussion, we will consider only the static class $\tcS_1$ as the other case is similar. 

Definition~\ref{defn:AM} ensures that there is  a Lipschitz function $y = g(x) \in (-1, 1)$, $x \in \T$, such that $\{\chi(x, g(x)): x \in \T\} = \tcS_1(H, c)$. 
\[
	\tcS_1(H, c) = \{\chi(x, f(x)): x \in \T\} =: \{F_c(x) = (F_c^\theta(x), F_c^p(x)): x \in \T\} \subset \T^n \times \R^n. 
\]

Now suppose $c, c' \in \Gamma_*(H) \cap B_\sigma(c_*)$, first we have:
\begin{lem}
	There is $C>0$ such that 
	\[
		\sup_x \|F_c(x) - F_{c'}(x)\| \le C \|c - c'\|^{\frac12}. 
	\]
\end{lem}
\begin{proof}
	The proof is the same as the one in Lemma~5.8, \cite{BKZ} where the proof only used the symplecticity of the cylinder. The the main idea is that $\sup_x \|F_c(x) - F_{c'}(x)\|$ is $\frac12$ H\"older with respect to the area between the two invariant curves restricted to the cylinder, while the latter is equivalent to the symplectic area by assumption. A direct calculation shows the symplectic area is bounded by $\|c' - c\|$. 
\end{proof}

\begin{figure}[t]
	\centering
	\includegraphics[width=2.85in]{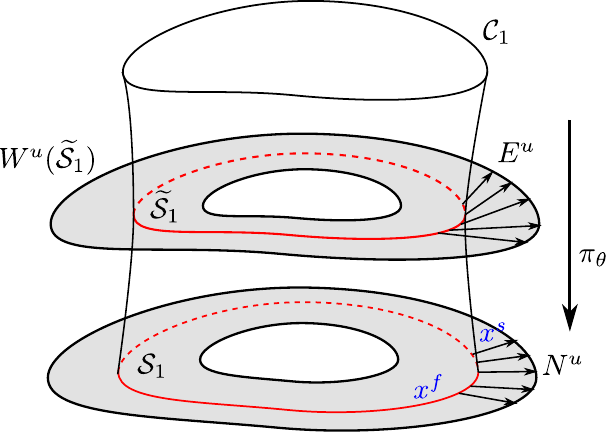}
	\caption{Local coordinates in the configuration space near $\cS_1$.} \label{fig:dr-coordinates}
\end{figure}

Let us now denote 
\[
	u_c(\theta) = h_{H \circ \Xi, \xi^*c}(\zeta_1, \theta), \quad \zeta_1 \in \cS_1(H, c). 
\]
Before moving forward, we define a convenient local coordinate system near $\cS_1$.  According to Proposition~\ref{prop:local-unst-wkam},  the graph $(\theta, \nabla u_c(\theta))$ locally coincides with $W^u(\tcS_1(H, c))$. The (strong) unstable bundle $E^u(z)$ of $z \in \tcS_1(H,c)$ is tangent to $W^u(\tcS_1(H,c))$ at every point, 
and transverse to the tangent cone of $\tcS_1(H,c)$. Since $W^u(\tcS_1(H, c))$ is a Lipschitz graph over $\T^n$,  the projection $N^u(z) : = d\pi_\theta E^u(z)$ forms a non-zero section of the normal bundle to $\cS_1(H,c) = \pi_\theta \tcS_1(H, c)$ within the configuration space $\T^n$. By choosing an orthonormal basis $e_1(z), \cdots, e_{n-1}(z)$, $N^u(z)$ naturally defines a coordinate system $(x^f, x^s) \in \T \times \R^{n-1}$ on the tubular neighborhood of $\cS_1$:
\[
	\iota(x^f, x^s) = F_c^\theta(x^f) + \sum_{i = 1}^{n-1} e_i(F_c^\theta(x^f)) x^s_i. 
\]
See Figure~\ref{fig:dr-coordinates}. We note that the coordinate system is only H\"older since the unstable bundle is only H\"older a priori. Therefore, we consider $C^\infty$ functions that approximate $F_c^\theta$ and $e_i$ in the $C^0$ sense, th new coordinate system is still well defined near $\cS_1$, and the $x^s$ coordinate projects onto the unstable direction.  In the sequel, we fix such a coordinate system using $W^u(\tcS_1(H_*, c_*))$. Due to semi-continuity, for $(H, c)$ close to $(H_*, c_*)$, the coordinate system is defined in a neighborhood of $\cS_1(H, c)$.

\begin{figure}[t]
	\centering
	\includegraphics[width=4.4in]{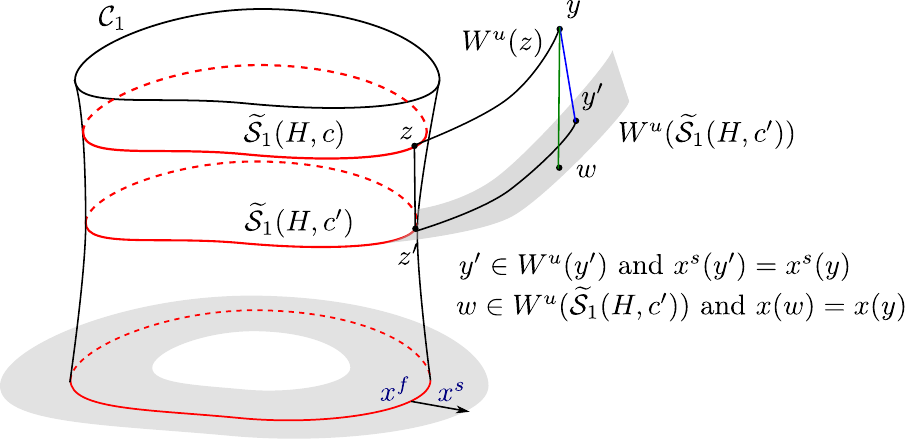}
	\caption{Proof of Lemma~\ref{lem:holder}} \label{fig:holder}
\end{figure}

By assumption, the foliation to $W^u(\tcS_1(H, c))$, $c \in \Gamma_*(H)$ by strong unstable manifold is $\beta_0$-H\"older for some $\beta_0 > 0$ (see \cite{PSW1997}). As a result, we obtain the following regularity:

\begin{lem}\label{lem:holder}
There is $\sigma>0$ such that for $H \in \cV_\sigma(H_*)$, $c, c' \in \Gamma_*(H) \cap B_\sigma(c_*)$,  $\theta \in B_\sigma(\cS_1(H, c_*))$, there is $\beta>0$,  $C_1 > 0$, $C_2 \in \R$ such that 
\begin{enumerate}
	\item 	$|  \nabla u_c(\theta) - \nabla u_{c'}(\theta) | \le C_1 \, \|c - c'\|^{\beta}$;
	\item $|u_c(\theta) - u_{c'}(\theta)- C_2| \le C_1 \, \|c - c'\|^{\beta}$. 
\end{enumerate}
Moreover, the same holds with $\cS_1$ replaced with $\cS_2$. 
\end{lem}
\begin{proof}
	The proof is essentially the same as Lemma~5.9 in \cite{BKZ}. Let $\sigma_1$ be small enough such that the local coordinates $(x^f, x^s)$  is defined for $|x^s| < \sigma$. We then consider the weak KAM solution $u_c \circ \iota(x^f, x^s)$ instead, which we still denote as $u_c(x^f, x^s)$, abusing the notation. 

	Let $|x^s| < \sigma_1$, let $y = (x^f, x^s, \nabla u_c(x^f, x^s))$, and let $z \in \cS_1(H, c)$ be such that $y \in W^u(z)$. We then define $z' \in \cS_1(H, c)$ be the unique such point with $x^f(z') = x^f(z)$. Finally define $y' \in W^u(z')$ be such that $x^s(y') = x^s(y)$, which is possible since $W^u(z')$ is a graph over the $x^s$ coordinates. See Figure~\ref{fig:holder}. 

	We note that within the center unstable manifold $W^u(\cC_1)$, the NHIC $\cC_1$ on one hand, and $x^s = x^s(y)$ on the other hand serves as two transversals to the strong unstable foliation $\{W^u(\cdot)\}$. Since the foliation is $\beta_0$-H\"older, there exists $C>0$ (throughout the proof, $C$ denotes a generic constant) such that 
	\[
		\|y - y'\| \le C \|z - z'\|^{\beta_0} \le C \|c - c'\|^{\frac{\beta_0}2}. 
	\]
	Denote $w = (x, \nabla u_{c'}(x))$, and noting $y' \in W^u(\cS_1(N,c')) = \{(x, \nabla u_{c'}(x))\}$ which is locally a $C^1$ graph, we get for $C > 0$
	\[
		\|w - y'\| \le C \|\pi_x(w) - \pi_x(y')\| = C \|\pi_x(y) - \pi_x(y')\| \le   C \|y - y'\|,
	\]
	therefore 
	\[
		\|\nabla u_c(x) - \nabla u_{c'}(x)\| \le \|w - y\| \le \|w -y'\| + \|y - y'\| \le C \|y - y'\| \le C \|c - c'\|^{\frac{\beta_0}{2}}. 
	\]
	We now revert the local coordinate $\iota$ to obtain item 1 with $\beta = \beta_0/2$, and possibly changing $\sigma$ and $C_1$. Item 2 is obtained from item 1 by direct integration. 
\end{proof}

\subsection{Bifurcation type}
\label{sec:bif-type-def}

We prove Theorem~\ref{thm:bifur-non-deg} in this section. Suppose $(H_*, c_*)$ is of bifurcation type, let $H^1, H^2$ be the Hamiltonians as defined in Definition~\ref{def:bif}, then there exists cylinders $\cC_1, \cC_2$ containing the local Aubry sets $\tcA_{H^1}(c)$ and $\tcA_{H^2}(c)$. 

Let us denote $\alpha^i_H(c) = \alpha_{H^i}(c)$, called the local alpha functions. Then $\tcA_H(c)$ has two static classes if and only if $\alpha^1_H(c) = \alpha^2_H(c)$. Moreover, for each $c$, the rotation vector $\rho_{H^i}(c)$ is uniquely defined, as a result, $\alpha^i_{H}$ is a $C^1$ function. 

\begin{prop}\label{prop:bif-gen}
Let $H_*, c_*$ be of bifurcation AM type. There is $\sigma>0$ such that for an open and dense subset of $H \in \cV_{\sigma}(H_*)$, such that there are at most finite many $c \in \overline{B_{\sigma}}(c_*)\cap \Gamma$ for which $\alpha^1_H(c) = \alpha^2_H(c)$. 

For each such (bifurcation) $c$, we have
$\tcA_{H^1}^0(c),\tcA_{H^2}^0(c)$ are both supported on hyperbolic periodic orbits. Moreover, $\tcN_H^0(c) \setminus \tcA_H^0(c)$ is a discrete set. 
\end{prop}
\begin{proof}
	Let $\sigma$ be as in Definition~\ref{def:bif}.  For the rest of the proof, we refer to $\overline{B_{\sigma}(c_*)} \cap \Gamma$ as $\Gamma$ to simplify notations. 

	Consider the family of Hamiltonians $H_\lambda = H(\theta, p, t) - \lambda f$ ($f$ is the mollifier function in Definition~\ref{def:bif}). Then we have
	\[
		\alpha^1_{H_\lambda}(c) = \alpha^1_H(c) + \lambda, \quad 
		\alpha^2_{H_\lambda}(c) = \alpha^2_H(c). 
	\]
	By Sard's theorem, there exists a full Lebesgue measure set $E$ of  regular values $\lambda$ of the function $\alpha^1_H - \alpha^2_H$, which implies for $\lambda \in E$, $0$ is a regular value of  $\alpha^1_{H_\lambda} - \alpha^2_{H_\lambda}(c)$ which implies the equation has only finitely many solutions on $\Gamma$, and at each solution, $\alpha^1$ and $\alpha^2$ has different derivatives when restricted to $\Gamma$. Note that this property is open in $H$, which implies the claim of our proposition is an open property. We only need to prove density.  

	Let $c_i(\lambda) \in \Gamma$, $i = 1, \cdots , N$ be the set on which $\alpha^1_{H_\lambda} = \alpha^2_{H_\lambda}$, then since $\frac{d}{d\lambda} \alpha^1_{H_\lambda}(c_i) \ne \frac{d}{d\lambda} \alpha^2_{H_\lambda}(c_i)$, each  $c_i(\lambda)$ is locally a monontone function in $\lambda$. 

	We now impose the assumption that $H$ is a Kupka-Smale system, namely, all periodic orbits are non-degenerate. Then in this setting, all invariant measures of rational rotation vectors are supported on hyperbolic periodic orbits, by a further perturbation, we can ensure for each $c$ there is only one minimal periodic orbit. Since the Aubry set is upper semi-continuous when there is only one static class, and the hyperbolic periodic orbit is structurally stable, we obtain that each hyperbolic periodic orbit is the Aubry set for an open set of $c$'s. Let us denote by $R(H, \Gamma)$ the set of all $c \in \Gamma$ such that $\tcA_{H}^0(c)$ is a hyperbolic periodic orbit, then $R(H^1, \Gamma)$ and $R(H^2, \Gamma)$ are both open and dense in $\Gamma$. Then for each $c_i$, there is $d_i >0$ and an open and sense subset of $\lambda \in B_{d_i}(0)$ such that $c_i(\lambda) \in R(H^1, \Gamma) \cap R(H^2, \Gamma)$. Let $d = \min d_i$, then there is an open and dense set of $\lambda \in B_d(0)$ such that for all $i = 1, \cdots, N$, $c_i(\lambda) \in R(\Gamma)$. For these $c_i$'s, $\tcA^0_{H^1}(c_i), \tcA^0_{H^2}(c_i)$ are both hyperbolic periodic orbits. Using the Kupka-Smale theorem again, it is an open and dense property such that the stable and unstable manifolds of $\tcA^0_{H^1}(c_i)$ and $\tcA^0_{H^2}(c_i)$ intersect transversally. Due to dimension considerations, the intersection is a discrete set.  Since each orbit in $\tcN_{H}(c_i) \setminus \tcA_{H}(c_i)$ is a heteroclinic orbit between $\tcA^0_{H^1}(c_i)$ and $\tcA^0_{H^2}(c_i)$, it is also a discrete set. 

	We have now proven the claim of our proposition holds for an arbitrarily small perturbation $H_\lambda$ of $H$, and hence is dense on $\cV_\sigma(H_*)$. 
\end{proof}

An analogous statement holds for the asymmetric bifurcation case:
\begin{prop}
\label{prop:asym-bif-gen}
Let $H_*, c_*$ be of asymmetric bifurcation type, then there is $\sigma > 0$ such that for an open and of $H \in \cV_{\sigma}(H_*)$, such that there a unique $c \in \overline{B_{\sigma}}(c_*)\cap \Gamma$ for which $\alpha^1_H(c) = \alpha^2_H(c)$. Moreover, at such (bifurcation) $c$, we have
$\tcA_{H^1}^0(c),\tcA_{H^2}^0(c)$ are both supported on hyperbolic periodic orbits. Moreover, $\tcN_H^0(c) \setminus \tcA_H^0(c)$ is a discrete set. 
\end{prop}
\begin{proof}
The proof is nearly identical to Proposition~\ref{prop:bif-gen} with the simplification that for $\sigma$ small enough,  $\tcA_{H^2}(c)$ is always the same hyperbolic periodic orbit, and $\alpha_{H^2}(c)$ is a linear function. Since $\alpha_{H^1}(c)$ is convex on $\Gamma$, there is at most one bifurcation. The rest of the proof is identical. 
\end{proof}

\begin{proof}[Proof of Theorem~\ref{thm:bifur-non-deg}]
	\emph{Case 1}:
	Let $H_*, c_*$ be of bifurcation AM type, and let $\sigma>0$ be such that Proposition~\ref{prop:bif-gen} holds on a open and dense subset $\cR_1$ of $\cV_{\sigma}(H_*)$. Then on $\Gamma \cap B_\sigma(c_*)$ there are at most finitely many bifurcations $c_i$, $i = 1, \cdots, N$. Moreover, for each bifurcation value $c_i$, there is $\sigma_i>0$ such that for $c \in \overline{B_{\sigma_i}(c_i)}$, the local Aubry sets are hyperbolic periodic orbits. This means that for each (non-bifurcation) $c \in \overline{B_{\sigma_i}(c_i)} \setminus c_i$, the sets $\tcN^0_H(c) = \tcA^0_H(c)$ are contractible (in fact finite). At the bifurcation values, the set $\tcN^0_H(c_i) \setminus \tcA^0_H(c_i)$ is discrete. 

	We now consider the set 
	\[
		\left( \Gamma \cap  \overline{B_\sigma(c_*)} \right)  \setminus  \bigcup_{i=1}^N B_{\sigma_i}(c_i)
	\]
	which is compact with finitely many connected components. On each of the components $\tcA^0_H(c)$ has a unique static class. We apply Theorem~\ref{thm:AM-non-deg} to get there is $\sigma'>0$ such that the dichotomy of Theorem~\ref{thm:AM-non-deg} holds for a residual subset $\cR_2$ of $H \in \cV_{\sigma'}(H_*)$ on each of the connected components. The theorem follows by taking the intersection of $\cR_1$ and $\cR_2$ as well as the smaller value of $\sigma$ and $\sigma'$. 

	\emph{Case 2}: If $(H_*, c_*)$ is of asymmetric bifurcation type, the proof is the same with Proposition~\ref{prop:asym-bif-gen} replacing Proposition~\ref{prop:bif-gen}. 
\end{proof}


\section{Aubry-Mather type at the single resonance}
\label{sec:sr-AM}

\subsection{Normally hyperbolicity and localization of Aubry/Ma\~ne sets in the single maximum case}

In this section we consider the Hamiltonian system 
\[
	N_\epsilon = H_0(p) + \epsilon Z(\theta^s, p) + \epsilon R(\theta, p, t), 
\]
where $\theta = (\theta^s, \theta^f) \in \T^{n-1} \times \T$ and $p =(p^s, p^f) \in \R^{n-1}\times \R$. Consider the resonant curve
\[
	\Gamma = \{ p \in B^n: \quad \partial_{p^s}H_0(p) = 0\} = 
	\{p_*(p^f) = (p^s_*(p^f), p^f):\quad p^f \in [a_-, a_+]\}. 
\]

We first consider the case where $p \in \Gamma$ satisfy the condition $[SR1_\lambda]$, namely for each $p \in B_\lambda (p_0) \cap \Gamma$,  $Z(\theta^s, p_*(p^f))$ has a unique global maximum at $\theta^s_*(p^f)$, and
\[
	D^{-1} I \le \partial^2_{pp}H_0 \le DI, \quad \|Z\|_{C^3}\le 1, \quad \lambda I \le - \partial^2_{\theta^s \theta^s} Z(\theta^s_*(p^f), p_*(p^f)) \le \lambda I, 
\]
and that for $K > 0$,  and $p_0 = p_*(a_0) \in \Gamma$, 
\[
	\|R\|_{C^2_I(\T^n \times B_{K \sqrt{\epsilon}}(p_0)\times \T)} \le \delta,
\]
note that we are using the rescaled norm $C^2_I$ (see \eqref{eq:rescaled-norm}). 
We then set $K_1 = K/D$ and  consider the local segment
\[
	\Gamma(\epsilon, p_0) = \{p_*(p^f): \quad p^f \in [a_0 - K_1\sqrt{\epsilon}, a_0 + K_1\sqrt{\epsilon}] \} \subset \Gamma \cap B_{K\sqrt{\epsilon}}(p_0), 
\]
which is contained in $B_\lambda(p_0) \cap \Gamma$ if $\epsilon_0$ is small enough depending on $K$. 
Throughout this section, we write $f = O(g)$ if $|f| \le C |g|$ for $C>0$ that may depend only on $D$, $\lambda$ and $n$. 

\begin{thm}[Proof in Section~\ref{sec:sr-hyp-coord} and Section~\ref{sec:sr-nhic-proof}, see also \cite{BKZ}, Theorem 3.1] \label{thm:nhic-sr}
Assume that $K_1 > 2$.  There is $\delta_0, \epsilon_0 > 0$ and $C = C(D, \lambda, n) >1$, such that if  $0 < \delta < \delta_0$ and $0 < \epsilon < \max\{\epsilon_0, \sqrt{\delta}\}$, there is a $C^2$ map 
\[
	(\Theta^s, P^s)(\theta^f, p^f, t): \T \times [a_0 - K_1\sqrt{\epsilon}/2, a_0 + K_1\sqrt{\epsilon}/2] \times \T \to \T^{n-1} \times \R^{n-1}, 
\]
such that $\cC = \{(\theta^s, p^s) = (\Theta^s, P^s)(\theta^f, p^f, t)\}$ is weakly invariant in the sense that the vector field is tangent to $\cC$. $\cC$ is contained in the set 
\[
	V = \{(\theta, p, t):\quad \|\theta^s - \theta^s_*(p^f)\| \le C^{-1}, \quad \|p^s - p^s_*(p^f)\| \le C^{-1} \}
\]
and it contains all the invariant set contained in $V$. Moreover, we have 
\[
	\left\| \Theta^s(\theta^f, p^f, t) - \theta^s_*(p^f) \right\| \le C \delta, \quad
	\left\|  P^s(\theta^f, p^f, t) - p^s_*(p^f) \right\| \le C \delta \sqrt{\epsilon}, 
\]
\[
	\|\partial_{p^f}\Theta^s\| \le C \sqrt{\delta/\epsilon}, \quad \|\partial_{(\theta^f, t)}\Theta^s\| \le C \sqrt{\epsilon}, \quad
	\|\partial_{p^f}P^s\| \le C, \quad \|\partial_{(\theta^f, t)} P^s\| \le C \sqrt{\epsilon}. 
\]
The cylinder $\cC$ is normally hyperbolic with its stable/unstable bundle projects onto the $\theta^s$ direction. 
\end{thm}

\begin{thm}[Proof is in Section~\ref{sec:local-aubry}, see also \cite{BKZ}, Theorem 4.1, Theorem 4.2]\label{thm:var-local-sr}
There is $\delta_0 = \delta_0(\lambda, n, D) >0$  and $\epsilon_0 = \epsilon_0(\lambda, n, D, \delta)$  such that if $0 < \delta < \delta_0$ and $0 < \epsilon < \epsilon_0$,  the Ma\~ne set of the cohomology $c$ satisfies
\[
	\tcN_{N_\epsilon}(c) \subset B_{\delta^{1/5}}(\theta_*^s) \times \T \times B_{\sqrt{\epsilon}\cdot \delta^{1/16}} \times \T \subset \T^{n-1} \times \T \times \R^n \times \T. 
\]
If $u$ is a weak KAM solution of $N_\epsilon$ at $c$, then the set $\tcI(u, c) \subset \T^n \times \R^n$ is contained in a $18\sqrt{D\epsilon}-$Lipschitz graph above $\T^n$. 
\end{thm}

The theorems as stated are analogous to the cited theorems in \cite{BKZ}. The main difference is that we now assume the much weaker assumption $\|R\|_{C^2_I} \le \delta$. Nevertheless, we now check that the method in \cite{BKZ} applies in the same way, and leading to the estimates as stated. 

\subsection{Aubry-Mather type at single resonance}

\begin{thm}\label{thm:sr-AM-type-single}
Let $c = p_*(p^f)$ with $p^f \in [a_0 - K_1\sqrt{\epsilon}/4, a_0 + K_1\sqrt{\epsilon}/4]$, there there is $\epsilon_0, \delta_0>0$ such that for  $N_\epsilon = H_0 + \epsilon Z + \epsilon R$ with $0 < \epsilon < \epsilon_0$ and  $\delta < \delta_0$, then $N_\epsilon, c$ is of Aubry-Mather type with the hyperbolic cylinder given by the embedding
\[
	\chi(\theta^f, p^f) = (\Theta^s, P^s)(\theta^f, p^f, 0)
\]
where $\Theta^s, P^s$ is from Theorem~\ref{thm:nhic-sr}. 
\end{thm}

We need the following statement. 
\begin{lem}[\cite{BKZ}, Proposition 4.11]\label{lem:graph-sr}
Let $u$ be a (discrete) weak KAM solution for $N_\epsilon$ at cohomology $c$, let
\[
	\cG_{c, u} = \{(\theta, c + \nabla u(\theta))\}
\]
be the associated pseudograph. Then there is $C>0$ depending only on $n, D$ such that  for any $k \ge 1/\sqrt{\epsilon}$,  $(\theta_1, p_1), (\theta_2, p_2) \in \phi^{-k}\cG_{c, u}$, we have 
\[
	\|p_2 - p_1\| \le C \sqrt{\epsilon}\|\theta_2 - \theta_1\|. 
\]
\end{lem}

We prove Theorem~\ref{thm:sr-AM-type-single} assuming Theorem~\ref{thm:nhic-sr} and \ref{thm:var-local-sr}. 
\begin{proof}
	[Proof of Theorem~\ref{thm:sr-AM-type-single}]
	Let us denote $\cC$ the cylinder in Theorem~\ref{thm:nhic-sr} and $\cC^0 = \cC \cap \{t =0\}$. 

	Let $\delta_0, \epsilon_0$ be small enough such that Theorem~\ref{thm:nhic-sr} and \ref{thm:var-local-sr} applies for $0 < \epsilon < \epsilon_0$ and $0 < \delta < \delta_0$. In particular, these statements holds on a open set $N_\epsilon$ and $c$, as required by Definition~\ref{defn:AM}. The embedding as described is a weakly normally hyperbolic invariant cylinder. Moreover, Theorem~\ref{thm:nhic-sr} implies that for any two tangent vectors $v, v' \in T_z\cC^0$, we have 
	$|d\Theta^s \wedge dP^s (v, v')| \le C \sqrt{\delta} |d\theta^f \wedge dp^f(v, v')|$, and therefore for $\delta$ small enough
	\[
		\left| (d\Theta^s \wedge dP^s + d\theta^f \wedge dp^f) (v, v') \right| \ge (1 - C\sqrt{\delta}) \ge \frac12 |d\theta^f \wedge dp^f(v, v')|. 
	\]
	Since the form $d\theta^f \wedge dp^f$ is non-degenerate on $\cC^0$, the symplectic form $d\theta^s \wedge dp^f + d\theta^f \wedge dp^f$ is non-degenerate when restricted to $\cC^0$. 

	Theorem~\ref{thm:var-local-sr} implies that for $\delta$ small enough, the Ma\~ne set $\tcN_{N_\epsilon}(c)$ is contained in the neighborhood $V$ described in Theorem~\ref{thm:nhic-sr}. Since the Ma\~ne set is invariant, it must be contained in the cylinder $\cC$. Therefore $\tcA^0_{N_\epsilon}(c) \subset \tcN^0_{N_\epsilon}(c) \subset \cC^0$. 

	We now show thta $\tcA^0_{N_\epsilon}$ is a Lipschitz graph over $\theta^f$. Let $(\theta_1, p_1), (\theta_2, p_2) \in \tcA^0_{N_\epsilon}(c)$, then by Theorem~\ref{thm:var-local-sr} we have 
	\[
		\|p_2 - p_1\| \le 18D \sqrt{\epsilon} \|\theta_2 - \theta_1\| \le 18D\sqrt{\epsilon}\left( \|\theta^s_2 - \theta^s_1\|  + \|\theta^f_2 - \theta^f_1\|\right). 
	\]
	By Theorem~\ref{thm:nhic-sr}, 
	\begin{equation}
		\label{eq:theta-s-lip}
		\|\theta_2^s - \theta_1^s\| \le C( 1 + \sqrt{\delta/\epsilon}) \left( \|\theta^f_2 - \theta^f_1\| + \|p_2 - p_1\| \right).  
	\end{equation}
	Combine everything, we get for some constant $C_1$ depending on $n, D$, 
	\[
		\left(  1 - C_1(\sqrt{\epsilon} + \sqrt{\delta}) \right) \|p_2 - p_1\| \le C_1(\sqrt{\epsilon} + \sqrt{\delta}) \|\theta_2^f - \theta_1^f\|. 
	\]
	Which implies $\|p_2 - p_1\| \le \|\theta^f_2 - \theta^f_1\|$ if $\delta, \epsilon$ is small enough depending only on $C_1$. Combine with \eqref{eq:theta-s-lip} we get $\|(\theta_2, p_2) - (\theta_1, p_1)\| \le 2 \|\theta^f_2 - \theta^f_1\|$ if $\epsilon, \delta$ is small enough. 

	Finally, let $u$ be a weak KAM solution for $N_\epsilon$ at cohomology $c$, and assume that $\tcA^0_{N_\epsilon}(c)$ projects onto $\theta^f$ component. Since the strong unstable manifolds depends $C^1$ on the base point, the unstable manifold $W^u(\tcA^0_H(c))$ is a Lipschitz manifold. Moreover, since the strong unstable direction projects onto the $\theta^s$ direction, $W^u(\tcA^0_H(c))$ is locally a graph over $\theta = (\theta^f, \theta^s)$. This verifies (2)(b) of Definition~\ref{defn:AM}. 
	We have verified all conditions in Definition~\ref{defn:AM} and therefore $N_\epsilon, c$ is of Aubry-Mather type. 
\end{proof}

\subsection{Bifurcations in the double maxima case}

In this section we assume that for the Hamiltonian
\[
	N_\epsilon = H_0 + \epsilon Z + \epsilon R, 
\]
where $Z$ satisfies the condition $[SR2_\lambda]$, namely for all  $p\in B_\lambda(p_0) \cap \Gamma$,
there exists two local maxima
$\theta^s_1(p)$ and 
$\theta^s_2(p)$ of the function $Z(.,p)$ 
in $ \T^{n-1}$ satisfying  
\begin{align*}
	\partial^2_{\theta^s} Z(\theta^s_1(p),p)< \lambda I
	\quad,\quad 
	\partial^2_{\theta^s} Z(\theta^s_2(p),p)< \lambda I,\qquad \qquad \qquad \\
	Z(\theta^s,p) < \max \{Z(\theta^f_1(p),p),Z(\theta^f_2(p),p)\}- \lambda
	\big(\min\{d(\theta^s-\theta^s_1), d(\theta^s-\theta^s_2)\}\big)^2. 
\end{align*}
Let $f: \T^{n-1} \to \R$ be a bump function satisfying the following conditions:
\[
	f|_{B_\lambda(\theta^s_1(p_0))} = 0, \quad f|_{B_\lambda(\theta^s_2(p_0))} = 1, 
\]
and $0 \le f \le 1$ otherwise. Define 
\[
	Z_1 = Z - f, \quad Z_2 = Z - (1- f), 
\]
then $Z_1(\cdot, p_0)$ has a unique maximum at $\theta_1^s(p_0)$, while $Z_2(\cdot, p_0)$ has a unique maximum at $\theta^s_2(p_0)$. There is $\kappa > 0$ depending only on $\lambda$ such that the same holds for $p \in B_\kappa(p_0)$. To the Hamiltonian 
\[
	N_\epsilon^i = N + \epsilon Z_i + \epsilon R, \quad i = 1,2 
\]
we may apply Theorem~\ref{thm:nhic-sr} and \ref{thm:var-local-sr} to obtain existence of the NHIC $\cC_i$ which contains the local Aubry sets $\cA^i(c)$ for $c \in B_{\kappa/2}(p_0)$. Morevoer, we may define the local alpha functions $\alpha^i(c) = \alpha_{N^i_\epsilon}(c)$ similar to Section \ref{sec:bif-type-def}. The cohomology $c \in \Gamma \cap B_{\kappa/2}(p_0)$ is of bifurcation type if and only if $\alpha^1(c) = \alpha^2(c)$, and is of Aubry-Mather type if and only if $\alpha^1(c) \ne \alpha^2(c)$. 

Moreover, in this case we have the following analog of Theorem~\ref{thm:var-local-sr}. 
\begin{thm}[Proof is in Section~\ref{sec:local-aubry}, see also \cite{BKZ}, Theorem 4.5] \label{thm:sr-double-loc}
If $\delta = \delta(\lambda, n, D)>0$ is small enough and if $\epsilon < \epsilon_0 = \epsilon_0(\lambda, n, D, \delta)$ is small enough,  for $c \in B_{\kappa/2}(p_0)\cap \Gamma$ the Aubry set at cohomology $c$ of the Hamiltonian $N_{\epsilon}$ satisfies
\[
	\tcA(c)\subset \big( B(\theta^s_1, \delta^{1/5})\cup B(\theta^s_2, \delta^{1/5}) \big)
	\times \T \times B(c, \sqrt{\epsilon}\delta^{1/16})\times \T
	\subset 
	\T^{n-1}\times \T \times \R^n \times \T.
\]
If, moreover, the projection $\theta^s(\cA^0(c)) \subset \T^{n-1}$ 
is contained in one of the (disjoint) balls $B(\theta^s_i,\delta^{1/5})$, then 
the projection  $\theta^s(\cN^0(c)) \subset \T^{n-1}$  of the Ma\~né set is contained in the same ball
$B(\theta^s_i,\delta^{1/5})$.
\end{thm}

\begin{thm}\label{thm:sr-bif-type}
Suppose $Z$ satisfies condition $SR2_\lambda$ at $p = c_*$. Then there exists $\epsilon_0, \delta_0>0$, such that if $0 < \epsilon < \epsilon_0$ and $0 < \delta < \delta_0$, such that if $N_\epsilon = H_0 + \epsilon Z + \epsilon R$ with $\|R\|_{C^2_I} < \delta$, $N_\epsilon, c_*$ is of bifurcation Aubry-Mather type. 
\end{thm}
\begin{proof}
	Let $V_1 = B_\lambda(\theta^s_1(p_0)) \times \T \subset \T^{n-1} \times \T$, and $V_2 = B_\lambda(\theta^s_2(p_0))\times \T$. We check that the functions $N_\epsilon^1$ and $N_\epsilon^2$ both satisfies the conditions of Theorem~\ref{thm:sr-AM-type-single}, and as a result, conditions (1) and (2) of Definition~\ref{def:bif} are satisfied. It suffices to check item (3). 
\end{proof}

\subsection{Hyperbolic coordinates}
\label{sec:sr-hyp-coord}

The Hamiltonian flow admits the following equation of motion :
\begin{equation}\label{eq:perturbed}
\begin{cases}
	\dot{\theta}^s = \partial_{p^s}H_0 + \epsilon\partial_{p^s}Z + 
	\epsilon\partial_{p^s} R \\
	\dot{p}^s = -\epsilon \partial_{\theta^s}Z - \epsilon\partial_{\theta^s} R \\
	\dot{\theta}^f = \partial_{p^f}H_0 + \epsilon\partial_{p^f}Z +  
	\epsilon\partial_{p^f} R \\
	\dot{p}^f = -\epsilon \partial_{\theta^f} R \\
	\dot{t}=1
\end{cases}.
\end{equation}
The Hamiltonian structure of the flow 
is not used in the following proof.

The system \eqref{eq:perturbed} is a perturbation of 
\[
	\dot{\theta}^s = \partial_{p^s}H_0,\quad
	\dot{p}^s = -\epsilon \partial_{\theta^s}Z, \quad
	\dot{\theta}^f = \partial_{p^f}H_0,\quad
	\dot{p}^f = 0, \quad
	\dot{t}=1
\]
which admits a normally hyperbolic invariant cylinder 
\[
	\{( \theta^s_*(p^f), \theta^f, p^s_*(p^f), p^f, t): \quad \theta^f, t \in \T, \}.
\]
Set 
\[
	B(p^f):=\partial^2_{p^sp^s}H_0(p_*(p^f)), \quad A(p^f):=- \partial^2_{\theta^s\theta^s}Z(\theta^s_*(p^f),p_*(p^f)), 
\]
then as in \cite{BKZ}, there is a positive definite matrix $T(p^f)$ such that 
\[
	\Lambda(p^f) = T(p^f)A(p^f)T(p^f)=T^{-1}(p^f)B(p^f)T^{-1}(p^f).
\]
We lift the equation to the universal cover, and consider the change of variable
\begin{equation}
	\label{eq:sr-hyp-coord}
	\begin{aligned}
		& x=T^{-1}(p^f) (\theta^s-\theta^s_*(p^f))+ \epsilon^{-1/2}T(p^f)(p^s-p^s_*(p^f)) \\
		& y=T^{-1}(p^f) (\theta^s-\theta^s_*(p^f))- \epsilon^{-1/2}T(p^f)(p^s-p^s_*(p^f)), \\
		& I=\epsilon^{-1/2}(p^f - a_0),\quad \Theta=\gamma \theta^f,
	\end{aligned}
\end{equation}
where $0 < \gamma < 1$ is a parameter to be determined. 

\begin{lem}[\cite{BKZ}, Lemma 3.1, Lemma 3.2]\label{lem:all-estimates}
For each $p^f \in [a_-, a_+]$, we have $\Lambda(p^f) \ge \sqrt{\lambda/D}I$. We also have the following estimates:
\[
	\|T\|_{C^2}, \|T^{-1}\|_{C^2},  \|\partial_{p^f}\theta^s_*\|_{C^0}, \|p^s_*\|_{C^2} = O(1), \quad
	\|\theta^s - \theta_*^s\|_{C^0} \le O(\rho), \quad \|p^s - p^s_*\|_{C^0} \le O(\sqrt{\epsilon} \rho), 
\]
where $\rho = \max\{\|x\|, \|y\|\}$. 
\end{lem}
Note   we are using the regular $C^2$ norm in Lemma~\ref{lem:all-estimates} since $T$ depends only on $H_0$ and $Z$ which are bounded in the regular norm.

\begin{lem}\label{lem:iso-block}
The equation of motion in the new variables takes the form
\[
	\dot{x} = \sqrt{\epsilon}\Lambda(a_0 + \sqrt{\epsilon} I) x + \sqrt{\epsilon}\,  O(\delta + \rho^2), \quad
	\dot{y} = - \sqrt{\epsilon}\Lambda(a_0 + \sqrt{\epsilon} I) y + \sqrt{\epsilon} \, O(\delta + \rho^2),
\]
and $\dot{I} = O(\sqrt{\epsilon}\delta)$.
\end{lem}
\begin{proof}
	The last equation is straight forward. 
	We only prove the equation for $\dot{x}$ as the calculations for $\dot{y}$ are the same. In the original coordinates, we have
	\[
		\begin{aligned}
			& \dot \theta^s=B(p^f)(p^s-p^s_*(p^f))+O(\|p^s-p^s_*(p^f)\|^2)+ O(\delta\sqrt{\epsilon}), \\
			& \dot p^s=\epsilon A(p^f)(\theta^s-\theta^s_*(p^f))+ O(\epsilon\|\theta^s-\theta^s_*(p^f)\|^2)+O(\epsilon \delta),
		\end{aligned}
	\]
	where we used $\|\partial_{p}R\| = \epsilon^{-\frac12}\|\partial_I R\| = O(\epsilon^{-\frac12}\delta)$. Differentiate \eqref{eq:sr-hyp-coord},  use $\dot{p}^f = O(\epsilon \delta)$, Lemma~\ref{lem:all-estimates}, we get 
	\[
		\begin{aligned}
			& \dot{x} = T^{-1} \dot{\theta}^s + \epsilon^{-\frac12} T \dot{p}^s + O(\sqrt{\epsilon}\delta) \\
			& = \sqrt{\epsilon} T^{-1} B T^{-1} \cdot \epsilon^{-\frac12} T (p^s - p^s_*(p^f)) + \sqrt{\epsilon}  T A T \cdot  T^{-1} (\theta^s - \theta^s_*) + O(\sqrt{\epsilon} \delta + \sqrt{\epsilon} \rho^2) \\
			& = \sqrt{\epsilon}\Lambda(p^f) x + \sqrt{\epsilon}\,  O(\delta + \rho^2).   
		\end{aligned}
	\]
\end{proof}

\begin{lem}\label{lem:linear-matrix}
Suppose $\sqrt{\epsilon} \le \delta$, then in  the new coordinates $(x, y, \Theta, I, t)$ the linearized system is given by 
\[
	L=\begin{bmatrix}
	\sqrt{\epsilon}\Lambda &	0&		0 & 0  &0\\
	0&-\sqrt{\epsilon}\Lambda &		0 &		0&0\\
	0 &0 & 0 &		0  & 0\\
	0 & 0 & 0 & 0 & 0\\
	0 & 0 & 0 & 0 & 0
\end{bmatrix}+ \sqrt{\epsilon} O(\delta\gamma^{-1} + \rho + \gamma),
\]
where $\rho = \max\{\|x\|, \|y\|\}$. 
\end{lem}
\begin{proof}
	In the original coordinates, the linearized equation is given by the matrix
	\[
		\tilde{L} = 
		\begin{bmatrix}
			O(\sqrt{\epsilon}\delta) & B + O(\delta) & O(\sqrt{\epsilon}\delta) & \partial^2_{p^f p^s}H_0 + O(\delta) & 0 \\
			\epsilon A + O(\epsilon \rho) & O(\sqrt{\epsilon}\delta) & 0 & O(\sqrt{\epsilon}\delta) & 0 \\
			O(\sqrt{\epsilon}\delta) & O(1) & O(\sqrt{\epsilon}\delta) & O(1) & 0 \\
			0 & O(\sqrt{\epsilon}\delta) & 0 & O(\sqrt{\epsilon}\delta) & 0 \\
			0 & 0 & 0 & 0 & 0 
		\end{bmatrix} + O(\epsilon \delta). 
	\] 
	The coordinate change matrix is
	\[
		\left[ \frac{\partial(\theta^s, p^s, \theta^f, p^f, t)}{\partial(x, y, \Theta, I, t)} \right] = 
		\begin{bmatrix}
			T/2 &T/2   & 0 &O(\sqrt{\epsilon})& 0 \\
			\sqrt{\epsilon}T^{-1}/2 &-\sqrt{\epsilon}T^{-1}/2 & 0 & \sqrt{\epsilon}\partial_{p^f}p^s_*+O(\epsilon\rho)&0\\
			0 &0  & \gamma^{-1} & 0 & 0\\
			0 & 0 & 0 & \sqrt{\epsilon} & 0\\
			0 & 0 & 0 & 0 & 1
		\end{bmatrix}.
	\]
	The product is
	\[  
		\begin{aligned}
			& 		\tilde{L} \left[ \frac{\partial(\theta^s, p^s, \theta^f, p^f, t)}{\partial(x, y, \Theta, I, t)} \right] =  O(\epsilon \delta \gamma^{-1} + \epsilon \rho) +  \\
			& 		
			\begin{bmatrix}
				\sqrt{\epsilon} B T^{-1}/2 + O(\sqrt{\epsilon}\delta) & - \sqrt{\epsilon}BT^{-1}/2 + O(\sqrt{\epsilon}\delta) & O(\sqrt{\epsilon}\delta \gamma^{-1}) & O(\sqrt{\epsilon}\delta) & 0 \\
				\epsilon A T/2  &  \epsilon A T/2 &  0  & 0 & 0 \\
				O(\sqrt{\epsilon}) & O(\sqrt{\epsilon}) & O(\sqrt{\epsilon}\delta\gamma^{-1}) & O(\sqrt{\epsilon}) & 0 \\
				0 & 0 & 0 & 0 & 0 \\
				0 & 0 & 0 & 0 & 0
			\end{bmatrix}. 
		\end{aligned}
	\]
	Most of the computations are straightforward, with the exception of the fourth row, first column, which contains the following cancellation:
	\[
		\partial^2_{p^sp^f}H_0 \partial_{p^f} p^s_* + \partial^2_{p^fp^f}H_0  = 
		\partial_{p^f}\left( \partial_{p^f} H_0(p_*(p^f)) \right) = 0
	\]
	since $\partial_{p^f} H_0(p_*(p^f)) = 0$ for every $p^f$ by definition. 

	The differential of the inverse coordinate change is
	\[
		\left[ \frac{\partial(x, y, \Theta, I, t)}{\partial(\theta^s, p^s, \theta^f, p^f, t)} \right] = \begin{bmatrix}
		T^{-1} & \epsilon^{-1/2}T   & 0 &O(\epsilon^{-1/2}\lambda^{-1/4})& 0 \\
		T^{-1} & -\epsilon^{-1/2}T  &0& O(\epsilon^{-1/2}\lambda^{-1/4})&0\\
		0 &0  & \gamma & 0 & 0\\
		0 & 0 & 0 & \epsilon^{-1/2}& 0\\
		0 & 0 & 0 & 0 & 1
	\end{bmatrix}. 
\]
Finally, the new matrix of the linearized equation is 
\[
	\begin{aligned}
		& L = \left[ \frac{\partial(x, y, \Theta, I, t)}{\partial(\theta^s, p^s, \theta^f, p^f, t)} \right]  \tilde{L} \left[ \frac{\partial(\theta^s, p^s, \theta^f, p^f, t)}{\partial(x, y, \Theta, I, t)} \right]  \\
		& = O(\sqrt{\epsilon} (\delta\gamma^{-1} + \rho)) +
		\begin{bmatrix}
			\sqrt{\epsilon} \Lambda & 0 & 0 & 0 & 0 \\
			0 & - \sqrt{\epsilon}\Lambda & 0 & 0 & 0 \\
			O(\sqrt{\epsilon}\gamma) & O(\sqrt{\epsilon}\gamma) & O(\sqrt{\epsilon}\gamma) & 0 & 0 \\
			0 & 0 & 0 & 0 & 0 \\ 
			0 & 0 & 0 & 0 & 0  
		\end{bmatrix}. 
	\end{aligned}
\]
\end{proof}

\subsection{Normally hyperbolic invariant cylinder}
\label{sec:sr-nhic-proof}

We state the following abstract statement for existence of normally hyperbolic invariant manifolds, given in \cite{BKZ}.

Let $F:\R^n \to \R^n$ be a $C^1$ vector field. 
We split the space $\R^n$ as $\R^{n_u}\times\R^{n_s}\times \R^{n_c}$,
and denote by $z=(u,s,c)$ the points of $\R^n$. We denote by
$(F_u,F_s,F_c)$ the components of $F$:
$$
F(x)=(F_u(z),F_s(z),F_c(z)).
$$
We study the flow of $F$ in the domain
$$\Omega=B^u\times B^s\times \Omega^c
$$
where $B^u$ and $B^s$ are  the open Euclidean balls of radius $r_u$  
and $r_s$ in $\R^{n_u}$ and $\R^{n_s}$, and $\Omega^c = \Omega^{c_1} \times \R^{c_2}$ is a convex 
open subset of $\R^{n_c}$. We denote by
$$
L(z)=dF(z)=\begin{bmatrix}
L_{uu}(z)&L_{us}(z)&L_{uc}(z)\\L_{su}(z)&L_{ss}(z)&L_{sc}(z)\\L_{cu}(z) &L_{cs}(z)&
L_{cc}(z)
\end{bmatrix}
$$
the linearized vector field at point $z$. We  assume that $\|L(z)\|$ is bounded
on $\Omega$, which implies that each trajectory of $F$ is defined until it leaves $\Omega$.
We denote by $W^c(F, \omega)$ the union of all full orbits contained in $\Omega$, 
$W^{sc}(F, \Omega)$ the set of points whose positive  orbit remains inside
$\Omega$, and by $W^{uc}(F, \Omega)$ the set of points whose negative  orbit remains inside
$\Omega$.

Let us further consider a positive parameter $b>0$, and consider the set $\Omega_b^{c_2} = B_{b}(\Omega^{c_2})$ and $\Omega_b^c = \Omega^{c_1} \times \Omega^{c_2}_b$ . 

\begin{prop}[\cite{BKZ}, Proposition A.6]\label{realNHI}
Let 
$
F: \R^{n_u}\times \R^{n_s}\times \Omega^c_b
\to
\R^{n_u}\times \R^{n_s}\times \R^{n_c}
$
be a $C^2$ vector field.
Assume that there exists $\alpha,m,\sigma>0$ such that
\begin{itemize}
	\item $F_u(u,s,c)\cdot u> 0$ on  $\partial B^u \times \bar B^s \times \bar \Omega^c_b$.
	\item $F_s(u,s,c)\cdot s< 0$ on  $ \bar B^u \times \partial B^s \times \bar \Omega^c_b$.
	\item
	$L_{uu}(z)\geq \alpha I , \quad L_{ss}(z)\leq -\alpha I$ for each $x\in \Omega_b$
	in the sense of quadratic forms.
	\item
	$
	\|L_{us}(z)\|+\|L_{uc}(z)\|+\|L_{ss}(z)\|+\|L_{sc}(z)\|+
	\|L_{cu}(z)\|+\|L_{cs}(z)\|+\|L_{cc}(z)\| \leq m
	$ for each $x\in \Omega_b$.
	\item
	$
	\|L_{us}(z)\|+\|L_{uc}(z)\|+\|L_{ss}(z)\|+\|L_{sc}(z)\|+
	\|L_{cu}(z)\|+\|L_{cs}(z)\|+\|L_{cc}(z)\|+
	2\|F_{c_2}(z)\|/b\leq m
	$ for each $z\in \Omega_{b}-\Omega$.
\end{itemize}
Assume furthermore that
$$
K:= \frac{m}{\alpha-2m}\leq \frac{1}{8},
$$
then there exist $C^2$ maps 
$$ w^{sc}:B^s \times \Omega^c_{b} \to B^u,
\quad  w^{uc}:B^u \times \Omega^c_{b} \to B^s,
\quad
w^c:\Omega^c_{b}\to B^u\times B^s
$$
satisfying the estimates
$$\|d w^{sc}\|\leq K,\quad \|d w^{uc}\|\leq K,
\quad \|d w^c\|\leq 2K,
$$
the graphs of which 
respectively contain $W^{sc}(F,\Omega), W^{uc}(F,\Omega), W^{c}(F,\Omega)$.
Moreover, the graphs of the restrictions of $w^{sc},w^{uc}$ and $w^c$ to, respectively,
$B^s\times \Omega^c$, $B^u\times \Omega^c$ and $\Omega^c$, are tangent to the flow.

There exists an invariant  $C^1$ foliation of the graph of $w^{uc}$ whose leaves are graphs of $K$-Lipschitz maps above $B^u$.
The set $W^{uc}(F,\Omega)$ is a union of leaves : it has the structure of an invariant $C^1$ lamination.
Two points $x,x'$ belong to the same leaf of this lamination if and only if $d(x(t),x'(t))e^{t\alpha/4}$ is bounded on $\R^-$.

If  in addition  there exists a group $G$ of translations  
of $\R^{n_{c_1}}$
such that
$
F\circ (id\otimes id\otimes g\otimes id)=F
$
for each $g\in G$, 
then the maps $w^*$ can be chosen such that
\begin{equation}\label{eq-w}
w^{sc}\circ(id\otimes g\otimes id)=w^{sc},\quad
w^{uc}\circ(id\otimes g\otimes id)=w^{uc}, \quad
w^c\circ (g\otimes id)=w^c
\end{equation}
for each $g\in G$. The lamination is also translation invariant.
\end{prop}

\begin{proof}[Proof of Theorem~\ref{thm:nhic-sr}]
	Consider the equation of $N_\epsilon$ in the $(x, y, \theta^f, p^f, t)$ coordinates, and set $B^u = \{\|x\|< \rho\}$, $B^s = \{\|y\| < \rho\}$, $\Omega^{c_2} = \{\|I^f\|\le K_1/2\}$, $\Omega^{c_1} = \{(\theta^f, t)\}$ with the group translation $\Z^2$, $\alpha = \sqrt{\lambda/4D}$, $b = 1 < K_1/2$.  Note that our choice of $b$ implies $\Omega_b^{c_2} \subset \{\|I^f\| < K_1\}$. We fix $\gamma = \sqrt{\delta}$

	According to Lemma~\ref{lem:all-estimates} and Lemma~\ref{lem:iso-block}, we have for $\|x\| = \rho$, 
	\[
		\dot{x} \cdot x \ge \sqrt{\epsilon} \alpha \rho^2 - \rho \sqrt{\epsilon} O(\delta + \rho^2) = \sqrt{\epsilon}\rho \left(  \rho - O(\delta + \rho^2) \right) > 0
	\]
	as long as $\rho \ge C \delta$ for $C$ large enough. This verifies the first bullet point assumption. The second assumption is verified in the same way. 

	By Lemma~\ref{lem:linear-matrix}, we have for $C_1 >1$ depending on $\lambda, D, n$, 
	\[
		\begin{aligned}
			& 		\|L_{us}(z)\|+\|L_{uc}(z)\|+\|L_{ss}(z)\|+\|L_{sc}(z)\|+\|L_{cu}(z)\|+\|L_{cs}(z)\|+\|L_{cc}(z)\|  \\
			& = \sqrt{\epsilon} O(\delta \gamma^{-1} + \rho + \gamma) \le C_1 \sqrt{\epsilon}(\sqrt\delta + \rho)=:m/2,
		\end{aligned}
	\]
	and $L_{uu} \ge (2\alpha - m)I \ge \alpha I$ as long as $m < \alpha$. This is possible as long as $\epsilon, \delta, \rho < C_1^{-1} < \alpha$. Finally, we check that if $z \in \Omega_b \setminus \Omega$, $\|F_{c_2}(z)\| = O(\sqrt{\epsilon}\delta) \le C_1\sqrt{\epsilon}(\sqrt\delta + \rho)/2  = m/4$ if $C_1$ is chosen large enough. These estimates ensure all the bullet point conditions are satisfied. By choosing $C_1$ even smaller, we can ensure $m < \alpha/10$ which implies $K = m/(\alpha - 2m) \le 1/8$. 

	To summarize, we have shown all conditions of Proposition~\ref{realNHI} is satisfied if $0 < \epsilon < C_1^{-1}$, $0 < \delta < C_1^{-1}$, $C \delta \le  \rho \le C_1^{-1}$. We apply the Proposition twice, once for $\rho = C_1^{-1}$ and once for $\rho = C\delta$. The first application shows the invariant cylinder is the maximal invariant set in the set $\Omega^c = \{\|x\|, \|y\| \le C_1^{-1}\}$. The second application shows that the cylinder is in fact contained in the set $\{\|x\|, \|y\| \le C\delta\}$. Moreover, in the second application, we get the estimate $m = O(\sqrt{\epsilon\delta})$, which allows the estimate $K = O(\sqrt{\epsilon\delta})$. 

	We now return to the original coordinates using the formula
	\begin{equation}
		\label{eq:graph-explicit}
		\begin{aligned}
			\Theta^s(\theta^f,p^f,t)&=\theta^s_*(p^f)+\frac{1}{2}T(p^f)
			\cdot(w^c_u+w^c_s)(\gamma \theta^f, a_0 + \epsilon^{-1/2}p^f,t)\\
			P^s(\theta^f,p^f,t)&=p^s_*(p^f)+\frac{\sqrt{\epsilon}}{2}T^{-1}(p^f)\cdot
			(w^c_u-w^c_s)(\gamma \theta^f,a_0 + \epsilon^{-1/2}p^f,t).
		\end{aligned}
	\end{equation}
	All the estimates stated in Theorem \ref{thm:nhic-sr} follow directly
	from these expressions, and from the fact that $\{\|x\|, \|y\| \le C\delta\}$,  $\|dw^c\|\leq 2K = O(\sqrt{\epsilon\delta})$.
\end{proof}

\subsection{Localization of the Aubry and Ma\~ne sets}
\label{sec:local-aubry}

Let $D_1 = 2D$,  the Hamiltonian $N_\epsilon$ satisfies the following estimates if $\delta$ is small enough depending only on $D$:
\[
	D_1^{-1} I \le \partial^2_{pp} N_\epsilon  \le D I, \quad
	\|\partial^2_{\theta p} N_\epsilon\| \le 2\sqrt{\epsilon},\quad 
	\|\partial^2_{\theta\theta}N_\epsilon\| \le 3 \epsilon. 
\]
Let $L(\theta, p, t)$ be the Lagrangian associated to $N_\epsilon$, and $L_0$ the Lagrangian of $H_0$.

\begin{lem}\label{lem:sr-L}
The following estimates holds for the Lagrangian $L$. 
\begin{enumerate}
	\item  (\cite{BKZ}, Lemma 4.1) For $K>0$, the image of the set $\T^n \times B_{K \sqrt{\epsilon}} \times \T$ under the diffeomorphism $\partial_p N_\epsilon$ contains the set $\T^n \times B_{K_1 \sqrt{\epsilon}}(c) \times \T$, where $K_1 = K/(4D_1)$. 
	\item (\cite{BKZ}, Lemma 4.2) The estimates
	\[
		\|\partial^2_{\theta v}L\|_{C^0} \le 2 D_1 \sqrt{\epsilon}, \quad \|\partial^2_{\theta\theta} L \|_{C^0} \le 3\epsilon
	\]
	holds on $\T^n \times B_{K_1 \sqrt{\epsilon}}(c) \times \T$. 

	\item (\cite{BKZ}, Lemma 4.3) For $v \in B_{K_1 \sqrt{\epsilon}}(c)$, we have 
	\[
		|L(\theta, v, t) - \left( L_0(v) - \epsilon Z(\theta^s, c) \right)| \le 2\epsilon \delta. 
	\]
	\item (\cite{BKZ}, Lemma 4.4) The alpha function $\alpha_{N_\epsilon}(c)$ satisfies
	\[
		|\alpha_{N_\epsilon}(c) - (H_0(c) + \epsilon \max Z(\cdot, c))| \le 2\epsilon \delta. 
	\]
	\item (\cite{BKZ}, Lemma 4.5) There is $C > 1$ depending only on $D$ such that if $\epsilon < C^{-1} \delta$, we have the estimates
	\begin{align}
		L(\theta,v,t)-c\cdot v +\alpha(c)
		&\geq \|v-\partial H_0(c)\|^2/(4D_1)-\epsilon \hat Z_c(\theta^s)-4\epsilon \delta\\
		L(\theta,v,t)-c\cdot v +\alpha(c)
		&\leq
		D_1\|v-\partial H_0(c)\|^2 -\epsilon \hat Z_c(\theta^s)+4\epsilon \delta
	\end{align}
	for each $(\theta, v,t)\in \T^n\times \R^n \times \R$,
	where
	$\hat Z_c(\theta^s):=Z(\theta^s,c)-\max_{\theta^s}Z(\theta^s,c)$.
\end{enumerate}
\end{lem}
\begin{proof}
	We remark that in \cite{BKZ} it is assumed that $\|\partial^2_{\theta p} N_\epsilon\| \le 2\epsilon$ instead of $2 \sqrt{\epsilon}$ as we assumed. However the same calculations as given in the cited lemmas prove Lemma~\ref{lem:sr-L}, once appropriate changes are made. 
\end{proof}

\begin{prop}\label{prop:near-auto-semi-concave-esp}
For each $c \in \R^n$, the weak KAM solution $u$ of cohomology $c$ is $3\sqrt{D_1\epsilon}/2-$semi-concave. 
\end{prop}

\begin{proof}
	[Proof of Theorem~\ref{thm:var-local-sr}] The proofs of Theorem 4.1 and 4.2 in \cite{BKZ} replies on the estimates (3)(4)(5) in Lemma~\ref{lem:sr-L}, as well as Proposition~\ref{prop:near-auto-semi-concave-esp}. We have arrived at the same estimates using weaker assumption, the only change is we replaced the constant $D$ with $D_1 = 2D$. The rest of the proofs are exactly identical. 
\end{proof}

\begin{proof}
	[Proof of Theorem~\ref{thm:sr-double-loc}] This is similar to Theorem~\ref{thm:var-local-sr}. The proof of Theorem 4.5, \cite{BKZ} applies once we taking into account Lemma~\ref{lem:sr-L} and Proposition~\ref{prop:near-auto-semi-concave-esp}. 
\end{proof}


\newpage




\section{Normally hyperbolic cylinders at double resonance}
\label{sec:NHIC-DR-proof}

We prove Theorem~\ref{thm:NHIC-DR} in these sections. Recall that the theorem deals with two cases:
\begin{enumerate}
	\item (Simple critical homology) In this case we show the homoclinic curve $\eta_{\pm h}^0$ can be extended to periodic orbits both in positive and negative energy. The union of these periodic orbits form a $C^1$ normally hyperbolic invariant manifold. 
	\item (Non-simple homology) In this case we show that for positive energy, there exists periodic orbits shadowing $\eta_{h_1}^0$ and $\eta_{h_2}^0$ in a particular order. 
\end{enumerate}
Our strategy is to prove existence of these periodic orbits as hyperbolic fixed points of composition of local and global maps. A main technical tool to prove existence and uniqueness of these fixed points is the Conley-McGehee isolation block (\cite{McG73}). The plan of this section is as follows. 

In section \ref{sec:system-normal-form} we state a standard normal from near the hyperbolic fixed point.

In  section \ref{sec:Shilnikov-bvp} we study the property of the Shil'nikov boundary value problem. 

In section \ref{sec:local}, we apply results of Section~\ref{sec:Shilnikov-bvp} to establish strong hyperbolicity of the local map $\Phi_{loc}^*$
as well as existence of unstable cones. Since the global maps $\Phi_{glob}^*$ have bounded time, they have bounded norms and the linearization of the proper compositions $\Phi_{glob}^*\Phi_{loc}^*$ are dominated by the local component.

In section \ref{sec:Conley-McGehee} we give definition and derive simple properties of isolating blocks  of Conley-McGehee \cite{McG73}.

In section \ref{single-leaf}, under non-degeneracy conditions $[DR1^c]-[DR4^c]$, we construct isolating blocks for the proper
compositions of $\Phi_{glob}^*\Phi_{loc}^*$. 

In section \ref{double-leaf} we extend this analysis to $\Phi_{glob}^*\Phi_{loc}^* \cdots \Phi_{glob}^*\Phi_{loc}^*$. This would imply existence of families of shadowing orbits in non-simple case. 

In section \ref{NHIC-isolating-block} we complete the proof of Theorem \ref{thm:NHIC-DR} by showing that periodic orbits constructed in the previous  two sections, forms a normally hyperbolic invariant cylinder. Moreover, they  coincide with the shortest geodesics for the Jacobi metric. 

In Section~\ref{sec:homology} we prove Lemma~\ref{lem:homotopy} which describes in the non-simple case, the order at which the simple periodic orbits are shadowed.

\subsection{Normal form near the hyperbolic fixed point}
\label{sec:system-normal-form}

We describe a normal form near the hyperbolic fixed point (assumed to be $(0,0)$) of the slow Hamiltonian $H^s:\T^2\times \R^2 \to \R$. For the rest of this section, we drop the superscript $s$ to abbreviate notations.
In a neighborhood of the origin, there exists a a symplectic linear
change of coordinates under which the system has the normal form
$$ H(u_1, u_2, s_1, s_2) = \lambda_1 s_1 u_1 + \lambda_2 s_2 u_2 + O_3(s,u).$$
Here $s=(s_1,s_2)$, $u=(u_1, u_2)$, and $O_n(s,u)$ stands for a
function bounded by $C|(s,u)|^n$. By taking a standard straightening coordinate change, we get:

\begin{lem}
After an $C^{r-1}$ symplectic coordinate change $\Phi$, the Hamiltonian takes the form 
\[
	N = H \circ \Phi = \lambda_1 s_1 u_1 + \lambda_2 s_2 u_2 + \sum_{i, j = 1, 2} s_i u_j O_1(s, u), 
\]
and the equation is 
\begin{equation}
  \label{eq:norm-form2}
  	\bmat{\dot{s} \\ \dot{u}} = \bmat{-\Lambda s + sO_1(s, u) \\ \Lambda u + u O_1(s, u)},
\end{equation}
where $\Lambda = diag\{\lambda_1, \lambda_2\}$. 
\end{lem}
\begin{proof}
    Since $(0,0)$ is a hyperbolic fixed point, for sufficiently small
    $r>0$, there exists stable manifold $W^s=\{(u=U(s), |s|\le r\}$ and
    unstable manifold $W^u=\{s=S(u), |u|\le r\}$ containing the
    origin. All points on $W^s$ converges to $(0,0)$ exponentially in
    forward time, while all points on $W^u$ converges to $(0,0)$
    exponentially in backward time. These manifolds are Lagrangian; as a
    consequence, the change of coordinates $s'=s-S(u)$,
    $u'=u-U(s')=u-U(s-S(u))$ is symplectic. Under the new coordinates, we
    have that $W^s=\{u'=0\}$ and $W^u=\{s'=0\}$. We abuse notation and
    keep using $(s,u)$ to denote the new coordinate system.

    Under the new coordinate system, the Hamiltonian has the form
    $$ H(s,u) = \lambda_1 s_1 u_1 + \lambda_2 s_2 u_2 + H_1(s,u), $$
    where $H(s,u) = O_3(s,u)$ and $H_1(s,u)|_{s=0} =
    H_1(s,u)|_{u=0}=0$. The Hamiltonian and the vector field takes the desired form under this coordinate. 
\end{proof}

\subsection{Shil'ni\-kov's boundary value problem}\label{sec:Shilnikov-bvp}

Recall the definition of the local map (Section~\ref{sec:intro-local-map}). We will define the sections $\Sigma^s_\pm$ to be a subset of the section $s_1 = \pm \delta$, and the sections $\Sigma^u_\pm$ to be contained in $u_1 = \pm \delta$. In this section, we study the properties of the local map via the Shil'nikov boundary value problem. 

\begin{prop}[Shil'nikov \cite{Shil67}, Lemma 2.1, 2.2, 2.3]\label{prop:bvp}
There exists $C>0$, $\delta_0>0$ and $T_0 > 0$ such that for each $0 < \delta < \delta_0$,  any $s^{in}=(s_1^{in},s_2^{in})$, $u^{out}=(u_1^{out},  u_2^{out})$ with $\|(s, u)\| \le \delta$ and any large $T>T_0$, there  exists a unique solution $(s^T,u^T):[0,T]\to B_{C\delta}$ of the system  (\ref{eq:norm-form2}) with the property $s^T(0)=s^{in}$ and   $u^T(T)=u^{out}$. Moreover, we have 
\[
\|s^T(t)\| \le C \delta e^{-\lambda_1 t/2},\quad \|u^T(t)\| \le C \delta e^{-\lambda_1 (T - t)/2}, 
\]
and 
\[
\left\| \frac{\partial s^T(t)}{\partial s^{in}} \right\| + \left\| \frac{\partial u^T(t)}{\partial s^{in}} \right\| \le C e^{-\lambda_1 t/2}, \quad
\left\| \frac{\partial s^T(t)}{\partial u^{out}} \right\| +	\left\| \frac{\partial u^T(t)}{\partial u^{out}} \right\| \le C e^{-\lambda_1 (T-t)/2}.
\]
\end{prop}

Consider the space of all smooth curves $(s, u): [0, T] \to B_{C \delta}$, define a map

\begin{cor}\label{cor:bvp-T-der}
Let $(s^T, u^T)$ be the solution in Proposition~\ref{prop:bvp} with fixed boundary values $s^{in}, u^{out}$. Then
\begin{equation}
  \label{eq:bvp-T-der-1}
  \left\|  \frac{d}{dT} (s^T(t), u^T(t)) \right\| \le C e^{-\lambda_1 (T - t)/2},\quad
\left\|  \frac{d}{dT} (s^T(T- t), u^T(T - t)) \right\| \le C e^{-\lambda_1 t/2},
\end{equation}
and 
\begin{equation}
  \label{eq:bvp-T-der-2}
  \begin{aligned}
 & 	\frac{d}{d\tau}\Bigr|_{\tau =T}(s^\tau(T), u^\tau(T)) = - X_N(s^T(T), u^T(T)) + O(e^{-\lambda_1 T/2}),  \\
& \frac{d}{d\tau}\Bigr|_{\tau =T}(s^\tau(T-\tau), u^\tau(T-\tau)) = X_N(s^T(0), u^T(0)) + O(e^{-\lambda_1 T/2}), 
\end{aligned}
\end{equation}
where $X_N$ denote the Hamiltonian vector field of $N$.
\end{cor}
\begin{proof}
Consider two solutions $(s^{T_1}, u^{T_1})$ and $(s^{T_2}, u^{T_2})$ of the boundary value problem, we extend the definition of the solutions to $\R$ by solving the ODE. We note that $(s^{T_1}, u^{T_1}): [0, T_2] \to \R^4$ satisfies the boundary condition $s^{in}$ and $u^{T_1}(T_2)$ on $[0, T_2]$, and by Corollary~\ref{cor:bvp-T-der}, 
\[
	\left\| (s^{T_1}, u^{T_1})(t) - (s^{T_2}, u^{T_2})(t) \right\| \le C e^{\lambda_1(T_2- t)/2}\|u^{T_1}(T_2) - u^{out}\|. 
\]
Note that
\[
	\lim_{T_2 \to T_1} \frac{u^{T_1}(T_2) - u^{out}}{T_2 - T_1} = \lim_{T_2 \to T_1} \frac{u^{T_1}(T_2) - u^{T_1}(T_1)}{T_2 - T_1} = \frac{d}{dt}\Bigr|_{t = T_1}u^{T_1}(t), 
\]
and as a result, 
\[
	\left\| \frac{d}{dT}\Bigr|_{T = T_1} (s^T, u^T)(t) \right\| \le C e^{\lambda_1(T_1- t)}. 
\]
This proves the first half of \eqref{eq:bvp-T-der-1} while the second half is similar. 

For \eqref{eq:bvp-T-der-2}, note 
\[
	\frac{d}{d\tau}\Bigr|_{\tau =T}(s^\tau(T), u^\tau(T)) = \frac{d}{d\tau}\Bigr|_{\tau = T}(s^\tau(\tau), u^\tau(\tau))
	 - \frac{d}{d\tau}\Bigr|_{\tau = T}(s^T(\tau), u^T(\tau)), 
\]
the first half of \eqref{eq:bvp-T-der-2} follows. The second half is similar. 
\end{proof}

Let $(v_{s_1}, v_{s_2}, v_{u_1}, v_{u_2})$ denote the coordinates for the tangent space induced by $(s_1, s_2, u_1, u_2)$. For $K>0$ and $x\in B_r$, we define the \emph{strong unstable cone} by 
\begin{equation}
\label{eq:unstable-cone-DR}
C^u_K(x)=\{K^2|v_{u_2}|^2 > |v_{u_1}|^2 + |v_{s_1}|^2 + |v_{s_2}|^2\}
\end{equation}
and the \emph{strong stable cone} to be
\[
C^s_K=\{K^2 |v_{s_2}|^2 > |v_{s_1}|^2 + |v_{u_1}|^2 + |v_{u_2}|^2\}.
\]

The following statement follows from the hyperbolicity of the fixed point $O$ via standard techniques. 
\begin{lem}[See for example, \cite{HPS77}]
\label{lem:local-cone} 
For any
$0<\kappa<\lambda_2-\lambda_1$, there exists $\delta=\delta(\kappa, K)$ and $C = C(\kappa, K)> 1$ such
that the following holds:
\begin{itemize}
	\item If $\phi_t(x)\in B_{4\delta}(O)$ for $0\le t \le t_0$, then $D\phi_t(C^{u}_K(x))\subset C^{u}_K(\phi_t(x))$ for all $0\le t\le t_0$. Furthermore, for any $v\in C^{u}_K(x)$,
	\[
	|D\phi_t(x)v| \ge C^{-1} e^{(\lambda_2-\kappa)t}, \quad 0 \le t\le t_0.
	\]
	\item If $\phi_{-t}(x)\in B_{4\delta}(O)$ for $0\le t \le t_0$, then $D\phi_{-t}(C^{s}_K(x))\subset C^{s}_K(\phi_{-t}(x))$ for all $0\le t\le t_0$. Furthermore, for any $v\in C^{s}_K(x)$,
	\[
	|D\phi_{-t}(x)v| \ge C^{-1} e^{(\lambda_2-\kappa)t}, \quad 0 \le t\le t_0.
	\]
\end{itemize}
\end{lem}

\begin{lem}\label{lem:bvp-cone}
For any $\kappa > 0$ there is $\delta_0 >0$ such that for any $0 < \delta \le \delta_0$, $|s_1^{in}| = |u_1^{out}| = \delta$, and $|s_2^{in}|, |u_2^{out}| \le \kappa\lambda_1 \delta/(2\lambda_2)$, we have 
\[
 \left| \frac{d}{dT}u_2^T(0) \right| \le \kappa \left| \frac{d}{dT} u_1^T(0) \right|, \quad 
 \left| \frac{d}{dT} s_2^T(T) \right| \le \kappa \left|  \frac{d}{dT} s_1^T(T) \right|. 
\]
Moreover, by integrating in $T$, we get 
\[
|u_2^T(0)| \le \kappa |u_1^T(0)|, \quad |s_2^T(T)| \le \kappa |s_1^T(T)|. 
\]
\end{lem}
\begin{proof}
Given $\kappa > 0$ we can choose $\delta_0$ small enough such that the backward flow $D\phi_{-t}$ preserves the cone 
\[
	(C_K^u)^c = \{ \|v_{u_2}\| \le K^{-1} \|(v_{s_1}, v_{s_2}, v_{u_1})\|\}, 
\]
where $K = \kappa^{-1}$. Note that 
\[
	\frac{d}{dT}(s^T(0), u^T(0)) = D\phi_{-T} \frac{d}{d\tau}\Bigr|_{\tau =T}(s^\tau(T), u^\tau(T)). 
\]
Since the vector 
\[
\begin{aligned}
 & 	\frac{d}{d\tau}\Bigr|_{\tau =T}(s^\tau(T), u^\tau(T)) = - X_N(s^T(T), u^T(T)) + O(e^{-\lambda_1 T/2}) \\
 &=    ( \lambda_1 s_1^T(T),  \lambda_2 s_2^T(T), \lambda_1 u_1^T(T), \lambda_2 u_2^T(T)) + O(\delta^2) + O(e^{-\lambda_1 T/2})  \\
& = (0, 0, \lambda_1 \delta, \lambda_2 u_2^{out}) + O(\delta^2) + O(e^{-\lambda_1 T/2}) \in (C_K^u)^c
\end{aligned}
\]
if $|u_2^{out}| \le \kappa\lambda_1 \delta/(2\lambda_2)$, $\delta$ is small and $T$ large enough. It follows that $\frac{d}{dT}(s^T(0), u^T(0)) \in (C_K^u)^c$. Keep in mind that $\frac{d}{dT}s^T(0) = 0$, the first half of our estimate follows. The second half is proven in a symmetric way.
\end{proof}

\subsection{Properties of the local maps}
\label{sec:local}

Recall that $\gamma^+$ is a non-degenerate homoclinic orbit, and $\gamma^-$ is its time-reversal. Let us define 
\begin{equation}
\label{eq:Sigma-su}
\begin{aligned}
&   	\Sigma^s_+ = \{(\delta, s_2, u_1, u_1): \,  |s_2|, |u_1|, |u_2| \le \delta\}, \\
&	\Sigma^u_+ = \{(s_1, s_2, \delta, u_1): \,  |s_2|, |u_1|, |u_2| \le \delta\}. 
\end{aligned}
\end{equation}
Denote  $q^+ = \gamma^+ \cap \{s_1 = \delta\} = \{( \sigma, s_2^+, 0, 0 )\}$, and $p^+ = \gamma^+ \cap \{u_1 = \delta\} = \{(0, 0, \delta, u_2^+)\}$. Let $\kappa>0$ be a parameter to be defined in Proposition~\ref{prop:local-rect}. Since $\gamma^+$ is tangent to the $s_1, u_1$ axes, we can choose $\delta>0$ such that $|s_2^+|, |u_2^+| <  \kappa \delta$, therefore $q^+, p^+$ are contained $\Sigma^s_+$ and $\Sigma^u_+$ respectively. Define 
\[
\begin{aligned}
& 	l^s = \{(\delta, s_2, 0, 0): \, -\kappa \delta \le s_2 \le \kappa \delta\} \subset W^s(O)\cap \Sigma^s_+, \\
&	l^u = \{(0, 0, \delta, u_2): \, - \kappa \delta \le u_2 \le \kappa \delta\} \subset W^u(O) \cap \Sigma^u_+. 
\end{aligned}
\]

Let 
\[
\Sigma^{s, E}_+ = \Sigma^s_+ \cap \{N(s, u) = E\}, \quad 	\Sigma^{u, E}_+ = \Sigma^u_+ \cap \{N(s, u) = E\} 
\]
be restriction of the sections to an energy $E$ close to $0$. We would like to study the domain of the restricted local map $\Phl^{++}: \Sigma^{s, E}_+ \to \Sigma^{u, E}_+$. 

A rectangle $R$ is a diffeomorphic image of the Euclidean rectangle in $\R^2$. Let us label the vertices by $1,2,3, 4$ in clockwise order, and call the four sides $l_{12}, l_{34}, l_{14}, l_{23}$. 

\begin{figure}[t]
\centering 
\includegraphics[width=3in]{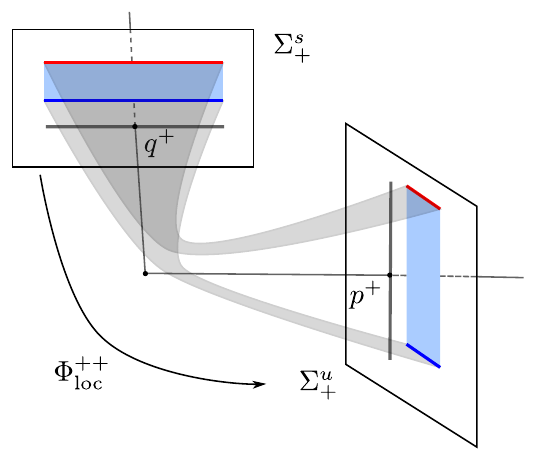} 
\caption{Rectangles mapped under $\Phl^{++}$}
\label{fig:rectangles}
\end{figure}

\begin{prop}
\label{prop:local-rect}
There is $e, \kappa>0$ such that for each $0 < E < e$, there is a rectangle $R^{++}(E) \subset \Sigma^{s, E}_+$ with sides $l_{ij}(E)$, such that:
\begin{enumerate}
	\item  $\Phl^{++}$ is well defined on $R^{++}(E)$, and its image $\Phl^{++}(R^{++}(E))$ is also a rectangle, it's four sides denoted $l_{ij}'(E)$. 
	\item As $E \to 0$, both $l_{12}(E)$ and $l_{34}(E)$ converges to $l^u$ in Hausdorff metric. Similarly, both $l_{14}'(E)$ and $l_{23}'(E)$ converges  to $l^u$. 
	\item If $(s,u) \in \Sigma^{s, E}_+$, $(s', u')\in \Sigma^{u, E}_+$ is such that  $\Phl^{++}(s, u) = (s', u')$, with $(s, u) \in B_{\kappa\delta}(q^+)$ and $(s', u') \in B_{\kappa\delta}(p^+)$, then $(s, u) \in R^{++}(E)$.
\end{enumerate}
\end{prop}
See Figure~\ref{fig:rectangles}. 

In order to do this, we prove a version of Proposition~\ref{prop:bvp} using $E$ as a parameter. 
\begin{prop}
\label{prop:bvp-E}
There is $e, \kappa > 0$ such that for each $0 < E < e$, $s^{in}$ and $u^{out}$ satisfies 
\[
	|s^{in}_2|, |u^{out}_2| \le \frac{\lambda_1}{2\lambda_2}\kappa  \delta, \quad  s_1^{in} = u_1^{out} = \delta, 
\]
then there is $T = T_E > 0$, and a unique orbit $(s^E, u^E): [0, T_E] \to B_{4\delta}(O)$, such that 
\[
s^E(0) = s^{in}, \quad u^E(T_E) = u^{out}. 
\]
\end{prop}

\begin{proof}
[Proof of Proposition~\ref{prop:bvp-E}]
Let $S_E$ denote the energy surface $\{N(s, u) = E\}$. 
Given $s^{in}$ and $s^{out}$ and $T>  T_0$, Proposition~\ref{prop:bvp} implies the existence of a solution $(s^T, u^T)$ solving the boundary value problem. Moreover, as $T \to \infty$, $(s^T, u^T)(0) \to (s^{in}, 0) \in S_0$. Writing $E(T)$ the energy of the orbit $(s^T, u^T)$, we have $E(T) \to 0$ as $T \to \infty$. \\

\textbf{Claim}: (1) $E(T)$ is positive and strictly monotone, therefore $E(T)$ is one-to-one and onto for $T \in [T_0, \infty)$; (2) There is $e>0$ such that uniformly over all $s^{in}, u^{out}$, we have $0 < E(T) < e$ if $T \ge T_0$, and therefore the inverse $T_E$ is well defined for $E \in (0, e]$.  

By Lemma~\ref{lem:bvp-cone} $|u_2^T(0)| \le \kappa |u_1^T(0)|$. We first show $N(s^T(0), u^T(0)) > 0$. 

Let us consider the energy function $N(s_1, u_1, s_2, u_2) = \lambda_1 s_1 u_1 + \lambda_2 s_2 u_2 + H_1(s, u)$ with $H_1|_{s=0} = H_1|_{u=0} = 0$. Recall that the section $\Sigma^+ = \{(s, u): \, s_1 = \delta, |(s, u)| \le \delta\}$. We consider the submanifold $\Sigma \cap \{N(s,u) = 0\}$. Since
\[
\partial_{u_1} N(\delta, u_1, s_2, u_2) = \lambda_1 \delta  + \partial_{u_1} H_1 = \lambda_1 \delta + O(\delta^2) > 0, 
\]
by implicit function theorem, for $|(s, u)| \le \delta$ the surface $\Sigma \cap \{N =0\}$ is given by a graph $u_1 = u_1^c(s_2, u_2)$. We conclude that $\Sigma$ is divided by the graph $u_1 = u_1^c(s_2, u_2)$ into two components, and:
\begin{itemize}
 \item  $N(s, u) >0$ whenever $u_1 > u_1^c(s_2, u_2)$, and $N(s, u) < 0$ when $u_1 < u_1^c(s_2, u_2)$.
\end{itemize}
 Moreover, $u_1^c|_{u_2 = 0} = 0$, and by differentiating the implicit function, we get $|\partial_{u_2}u_1^c| \le C$ where $C$ is a constant depending on $\|N\|_{C^2}$. Therefore $|u_1^c(s_2, u_2)| \le C |u_2|$. This implies if $(\delta, u_1, s_2, u_2) \in \Sigma$ satisfies $u_1 > C|u_2| \ge |u_1^c|$, we have $N(\delta, u_1, s_2, u_2) >0$. In particular, this is satisfied for $(s^T(0), u^T(0))$ if $\kappa < C^{-1}$.

$E(T)$ is a continuous function defined on $T \ge T_0$, and $\lim_{T \to \infty}E(T) = 0$. We now prove that it is monotone. Using Lemma~\ref{lem:bvp-cone}, 
\[
\begin{aligned}
&	\frac{d}{dT}E(T) = \frac{d}{dT}N(s^T(0), u^T(0)) = \nabla N \cdot \left( 0, 0, \frac{d}{dT} u_1^T(0), \frac{d}{dT} u_2^T(0) \right) \\
& = (\lambda_1 s_1^{in} + O_2(s, u) )\frac{d}{dT} u_1^T(0) + (\lambda_2 s_2^{in} + O_2(s, u))\frac{d}{dT} u_2^T(0) \ne 0
\end{aligned}
\]
if $\kappa$ is sufficiently small. Since we have proved $E(T)>0$, claim (1) follows. Finally, Proposition~\ref{prop:bvp} implies
\[
0 < E(T) \le C \delta^2 e^{-\lambda_1 T/2}, 
\]
claim (2) follows. 
\end{proof}

\begin{proof} [Proof of Proposition~\ref{prop:local-rect}]
Let $\kappa>0$ be small enough such that Proposition~\ref{prop:bvp-E} applies, and let us rename $\lambda_1 \kappa/(2\lambda_2)$ into $\kappa$. Consider two parameters $a \in [-\kappa\delta, \kappa\delta]$ and $b \in [-\kappa\delta, \kappa\delta]$, and then there exist unique orbits $s^E_{a, b}$ solving the boundary value problem
\[
s^E_{a, b}(0) = (\delta, a), \quad  s^E_{a, b}(T_{a, b}(E)) = (\delta, b), 
\]
then 
\[
\begin{aligned}
R^{++}(E) & = \left\{ 
(s^E_{a, b}(0), u^E_{a, b}(0)): \,  a, b  \in [\delta/2, \delta/2] 	\right\}, \\
\Phl^{++}\left( R^{++}(E) \right) &= \left\{ 
(s^E_{a, b}(T_{a, b}(E)), u^E_{a, b}(T_{a, b}(E))): \,  a, b  \in [\delta/2, \delta/2] 	\right\} . 
\end{aligned}
\]

Note that two sides of the rectangle $R^{++}(E)$ are graphs over $s_2 \in [-\kappa\delta, \kappa \delta]$, since the $u_1, u_2$ component converge to $0$ exponentially fast as $T \to \infty$ (and $E \to 0$), we conclude that these two sides converge to $l^s$. The same can be said about the rectangle $\Phl^{++}\left( R^{++}(E) \right)$ and $l^u$. 

Finally, note that by definition, $R^{++}(E)$ contains the initial point of all orbits $(s, u): [0, T] \to S_E$ such that $s_1(0) = \delta$, $|s_2(0)|\le \kappa \delta$, $u_1(T) = \delta$, $|u_2(T)| \le \kappa \delta$, which contains the orbits such that $(s, u)(0) \in B_{\kappa \delta}(p)$ and $(s, u)(T) \in B_{\kappa \delta}(p^+)$. 
\end{proof}

For each of the other symbols $--$, $+-$ and $-+$, analogous statements to Proposition~\ref{prop:local-rect} hold, after making appropriate changes. We only make the remark that for symbol $--$, the range is energy is also $0 < E < e$, while for the symbol $+-$ and $-+$, the energy for Proposition~\ref{prop:local-rect} is $-e < E < 0$. 

\begin{figure}[t]
\centering
\includegraphics[width=5in]{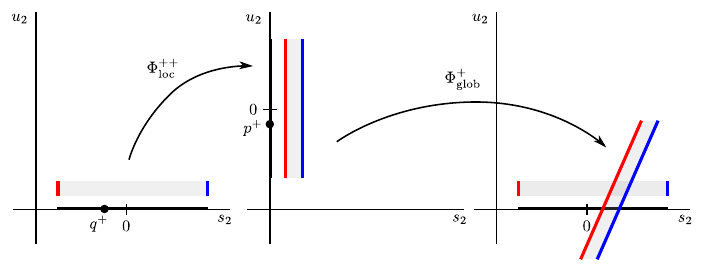}
\caption{Local map and global map in $s_2, u_2$ coordinates}
\label{fig:localmap}
\end{figure}

We now turn to the global map. We have $p^+ \in l^u$, $q^+ \in l^s$, and $\Phg^+(p^+) = q^+$. Moreover, condition $[DR4]^c$ imply that $\Phg^+(l^u)$ intersects $l^s$ transversally in the energy surface $S_E$. Let us note that the return time for the Poincar\'e map $\Phg^+$ is uniformly bounded, the restricted map $\Phg^+|S_E$ depends smoothly on $E$. As a result, we have the following corollary. 
\begin{cor}\label{trans-rect}
There exists $e>0$ such that, the following hold.  
\begin{enumerate}
	\item  For $0 < E < e$,  $\Phg^+ \circ \Phl^{++}(R^{++}(E))$ intersects $R^{++}(E)$ transversally.  This means, the  images of $l_{12}$ and $l_{34}$ intersect $l_{12}$ and $l_{34}$ transversally, and the images of $\gamma_{14}$ and $\gamma_{23}$ does not intersect $R^{++}(T)$.
	\item For $0 < E < e$,  $\Phg^-\circ \Phl^{--}(R^{--}(T))$ intersects $R^{--}(T)$ transversally.
	\item For $-e < E < 0$, $\Phg^- \circ \Phl^{+-} (R^{+-}(E))$  intersect $R^{-+}(E)$ transversally, and $\Phg^+\circ \Phl^{-+}(R^{-+}(E))$ intersect $R^{+-}(E)$ transversally.
\end{enumerate}
\end{cor}
See Figure~\ref{fig:localmap} for a demonstration. 

\subsection{Conley-McGehee isolating blocks}
\label{sec:Conley-McGehee}

The rectangles as constructed form isolation blocks of Conley and McGehee (\cite{McG73}). 

A rectangle $R = I_1\times I_2\subset \R^d\times \R^k$,
$I_1=\{\|x_1\|\le 1\}$, $I_2=\{\|x_2\| \le 1\}$ is called
{\it an isolating block} for the $C^1$ diffeomorphism $\Phi$, if the following hold:
\begin{enumerate}
	\item The projection of $\Phi(R)$ to the first component covers $I_1$.
	\item $\Phi|I_1 \times \partial I_2$ is homotopically equivalent to
	the identity restricted on the set $I_1 \times (\R^k\setminus \ \text{int}\ I_2)$.
\end{enumerate}
If $R$ is an isolating block of $\Phi$, then the set
$$W^+=\{x\in R: \Phi^k(x)\in R,\ k\ge 0\}\quad
$$
$$ \text{ (resp. }
W^-=\{x\in R: \Phi^{-k}(x)\in R,\ k\ge 0\})
$$
projects onto $I_1$ (resp. onto $I_2$) (see \cite{McG73}). If some
additional cone conditions are satisfied, then $W^+$ and $W^-$ are
in fact $C^1$ graphs. Note that in this case, $W^+\cap W^-$ is the unique
fixed point of $\Phi$ on $R$.

As usual, we denote by $C^u_K(x)=\{c\|v_1\| \le \|v_2\|\}$ the unstable cone at $x$.
We denote by $\pi C^u_K(x)$ the set $x+ C^u_K(x)$, which corresponds to
the projection
of the cone $C^u_K(x)$ from the tangent space to the base set. The stable
cones are defined similarly.  Let $U\subset \R^d\times \R^k$ be an open set and
$\Phi:U\to \R^d\times \R^k$ a $C^1$ map. Denote $D\Phi_x$ the linearization of
$\Phi$ at $x$.
\begin{enumerate}
	\item[C1.] $D\Phi_x$ preserves the cone field $C^u_K(x)$, and there exists
	$\Lambda>1$ such that $\|D\Phi(v)\|\ge \Lambda \|v\|$ for any $v\in C^u_K(x)$.
	\item[C2.] $\Phi$ preserves the projected restricted cone field $\pi C^u_K$, i.e.,
	for any $x\in U$,
	$$\Phi(U\cap \pi C^u_K(x)) \subset C^u_K(\Phi(x))\cap \Phi(U).$$
	\item[C3.] If $y\in \pi C^u_K(x)\cap U$, then $\|\Phi(y)-\Phi(x)\| \ge \Lambda \|y-x\|$.
\end{enumerate}
The unstable cone condition guarantees that any forward invariant set is contained in
a Lipschitz graph.
\begin{prop}[See \cite{McG73}]\label{Lipschitz}
Assume that $\Phi$ and $U$ satisfies [C1]-[C3], then any forward invariant set
$W\subset U$ is contained in a Lipschitz graph over $\R^k$ (the stable direction).
\end{prop}
Similarly, we can define the conditions [C1]-[C3] for the inverse map and the stable cone,
and refer to them as ``stable [C1]-[C3]'' conditions.  Note that if $\Phi$ and $U$ satisfies
both the isolating block condition and the stable/unstable cone conditions, then $W^+$
and $W^-$ are transversal Lipschitz graphs. In particular, there exists a unique
intersection, which is the unique fixed point of $\Phi$ on $R$.  We summarize as follows.

\begin{cor}\label{unique-fixed-point}
Assume that $\Phi$ and $U$ satisfies the isolating block condition,
and that $\Phi$ and $U$ (resp. $\Phi^{-1}$ and $U \cap \Phi(U)$)
satisfies the unstable (resp. stable) conditions [C1]-[C3].
Then $\Phi$ has a unique hyperbolic fixed point in $U$.
\end{cor}

\subsection{Periodic orbit in simple homologies}
\label{single-leaf}

We now apply the isolating block construction to the maps and rectangles obtained
in Corollary~\ref{trans-rect}.

\begin{prop}\label{prop:fixedpt}
There exists $e>0$ such that the following hold.
\begin{itemize}
	\item For $0 < E < e$, $\Phg^+\circ \Phl^{++}$ has a unique fixed point $p^+(E)$
	on $\Sigma^s_+\cap R^{++}(E)$;
	\item For $0 < E < e$, $\Phg^-\circ \Phl^{--}$ has a unique fixed point $p^-(E)$
	on $\Sigma^s_-\cap R^{--}(E)$;
	\item  For $-e < E < 0$: $\Phg^+ \circ \Phl^{-+}\circ \Phg^-\circ \Phl^{+-}$
	has a unique fixed point $p^c(E)$ on $R^{+-}(E)\cap (\Phg^-\circ \Phl^{+-})^{-1}(R^{-+}(E))$.
\end{itemize}
\end{prop}

To prove Proposition~\ref{prop:fixedpt}, we notice that the rectangle $R^{++}(T)$ has
$C^1$ sides, and there exists a $C^1$ change of coordinates turning it to a standard
rectangle. It's easy to see that the isolating block conditions are satisfied for
the following maps and rectangles:
\[
\Phg^+\circ \Phl^{++}\quad \text{and}\quad R^{++}(E), 
\qquad \Phg^-\circ \Phl^{--}\quad \text{and}\quad R^{--}(E),
\]
\[
\Phg^+\circ \Phl^{-+} \circ \Phg^- \circ \Phl^{+-} \quad \text{ and } 
\quad (\Phg^- \circ \Phl^{+-})^{-1}R^{-+}(E)\cap R^{+-}(E).
\]
It suffices to prove the stable and unstable cone conditions [C1]-[C3] for the corresponding return map and rectangles. Given $z = (s,u) \in S_E$, let us define the restricted cones
\[
C^{u, E}_K(z) = C^u_K \cap T_z S_E, \quad C^{s, E}_K = C^s_K \cap T_z S_E. 
\]

\begin{lem}
For each orbit $(s, u)(t)$ associated to orbits of the local map, the restricted cones $C^{u, E}_K$ and $C^{s, E}_K$ are non-empty. 
\end{lem}
\begin{proof}
We only prove the statement for $C^{u, E}_K$ since the stable one is completely symmetric. Moreover, observe that if we prove the cone $C^{u, E}_K$ is non-empty at $(s, u)(0)$, the same holds for the entire orbit, since the unstable cone is invariant under the forward dynamics. 

Note that
\[
\nabla N = (\lambda_1 u_1 + u O_1, \lambda_2 u_2 + u O_1,   \lambda_1 s_1 + s O_1, \lambda_2 s_2 + s O_1),
\]
and $\|u(0)\|$ is exponentially small in $T$, we have $\nabla N \sim (0, 0,   \lambda_1 s_1, \lambda_2 s_2)$. Since $|s_2|\le \delta = |s_1|$ on $\Sigma^s_+$, we have the angle between $\nabla N$ and $u_1$ axis is bounded from below. As a consequence, there exists $K>0$,
such that $C^{u,K}$ has nonempty intersection with the tangent
direction of $S_E$ (which is orthogonal to $\nabla N$). The lemma follows.
\end{proof}

We will only prove the [C1]-[C3] conditions conditions for the unstable cone $C_{E}^{u,\,c}$, the map $\Phg^+\circ \Phl^{++}$ and the rectangle $R^{++}(T)$; the proof for the other cases  can be obtained by making obvious changes to the case covered.

\begin{lem}\label{local-global-cone}
There exists $T_0>0$ and $c>0$ such that the following hold. Assume that $U\subset \Sigma^s_+\cap B_r$ is a connected open set on which the local map $\Phl^{++}$ is defined, and for each $x\in U$, 
\[
\inf\{ t\ge 0: \varphi_t(x)\in \Sigma^u_+\} \ge T_0.
\]
Then the map $D(\Phg^+\circ \Phl^{++})$ preserves the cone field  $C^u_K$, and the inverse $D(\Phg^+\circ \Phl^{++})^{-1}$ preserves  the non-empty $C^s_K$. Moreover, the projected cones $\pi C^u_K\cap U$ and $\pi C^s_K\cap V$ are preserved by  $\Phg^+\circ \Phl^{++}$ and its inverse, where $V=\Phg^+\circ \Phl^{++}(U)$.

The same set of conclusions hold for the restricted version. Namely,  we can replace $C^u_K$ and $C^s_K$  with $C^{u, E}_K$ and $C^{s, E}_K$, and $U$  with $U \cap S_E$.
\end{lem}

Recall that $l^u \ni p^+$ is the intersection of the unstable manifold $W^u(O)$ with the section $\Sigma^u_+$. Let $T^u$ be the tangent vector to $l^u$ at $p^+$, and $T^s$ the tangent vector to $T^s$ at $q^+$. We will show that if $\Phl^{++}|_{S_E}(x) = y$, then the image of the unstable cone $D\Phl^{++}(x) C^u_K$ is very close to $T^u$. This happens because the flow of tangent vector is very close to that of a linear flow. 

Assume that $\phi_t$ is a flow on $\R^d\times \R^k$, and $x_t$ is a trajectory of the flow. Let $v(t)=(v_1(t), v_2(t))$ be a solution of the variational equation, i.e. $v(t) = D\phi_t(x_t) v(0)$. 

\begin{lem}\label{cone-width} With the above notations assume that there exists $b_2>0$, $b_1 < b_2$ and $\sigma,\delta>0$ such that the variational equation 
\[
\dot{v}(t) =
\begin{bmatrix}
A(t) & B(t) \\
C(t) & D(t)
\end{bmatrix}
\begin{bmatrix}
v_1(t) \\ v_2(t)
\end{bmatrix}
\]
satisfy $A \le b_1 I$ and $D\ge b_2 I$ as quadratic forms, and $\|B\|\le \sigma$,
$\|C\|\le \delta$.

Then for any $c>0$ and $\kappa>0$, there exists $\delta_0>0$ such that
if $\ 0< \delta, \sigma< \delta_0$, we have
$$ (D\phi_t)\,C^u_K \subset C^{u,\beta_t}, \quad
\beta_t = ce^{-(b_2-b_1-\kappa)t} + \sigma/(b_2-b_1-\kappa). $$
\end{lem}
\begin{proof}
Denote $\gamma_0 = c$. The invariance of the cone field is equivalent to
$$ \frac{d}{dt} \left( \beta_t^2\langle v_2(t), v_2(t)\rangle - \langle v_1(t),
v_1(t)\rangle \right)\ge 0.$$
Compute the derivatives using the variational equation,  apply the norm bounds and
the cone condition, we obtain
$$
2\beta_t \left(  \beta_t' + (b_2 - \delta \beta_t -b_1) \beta_t - \sigma \right)\|v_2\|^2
\ge 0.$$
We assume that $\beta_t \le 2\gamma_0$, then for sufficiently small $\delta_0$,
$\delta \beta_t \le \kappa$. Denote $b_3 = b_2 - b_1 - \kappa$ and let
$\beta_t$ solve the differential equation
$$ \beta_t' = - b_3 \beta_t + \sigma. $$
It's clear that the inequality is satisfied for our choice of $\beta_t$.
Solve the differential equation for $\beta_t$ and the lemma follows.
\end{proof}

\begin{proof}[Proof of Lemma~\ref{local-global-cone}]
We will only prove the unstable version. By Assumption 4, there exists $c>0$
such that $D\Phg^{+}(q^+)T^{uu}(q^+)\subset C^u_K(p^+)$. Note that as
$T_0 \to \infty$, the neighborhood $U$ shrinks to  $p^+$ and $V$ shrinks
to $q^+$. Hence there exists $\beta>0$ and  $T_0>0$ such that
$D\Phg^+(y)C^{u,\,\beta}(y)\subset C^u_K$ for all $y\in V$.

Let $(s,u)(t)_{0\le t\le T}$ be the trajectory from $x$ to $y$.
By Proposition~\ref{prop:bvp}, we have $\|s\|\le e^{-\lambda_1T/4}$
for all $T/2 \le t\le T$. It follows that the matrix for the  variational
equation
\be \label{variational-matrix}
\begin{bmatrix}
A(t) & B(t) \\
C(t) & D(t)
\end{bmatrix}=
\begin{bmatrix}
-\diag\{\lambda_1, \lambda_2\} +O(s) & O(s) \\
O(u) & \diag\{\lambda_1, \lambda_2\} + O(u)
\end{bmatrix}
\ee
satisfies $A\le -(\lambda_1 -\kappa)I$, $D\ge (\lambda_1-\kappa)I$,
  $\|C\| = O(\delta)$ and 
  $\|B\| = O(e^{-(\lambda_1-\kappa)T/2})$.
   As before $C^u_K(x) = \{\|v_s\| \le c \|v_u\|\}$,
   Lemma~\ref{cone-width} implies
  $$ D\phi_T(x)C^u_K(x) \subset C^{u,\beta_T}(y),$$
  where $\beta_T = O(e^{-\lambda'T/2})$ and $\lambda' =
  \min\{\lambda_2-\lambda_1-\kappa, \lambda_1-\kappa\}$. Finally, note
  that $D\phi_T(x)C^u_K(x)$ and $D\Phl^{++}(x)C^u_K(x)$ differs
  by the differential of the local Poincar\'e map near $y$. Since near $y$
  we have $|s|=O(e^{-(\lambda_1-\kappa)T})$, using the equation of motion,
  the Poincar\'e map is exponentially close to identity on the $(s_1,s_2)$
  components, and is exponentially close to a projection to $u_2$ on the
  $(u_1,u_2)$ components. It follows that the cone $C^{u,\,\beta_T}$
  is mapped by the Poincar'e map into a strong unstable cone with
  exponentially small size. In particular, for $T\ge T_0$, we have
  $$ D\Phl^{++}(x)C^u_K(x) \subset C^{u,\,\beta}(y), $$
  and
the first part of the lemma follows. To prove the restricted version
we follow the same arguments.
\end{proof}

Conditions [C1]-[C3] follows, and this concludes the proof of Proposition~\ref{prop:fixedpt}.

\subsection{Periodic orbits for non-simple homology}
\label{double-leaf}

In the case of the non-simple homology, there exist two rectangles $R_1$ and $R_2$, whose images under $\Phg\circ \Phl$ intersect themselves transversally, providing a ``horseshoe'' type picture.

\begin{prop}\label{prop:non-simple-per-pt}
There exists $e>0$ such that the following hold:
\begin{enumerate}
	\item For all $0< E \le e$, there exist rectangles
	$R_1(E),R_2(E)\in \Sigma^{s,E}_+$ such that for $i=1, 2$,
	$\Phg^i \circ \Phl^{++} (R_i)$ intersects both $R_1(E)$ and $R_2(E)$ transversally.
	\item Given $\sigma=(\sigma_1, \cdots, \sigma_n)$, there exists a unique fixed point $p^{\sigma}(E)$ of
	\begin{equation}
	\label{eq:non-simple-composed-map}
	\prod_{i=n}^1 \left(\Phg^{\sigma_i} \circ  \Phl^{++} \right)|_{R_{\sigma_i}(E)}
	\end{equation}
	on the set $R_{\sigma_1}(E)$.
	\item The curve $p^{\sigma}(E)$ is a $C^1$ graph over the $u_1$ component
	with uniformly bounded derivatives. Furthermore, $p^{\sigma}(E)$  approaches
	$p^{\sigma_1}$ and for each $1 \le j\le n-1$,
	$$ \prod_{i=j}^1 \left(\Phg^{\sigma_i} \circ  \Phl^{++} \right) (p^\sigma(E)) $$
	approaches $p^{\sigma_{j+1}}$ as $E\to 0$.
\end{enumerate}
\end{prop}

The proof is analogous to that of Proposition~\ref{prop:fixedpt} and we omit it. 

\subsection{Normally hyperbolic invariant cylinders for the slow mechanical system}
\label{NHIC-isolating-block}

\begin{figure}[t]
\centering
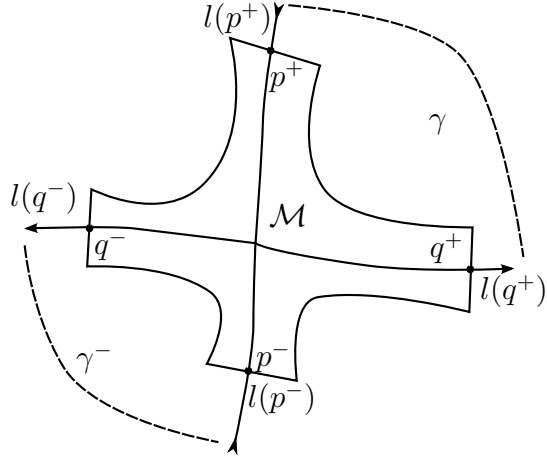
\caption{Invariant manifold $\mM$ near the origin}
\label{fig:mnfd}
\end{figure}

We prove Theorem~\ref{thm:NHIC-DR} in this section. 
\begin{proof}
\emph{Non-simple case}. 
If the homology $h$ is non-simple, then $h = n_1 h_1 + n_2 h_2$ with $h_1, h_2$ being simple homologies. Let $(\sigma_1, \dots, \sigma_n)$ be the sequence determined by Lemma~\ref{lem:homotopy}. Apply Proposition~\ref{prop:non-simple-per-pt}, we obtain the fixed points $p^\sigma(E)$ for all $0 < E < e$. The fixed points corresponds to hyperbolic periodic orbits that we call $\eta_h^E$. Let $\gamma_h^E$ be the projection of $\eta_h^E$ to the configuration space, we now prove that they must be identical to the shortest curve in Jacobi metric $g_E$, after a reparametrization. According to the condition $[DR3^c]$, $\gamma_h^0$ is the unique shortest curve for the Jacobi metric $g_0$, and there exist $c_0$ such that $\eta_h^0 = \tcA_{H^s}(c_0) = \tcN_{H^s}(c_0)$. Any $g_E$ shortest curve $\gamma'_E$ corresponds to Aubry set of cohomology $c_E$, which lifts to an orbit $\eta_E'$ in phase space. Using semi-continuity, $\eta_E'$   must be contained in a neighborhood of $\eta_h^0$ in the phase space. In particular, it must intersect the sections $\Sigma^s_+$ and $\Sigma^u_+$ sufficiently close to $q^+$ and $p^+$. According to Proposition~\ref{prop:bvp-E}, item 3, the intersection with $\Sigma^s_+$ must be contained in the rectangle $R^{++}(E)$. Since the fixed point $p^\sigma(E)$ is the unique fixed point for the map \eqref{eq:non-simple-composed-map}, we conclude that $\eta_E' = \eta_h^E$.

\emph{Simple critical case}. The existence  the periodic orbits follows from Proposition~\ref{prop:fixedpt} in the same way as the non-simple case. Also, by the same reasoning, we know that the orbits $\eta_h^E$, $\eta_{-h}^E$ must coincide with the minimal geodesics of the Jacobi metric. It suffices to show that
\[
\cM = 	\bigcup_{0 < E < e} \left( \eta_h^E \cup  \eta_{-h}^E \right) \cup \eta_h^0 \cup \eta_{-h}^0 \cup \bigcup_{-e < E < 0} \eta_c^E
\]
form a $C^1$ normally hyperbolic invariant cylinder. 

Denote
\[
l^+(p^+)=\{p^+(E)\}_{0<E\le E_0}, \quad l^+(p^-)=\{p^-(E)\}_{0<E\le E_0},
\]
$l^+(q^+) = \Phl^{++}(l^+(p^+))$ and $l^+(q^-) = \Phl^{--}(l^+(q^-))$. Note that the superscript of $l$ indicates positive energy instead of the signature of the homoclinics. We denote
\[
l^-(p^+) = \{p^c(E)\}_{-E_0 \le E <0}
\]
$l^-(q^-) = \Phl^{+-}(l^-(p^+))$, $l^-(p^-) = \Phg^-(l^-(q^-))$ and $l^-(q^+) = \Phl^{-+}(l^-(p^-))$. An illustration of $\mM$ the curves $l^\pm$ are included in Figure~\ref{fig:mnfd}. 

$l(p^+)$ and $l(p^-)$ are both graphs over the $u_1$ variable. By considering the image $\bigcup_{t >0}\phi_t (l(p^+))$ and $\bigcup_{t >0}\phi_t (l(q^+))$, using Lemma~\ref{lem:bvp-cone}, the projection of $\cM$ to $s_1u_1$ plane contains a neighborhood of $0$. Standard cone arguments implies it must be a Lipschitz graph over $s_1u_1$ near $0$. Finally, standard arguments in partial hyperbolicity (see \cite{HPS77}) implies the manifold $\cM$ is $C^{1+\alpha}$, where $\alpha$ depends on the ratio between the central and hyperbolic exponents. 

Finally, note that all the arguments applies to the normal form system \eqref{eq:norm-form2}. The same conclusions hold for small $C^2$ perturbation to the Hamiltonian (which leads to a small $C^1$ perturbation of the normal form). 
\end{proof}

\subsection{Cyclic concatenations of simple geodesics}
\label{sec:homology}

We provide the proof of the auxilliary result Lemma~\ref{lem:homotopy} before proceeding to the next section.

Denote $\gamma_1=\gamma^{h_1}_0$ and $\gamma_2=\gamma^{h_2}_0$ and $\gamma=\gamma_h^0$. Recall that $\gamma$ has homology class $n_1 h_1 + n_2 h_2$ and is the concatenation of $n_1$ copies of $\gamma_1$ and $n_2$ copies of $\gamma_2$. Since $h_1$ and $h_2$ generates $H_1(\T^2, \Z)$, by introducing a linear change of coordinates, we may assume $h_1=(1,0)$ and $h_2=(0,1)$.

Given $y\in \T^2 \setminus \gamma\cup \gamma_1 \cup \gamma_2$, the fundamental group of $\T^2\setminus \{y\}$ is a free group of two generators, and in particular, we can choose $\gamma_1$ and $\gamma_2$ as generators. (We use the same notations for the closed curves $\gamma_i$, $i=1,2$ and their homotopy classes). The curve $\gamma$ determines an element
$$\gamma = \prod_{i=1}^n \gamma_{\sigma_i}^{s_i}, \quad \sigma_i \in \{1,2\}, \, s_i \in \{0,1\} $$
of this group. Moreover, the translation $\gamma_t(\cdot):=\gamma(\cdot +t)$ of $\gamma$ determines a new element by cyclic translation, i.e.,
$$ \gamma_t = \prod_{i=1}^n \gamma_{\sigma_{i+m}}^{s_{i+m}}, \quad m\in \Z,$$
where the sequences $\sigma_i$ and $s_i$ are extended periodically. We claim the following:

There exists a unique (up to translation) periodic sequence $\sigma_i$ such that $\gamma=\prod_{i=1}^n \gamma_{\sigma_{i+m}}$ for some $m\in \Z$, independent of the choice of $y$. Note that in particular, all $s_i=1$.

The proof of this claim is split into two steps.

\emph{Step 1.} Let $\gamma_{n_1/n_2}(t)=\{\gamma(0)+(n_1/n_2, 1)t, \, t\in \R\}$. We will show that $\gamma$ is homotopic (along non-self-intersecting curves) to $\gamma_{n_1/n_2}$. To see this, we lift both curves to the universal cover with the notations $\tilde{\gamma}$ and $\tilde\gamma_{n_1/n_2}$. Let $p.q\in \Z$ be such that $ pn_1 -qn_2 =1$ and define
$$ T\tilde\gamma(t) = \tilde\gamma(t) + (p,q). $$
As $T$ generates all integer translations of $\tilde\gamma$, $\gamma$ is non-self-intersecting if and only if $T\tilde\gamma\cap \tilde\gamma = \emptyset$. Define the homotopy $\tilde\gamma_\lambda = \lambda \tilde\gamma +(1-\lambda) \tilde\gamma_{n_1/n_2}$, it suffices to prove $T\tilde\gamma_\lambda \cap \tilde\gamma_\lambda = \emptyset$. Take an additional coordinate change
$$
\begin{bmatrix}
x \\ y
\end{bmatrix} \mapsto
\begin{bmatrix}
n_1 & p \\ n_2 & q
\end{bmatrix}^{-1}
\begin{bmatrix}
x \\ y
\end{bmatrix},
$$
then under the new coordinates $T\tilde\gamma(t) = \tilde\gamma(t) +(1,0)$.

Under the new coordinates, $T\tilde\gamma\cap \tilde\gamma = \emptyset$ if and only if any two points on the same horizontal line has distance less than $1$. The same property carries over to $\tilde\gamma_\lambda$ for $0\le \lambda <1$, hence $T\tilde\gamma_\lambda \cap \tilde\gamma_\lambda = \emptyset$.

\emph{Step 2.} By step 1, it suffices to prove that  $\gamma = \gamma_{n_1/n_2}$ defines  unique sequences $\sigma_i$ and $s_i$. Since $\tilde\gamma_{n_1/n_2}$ is increasing in both coordinates, we have $s_i=1$ for all $i$. Moreover, choosing a different $y$ is equivalent to shifting the generators $\gamma_1$ and $\gamma_2$. Since the translation of the generators is homotopic to identity, the homotopy class is not affected. This concludes the proof of Lemma~\ref{lem:homotopy}.



\newpage


\section{Aubry-Mather type at the double resonance}
\label{sec:AM-double}

We consider the system 
\[
	H^s(\varphi, I) = K(I) - U(\varphi). 
\]
Given the homology $h \in H_1(\T^2, \Z)$, there is a curve $\bar{c}_h: (0, \bar{E}] \to H^1(\T^2, \R)$ such that $\bar{c}_h(E)$. There are two regimes, the ``high energy regime'', where we consider the cohomologies $\bar{c}_h(E)$ with $e \le E \le  \bar{E}$ where $e$ is a small parameter; and the critical regime, where $c$ is in a small neighborhood of $\bar{c}_h(0)$. 

\subsection{High energy case}

First we consider the ``non-critical energy case'', and show that the cohomologies as chosen are of Aubry-Mather type. For each $E> 0$, there are two possible behavior: 
\begin{enumerate}
	\item The Aubry set $\cA_{H^s}(\bar{c}_h(E)) = \gamma_h^E$, where $\gamma_h^E$ is the unique shortest geodesic in homology $E$. Let us denote the corresponding Hamiltonian orbit $\eta_h^E = (\varphi, I): [0, T_E] \to \T^2 \times \R^2$. 
	\item (Bifurcation) The Aubry set $\cA_{H^s}(\bar{c}_h(E)) = \gamma_h^E \cup \bar{\gamma}_h^E$, where $\gamma_h^E, \bar{\gamma}_h^E$ are the two shortest geodesic in homology $E$. Let us denote the corresponding Hamiltonian orbit $\eta_h^E, \bar{\eta}_h^E$. 
\end{enumerate}
\begin{thm}
	\label{thm:dr-high-AM}
	Given any $e>0$, there is $\epsilon_0, \delta > 0$ depending only on $H^s$ and $e$ such that the following holds for all $0 < \epsilon < \epsilon_0$ and all $U' \in \cV_\delta(U)$ (in the space $C^r(\T^2)$), the Hamiltonian  
	\[
		H^s_\epsilon(\varphi, I, \tau) = K(I) - U'(\varphi) + \sqrt{\epsilon} P, \quad \|P\|_{C^2} \le 1, 
	\]
	satisfies the following properties: 
	\begin{enumerate}
		\item Suppose $E \ge e$ is such that $\cA_{H^s}(\bar{c}_h(E))$ is a unique hyperbolic orbit, then $(H^s_\epsilon, \bar{c}_H(E))$ is of Aubry-Mather type. 
		\item Suppose $E \ge e$ is such that $\cA_{H^s}(\bar{c}_h(E))$ is the union of two hyperbolic orbits, then $(H^s_\epsilon, \bar{c}_h(E))$ is of bifurcation Aubry-Mather type. 
	\end{enumerate}
\end{thm}

We first discuss the non-bifurcation case, and the bifurcation case will be a simple corollary. 

Given $E_0 \ge 0$, denote $c_0 = \bar{c}_h(E_0)$. Suppose $\cA_{H^s}(c_0) = \gamma_h^E$ consists of a unique shortest curve, then there exists $E_1 < E_0 < E_2$ such that each $\eta_h^E$ for $E \in (E_1, E_2)$ is a hyperbolic periodic orbit. Then 
\[
	\cC_0 = \bigcup_{E \in (E_1, E_2)} \eta_h^E \subset \T^2 \times \R^2
\]
is a normally hyperbolic invariant cylinder. In order to give a proper parametrization for $\cC_0$, we consider a transversal section to $\eta_h^{E_0}$, which is also transversal to $\eta_h^E$, $E \in (E_1, E_2)$ if $E_1, E_2$ are close enough to $E_0$. Denote $z^E = \eta_h^E  \cap \Sigma$, then $Z^E$, $E \in (E_1, E_2)$ is a smooth function of $E$. We now define
\begin{equation}
	\label{eq:chi-0}
	\chi: \T \times (E_1, E_2) \to \T^2 \times \R^2, \quad  \chi(s, E) = \phi_{H^s}^{s T_E}(z^E),
\end{equation}
where $\phi_{H^s}^t$ is the Hamiltonian flow of $H^s$ and $T_E$ is the period of $\eta_h^E$.

\begin{lem}\label{lem:hyp-orbit}
	There is $\delta > 0$ depending on $K, U, E_1, E_2$ such that for all $U' \in \cV_\delta (U)$ and $H^s = K - U'$, 
	$\chi$ is a smooth embedding and $\chi(\T \times (E_1, E_2)) = \cC_0$. $\cC_0$ is normally hyperbolic, and for each $E$, $W^u(\eta_h^E)$ is a smooth graph over $\theta$ component on a neighborhood of $\gamma_h^E$. 
\end{lem}
\begin{proof}
	The fact that $\gamma_h^E$ is a unique non-degenerate shortest geodesic is robust, therefore $\chi$ is well defined for all $U' \in \cV_\delta (U)$, where $\delta$ depends on $K, U, E_1, E_2$. 

	The fact that $\chi$ is smooth and that the image is $\cC_0$ follows directly from definition. 
	Since $H^s(\chi(s, E)) = E$ by definition, have $\nabla H^s(\chi(s,E)) \cdot \partial_E \chi(s, E) = 1$. On the other hand, $\partial_s \chi(s, E) = X_{H^s}(\chi(s, E))$, where $X_{H^s}$ is the Hamiltonian vector field of $H^s$. Since $\nabla H^s \cdot X_{H^s} = 0$, we conclude that $\partial_s \chi, \partial_E \chi$ are linearly independent over $(s, E)\in \T \times [E_1, E_2]$.

	To see that the local stable-unstable manifolds are graphs over the $\theta$ component, we invoke the concept of Green bundles $\cG_\pm$. These  invariant bundles are defined for all orbits of the Aubry sets, and are Lipschitz graphs over the $\theta$ components. The Lipschitz constant is uniform over all $E> 0$.  According to \cite{Arn2010}, for hyperbolic periodic orbits, the bundle $\cG_-$ is the sum of the vector field direction and the unstable direction. As a result, the unstable bundle of $\eta_h^E$ projects onto a bundle transversal to $\gamma_h^E$, and the unstable manifold projects onto the $\theta$ component. 
\end{proof}

We now consider the Hamiltonian 
\begin{equation}
	\label{eq:Hs-eps}
	H^s_\epsilon(\varphi, I, \tau) = H^s(\varphi, I) + \sqrt{\epsilon} P(\varphi, I, \tau), \quad
	\varphi \in \T^2, \, I \in \R^2, \, \tau \in \T_{\sqrt{\epsilon}}. 
\end{equation}
\begin{prop}\label{prop:dr-highE-cyl}
Given any $\kappa>0$, $0 < e < \bar{E}$ there is $\epsilon_0 >0$ depending on $H^s, e, \kappa$ and $\kappa >0$ depending only on $H^s, e$; such that for all $E_0 \in [e, \bar{E}]$, there is $E_1 < E_0 < E_2$, such that for all $0 < \epsilon < \epsilon_0$, there is an embedding
\[
	\chi_\epsilon(x, y, \tau): \T \times (E_1, E_2) \times \T_{\sqrt{\epsilon}} \to \T^2 \times \R^2 \times \T_{\sqrt{\epsilon}}, \quad \pi_\tau \chi_\epsilon = id
\]
such that $C_\epsilon$ is a weakly invariant normally hyperbolic cylinder for the Hamiltonian flow of $H^s_\epsilon$ (see \eqref{eq:Hs-eps}). Moreover, we have 
\[
	\|\chi_\epsilon(\varphi, I, \tau) - \chi_0(\varphi, I)\|_{C^1} < \kappa, 
\]
and $\cC_\epsilon$ contains all the invariant sets in 
\[
	V = B_\kappa(\cC_0) \times \T_{\sqrt{\epsilon}}. 
\]
\end{prop}
Since $\| \sqrt{\epsilon} P\|_{C^2} \le \sqrt{\epsilon}$, this is a \emph{regular perturbation} of the vector field, compared to the \emph{singular perturbation} we see in Section~\ref{sec:sr-AM}. The proof follows from standard theory, see for example \cite{PSW1997}, \cite{DLS}. 

We will denote by $\chi_\epsilon^0(x, y) = \chi_\epsilon(\varphi, I, 0)$ the zero-section of the embedding, and 
\[
	\chi_\epsilon^0\left( \T \times (E_1, E_2) \right) = \cC_\epsilon^0
\]
which is invariant under the time-$\sqrt{\epsilon}$ map $\phi = \phi_{H^s_\epsilon}^{\sqrt{\epsilon}}$.

We have the following Lipschitz estimate for weak KAM solutions of $H_\epsilon^s$. 
\begin{prop}[Follows from Theorem~\ref{thm:energy-lip}]\label{prop:dr-highE-lip}
Given $R>0$,  there is $C>0$ depending only on $\|\partial^2 K\|$, $\|H^s\|_{C^2}$ and $R$, such that if $u$ is any weak KAM solution of $H^s_\epsilon$ at cohomology $c$ with $\|c\| \le R$, then for any 
\[
	(\varphi_1, I_1), \quad (\varphi_2, I_2) \in \tcI(c, u) =   \bigcap_{n \in \N}\phi^{-n}   \left(  \cG_{c, u} \right), 
\]
we have 
\[
	|H^s(\varphi_1, I_1) - H^s(\varphi_2, I_2)| \le C \epsilon^{\frac14} \|\varphi_1 - \varphi_2\| . 
\]
As a result, the same Lipschitz property holds on the sets $\tcI(c, u)$, $\tcA^0_{H^s_\epsilon}$. 
\end{prop}

\begin{proof}
	[Proof of Theorem~\ref{thm:dr-high-AM}]
	\emph{Part (1)}: 
	Let $c_0 = \bar{c}_h(E_0)$, and $\gamma_h^{E_0}$ is the unique shortest curve,  we show $(H_s^\epsilon, c_0)$ is of AM type. By Proposition~\ref{prop:dr-highE-cyl}, for sufficiently small $\epsilon$, we have $\tcN_{H^s_\epsilon}(c_0) \subset B_\kappa(\cC_0) \times \T_{\sqrt{\epsilon}}$. Due to semi-continuity, the same holds for $c \in \overline{B_\sigma(c_0)}$ for some $\sigma >0$. As a result, $\tcA_{H^s_\epsilon}(c) \subset \tcA_{H^s_\epsilon}(c) \subset \cC_\epsilon$. 

	We now prove the graph property. Suppose 
	\[
		(x_1, y_1), \quad(x_2, y_2) \quad \in (\chi_\epsilon^0)^{-1}(\tcA_{H^s_\epsilon}^0(c))
	\]
	where $c \in \overline{B_\sigma(c_0)}$. Observe that $H^s \circ \chi_0(x, y) = y$. Then due to Proposition~\ref{prop:dr-highE-lip}, 
	\[
		\|H^s \circ \chi_\epsilon^0(x, y) - y\|_{C^1} = \|H^s \circ \chi_\epsilon^0 - H^s \circ \chi_0\| \le \|H^s\|_{C^2} \|\chi_\epsilon^0 - \chi_0\|_{C^1} \le  C \delta,  
	\]
	where $C$ will be used to denote a generic constant. We have
	\begin{equation}
	  \label{eq:lip-double-res}
	  		\begin{aligned}
			& \|y_2 - y_1\| \le \|H^s \circ \chi_\epsilon^0(x_2, y_2) - H^s \circ \chi_\epsilon^0(x_1, y_1)\| 
			+  C \delta \|(x_2 - x_1, y_2 - y_1)\|  \\
			& \le C \epsilon^{\frac14} \| \pi_\varphi \circ \chi_\epsilon^0 (x_2, y_2) -  \pi_\varphi \circ \chi_\epsilon^0(x_1, y_1)\| + C \delta   \|(x_2 - x_1, y_2 - y_1)\|  \\
			& \le C (\delta + \epsilon^{\frac14}) \|(x_2 - x_1, y_2 - y_1)\| 
		\end{aligned}
	\end{equation}
	For $\epsilon, \delta$ small enough, we have $\|y_2 - y_1\| \le 2C(\delta + \epsilon^{\frac14}) \|x_2 - x_1\|  $. 

	Finally, assume that  $(\chi_\epsilon^0)^{-1}\tcA_{H^s_\epsilon}^0(c)$ project onto the $x$ component, we will show that $W^u(\tcA^0_{H^s_\epsilon}(c))$ is a graph over the $\theta$ component. According to Lemma~\ref{lem:hyp-orbit}, $W^u(\tcA^0_{H^s}(c))$ is a smooth graph over $\theta$ component, and the projection of the unstable bundle is normal to $\gamma_h^E$. Let us note that \eqref{eq:lip-double-res} implies $\tcA_{H^s_\epsilon}^0$ is close to $\tcA_{H^s}^0(c)$ in Lipschitz norm, and the unstable bundle depends smoothly on perturbation. Therefore the unstable bundle of $\tcA_{H^s_\epsilon}^0(c)$ is also transversal to $\cA_{H^s_\epsilon}(c)$, implying (2)(b) of Definition~\ref{defn:AM}. 

	\emph{Part (2)}: Suppose $c_0 = \bar{c}_h(E_0)$ is such that $\cA_{H^s}(c) = \gamma_h^{E_0} \cup \bar{\gamma}_h^{E_0}$ are two shortest curves. Let $V_1, V_2$ be neighborhoods of $\gamma^{E_0}$, $\bar{\gamma}^{E_0}$ such that $\pi_\varphi \eta^E \subset V_1$, $\pi_\varphi \bar{\eta}^E \subset V_2$ for all $E \in (E_0 - \delta, E_0 + \delta)$, and let 
\[
	f: \T^n \to [0, 1], \quad f|_{V_1} = 0, \quad f|_{V_2} = 1,
\]
be a smooth bump function. Then for 
\[
H^1 = H^s - f, \quad H^2 = H^s - (1 - f), 
\]
$\tcA_{H^1}(c_0) = \gamma^{E_0}$ and $\tcA_{H^2}(c_0) = \bar\gamma^{E_0}$. We then apply part (1) to obtain that for each $i = 1, 2$,  $H^i + \sqrt{\epsilon} P, c_0$ is of Aubry-Mather type. This implies $H^s_\epsilon, c_0$ is of bifurcation Aubry-Mather type (see Definition~\ref{def:bif}). 
\end{proof}

\subsection{Simple non-critical case}

Suppose $h$ is a simple non-critical homology, which means that for the energy $E = 0$, there is a unique shortest curve $\gamma_h^0$ in the homology $h$ corresponding to a hyperbolic periodic orbit $\eta_h^0$ of the Hamiltonian system. In this case, however, we have, for $c_0 = \bar{c}_h(0)$
\[
\cA_{H^s}(c_0) = \gamma_h^0 \cup O,
\]
where $O$ is the origin (which is where $U$ attains its minimum). If we consider $V_1, V_2$ disjoint open sets containing $\gamma_h^0$ and $0$ respectively, and define the local Aubry sets $\cA_{H^1}(c)$ and $\cA_{H^2}(c)$. It follows that $\cA_{H^1}(c) = \gamma_h^0$ and $\cA_{H^2}(c) = O$. 

\begin{thm}
\label{thm:dr-asym-bif}
Suppose $h$ is simple non-critical, and let $c_0 = \bar{c}_h(0)$. Then there is $\epsilon_0, \delta >0$ depending only on $K, U$ such that for $0 < \epsilon < \epsilon_0$ and $U' \in \cV_\delta(U)$, for
\[
H_\epsilon^s = K(I) - U'(\varphi) + \sqrt{\epsilon} P, \quad \|P\|_{C^2} \le 1, 
\]
the pair $(H^s_\epsilon, c_0)$ is of asymmetric bifurcation type. 
\end{thm}
\begin{proof}
The fact that $H^1, c_0$ is of Aubry-Mather type follows the same proof as Theorem~\ref{thm:dr-high-AM}. On the other hand, $\tcA_{H^s}(c_0)$ is a hyperbolic periodic orbit, which is robust under perturbation. Therefore Definition~\ref{def:asym-bif} is satisfied. 
\end{proof}

\subsection{Simple critical case}

\subsubsection{Proof of Aubry-Mather type using local coordinates}

\begin{thm}\label{thm:dr-am-crit}
Suppose $h \in H_1(\T^2, \Z)$ is a simple homology for $H^s$, and consider the cohomology class $c_0 = \bar{c}_h(0)$. Then there exists $\epsilon_0, e, \delta > 0$ depending only on $H^s, h$ such that for each $0 < \epsilon < \epsilon_0$, $U' \in \cV_\delta(U)$ and $c \in B_e(c_0)$, the pair $(H^s_\epsilon, c)$ is of Aubry-Mather type. 

Moreover, the same holds for $(H^s, \lambda c)$ for all $0 \le \lambda \le 1$. 
\end{thm}

We have the following (See Section~\ref{sec:NHIC-DR-proof}):
\begin{itemize}
 \item $\eta_h^0 = \tcA_{H^s}(c_0)$ contains the hyperbolic fixed point $(0,0)$, and is a homoclinic orbit to $(0, 0)$. 
 \item $(0, 0)$ admits eigenvectors $-\lambda_2 < - \lambda_1 < \lambda_1 < \lambda_2$.  Let $v_1^{s/u}$ and $v_2^{s/u}$ denote the eigendirections of the eigenvalues $\pm \lambda_1, \pm \lambda_2$.  Let $\Inv(\varphi, I) = (\varphi, - I)$ denote the involution of the Hamiltonian system. Since the flow is time-reversible, we have $\Inv(v_i^s) = \pm v_i^u$, $i = 1,2$. Without loss of generality, we assume $\Inv(v_i^s) = v_i^u$. As result, there exists $v_i, w_i \in \R^2$ such that 
 \[
 	v_i^s = (v_i, w_i), \quad v_i^u = (v_i, - w_i), \quad i = 1, 2. 
 \]

 \item There is a $C^1$ normally hyperbolic invariant manifold $\cM$ containing $\eta_h^0$. In particular, $\cM$ must contain $(0, 0)$ and it's tangent to the plane $\Span\{v_1^s, v_1^u\} = \R v_1 \oplus \R w_1 \subset \R^2 \times \R^2$ at $(0, 0)$. 

 The projection $\pi_\varphi \eta_h^0 =  \gamma_h^0$ then is a $C^1$ curve in $\T^2$, since $0 \in \T^2$ is the only possible discontinuity of the tangent direction, but at $0$ the curve is tangent to $v_1 = \pi_\varphi v_1^{s} = \pi_\varphi v_1^u$. Let $\pi_{v_1}: \R^4 \to \R$ and $\pi_{w_1} : \R^4 \to \R$ be the orthogonal projections to $v_1, w_1$ directions. 
\end{itemize}

We also need the following analog of Lemma~\ref{lem:hyp-orbit}. 
\begin{lem}
\label{lem:homoclinic-unst}
$W^u(\eta_h^0)$ is a Lipschitz graph over $\theta$ component on a neighborhood of $\gamma_h^0$. 
\end{lem}
\begin{proof}
In the proof of Lemma~\ref{lem:hyp-orbit}, the Lipschitz constant of the Green bundles is uniform over all energy $E$. The lemma follows by taking limit $E \to 0$ in the space of Lipschitz graphs. 
\end{proof}

\begin{figure}[t]
  \centering
  \includegraphics[width=2.35in]{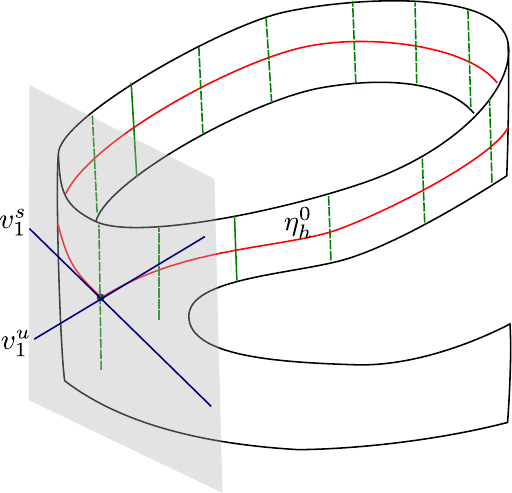}
  \caption{Coordinates near the homoclinic orbit: green lines indicate the level curves of the $x$ coordinate.}
  \label{fig:crit-coord-glob}
\end{figure}

We require a suitable parametrization of the cylinder $\cC(H^s)$ near the homoclinic $\eta_h^0$. An illustration of the parametrization is given in Figure~\ref{fig:crit-coord-glob}. 
\begin{prop}\label{prop:dr-crit-param}
For each $\kappa > 0$ there exists $\delta_1, \delta_2 > 0$, and a smooth embedding 
\[
	\chi_0= \chi_0(x, y): \quad \T \times (-\delta_2, \delta_2) \to \cM, 
\]
such that the following hold:
\begin{enumerate}
 \item $\cC(H^s): = \chi_0(\T \times (-\delta_2, \delta_2) ) \supset \eta_h^0$. We use the notation $\cC$ since the image of $\chi_0$ is a cylinder. 
 \item The cylinder $\cC(H^s)$ is symplectic. 
 \item $\chi_0$ is ``almost vertical'' near $x =0$, namely  $\pi_\varphi \circ \chi_0(x, y)$ is $\kappa-$Lipschitz in $y$ for all $|x| < \delta_1$. 
 \item The vertical coordinate is given by the energy function away from $x = 0$, namely for all $|x| \ge \delta_1/2$ we have $H^s(\chi_0(x,y)) = y$. 
\end{enumerate}
\end{prop}

We now prove Theorem~\ref{thm:dr-am-crit} assuming Proposition~\ref{prop:dr-crit-param}. 

\begin{proof}
[Proof of Theorem~\ref{thm:dr-am-crit}]
Let us consider the system 
\[
H^s_\epsilon = K(I) + U'(\varphi) + \sqrt{\epsilon} P. 
\]
where $\|P\|_{C^2} \le 1$ and $\|U - U'\|_{C^2} \le \delta$. Then for $\epsilon_0, \delta$ small enough depending only on $K, U$, standard perturbation theory implies $H_\epsilon^s$ admits a normally hyperbolic weakly invariant manifold $\cC(H^s_\epsilon)$, let us denote by $\cC_\epsilon$ its zero section. Then $\cC_\epsilon$ is invariant under the time-$\sqrt{\epsilon}$ map $\Phi = \Phi_{H^s_\epsilon}^{\sqrt{\epsilon}}$, and it admits a parametrization $\chi_\epsilon: \T \times (-\delta_2, \delta_2) \to \cC_\epsilon$, and $\|\chi_\epsilon - \chi_0\|_{C^1} = o(1)$ as $\epsilon, \delta \to 0$. The cylinder is symplectic since symplecticity is open under perturbations. Moreover, for $e>0$ small enough, the set $\tcN_{H^s}(c) \subset \cC(H^s)$ for all $c \in B_e(c_0)$, and since $\tcN$ is upper semi-continuous under Hamiltonian pertubations, $\tcN_{H^s_\epsilon}(c)$ is close to $\cC(H^s_\epsilon)$ for small $\epsilon$. Since $\cC(H^s_\epsilon)$ contains all the invariant sets in its neighborhood, we conclude that $\tcA_{H^s_\epsilon}(c) \subset \tcN_{H^s_\epsilon}(c) \subset \cC(H^s_\epsilon)$. 

We now show $\chi_\epsilon^{-1} \tcA_{H^s_\epsilon}(c)$ is a Lipschitz graph for $c \in B_e(c_0)$. Let $(\varphi_i, I_i) = \chi_\epsilon(x_i, y_i)$, $i = 1,2$ be two points in $\tcA_{H^s_\epsilon}(c)$. 

The proof consists of two cases. In the first case we use the almost verticality of the cylinder, and the idea is similar to the proof of Theorem~\ref{thm:sr-AM-type-single}. In the second case we use the strong Lipschitz estimate for the energy $H^s$, and the idea is similar to the proof of Theorem~\ref{thm:dr-high-AM}. 

\emph{Case 1}. $|x_1|, |x_2| < \delta_1$. In this case, we apply the \emph{a priori}  Lipschitz estimates for the Aubry sets: there is $C >0$  depending only $\partial^2K$, $\|H^s\|_{C^2}$ such that 
\[
	\|I_2 - I_1\| \le C\|\varphi_2 - \varphi_1\|. 
\]
Let $\kappa$ be as in Proposition~\ref{prop:dr-crit-param}, and let $\epsilon_0$ be small enough such that $\pi_\varphi \chi_\epsilon(x, y)$ is $2\kappa$-Lipschitz in $y$. Since $\chi_0$ is an embedding, there is $C>1$ depending only on $H^s$ such that for all $\epsilon_0$ small enough, 
\[
	\|\varphi_2 - \varphi_1\| + \|I_2 - I_1\| \ge C^{-1}\left( \|x_2 - x_1\| + \|y_2 - y_1\|\right). 
\]
Let $(\varphi_3, I_3) = \chi_\epsilon(x_2, y_1)$, then 
\[
	\|\varphi_3 - \varphi_2\| \le \| \chi_\epsilon(x_2, y_1) - \chi_\epsilon(x_2, y_2) \|\le 2\kappa \|y_2 - y_1\|, 
\]
\[
	\|\varphi_3 - \varphi_1\| \le \|\chi_\epsilon(x_2, y_1) - \chi_\epsilon(x_1, y_1)\| \le C\|x_2 - x_1\|, 
\]
combine all estimates, we get 
\[
\begin{aligned}
  & 	C^{-1}\left( \|x_2 - x_1\| + \|y_2 - y_1\|\right) \le (1 + C)\|\varphi_2 - \varphi_1\| \\
  & \quad \le 	(1 + C)\left(  2\kappa\|y_2 - y_1\|  + C \|x_2 - x_1\|\right),
\end{aligned}
\]
or 
\[
    \left( C^{-1} - 2(1+ C)\kappa \right) \|y_2 - y_1\| \le \left( C^{-1} + (1+ C)C \right) \|x_2 - x_1\|
\]
which is what we need if $\kappa < 1/(4C(1+C))$. 

\emph{Case 2}. $|x_1|, |x_2| > \delta_1/2$. In this case we apply Proposition~\ref{prop:dr-highE-lip}, to get 
\[
	\|H^s(\varphi_2, I_2) - H^s(\varphi_1, I_1)\| \le C \epsilon^{\frac14} \|\varphi_2 - \varphi_1\|. 
\]
Assume that $\epsilon$ is small enough such that $\|\chi_\epsilon - \chi_0\|_{C^1} < \kappa$. Then a computation identical to \eqref{eq:lip-double-res} implies $\|y_2 - y_1\| \le 2C(\kappa + \epsilon^{\frac14})\|x_2 - x_1\|$ if $\epsilon, \kappa$ is small enough depending only on $C$. 

We obtain the Lipschitz property of $\chi_\epsilon^{-1} \tcA_{H^s_\epsilon}^0(c)$ after combining the two cases. 

Finally, we show $W^u(\tcA_{H^s_\epsilon}^0(c))$ is a graph over $\varphi$ when $\chi_\epsilon^{-1}\tcA_{H^s_\epsilon}^0(c)$ projects onto the $x$ component. It suffices to check that the unstable bundle $E^u$ is uniformly transverse to the projection $\cA_{H^s_\epsilon}^0(c)$. In case 1, the almost verticality of the cylinder implies $\cA_{H^s_\epsilon}^0(c)$ differs from $\cA_{H^s}(c_0) = \gamma_h^0$ by $O(\kappa)$ in Lipschitz norm. Given that the tangent vector of $\gamma_h^0$ is close to the weak directions $v_1^{s/u}$ in case 1, while the projection of $E^u$ is close to $v_2^{s/u}$, the desired transversality holds. In case 2, similar to the proof of Theorem~\ref{thm:dr-high-AM}, $\tcA_{H^s_\epsilon}(c)$ is also close to $\tcA_{H^s}(c_0)$ in Lipschitz norm (for a different reason, i.e \eqref{eq:lip-double-res}). The claim follows, similar to the proof of Theorem~\ref{thm:dr-high-AM}, using Lemma~\ref{lem:homoclinic-unst}. 
\end{proof}

\subsubsection{Construction of the local coordinates}

 This is done separately near the hyperbolic fixed point (local) and away from it (global). Furthermore, the local coordinate requires a preliminary step. 
\begin{lem}\label{lem:crit-cor-pre}
For each $\kappa > 0$ there is $\delta > 0$, and a smooth embedding 
\[
	\chi_{\pre} = \chi_{\pre}(x, y): (-\delta, \delta) \times (-\delta, \delta) \to \cM \cap B_{2\delta}(0, 0) \subset \T^2 \times \R^2, 
\]
satisfying 
\begin{equation}
  \label{eq:almost-vertical}
  	\|\chi_{\pre} \circ (\pi_{v_1}, \pi_{w_1}) - \id\|_{C^1} < \kappa. 
\end{equation}
There is $0 < \delta_1 < \delta$ such that the curve $\chi_{\pre}^{-1} \eta_h^0 \cap \pi_x^{-1}\{(-\delta_1, \delta_1)\}$ is given by a Lipschitz graph $\{(x, g(x)): x \in (-\delta_1, \delta_1)\}$. 
\end{lem}
\begin{proof}
The existence of the local coordinate follows from $T_{(0, 0)}\cM = \R v_1 \oplus \R w_1$, the fact that $\cM$ is $C^1$, and the implicit function theorem. Since $\gamma_h^0 = \pi_\varphi \eta_h^0$ is a $C^1$ curve tangent to $v_1$ at $0 \in \T^2$, $\gamma_h^0$ can be reparametrized using its projection to $v_1$ direction. Since $\eta_h^0$ is a Lipschitz graph over $\gamma_h^0$, the second claim follows. 
\end{proof}

\begin{lem}
For each  $\kappa>0$ there is $\delta_1 > \delta_2 >0 $,  a smooth embedding
\[
	\chi_\loc = \chi_\loc(x, y) : \quad (-\delta_1, \delta_1) \times (-\delta_1, \delta_1) \to \cM \cap B_{2\delta}(0, 0), 
\]
a neighborhood $V$ of the local homoclinic $\eta_h^0 \cap \pi_x^{-1}\{(-\delta_1, \delta_1)\}$ on which the cylinder is ``almost vertical'' in the sense that there is $C>0$ such that 
\[
	\pi_\varphi \chi_\loc(x, y) \text{ is } C\kappa-\text{Lipschitz in } y. 
\]
The pull back $\chi_{\loc}^{-1} \eta_h^0 \cap \pi_x^{-1}\{(-\delta_1, \delta_1)\}$ is a Lipschitz graph over $x$,  and in addition, 
\begin{equation}
  \label{eq:almost-vertical-loc}
  	H^s(\chi_\loc(x, y)) = y, \quad \text{ for all } \quad \delta_1/2 < |x| < \delta_1, \, |y| < \delta_2.
 \end{equation}
\end{lem}
The manifold $\chi_\loc\left(-\delta_1, \delta_1) \times (-\delta_1, \delta_1)\right)$ is symplectic. 
\begin{proof}
First we show the image of $\chi_\pre$ is symplectic. 
To see this, note that $H^s(\varphi, I) = \frac12 A I \cdot I -  \frac12 B \varphi \cdot \varphi + O_3(I, \varphi)$, where $A = \partial^2_{II}K$ and $B = \partial^2_{\varphi\varphi}U(0)$ are both positive definite. Moreover, we have $\lambda_1 w_1 = A v_1$ and $\lambda_1 v_1 = B w_1$.
From \eqref{eq:almost-vertical}, we have 
\[
	\partial_x \chi_\pre(x, y) = v_1 + O(\kappa), \quad\partial_y \chi_\pre(x, y) = w_1 + O(\kappa),
\]
 As a result, let $\omega$ be the standard symplectic form, we have 
\[
	\omega(\partial_x \chi_\pre, \partial_y \chi_\pre) = \omega( (v_1, 0), (0, w_1)) + O(\kappa) = \lambda_1 v_1 \cdot A v_1 + O(\kappa) 
\]
is uniformly bounded away from $0$ if $\kappa$ is small enough.

Let $\delta, \delta_1$, $\chi_\pre$ and $g$ be from Lemma~\ref{lem:crit-cor-pre}. We claim that after possibly shrinking $\delta_1$,  there is $\delta_2 > 0$ and $C>1$ such that for all $(x, g(x) + y) \in \R^2$, $|y|< \delta_2$, $\delta_1/2 <|x| < \delta_1$ we have 
\[
	C \delta_1 > |\partial_y\left(  H^s \circ \chi_\pre(x, g(x) + y) \right)| > C^{-1} \delta_1 > 0. 
\]
We will only prove it for the case $\delta_1/2 < x < \delta$ as the other half is symmetric.

 Since $\eta_h^0$ is tangent to the stable/unstable vectors $v_1^{s/u}$, we have
\[
\begin{aligned}
 & 	\chi_\pre(x, g(x)) = x v_1^u + O(x^2) = x (v_1, w_1) + O(x^2), \\
 & \chi_\pre(x, g(x) + y) = x(v_1, w_1) + O(x^2) + O(y) , \quad x > 0. 
\end{aligned}
\]
We have 
\[
\begin{aligned}
 &    \partial_y \left( H^s \circ \chi_\pre(x, g(x) + y) \right) = 
   \left( x (Bv_1, A w_1) + O(x^2) + O(y) \right) \cdot \left( w_1 + O(\kappa) \right)  \\
   & = x A w_1 \cdot w_1 + O(x^2) + O(y) + O(\kappa (|x| + |y|)) \\
   & \ge  4C^{-1} \delta_1 + O(x^2) + O(y) + O(\kappa(|x| + |y|)) > C^{-1} \delta_1,
\end{aligned}
\]
if $8C^{-1} = \|A w_1 \cdot w_1\|$,  $\delta_1/2 < x < \delta_1$, $\delta_1 < C^{-1}$, $|y| < \delta_2 < C^{-1} \delta_1$, and $\kappa < C^{-1}/2$. The upper bound can be obtained similarly. This proves our claim.

We consider the function 
\[
	F = F(x, y) = \frac{\delta_1^{-1}}{\partial_y \left( H^s \circ \chi_\pre(x, g(x) + y) \right)}, \quad \delta_1/2 < |x|  < \delta_1, \, |y| < \delta_2. 
\]
and for each fixed $x$, let $Y_F(x, y)$ denote the solution to the ODE
\[
	\frac{d}{dy} Y_F = F(x, y), \quad Y_F(0) = g(x), 
\]
then $\partial_y H^s\circ \chi_\pre (x, Y_F(x, y)) = 1$, therefore $H^s\circ \chi_\pre (x, Y_F(x, y)) = y$.

Finally, let us define the vector field $G(x, y)$ via 
\[
G(x, y) = 
\begin{cases}
(0, F(x, y)),  &  \delta_1/2 < |x|  < \delta_1, \, |y| < \delta_2; \\
(0, 1), & |x| < \delta_1/4, \text{ or } |y| > 2\delta_2; 
\end{cases}
\]
and smoothly interpolated (keeping the first coordinate $0$) in between. Since $C^{-1} < |F_1| < C$, this can be done keeping $\|G\|_{C^1} = O(\delta_1^{-1})$.  Let $g_1(x)$ be a mollified version of $g$ such that $g_1(x) = g(x)$ for all $|x| \ge \delta_1/4$, and $|g_1(x) - g(x)|\le \delta_1/2$. Finally define
\[
	\chi_\loc(x, y) = \chi_\pre \circ \Phi_G^y(x, g_1(x)), 
\]
where $\Phi_G(x_0, y_0)$ is the time-$y$-flow of $G(x, y)$. The modification $g_1$ is to ensure the coordinate system is smooth. See Figure~\ref{fig:crit-local}. 

\begin{figure}[t]
  \centering
  \includegraphics[width=3in]{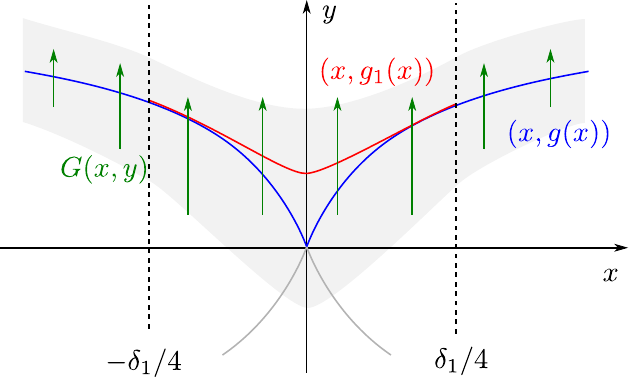}
  \caption{Construction of local coordinate}
  \label{fig:crit-local}
\end{figure}

 We have:
\begin{itemize}
 \item $\chi_\loc(x,y) = \chi_\pre(x,g_1(x) + y)$ when $|x| < \delta_1/4$ and $|y| > 2\delta_1$.
 \item $H^s \circ \chi_\loc(x, y) = H^s \circ \chi_\pre(x, Y_F(x, y)) = y$  when $\delta_1/2 < |x| < \delta_1$ and $|y| < \delta_2$. 
 \item The function $\chi_\pre$ is $\kappa$-Lipschitz in $y$ (see \eqref{eq:almost-vertical}). Therefore \eqref{eq:almost-vertical-loc} holds when $|x| < \delta_1/4$ and $|y| > 2\delta_1$.

   Moreover, given that the flow $\Phi_G^y$ is vertical: $\pi_x \Phi_G^y(x_0, y_0) = x_0$, and $|\partial_{y_0} \Phi_G^y| \le 1 +  |y| \cdot \|G\|_{C^1} = 1+  O(\delta_1^{-1} |y|)$. Therefore 
 \[
 	\left| \partial_y( \pi_\varphi\chi_\pre) \circ \Phi_G^y(x, g_1(x)) \right| \le \|\partial_y( \pi_\varphi\chi_\pre)\| \cdot \left( 1 + O(\delta_1^{-1} |y|) \right) \le 2\kappa
 \]
 if $|y| < \delta_2$ is sufficiently small. As a result, \eqref{eq:almost-vertical-loc} holds on the set $V: = \{|y| < \delta_2\}$ which is the gray area in Figure~\ref{fig:crit-local}. 
 \item The curve $\eta_h^0 = \left\{ \chi_\pre(x, g(x)): \, |x| < \delta_1 \right\}$ coincide with $\chi_\loc(x, 0)$ when $|x| \ge \delta_1/4$ and coincides with $\chi_\loc(x, y - g_1(x) + g(x))$ when $|x| \le \delta_1/4$, and therefore is a Lipschitz graph under $\chi_\loc^{-1}$ . 
\end{itemize}
\end{proof}

\begin{proof}
[Proof of Proposition~\ref{prop:dr-crit-param}]
Consider the $\chi_\loc$ coordinate system as constructed. By construction, the sections $\Sigma_\pm = \chi_\loc(\{|y| \le \delta_2, x = \pm \delta_1\})$ are transversal to the Hamiltonian flow, and therefore there is a Poincare map $\Phi: \Sigma_+ \to \Sigma_-$. Moreover, we must have $\chi_\loc^{-1} \circ \Phi \circ \chi_\loc(\delta_1, y) = (- \delta_1, y)$, since the flow preserves energy. Let $T(y)$ denote the time it takes for $\phi_{H^s}^t$ to flow from $\chi_\loc(\delta_1, y)$ to $\chi_\loc(-\delta_1, y)$.

We now define 
\[
	\chi(x, y) = 
	\begin{cases}
	\chi_\loc(x, y), & 0 \le x < \delta_1;\\
	\chi_\loc(x - 1, y), & 1 - \delta_1  \le x \le 1; \\
	\phi_{H^s}^{s(x, y)}(\delta_1, y), \quad s(x,y) = \frac{x - \delta_1}{1 - 2\delta_1}T(y), & \delta_1 \le x \le 1 - \delta_1.
	\end{cases}
\]
Then $\chi(x,y)$ is an embedding $\T \times (-\delta_2, \delta_2) \to \T^2 \times \R^2$ satisfying item (1)(3)(4) of Proposition~\ref{prop:dr-crit-param}. 

It remains to prove (2), namely symplecticity. For $|x| \le \delta_1$, this is covered in Proposition~\ref{prop:dr-crit-param}. For $|x| \ge \delta_1$, note that the tangent plane to the cylinder is spanned by two vector fields, one being $X_{H^s}$ and the other is $\partial_y \chi$. We have $\omega(X_H, \partial_y \chi) = \nabla H \cdot \partial_y \chi = 1$ by construction, therefore the manifold is symplectic. 
\end{proof}



\section{Forcing equivalence between kissing cylinders}
\label{sec:forcing-for-kissing}

In this section we prove Theorem~\ref{thm:forcing-jump}. We assume that $h$ is a non-simple homology class which satisfies $h = n_1 h_1 + n_2 h_2$, where $h_1, h_2$ are simple homologies. They are associated with curves $c_h(E)$ and $c_{h_1}^\mu(E)$ in the corresponding channels, where we assume 
\[
   \left\|  c_h(0) - c_{h_1}^\mu(0) \right\| < \mu. 
\]
Our goal is to prove forcing equivalence of the cohomologies $\Phi_L^*(c_h(E_1))$ and $\Phi^*_L(c_{h_1}^\mu(E_2))$ for some $E_1, E_2 \in (e, e+ \mu)$, where $e>0$ is sufficiently small. 

We will construct a variational problem which proves forcing equivalence for the original Hamiltonian $H_\epsilon$ using definition (Definition~\ref{defn:forcing}). The proof consists of four steps.

\begin{enumerate}
\item We construct a special variational problem for the slow
mechanical system $H^s$. A solution of this variational problem is
an orbit ``jumping'' from one homology class $h$ to the other $h_1$.
The same can be done with $h$ and $h_1$ switched.

\item We modify this variational problem for the fast
time-periodic perturbation of $H^s$, i.e. for the perturbed
slow system $H_\eps^s(\varphi^s, I^s,\tau)=K(I^s) - U(\varphi^s)+
\sqrt \eps P(\varphi^s, I^s,\tau)$ with $\tau \in \sqrt \eps \,\T$. This is achieved by applying the perturbative results established in Section~\ref{sec:pert-weak-kam}. 

 Recall the original Hamiltonian system $H_\eps$ near a double
resonance can be brought to a normal form $N_\eps=H_\eps\circ \Phi_\eps$
and this normal form, in turn, is related to the perturbed slow system
through coordinate change and energy reduction (see section \ref{slow-fast-section}). The variational problem for $H_\epsilon^s$ can then be converted to a variational problem for the original $H_\epsilon$. 

\item Using this variational problem we prove forcing relation
between $ c_h(E_1)$ and $c_{h_1}^\mu(E_2)$. 
\end{enumerate}

\subsection{Variational problem for the slow mechanical system}
\label{sec:var-slow-mech-system}

\begin{figure}[t]
  \centering
  \includegraphics[width=4in]{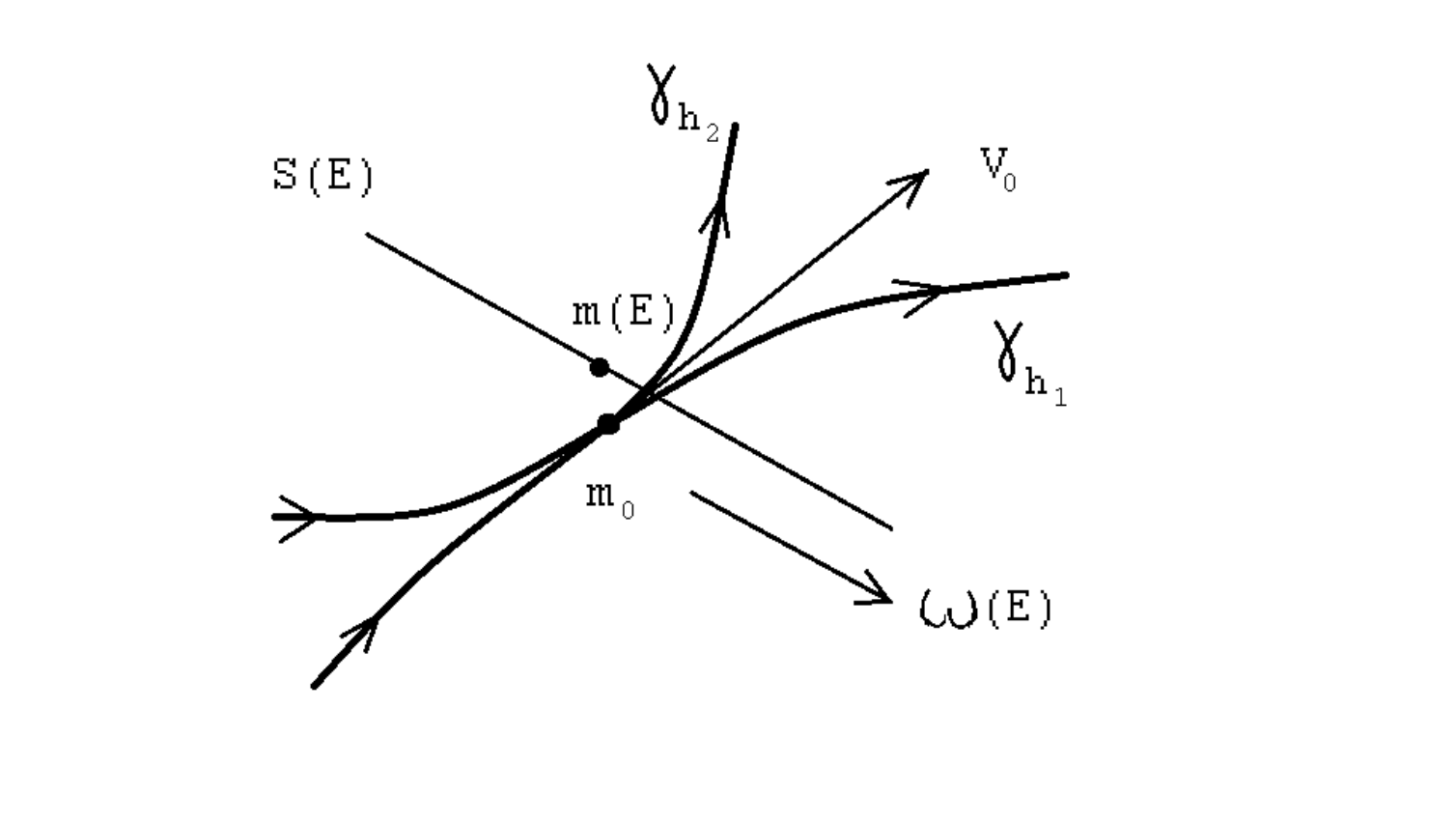}
  \caption{Jump from one cylinder to another in the same homology}
  \label{fig:local-jump-picture}
\end{figure}

The slow system is given by $H^s(\varphi, I) = K(I) - U(\varphi)$ where we assumed that the minimum of $U$ is achieved at $0$. Given $m \in \T^2$, $a > 0$ and a unit vector $\omega \in \R^2$, define 
\[
	S(m, a,\omega) = \{m + \lambda \omega: \lambda \in (-a, a)\}.
\]
$S(m, a,\omega)$ is a line segment in $\T^2$ and we will refer to it as a \emph{section} (see Figure \ref{fig:local-jump-picture}).

Given $c_1, c_2 \in \R^2$, we say that $c_1, c_2$ has a non-degenerate connection for $H^s$,  at the section $S(m, a, \omega) \in \T^2$, such that the following conditions hold.
\begin{itemize}
\item[{[N1A]}]  
\footnote{The ``A'' in [N1A]-[N3A] stands for ``autonomous''.}
\[
	(c_1 - c_2) \perp S , \quad \alpha_H(c_1) = \alpha_H(c_2). 
\]
\item[{[N2a]}] There exists a compact set $K \subset S_0$ such that for all $x \in \cA_H(c_1)$ and $z \in  \cA_H(c_2)$, the minimum of the variational problem 
\[
	\min_{y \in \overline{S}}\{ h_{H^s, c_1}(x,y) + h_{H^s, c_2}(y,z)\}
\]
is never achieved outside of $K$. 
\item[{[N3A]}]  Suppose the above minimum is achieved at $y_0$, let $p_1 - c_1$ be any super-differential of $h_{H^s,c_1}(x, \cdot)$ at $y_0$, and $-p_2 + c_2$ a super-differential of $h_{H^s,c_2}(\cdot, z)$ at $y_0$, then
\[
	\partial_p H^s(y_0, p_i) \cdot S^\perp,  \quad i=1,2
\]
have the same signs, here $S^\perp$ denote a normal vector to $S$.
\end{itemize}

\begin{rmk}
Conditions of the type [N1A]-[N3A] are common in variational construction of shadowing orbits, called the ``no corners'' conditions (see for example \cite{Bs1}) . They imply that the minimizers of the $h_{H^s, c_1}(x, y)$ and $h_{H^s, c}(y, z)$ concatenates to a smooth trajectories of the Euler-Lagrange equation. We take advantage of this fact to prove forcing relation using the defintion. 
\end{rmk}

Recall that the cohomology classes $\barc_h(E)$, $\barc_{h_1}^\mu(E)$ are chosen in the channels of $h$ and $h_1$, i.e. $\LF(\lambda_h^E h)$ and $\LF(\lambda_{h_1}^E h_1)$. Since $h$ is non-simple Proposition~\ref{prop-E0} implies the channel will pinch to a point, and $\barc_h(0)$ is uniquely chosen. The channel of $h_1$ as positive width at $E=0$, and $\barc_h(0)$ is at the boundary of this segment. In particular, since the channel for $h_1$ is parallel to $h_1^\perp$, we always have 
\begin{equation}
  \label{eq:ch1-orth}
  \barc_{h_1}^\mu(0) - \barc_h(0) \perp h_1. 
\end{equation}

We will choose $\barc_{h_1}(0)$ sufficiently close to $\barc_h(0)$ according to the following proposition. 
\begin{prop}\label{prop:critical-jump}
Suppose the slow mechanical system $H^s$ satisfies conditions [A0]-[A4] Then there is $\mu>0$ depending only on $H^s$ such that if $|\barc_{h_1}^\mu(0) - c_h(0)| < \mu$, we have: There exists a section $S(0) = S(m(0), a(0), \omega(0))$, such that conditions [N1A]-[N3A] is satisfied for $\barc_{h_1}(0)$, $\barc_h(0)$ and section $S(0)$. 

Moreover, the same holds with $\barc_{h_1}(0)$, $\barc_h(0)$ switched. 
\end{prop}
\begin{proof}
We first show that [N1A]-[N3A] are satisfied  when $c_1 = c_2 = \barc_{h}(0)$, then use perturbation arguments. Note that [N1A] is trivially satisfied when $c_1 = c_2$. 

Recall that $\cA_{H^s}(\bar{c}_h(0)) = \gamma_{h_1}^0 \cup \gamma_{h_2}^0$, and the curves $\gamma_{h_1}^0$ and $\gamma_{h_2}^0$ are tangent to a common direction
at $m_0$, which we will call $v_0$.
By the choice of $h_1$, $v_0$ is not parallel to $h_1$. We now choose $\omega(0) = \frac{h_1}{\|h_1\|}$,    $m(0)$ sufficiently close to $0$,  and $a(0)=a>0$, such that $S(m(0), \omega(0), a(0))$ intersect $\cA_{H^s}(\bar{c}_h(0))$ transversally and is disjoint from $0$ (see Figure~\ref{fig:local-jump-picture}). 

The Aubry set $\tilde\cA_{H^s}(\bar{c}_h(0))$ supports a unique minimal measure, and therefore has a unique static class. Since for any $x, z \in \cA_{H^s}(\bar{c}_h(0))$, the function 
\[
	h_{\bar{c}_h(0)}(x, \cdot) + h_{\bar{c}_h(0)}(\cdot, z)
\]
reaches its global minimal at $\cN_{H^s}(\bar{c}_h(0))$, which coincides with $\cA_{H^s}(\bar{c}_h(0))$. Hence the minimum in 
\begin{equation}
  \label{eq:min-ch0}
  \min_{y\in S(0)}\  \left\{h_{\bar{c}_h(0)}(x, y) + h_{\bar{c}_h(0)}(y,z)
\right\}
\end{equation}
is reached at $S(0)\cap \cA_{H^s}(\bar{c}_h(0))$, which is compactly contained in $S(0)$. This implies [N2a] is satisfied for $c_1 = c_2 = \bar{c}_h(0)$ along the section $S(0)$. 

Moreover, for any $y_0$ reaching the minimum in \eqref{eq:min-ch0}, then according to the above analyis $y_0$ also reaches the global minimum of $h_{\bar{c}_h(0)}(x, \cdot) + h_{\bar{c}_h(0)}(\cdot, z)$. Then using Proposition~\ref{prop:action-concave}, the associated super-gradients $p_1 - \bar{c}_h(0)$ and $-p_2 + \bar{c}_h(0)$ (in [N3A]) satisfies $p_1 = p_2$, and $\partial_p H^s(y_0, p_1)=\partial_p H^s(y_0, p_2)$ is the velocity of the unique backward minimizer at $y_0$. To show that [N3A] holds, we only need to show $\partial_p H^s(y_0, p_1) \ne 0$. This is the case because  $0$ is the only equilibrium in $\cA_{H^s}(\bar{c}_h(0))$, and $S(0)$ is disjoint from $0$. 

We now perturb $c_1$ away from $\bar{c}_h(0)$ while keeping $c_2 = \bar{c}_h(0)$. We choose $c_1 = \barc_{h_1}^\mu(0)$ at the bottom of the $h_1$ channel, then \eqref{eq:ch1-orth} and $\omega(0) = h_1$ means [N1A] is still satisfied. Proposition~\ref{unif-barrier-conv} implies that [N2a] and [N3A] are both robust under the perturbation of $c_1$, therefore [N1A]-[N3A] holds if $\| \barc_{h_1}(0) - \barc_h(0)\|$ is small enough. 

To prove the same with $\barc_{h_1}(0)$, $\barc_h(0)$ switched, we perform the same perturbation argument keeping $c_1 = \bar{c}_h(0)$ and $c_2 = \barc_{h_1}(0)$. 
\end{proof}

We now perform one more step of perturbation by taking $E>0$. 

\begin{prop}\label{var-mech} Suppose the slow mechanical system
$H^s$ satisfies conditions $[DR1^c] - [DR4^c]$, then there exists $e>0$ such
that the following hold. For each $0 \le E \le e$,
there exists a section $S(E):=S(m(E), a(E), \omega(E))$, with $S(E) \to S(0)$ in Hausdorff distance, 
and the condition [N1A]-[N3A] are satisfied for $\barc_{h_1}^\mu(E), \barc_h(E)$ at $S(E)$. 

   Moreover, the same conditions are satisfied with
   $\bar{c}_h(E)$ and $\bar{c}_{h_1}^\mu(E)$ switched.
\end{prop}

\begin{proof}[Proof of Proposition~\ref{var-mech}]
Since the Aubry sets $\cA_{H^s}(\barc_{h_1}(0))$ and $\cA_{H^s}(\barc_h(0))$ has a unique static class,  Proposition~\ref{unif-barrier-conv} implies that [N2a] and [N3A] are both robust under the perturbation of $c_1$ and $c_2$. As a result, it suffices to construct a section $S(E) = S(m(E), \omega(E), a(E))$ such that [N1A] holds, and $S(E) \to S(0)$ as $E \to 0$. 

To do this, we choose $\omega(E)$ to be a unit vector orthogonal to $\barc_{h_1}(E) - \barc_h(E)$, and by continuity of the functions $\bar{c}$,  $\omega(E) \to \omega(0)$ as $E \to 0$. Since $\alpha_{H^s}(\barc_{h_1}(E)) = \alpha_{H^s}(\barc_h(E)) = E$, [N1A] is satisfied. We then choose $m(E) = m(0)$ and $a(E) = a(0)$. Clearly $S(E)\to S(0)$ in Hausdorff metric. The proposition follows. 
\end{proof}

\subsection{Variational problem for original coordinates}
\label{sec:var-org}

The original Hamiltonian $H_\epsilon$ is reduced via coordinate change, and time change, to 
\[
	H_\epsilon^s(\varphi^s, I^s, \tau) =
	 K(I^s) - U(\varphi^s) +
	\sqrt{\epsilon}\, P(\varphi^s, I^s, \tau),
\]
with $\|P\|_{C^2} \le C_1$, see Section~\ref{sec:aff-resc}. 
The system $H_\epsilon^s$ is defined on $\T^2 \times \R^2 \times \R$, but is $\sqrt{\epsilon}$ periodic in $\tau$, i.e. it is a periodic Tonelli Hamiltonian introduced in Section~\ref{sec:intro-Tonelli}. Denote $\Te = \R/(\sqrt{\epsilon}\Z)$, then $H_\epsilon^s$ projects to $\T^2 \times \R^2 \times \Te$.

We now define a variational problem for the perturbed system. First, we will adjust the cohomologies so that they have the same alpha function. 

\begin{lem}\label{lem:pert-E}
Fix $e>0$. There exists $C>0$, and $\epsilon_0>0$ depending only on $K(I)$, $\|U\|_{C^0}$ and $e$ such that
for any $\frac{e}{3} \le E \le \frac{2e}{3}$ and
$0 \le \epsilon \le \epsilon_0$, there exists $0 < E^\epsilon < e$ such that
  $$ \alpha_{H_\epsilon^s}(\barc_h(E))= \alpha_{H_\epsilon^s}^\mu(\barc_{h_1}(E^\epsilon
  )), \quad |E-E^\epsilon| \le C \sqrt \epsilon. $$
\end{lem}
\begin{proof}
We note that $\alpha_{H^s}(c_h(E)) = E = \alpha_{H^s}(c_{h_1}^\mu(E)$. By Lemma~\ref{lem:pert-alpha}, 
\[
	\left\| \alpha_{H^s_\epsilon}(c_h(E)) - E \right\|, \quad  \left\| \alpha_{H^s_\epsilon}(c_{h_1}^\mu(E)) - E \right\| \le C\sqrt{\epsilon}. 
\]
The Lemma easily follows if we choose $\epsilon_0$ such that  $e > 3 C \sqrt{\epsilon_0}$. 
\end{proof}

We define a section $S^\epsilon(E) = S(m(E), a(E), \omega^\epsilon(E))$, by keeping $m(E)$, $a(E)$ the same as before, with $\omega^\epsilon(E)$ to be a unit vector orthogonal to $\barc_h(E) - \barc_{h_1}^\mu(E^\epsilon)$. It is natural to study the extended section $S^\epsilon \times \Te \subset \T^2 \times \Te$, and condition [N1] becomes 
\[
	\bmat{c_1 - c_2 \\ - \alpha_{H_\epsilon^s}(c_1) + \alpha_{H_\epsilon^s}(c_2)} \perp S^\epsilon \times \Te,
\]
and the variational problem for [N2] is: For $(x,0) \in \cA_{\Hse}(\barc_h(E))$ and $(z,0) \in \cA_{\Hse}(\barc_{h_1}^\mu(E^\epsilon))$, consider the minimization 
\begin{equation}
  \label{eq:var-Hse}
  	\min_{(y,t) \in S^\epsilon \times \Te} \left\{  h_{\Hse, \barc_h(E)}(x,0;y,t) + h_{\Hse, \barc_{h_1}^\mu(E^\epsilon)}(y,t; z, 0)  \right\}.
\end{equation}

We will not, however, study conditions [N1]-[N3] for $\Hse$ directly. Instead, we transform the cohomology $\barc_h(E)$ and $\barc_{h_1}(E^\epsilon)$,  the section $S^\epsilon \times \Te$, and the variational problem \eqref{eq:var-Hse} directly to the original system $H_\epsilon$.

Recall that the original Hamiltonian $H_\epsilon$ can be brought into
a normal form system $N_\epsilon$. $N_\epsilon$ is related to $H_\epsilon^s$
via the coordinate $\Phi_{L}$, and an energy reduction, see section~\ref{sec:aff-resc}.
Denote
\[
	\Phi_{L}^1(\varphi, \tau) = B^{-1} 
	\begin{bmatrix}
	\varphi \\ \tau/\sqrt{\epsilon}
	\end{bmatrix}, 
	\quad 
	\left(  \Phi_{L}^1 \right)^{-1} (\theta, t) = 
	\begin{bmatrix}
	1 & 0 \\ 0 & \sqrt{\epsilon}
	\end{bmatrix} B 
	\begin{bmatrix}
	\theta \\ t
	\end{bmatrix}, 
\]
   this is the angular component of the affine coordinate change $\Phi_{L}$ (see \eqref{eq:linear-rescale}).
   The $3 \times 3$ matrix $B$ is defined in \eqref{eq:int-matrix}.

Given $\bar{c} \in \R^2$, we define 
\begin{equation}
  \label{eq:c-eps}
  c^\epsilon = p_0 + B_0^T \bmat{ \sqrt{\epsilon} \bar{c} \\ - \epsilon \alpha_{H_\epsilon^s}(\bar{c})}
\end{equation}
where $B_0$ is the first two rows of $B$, given precisely by $k_1^T$, $k_2^T$. According to Proposition~\ref{prop:He-He-slow-var}, we have
\begin{equation}
  \label{eq:c-alpha-explicit}
	\bmat{ c^\epsilon - p_0  \\ - \alpha_{H_\epsilon}(c^\epsilon) + H_0(p_0)} = B^T \bmat{ \sqrt{\epsilon} \, \bar{c} \\ - \epsilon\,\alpha_{H^s_\epsilon}(\bar{c})} = \Phi_{L}^* (\bar{c}, \alpha_{H^s}(\bar{c})),  
\end{equation}
where $\Phi_{L}^*$ is from \eqref{eq:Phi-L-star}, and coincide with the action component of $\Phi_{L}$. In particular, \eqref{eq:c-eps} is the first row of \eqref{eq:c-alpha-explicit}. Let us denote the cohomologies $c^\epsilon_h(E)$ and $c^\epsilon_{h_1}(E)$ the image of $\bar{c} = \bar{c}_h(E)$ and $\bar{c} = \bar{c}_{h_1}^e(E)$ under \eqref{eq:c-eps}.
 
 We define a section $\Sigma = \Sigma(\theta_0, a, \Omega, l)  \subset \T^2 \times \T$ (for the original Hamiltonian $H_\epsilon: \T^2 \times \R^2 \times \T \to \R$) by 
\begin{equation}
  \label{eq:3d-section}
   \Sigma(\theta_0, a, \Omega, l) = \{(\theta_0 + \lambda \Omega + l t, t)
\in \T^2 \times \T: \, -a < \lambda < a,\ t\in \T \},
\end{equation}
where 
\[
	\Omega \in \R^3, \quad l \in \Z^3. 
\]
The section 
$S(m, a, \omega)\times \sqrt{\epsilon}\,\T \subset \T^2 \times \sqrt{\epsilon}\,\T$
is mapped under $\Phi_{L}^1$ to $\Sigma(\theta_0, a, \Omega, l)$ with 
\begin{equation}
  \label{eq:section-Phi-Ls}
  	\theta_0 = B^{-1}\bmat{m \\ 0} \in \T^3, \quad \Omega = B^{-1} \bmat{\omega \\ 0} \in \R^3, 
	\quad l = B^{-1} \bmat{0 \\ 0 \\ 1} \in \Z^3. 
\end{equation}
Note that for $(c_i^\epsilon, \alpha_{H_\epsilon}(c_i^\epsilon)) = \Phi_{L}^*(\barc_i, \alpha_{H^s_\epsilon}(\barc_i))$, $i = 1,2$, using \eqref{eq:c-eps} and $\Sigma = \Phi_{L}^1(S\times \Te)$, we have 
\begin{equation}
  \label{eq:orth-section}
  	\bmat{\barc_1 - \barc_2 \\ - \alpha_{H^s_\epsilon}(\barc_1) + \alpha_{H^s_\epsilon}(\barc_2)} \perp S \times \Te \Longleftrightarrow 
	\bmat{c_1^\epsilon - c_2^\epsilon \\ - \alpha_{H_\epsilon}(c_1^\epsilon) + \alpha_{H_\epsilon}(c_2^\epsilon)} \perp \Sigma
\end{equation}

We say that $c_1, c_2 \in \R^2$ has non-degnerate connection   along a section $\Sigma(\theta, a, \Omega, l)$, if the following conditions hold. 
\begin{enumerate}
 \item[{[N1]}] We have 
 \[
 	\bmat{ c_1 - c_2 \\ - \alpha_{H_\epsilon}(c_1) + \alpha_{H_\epsilon}(c_2)} \perp \Sigma. 
 \]
 \item[{[N2]}] There exists compact set $K \subset \Sigma$ such that: for each $(x,0) \in \cA_{H_\epsilon}(c_1)$ and $(z,0) \in \cA_{H_\epsilon(c_2)}$, the minimum in 
 \[
 	\min_{(y,t) \in \overline{\Sigma}} \left\{  h_{H_\epsilon, c_1}(x, 0; y, t) + h_{H_\epsilon, c_2}(y,t; z, 0)  \right\}
 \]
 is never achieved outside of $K$. 
 \item[{[N3]}]  Assume that the minimum [N2] is reached at $(y_0,t_0)$, and let
$p_1-c_1$ and $-p_2+c_2$ be any super-differentials of $h_{c_1}(x,0; \cdot, t)$
and $h_{c_2}(\cdot,t_0; z,0)$ respectively.  Then
  $$
    (\partial_pH(y_0, p_1, t_0), 1) \cdot \Sigma^\perp , \quad
    (\partial_p H(y_0, p_2, t_0), 1) \cdot \Sigma^\perp
    $$
 have the same signs, where $\Sigma^\perp$ is a normal vector to $\Sigma$. 
\end{enumerate}

\begin{prop}\label{int-min}
Consider $c_h^\epsilon(E), c_{h_1}^\epsilon(E)$ be defined from $\bar{c}_h(E), \bar{c}_{h_1}(E)$ using \eqref{eq:c-eps}. Let $E^\epsilon$ be as in Lemma~\ref{lem:pert-E}, and let the section $S^\epsilon(E)$ be as in \eqref{eq:var-Hse}, and the section $\Sigma^\epsilon(E)$ be obtained from $S^\epsilon(E) \times \sqrt{\epsilon}\T$ via \eqref{eq:section-Phi-Ls}. 

Then for each $0 < \epsilon < \epsilon_0$, we have
\[
	c^\epsilon_h(E), c_{h_1}^\epsilon(E^\epsilon), \Sigma^\epsilon(E)
\]
satisfies [N1]-[N3]. Moreover, the same holds with $c^\epsilon_h(E)$ and $c^\epsilon_{h_1}(E^\epsilon)$ switched. 
\end{prop}

The following statement is a direct application of Theorem~\ref{thm:jump} which holds for general Tonelli Hamiltonians. 
\begin{prop}\label{min-forcing}
Assume that the conclusions of Proposition~\ref{int-min} hold. In addition, assume that both $\cA_{H_\epsilon}(c^\epsilon_h(E))$ and $\cA_{H_\epsilon}({c}^\epsilon_{h_1}(E^\epsilon))$ admits a unique static class.  Then
  $$ c_{h}(E) \dashv \vdash c_{h_1}^{e, \mu}(E^\epsilon). $$
\end{prop}

\begin{proof}[Proof of Theorem~\ref{thm:forcing-jump}]
By  Proposition~\ref{int-min}, for the system $H_\epsilon$, which is a perturbation of $H_\epsilon$, all conditions of Proposition~\ref{min-forcing} are satisfied, except the condition of uniqueness of static classes. For this we consider again the residual condition that all rational minimal periodic orbits are hyperbolic, and has transversal homoclinic and heteroclinic intersections. Under this condition all cohomologies $c_h(E)$ and $c_{h_1}^{e, \mu}(E)$ has unique static class (using the fact that they are of Aubry-Mather type). 
\end{proof}

We prove  Proposition~\ref{int-min} in section~\ref{sec:scale-barrier} and  Proposition~\ref{min-forcing} in section~\ref{sec:jump-mech}.

\subsection{Scaling limit of the barrier function}
\label{sec:scale-barrier}

In this section we prove Proposition~\ref{int-min}.  Using the choice of the cohomology classes and the sections, together with \eqref{eq:orth-section}, we get [N1] is satisfied for $c_h^\epsilon(E)$, $c_{h_1}^\epsilon(E^\epsilon)$, and $\Sigma^\epsilon(E)$. It suffices to prove [N2] and [N3].  We will show that the variational problem  in [N2] is 
a scaling limit of the variational problem \eqref{eq:var-Hse}. 

\begin{prop}\label{scaling-limit}
The family of functions $ h_{H_\epsilon, {c}_h(E)}
(\mathcal \chi^{E,\epsilon},0; \cdot,t)/\sqrt{\epsilon}$ is uniformly semi-concave, and
  $$ \lim_{\epsilon\to 0+}
\ \sup_{(\theta,t)\in \T^2\times \T}\Bigl|
h_{H_\epsilon, c^\epsilon_h(E)}(\chi^{E,\epsilon},0; \theta,t)/\sqrt{\epsilon} -
h_{H^s, \bar{c}_{h}(E)}(x^E, B_0 \bmat{\theta \\ t})\Bigr|=0. $$
uniformly over
\[
	(\chi^{E,\epsilon},0)\in \cA_{H_\epsilon}({c}_h(E)), \quad 	(x^E,0) \in \cA_{H_\epsilon}(\bar{c}_h(E)). 
\]
 The same conclusions apply to the barrier function $h_{H_\epsilon, {c}_{h_1}^e(E)}(\cdot,\cdot; \xi^{E,\epsilon},0)/\sqrt{\epsilon}$.
\end{prop}
\begin{proof}
The uniform semi-concavity of barrier functions follows from Proposition~\ref{prop:semi-concave-near-int}. Moreover, according to Proposition~\ref{prop:He-He-slow-var}, item 2, the families of functions
\[
	h_{H_\epsilon, c^\epsilon_h(E)}(\chi^{E,\epsilon},0; \theta, t)/\sqrt{\epsilon}, \qquad
	h_{H_\epsilon^s, \barc_h(E)} \left(\Phi_{LS}^1(\chi^{E,\epsilon},0); \Phi_{LS}^1(\theta, t) \right)
\]
share the same limit points as $\epsilon \to 0$. Moreover, Proposition~\ref{prop:He-He-slow-var}, item 3 implies $\Phi_{LS}^1(\chi^{E,\epsilon},0) \in \cA_{H_\epsilon^s}(\barc_h(E))$. 

The functions $H_\epsilon^s$ form a uniform family of periodic Hamiltonians, and $H_\epsilon^s \to H^s$ in $C^2(\T^2 \times \R^2 \times \R)$. By construction, $\cA_{H^s}(\barc_h(E))$ has unique static class, therefore Proposition~\ref{unif-barrier-conv} applies, and for any $(x_E, 0) \in \cA_{H^s}(\barc_h(E))$,
\[
	h_{H_\epsilon^s, \barc_h(E)} \left(\Phi_{LS}^1(\chi^{E,\epsilon},0); \varphi, \tau) \right) \to h_{H^s, \barc_h(E)} \left(x^E,0; \varphi, \tau \right)
\]
uniformly. Finally, noticing 
\begin{equation}
  \label{eq:Phi-LS-lim}
  	\Phi_{L}^1(\theta, t) = \begin{bmatrix}
	\Id & 0 \\ 0 & \sqrt{\epsilon}
	\end{bmatrix} B 
	\begin{bmatrix}
	\theta \\ t
	\end{bmatrix} = 
	\begin{bmatrix}
	\Id & 0 \\ 0 & \sqrt{\epsilon}
	\end{bmatrix} 
	\bmat{ B_0 \\ k_3^T} 
	\begin{bmatrix}
	\theta \\ t
	\end{bmatrix} 
	 \to
	 \bmat{ B_0 \bmat{\theta \\ t} \\ 0} 
\end{equation}
as $\epsilon \to 0$, we obtain the desired limit. The case for $h_{H_\epsilon, {c}_{h_1}^e(E)}(\cdot,\cdot; \xi^{E,\epsilon},0)/\sqrt{\epsilon}$ is symmetric and we omit it. 
\end{proof}

\begin{proof}[Proof of Proposition~\ref{int-min}]
Condition [N1] is satisfied by the choice of section. Observe that 
\begin{align*}
 	B_0 \Sigma^\epsilon(E) & = \bmat{ \Id & 0 \\ 0 & 0} B \Sigma^\epsilon(E) = 
\bmat{ \Id & 0 \\ 0 & 0} B B^{-1} \bmat{\Id & 0 \\ 0 & \sqrt{\epsilon}^{-1}} (S^\epsilon(E) \times \Te)  \\
& =  \bmat{ \Id & 0 \\ 0 & 0} (S^\epsilon(E) \times \Te)  = S^\epsilon(E).
\end{align*}
Since $S^\epsilon(E) \to S(E)$ as $E \to 0$, combine with Proposition~\ref{var-mech} we get [N2] still holds for $0 \le \epsilon < \epsilon_0$, with $\epsilon_0$ depending only on $H^s$. [N3] follows from the same limiting argument and we omit it. 
\end{proof}

\subsection{The jump mechanism}
\label{sec:jump-mech}

In this section we show that [N1]-[N3] condition combined with uniqueness of static class imply $c_1 \vdash c_2$,  which is Proposition~\ref{min-forcing}. The discussions here apply to general Toneli Hamiltonians.  Denote 
\[
\tilde{A}_{H,c}(x, s, y, t) = A_{H, c}(x, s, y, t) + \alpha_H(c) (t-s). 
\]
The subscript $H$ may be omitted. 

\begin{thm}\label{thm:jump}
Assume that $c_1$, $c_2$ and $\Sigma$ satisfies the conditions [N1]-[N3], and in addition, 
\[
	\cA_{H}(c_1), \quad \cA_H(c_2)
\]
both have unique static class, then the following hold.
\begin{enumerate}
	\item (interior minimum) There exists $N<N', M < M' \in \NN$ and  a compact set $K' \subset \Sigma$, such that for
	any semi-concave function $u$ on $\T^2$, the minimum in  
	\begin{equation}\label{eq:min} 
	v(z):=\min \{ u(x)  + \tilde{A}_{c_1}(x,0; y, t+n) + \tilde{A}_{c_2}(y, t+n; z, n+m)\},
\end{equation}
where the minimum is taken in
\[
	x\in \T^2, \,(y,t)\in \Sigma, \, N \le n \le N', \, M \le m \le M',
\]
is never achieved for $(y,t)\notin K'$.
\item (no corner)
Assume the minimum in \eqref{eq:min} is achieved at $(y,t)=(y_0,t_0)$,  $(n,m)=(n_0,m_0)$, and the minimizing curves are  $\gamma_1:[0, t_0+n_0]\to \T^2$ and  $\gamma_2:[t_0+n_0,t_0 + n_0 +m_0]\to \T^2$. Then $\gamma_1$ and $\gamma_2$ connect to an orbit of  the Euler-Lagrange equation, i.e.
\[
	\dot\gamma_1(t_0+n_0) = \dot\gamma_2(t_0+n_0).
\]
\item (connecting orbits)
The function $v$ is semi-concave, and its associated pseudograph satisfies
\[
	\overline{\cG_{c_2, v}} \subset \bigcup_{0 \le t \le N'+M'}     \phi^t\cG_{c_1, u}. 
\]
As a consequence,
\[
	c_1 \vdash c_2.
\]
\end{enumerate}
\end{thm}

\begin{rmk}
We only need item 3 of the theorem for our purpose,  the first two items are stated to illuminate the idea. Item 1 can be seen as a finite time version of the variational problem in [N2]. It holds for sufficiently large $M$, $N$, due to uniform convergence of the Lax-Oleinik semi-group. Item 2 implies the minimizers concatenate to a real orbit of the system, which is crucial in proving item 3. 
\end{rmk}

\begin{lem}\label{uniform}
\begin{enumerate}
	\item Let $u$ be a continuous function on $\T^2$. The limit
	$$       \lim_{N \to \infty}\ \lim_{N' \to \infty}
	\ \min_{x\in
	\T^2,\, N \le n \le N'}\{u(x) + \tilde{A}_c(x,0;y,t+n)\} =
	$$
	$$
	=\min_{x\in D}\ \{u(x) + h_c(x,0;y,t)\} $$
	is uniform in $u$ and $(y,t)$.
	\item The limit
	$$ \lim_{N \to \infty}\ \lim_{N' \to \infty}\ \min_{ N \le n \le N'}
	\tilde{A}_c(y,t; z, n)=h_c(y,t; z,0)
	$$
	is uniform in $y,t,z$.
\end{enumerate}
\end{lem}

\begin{proof}
	The proof of the first item is similar to the proof of Proposition 6.3 of  \cite{Be} with some
	auxiliary facts proven in Appendix A there. The proof of the second item
	is similar to that of Proposition 6.1 from \cite{Be}.

	In both cases the action function, defined in (2.4) and (6.1) of \cite{Be}, is
	restricted to have integer time increment. For non-integer time increments
	the same argument applies.
\end{proof}

Using the representation formula (Proposition~\ref{prop:rep-formula}) and Lemma~\ref{unique-weak-kam}, we have the following characterization of the barrier functions.

\begin{lem}\label{barrier-split}
Assume that $\cA(c)$ has only one static class. For each point $(y,t)\in \T^2 \times \T$
and each $z\in \T^2$
\begin{enumerate}
	\item there exists $x_0 \in \T^2$ and $x_1\in \cA(c)$ such that
	$$ \min_{x\in \T^2}\ \{u(x) + h_c(x,0;y,t)\} =
	u(x_0) + h_c(x_0,0;x_1,0) + h_c(x_1, 0; y,t). $$
	\item there exists $z_1 \in \cA(c)$ such that
	$$ h_c(y,t; z,0) = h_c(y,t; z_1, 0) + h_c(z_1, 0; z,0).$$
\end{enumerate}
\end{lem}


\begin{proof}[Proof of Theorem~\ref{thm:jump}]
According to Lemma~\ref{uniform}, \eqref{eq:min} converges
uniformly as $N,M\to \infty$ to
$$\min_{x,y,t}\ \{u(x) +h_{c_1}(x,0;y,t) + h_{c_2}(y,t;z,0)\}, $$
which is equal to
\begin{multline*}
	\min_{(y,t)}\ \{u(x_0) + h_{c_1}(x_0,0; x_1,0)+
	h_{c_1}(x_1, 0;y,t) + h_{c_2}(y,t; z_1,0) + h_{c_2}(z_1,0;z,0) \} \\
	= \min_{(y,t)}\ \{ const + h_{c_1}(x_1, 0;y,t)
	+ h_{c_2}(y,t; z_1,0) + h_{c_2}(z_1,0;z,0)\}.
\end{multline*}
by Lemma~\ref{barrier-split}. Since the above variational problem
has a interior minimum due to  condition N2, by uniform convergence,
the finite-time variational problem \eqref{eq:min} also has
an interior minimum for sufficiently large $N,M$.

We now prove the second conclusion. Let $\gamma_1$ and $\gamma_2$ be
the minimizers for $\tilde{A}_{c_1}(x_0,0; y_0 , t_0 + n_0)$ and
$\tilde{A}_{c_2}(y_0, t_0 + n_0; z, n_0 + m_0)$, and let $p_1$ and $p_2$
be the associated momentum, we will show that
$$ p_1(t_0 +n_0) = p_2(t_0+n_0),$$
which implies $\dot\gamma_1(t_0+n_0)=\dot\gamma_2(t_0+n_0)$.
To abbreviate notations,  we write $p_1^0 = p_1(t_0 + n_0)$ and
$p_2^0 = p_2(t_0 + n_0)$for the rest of the proof. 

Note that
$$
u_1(x_0) + \tilde{A}_{c_1}(x_0, 0; y_0, t_0 + n_0) =
\min_{x\in \T^2}\{u_1(x) + \tilde{A}_{c_1}(x,0; y_0 , t_0 + n_0)\}.
$$
By semi-concavity, the function $u_1(x) + \tilde{A}_{c_1}(x,0; y_0 , t_0 + n_0)$
is differentiable at $x_0$ and the derivative vanishes. By Proposition
\ref{prop:action-concave} part 3,
\begin{equation}\label{eq:start-v} d_x u(x_0) = p_1(0) - c_1.
\end{equation}
By a similar reasoning, we have
\begin{equation}
	\label{eq:fin-interior-min}
	\begin{aligned}
		&         \tilde{A}_{c_1}(x_0, 0; y_0 , t_0 +n_0) +
		\tilde{A}_{c_2}(y_0, t_0 + n_0; z, n_0 + m_0) \\
		= & \min_{(y,t)\in \overline{\Sigma}} \{
			\tilde{A}_{c_1}(x_0, 0; y , t +n_0) + \tilde{A}_{c_2}(y, t + n_0; z, n_0 + m_0) \}.
		\end{aligned}
	\end{equation}

	By Proposition~\ref{prop:action-concave}, we know
	\[
		(p_1^0 - c_1,  \alpha_{H}(c_1) - H(y_0, t_0, p_0^1)), \quad
		(-p_2^0 + c_1,  -\alpha_{H}(c_2) + H(y_0, t_0, p_0^2)).
	\]
	are super-differentials of $ \tilde{A}_{c_1}(x_0, 0; y_0 , t_0 +n_0)$ and $\tilde{A}_{c_2}(y_0, t_0 + n_0; z, n_0 + m_0)$ at $(y_0, t_0 + n_0)$, since the minimum in \eqref{eq:fin-interior-min} is obtained inside of $\Sigma$, we have
	\[
		(p_1^0 - p_2^0, - H(y_0, t_0, p_1^0 + H(y_0, t_0, p_2^0))  + \left(  -c_1 + c_2, \alpha_H(c_1) - \alpha_H(c_2) \right)
	\]
	is orthogonal to $\Sigma$. Since the second term is orthogonal to $\Sigma$ by [N1], we obtain 
	\[
		(p_1^0 - p_2^0, - H(y_0, t_0, p_1^0) + H(y_0, t_0, p_2^0)) \in \R \Sigma^\perp. 
	\]
	We proceed to prove $p_1^0 = p_2^0$.

	Write $\Sigma^\perp = (v, w) \in \R^2 \times \R$, then there exists $\lambda_0$ such that 
	\[
		p_1^0 - p_2^0 = \lambda_0 v, \quad - H(y_0, t_0, p_1^0) + H(y_0, t_0, p_2^0) = \lambda_0 w. 
	\]
	Define 
	\[
		f(\lambda)  = H(y_0, t_0, p_1^0 + \lambda v) - H(y_0, t_0, p_1) + \lambda w, 
	\]
	then $f(\lambda_0) =0$. Since $f(\lambda)$ is a strictly convex on $\R$, there are at most two solutions, one of them is $\lambda=0$. Suppose there is a nonzero solution $\lambda_0$, then $f'(0)$ and $f'(\lambda_0)$ must have different signs. Since 
	\[
		f'(0) = (\partial_p H(y_0, t_0, p_1^0), 1) \cdot (v,w), \quad f'(\lambda_0) = (\partial_p H(y_0, t_0, p_2^0), 1) \cdot (v,w),
	\]
	have the same signs from [N3], we get a contradiction. As a result, $0$ is the only solution to $f(\lambda) =0$, indicating $p_1^0 = p_2^0$. Moreover, this implies $p_1^0=p_2^0$ is uniquely defined and the functions $\tilde{A}_{c_1}(x_0, 0; \cdot, \cdot)$, $A_{c_2}(\cdot, \cdot; z_0, n_0 + m_0)$ are differentiable at $(y_0, t_0 + n_0)$.

	As a consequence,  $(\gamma_1,p_1)$ and $(\gamma_2,p_2)$ connect as
	a solution of the Hamiltonian flow. Using (\ref{eq:start-v}), we have
	$$
	\phi^{n_0+m_0}(x_0, du_x(x_0) + c_1) = (z, p_2(n_0+m_0)).
	$$
	Note that $p_2(n_0+ m_0)-c_2$  is a super-differential to $v$ at $z$.
	If $v$ is differentiable at $z$, then $p_2 = dv(z) + c_2$.
	This implies
	$$
	\overline{\cG_{c_2, v}} \subset \bigcup_{0 \le t \le N'+M'}\phi^t\cG_{c_1, u}
	$$
	and the forcing relation.
\end{proof}




\newpage

\appendix


%
%
%

\section{Generic properties of mechanical systems on the two-torus}
\label{generic-mechanical-system}

Most of this section is devoted to proving Theorem~\ref{thm:DR-non-deg-HE}. At the end of the section
we prove Proposition~\ref{prop:DR-non-deg-Crit}. This states that the minimal shortest curves in a fixed homology class $h$ is generically hyperbolic, with finitely many bifurcations. To achieve this, we study generic properties of \emph{non-degenerate} orbits. It is well known that minimal non-degenerate orbits must be hyperbolic.

The proof of Theorem~\ref{thm:DR-non-deg-HE} consists of three parts.

\begin{enumerate}
	\item In section~\ref{sec:gen-per-orbits}, we prove a Kupka-Smale-like theorem
	about non-degeneracy of periodic orbits. For a fixed energy surface, generically,
	all periodic orbits are non-degenerate. This fails for an interval of energies.
	We show that while degenerate periodic orbits exists, there are only finitely many
	of them. Moreover, there could be only a particular type of bifurcation for
	any family of periodic orbits crossing a degeneracy.

	\item In section~\ref{sec:gen-min-per}, we show that a non-degenerate locally minimal orbit is always
	\emph{hyperbolic}. Using part I, we show that for each energy,  the globally minimal orbits is chosen from
	a finite family of \emph{hyperbolic} locally minimal orbits.

	\item In section~\ref{sec:proof-dr123}, we finish the proof by proving the finite local families obtained from part II are 	``in general position'', and therefore there are at most two global minimizers for each energy. 
\end{enumerate}

\subsection{Generic properties of periodic orbits}
\label{sec:gen-per-orbits}

We simplify notations and drop the supscript ``s'' from the notation of
the slow mechanical system. Moreover, we treat $U$ as a parameter, and write
\begin{equation}\label{eq:mech}
H^U(\varphi, I) = K(I) - U(\varphi), \quad
\varphi\in \T^2,\ I \in \R^2,\ U \in C^r(\T^2).
\end{equation}
We shall use $U$ as {\it an infinite-dimensional parameter}. As before $K$ is
a kinetic energy and it is fixed. Denote by $\cG^r=C^r(\T^2)$ the space
of potentials, $x$ denotes $(\varphi, I)$, $W$ denotes $\T^2\times \R^2$, and either
$\phi_t^U$ or $\Phi(\cdot, t, U)$ denotes the flow of (\ref{eq:mech}). We will
use $\chi^U(x)=(\partial K, \partial U)(x)$ to denote the Hamiltonian vector
field of $H^U$ and use $S_E$ to denote the energy surface $\{H^U=E\}$.
We may drop the superscript $U$ when there is no confusion.

By the invariance of the energy surface, the differential map $D_x\phi_t^U$
defines a map
$$
D_x\phi_t^U(x): T_xS_{H(x)} \to T_{\phi_t^U(x)}S_{H(x)}.
$$
Since the vector field $\chi(x)$ is invariant under the flow,
$D_x\phi$ induces a factor map
$$
\bar{D}_x\phi_t^U(x): T_xS_{H(x)}/\R \chi(x) \to
T_{\phi_t^U(x)}S_{H(x)}/\R \chi(\phi_t(x)).
$$
Given $U_0\in \cG^r$, $x_0\in W$ and $t_0 \in \R$, let
$$
\cV = V(x_0)\times (t_0-a, t_0+a)\times V(U_0)
$$
be a neighborhood of $(x_0,t_0, U_0)$, $V(x_0,t_0)$ of $(x_0,t_0)$, and
$V(\phi_{t_0}^{U_0}(x_0))$ a neighborhood of $\phi_{t_0}^{U_0}(x_0)$,  such that
$$
\phi_t^U(x)\in V(\phi_{t_0}^{U_0}(x_0)), \quad (x,t,U) \in \cV.
$$
By fixing the local coordinates on $V(x_0)$ and $V(\phi_{t_0}^U(x_x))$, we define
$$
\widetilde{D}_x\Phi: \cV \to Sp(2),
$$
where $\widetilde{D}_x\Phi(x,t,U)$ is the $2\times 2$ symplectic matrix associated to $\bar{D}_x\varphi_t^U(x)$ under the local coordinates. The definition depends
on the choice of coordinates.

Let $\{\phi_t^{U_0}(x_0)\}$ be a periodic orbit with minimal period $t_0$.
The periodic orbit is \emph{non-degenerate} if and only if $1$ is \emph{not}
an eigenvalue of $\widetilde{D}_x\Phi(x_0,t_0,U_0)$
\footnote{Note that we are interested in non-degeneracy for minimal period
of periodic orbits only. As the result eigenvalues given by $\exp(2\pi\,i\,p/q)$
with integer $p,q\ne 0$ are allowed}. Furthermore, we identify two types
of degeneracies:
\begin{enumerate}
	\item A degenerate periodic orbit $(x_0, t_0, U_0)$ is of \emph{type I} if $\widetilde{D}_x\Phi(x_0,t_0,U_0) =  Id$, the identity matrix;
	\item It is  of \emph{type II} if $\widetilde{D}_x\Phi(x_0,t_0,U_0)$ is conjugate
	to the matrix $[1, \mu; 0, 1]$  for $\mu \ne 0$.
\end{enumerate}
Denote
$$ N = \left\{\begin{bmatrix}
1 & \mu \\ 0 & 1
\end{bmatrix}: \mu \in \R\setminus \{0\} \right\},\quad
\mO(N) = \{ BAB^{-1}: A\in N, B\in Sp(2)\}. $$
Then $(x_0, t_0, U_0)$ is of type II if and only if
$\widetilde{D}_x\Phi(x_0,t_0,U_0)\in \mO(N)$.

\begin{lem} The set $\mO(N)$ is a $2$--dimensional submanifold of $Sp(2)$.
\end{lem}
\begin{proof}
Any matrix in $\mO(N)$ can be expressed by
$$ \begin{bmatrix}
a & b \\ c & d
\end{bmatrix}
\begin{bmatrix}
1 & \mu \\ 0  & 1
\end{bmatrix}
\begin{bmatrix}
d & -b \\ -c & a
\end{bmatrix} =
\begin{bmatrix}
1 -ac\mu & a^2\mu \\ -c^2\mu & 1-ac\mu
\end{bmatrix},
$$
where $ad-bc=1$ and $\mu \ne 0$. Write $\alpha = a^2\mu$ and
$\beta = ac\mu$, we can express any matrix in $\mO(N)$ by
\begin{equation}\label{eq:ON}
\begin{bmatrix}
1 - \beta & \alpha \\
1- \beta^2/\alpha & 1- \beta
\end{bmatrix}.
\end{equation}
\end{proof}

The standard Kupka-Smale theorem (see \cite{Oli08}, \cite{RiRu11})
{\it no longer holds for an interval of energies}. Generically, periodic
orbits appear in one-parameter families and may contain degenerate ones.
However, while degenerate periodic orbits may appear, generically,
a family of periodic orbits crosses the degeneracy ``transversally''.
This is made precise in the following theorem.

\begin{thm}\label{gen-per-orbit}
There exists residual subset of
potentials $\cG'$ of $\cG^r$, such that for all $U \in \cG'$, the following hold:
\begin{enumerate}
	\item The set of periodic orbits for $\phi_t^U$ form
	a submanifold of dimension $2$. Since a periodic orbit
	itself is a 1-dimensional manifold, distinct periodic orbits form
	one-parameter families.

	\item There is no degenerate periodic orbits of type I.

	\item The set of periodic orbits of type II form a 1-dimensional
	manifold. Factoring out the flow direction, the set of
	type II degenerate orbits are isolated.

	\item For $U_0 \in \cG^r$, let $\Lambda^{U_0} \subset W\times \R^+$
	denote the set of periodic orbits for $\phi^U_t$, and
	$\Lambda^{U_0}_N\subset \Lambda^{U_0}$ denote
	the set of type II degenerate ones. Then for any
	$(x_0, t_0)\in \Lambda_N^{U_0}$, let $V(x_0,t_0)$
	be a neighborhood of $(x_0,t_0)$. Then
	$$
	\widetilde{D}_x\Phi|_{U=U_0}: \Lambda^{U_0} \cap V(x_0,t_0) \to Sp(2)
	$$
	is transversal to $\mO(N) \subset Sp(2)$.
\end{enumerate}
\end{thm}

\begin{rmk}
Statement 4 of the theorem can be interpreted in the following way.
Let $A(\lambda)$ be the differential of the Poincare return map
on associated with a family of periodic orbits. Then if
$A(\lambda_0) \in \mO(N)$, then the tangent vector $A'(\lambda_0)$
is transversal to $\mO(N)$.
\end{rmk}
We can improve the set $\cG'$ to an open and dense set, if there is
a lower and upper bound on the minimal period.

\begin{cor}\label{gen-per-opendense}
\begin{enumerate}
	\item Given $0< T_0 < T_1$, there exists an open and dense subset
	$\cG' \subset \cG^r$, such that the set of periodic orbits with
	minimal period in $[T_0, T_1]$ satisfies the conclusions of Theorem~\ref{gen-per-orbit}.
	\item For any $U_0 \in \cG'$, there are at most finitely many type II
	degenerate periodic orbits. Moreover, there exists a neighbourhood
	$V(U_0)$ of\  $U_0$, such that the set of type II degenerate periodic
	orbits depends smoothly on $U$. (This means  the number of such periodic
	orbits is constant on $V(U_0)$, and each periodic orbit depends smoothly on $U$.)
\end{enumerate}
\end{cor}

We define
\begin{equation}\label{eq:map}
F: W \times \R^+ \times \cG^r \to W\times W,
\end{equation}
$$
F(x,t,U)=(x, \Phi(x, t, U)).
$$
$F$ is a $C^{r-1}-$map of Banach manifolds. Define the diagonal set by
$\Delta = \{(x,x)\}\subset W\times W$. Then $\{\phi_t^{U_0}(x_0)\}$
is an period orbit of period $t_0$ if and only if
$(x_0, t_0, U_0)\in F^{-1}\Delta$.
\begin{prop}\label{corank1}
Assume that $(x_0, t_0, U_0)\in F^{-1}\Delta$
or, equivalently, $x_0$ is periodic orbit of period $t_0$ for $H^U$
and that $t_0$ is the minimal period, then there exists a neighborhood
$\cV$ of $(x_0,t_0,U_0)$ such that the map
$$
d\pi^\perp_\Delta D_{(x,t,U)} F: T_{(x,t,U)}(W\times \R^+\times \cG^r)
\to T_{F(x,t,U)}(W\times W)/T\Delta
$$
has co-rank $1$ for each $(x,t,U)\in \cV$, where
$d\pi^\perp_\Delta T(W\times W) \to T(W\times W)/T\Delta$
is the standard projection along $T\Delta$.
\end{prop}
\begin{rmk}
If the aforementioned map has full rank, then the map is called
transversal to $\Delta$ at $(x_0,t_0,U_0)$. However, the transversality
condition fails for our map.
\end{rmk}
Given $\delta U\in \cG^r$, the directional derivative \ $D_U\Phi \cdot \delta U$ \
is defined as follows\
$\frac{\partial}{\partial \epsilon}|_{\epsilon=0}\Phi(x,t,U+\epsilon \delta U)$.
The differential $D_U\Phi$ then defines a map from $T\cG^r$ to $T_{\Phi(x,t,U)}W$.
The following hold for this differential map:
\begin{lem}\cite{Oli08}\label{trans-map}
Assume that there exists $\epsilon>0$ such that the orbit of $x$
has no self-intersection for the time interval $(\epsilon,\tau-\epsilon)$,
then the map
$$
D_U\Phi(x,\tau,U):\cG^r \to T_{\Phi(x,\tau,U)}W
$$
generates a subspace orthogonal to the gradient $\nabla H^U(\Phi(x,\tau,U))$
and the Hamiltonian vector field $\chi^U(\Phi(x,\tau,U))$ of $H^U$.
\end{lem}
\begin{proof}
We refer to \cite{Oli08}, Lemma 16 and 17. We note that while the proof
was written for a periodic orbit of minimal period $\tau$, the proof holds
for non-self-intersecting orbit.
	\end{proof}

	\begin{proof}[Proof of Proposition~\ref{corank1}]
	We note that if $\{\phi_t^U(x)\}$ is a periodic orbit of minimal period $\tau$,
	then the orbit  $\{\phi_t^{U'}(x')\}$ satisfies the assumptions of Lemma~\ref{trans-map}.
	It follows that the matrix
	$$
	d\pi^\perp_\Delta \circ D_UF =
	\begin{bmatrix}
	D_x\Phi -I & D_t \Phi & D_U\Phi
	\end{bmatrix}
	$$
	has co-rank $1$, since the last two component generates the subspace
	$$
	\text{Image }(D_U\Phi)+\R \chi^U,
	$$
	which is a subspace complementary to $\nabla H^U$.
	\end{proof}
	Proposition~\ref{corank1} allows us to apply the constant rank theorem
	in Banach spaces.
	\begin{prop}\label{mnfd-per}
	The set $F^{-1}\Delta$ as a subset of a Banach space is a submanifold
of codimension $2n-1$. If $r\ge 4$, then for generic $U\in \cG^r$,
$F^{-1}\Delta \cap \pi_U^{-1}\{U\}$ is a $2$-dimensional manifold.
\end{prop}
\begin{proof}
We note that the kernel and cokernel of the  map $d\pi \circ D_U F$
has finite codimension, hence the constant rank theorem
(see \cite{AMR88}, Theorem 2.5.15) applies. As a consequence, we may
assume that locally, $\Delta = \Delta_1 \times (-a, a)$ and that
the map $\pi_1\circ F$ has full rank. Since the dimension of
$\Delta_1$  is $2n-1$, $F^{-1}\Delta$ is a submanifold of
codimension $2n-1$. The second claim follows from Sard's theorem.
\end{proof}

Denote $\Lambda = F^{-1}\Delta$. On a neighbourhood $\cV$ of each
$(x_0,t_0,U_0) \in \Lambda$, we define the map
\begin{equation}\label{eq:matrix-map}
\widetilde{D}_x\Phi: \Lambda \cap \cV \to Sp(2), \quad \widetilde{D}_x\Phi(x,t,U) =
\widetilde{D}\phi_t^U(x).
\end{equation}

First we refer to the following lemma of Oliveira:
\begin{lem}[\cite{Oli08}, Theorem 18]\label{trans-det}
For each $(x_0,t_0, U_0)\in \Lambda$ such that $t_0$ is the minimal period,
let $\widetilde\cG$ be the set of tangent vectors in $T_{(x_0,t_0,U_0)}\Lambda$ of the form $(0,0,V)$.
Then the map
$$D_U\widetilde{D}_x\Phi: \widetilde\cG \to T_{\widetilde{D}_x\Phi(x_0,t_0,U_0)}Sp(2)$$  has full rank.
\end{lem}
\begin{cor}
The map \eqref{eq:matrix-map} is transversal to the submanifold $\{Id\}$ and  $\mO(N)$ of $Sp(2)$.
\end{cor}
Denote
$$\Lambda_{Id}= \Lambda \cap \widetilde{D}_x\Phi)^{-1}(\{Id\}) \text{ and }
\Lambda_N = \Lambda \cap \widetilde{D}_x\Phi)^{-1}(\mO(N)). $$
Note that the expression is well defined because both preimages are defined
independent of local coordinate changes.

\begin{proof}[Proof of Theorem~\ref{gen-per-orbit}]
The first statement of the theorem follows from Proposition~\ref{mnfd-per}.

As the subset $\{Id\}$ has codimension $3$ in $Sp(2)$, $\Lambda_{Id}$
has codimension $3$ in $\Lambda$, and hence has codimension $2n+2$ in
$W\times \R^+ \times \cG^r$. By Sard's lemma, for a generic $U\in \cG^r$,
the set $\Lambda_{Id}\cap \pi_U^{-1}$ is empty. This proves the second
statement of the theorem.

Since the set $\mO(N)$ has codimension $1$, $\Lambda_N$ has codimension $1$
in $\Lambda$, and hence has codimension $2n$ in $W\times \R^+\times \cG^r$.
As a consequence, generically, the set $\Lambda_N^U = \Lambda_N\cap \pi_U^{-1}$
has dimension $1$. This proves the third statement.

Fix $U_0\in \cG'$,  the set $\Lambda^{U_0} = \Lambda \cap \pi_U^{-1}(U_0)$
has dimension $2$, while $\Lambda_N^{U_0}\subset \Lambda^{U_0}$ has dimension $1$.
It follows that at any $(x_0,t_0) \in \Lambda_N^{U_0}$, there exists a tangent
vector
$$
(\delta x, \delta t) \in T_{(x_0,t_0)}\Lambda^{U_0}
\setminus T_{(x_0,t_0)}\Lambda_N^{U_0}
$$
such that
$$
(\delta x, \delta t, 0 ) \in T_{(x_0,t_0,U_0)}\Lambda \setminus T_{(x_0,t_0,U_0)}\Lambda_N. $$
It follows that
$$ D_{(x,t)}\widetilde{D}_x\Phi|_{U=U_0}(x_0,t_0) = D_{(x,t,U)}\widetilde{D}_x\Phi(x_0,t_0,U_0) \cdot (\delta x,
\delta t, 0) $$
is not tangent to $\mO(N)$. Since $\mO(N)$ has codimension $1$,
this implies that the map $\widetilde{D}_x\Phi|_{U=U_0}$ is transversal to $\mO(N)$.
This proves the fourth statement.
\end{proof}

\begin{proof}[Proof of Corollary~\ref{gen-per-opendense}]
	If a potential $U\in \cG'$, then by Theorem~\ref{gen-per-orbit} conditions
	1--4 are satisfied. In particular, all periodic orbits are either
	non-degenerate or degenerate satisfying conditions 3 and 4. Non-degenerate
	periodic orbits of period bounded both from zero and infinity form
	a compact set. Therefore, they stay non-degenerate for all potential
	$C^r$-close to $U$. By condition 3 degenerate periodic orbits are
	isolated. This implies that there are finitely many of them.
	Condition 4 is a transversality condition, which is $C^r$ open
	for each degenerate orbit.
	\end{proof}

	Fix $U \in \cG'$ as in Corollary~\ref{gen-per-opendense}. It follows
	that periodic orbits of $\phi_t^U$ for one-parameter families. We now
	discuss the generic bifurcation of such a family at a degenerate periodic or

	\begin{prop}\label{bif-deg}
	Let $(x_\lambda, t_\lambda)$ be a family of periodic orbits such that
	$(x_0, \lambda_0)$ is degenerate. The one side of $\lambda=\lambda_0$,
	the matrix $\widetilde{D}_x\phi_{t_\lambda}^U(x_\lambda)$ has a pair of
distinct real eigenvalues; on the other side of $\lambda=\lambda_0$,
it has a pair of complex eigenvalues.
\end{prop}
\begin{proof}
Write $A(\lambda) = \widetilde{D}_x\phi_{t_\lambda}^U(x_\lambda)$ for short. By choosing a proper local coordin
we may assume that $A(\lambda_0) = [1, \mu ; 0 , 1]$. The tangent space to $Sp(2)$ at $[1, \mu ; 0 , 1]$
is given by the set of traceless matrices $[a, b; c, -a]$. Using \eqref{eq:ON}, we have a basis of the tang
space
to $\mO(N)$ to $[1, \mu; 0, 1]$ is given by
$$ \begin{bmatrix}
0 & 1 \\ - \beta^2/\alpha^2 & 0
\end{bmatrix} \text{ and }
\begin{bmatrix}
-1 & 0 \\ - 2\beta/\alpha & -1
\end{bmatrix}. $$
An orthogonal matrix to this space, using the inner product $tr(A^TB)$, is given by $[0 , 0 ; 1, 0]$.
As a consequence, a matrix $[a, b; c, -a]$ is transversal to $\mO(N)$ if and only if $c\ne 0$.

The eigenvalues of the matrix
$$\begin{bmatrix}
1 + a h& \mu + bh\\ ch & 1 -ah
\end{bmatrix}$$
are given by $\lambda = 1 \pm \sqrt{a^2 h^2 - bch^2 -\mu c h}$. Using $\mu \ne 0$ and $c\ne 0$
we obtain that $a^2 h^2 - bch^2 -\mu c h$ changes sign at $h=0$. This proves our proposition.
\end{proof}

\subsection{Generic properties of minimal orbits}
\label{sec:gen-min-per}
In this section,  we analyze properties of families of minimal orbits. It is well known that a non-degenerate minimizing geodesic is hyperbolic. However, we have shown in the previous section that degenerate ones do exist. The main idea is when one extend a family of hyperbolic minimal orbits, the family must terminate at a degenerate orbit of type II. We then can ``slide'' different families against each other so that the degenerate orbit is never the shortest.

Let $d_E$ denote the metric derived from the Riemannian metric $g_E(\varphi, v) = 2(E + U(\varphi))K^{-1}(v)$.
We define the arclength of any continuous curve $\gamma:[t,s]\to \T^2$ by
$$ l_E(\gamma) = \sup \sum_{i=0}^{N-1}d_E(\gamma(t_i), \gamma(t_{i+1})),$$
where the supremum is taken over all partitions $\{[t_i,t_{i+1}]\}_{i=0}^{N-1}$ of $[t,s]$.
A curve $\gamma$ is called rectifiable if $l_E(\gamma)$ is finite. 

A curve $\gamma:[a,b]\to \T^2$ is called piecewise regular, if it is piecewise $C^1$ and
$\dot \gamma(t) \ne 0$ for all $t\in [a,b]$. A piecewise regular curve is always rectifiable.

A closed curve $\gamma$ is called a $(g_E, h)$ minimizer
We write
$$
l_E(h) = \inf_{\xi  \in \mC_h^E}l_E(\xi),
$$
where $\mC_h^E$ denote the set of all piecewise regular closed curves
with homology $h\in H_1(\T^2,\Z)$. A curve realizing the infimum is
the shortest geodesic curve in the homology $h$, which we will also refer
to as a global $(g_E,h)$-minimizer. 

It is well known that for any $E > -\min_\varphi U(\varphi)$,  a global $g_E-$minimizer is  a closed
$g_E-$geodesic. Hence, it corresponds to a periodic orbit of the Hamiltonian flow. 

It will also be convenient to consider the local minimizers. Let $V \subset \T^2$ be an open set, a closed continuous curve $\gamma: [a, b] \to V$ is a local minimizer if 
\begin{equation} \label{eq:local-min}
l_E(\gamma) = \inf_{\xi  \in \mC_h^E, \xi \subset \bar{V}}l_E(\xi). 
\end{equation}
Since $\gamma \Subset V$, there is an open set $V_1$ such that  $\gamma \Subset V_1  \Subset V_2$. Then by modifying the metric outside of $V_1$, we can ensure that $\gamma$ is a global minimizer of the new metric. 

\begin{prop}\label{nondeg-hyp}
Assume that $\gamma$ is a $(g_E, h)$ (local) minimizer, and assume that the associated Hamiltonian periodic orbit
$\eta$ is nondegenerate. Then $\eta$ is hyperbolic.
\end{prop}

An orbit $(\theta, p)(t), t \in \R$ of the Hamiltonian flow is called dis-conjugate if 
\[
D\phi_s \bV(x(t), p(t)) \pitchfork \bV(x(t+s), p(t+s)), \quad  t, s \in \R, \quad 
\bV(x, p) = \{0\} \times \R^2. 
\]
Let $\gamma: [0, T] \to \T^2$ be a $(g_E, h)$ minimizer, let $\eta: \R \to \T^2 \times \R^2$ be the associated periodic Hamiltonian orbit. Then:
\begin{lem}[\cite{DC, Arn2010}]\label{lem:dis-con-space}
$\eta$ is minimizing and disconjugate. Its differential map $D\phi_s(\eta(t))$ admits a $2-$dimensional invariant bundle (and the associated Poincar\'e map preserves a 1-dimensional bundle). 
\end{lem}

\begin{proof}[Proof of Proposition~\ref{nondeg-hyp}]
A non-degenerate periodic orbit $\eta:[0, T] \to \T^2 \times \R^2$ either is hyperbolic or the associated Poincare map has eigenvalues on the unit circle. In this case the Poincare map does not preserve any one-dimensional subspace, which contradicts with Lemma~\ref{lem:dis-con-space}. 
\end{proof}


\begin{thm}
\label{local-branch}
Given $0 < e < \bar{E}$, there exists an open and dense subset $\cG'$ of $\cG^r$,
such that for each $U\in \cG'$, the Hamiltonian
$H(\varphi, I) = K(I) - U(\varphi)$ satisfies
the following statements. There exists finitely many smooth families of
local minimizers
$$ \xi_j^E, \quad a_j -\sigma \le E \le b_j +\sigma, \, j = 1, \cdots, N, $$
and $\sigma>0$, with the following properties.
\begin{enumerate}
	\item All $\xi_j^E$ for $a_j - \sigma \le E \le b_j + \sigma$ are hyperbolic.
	\item $\bigcup_j [a_j, b_j] \supset [E_0, \bar{E}]$.
	\item For each $E_0 \le E \le \bar{E}$, any global minimizer is contained in the set of
	all $\xi_j^E$'s such  that  $E \in [a_j, b_j]$.
\end{enumerate}
\end{thm}

Proof of  Theorem~\ref{local-branch}, occupies the rest of this section.

\begin{lem}\label{prop-hyp-min}
Assume that $\gamma_{E_0}$ is a hyperbolic local $(\rho_{E_0}, h)-$minimizer. The following hold.
\begin{enumerate}
	\item There exists a neighbourhood $V$ of $\gamma_{E_0}$, such that $\gamma_{E_0}$ is the unique local
	$(g_E, h)$-minimizer on $V$.
	\item There exists $\delta>0$ such that for any $U'\in C^r(\T^2)$ with $\|U-U'\|_{C^2}\le \delta$ and
	$|E'-E_0|\le \delta$, the Hamiltonian $H'(\varphi, I) = K(I) - U'(\varphi)$ admits a hyperbolic local
	minimizer in $V$.
	\item There exists $\delta>0$ and a smooth family $\gamma_E\subset {V}$, $E_0 -\delta \le E \le E_0
	+\delta$, such that each of them is a hyperbolic local minimizer.
\end{enumerate}
\end{lem}
\begin{proof}
Let $\eta_{E_0}$ denote the Hamiltonian orbit of $\gamma_{E_0}$. Inverse function theorem implies that if $\eta_{E_0}$ is hyperbolic, then it is locally unique, and it uniquely extends to hyperbolic periodic orbit $\eta_{E, U'}$ if $E$ is close to $E_0$ and $U'$ is close to $U$. Let $\gamma_E$ be the projection of $\eta_E$, then they must be the unique local minimizers. 
\end{proof}

We now use the information obtained to classify the set of global minimizers.
\begin{itemize}
	\item Consider the Hamiltonian $H(\varphi, I) = K(I) - U(\varphi) +
	\min_\varphi U(\varphi)$. For $0 <E_0 < \bar{E}$, it is easy to see that
	any periodic orbit in the energy $E_0 \le E \le \bar{E}$ has a lower
	bound and upper bound on the minimal period, which depends only on $E_0$ and $\bar{E}$. Hence, Corollary~\ref{gen-per-opendense} applies. In particular, there are only finitely many degenerate periodic orbits and the number of them is constant under small perturbation of $U$. 

	\item By the previous item, all but finitely many global minimizers are nondegenerate, hence hyperbolic by Proposition~\ref{nondeg-hyp}. 

	\item Since a global minimizer is always a local minimizer, using
	Lemma~\ref{prop-hyp-min}, it extends to a smooth one parameter family of
	local minimizers. The extension can be continued until either the orbit is no-longer locally minimizing, or if the orbit becomes degenerate. Once a local minimizer becomes degenerate, this family can no longer be extended
	as local minimizers as by Proposition~\ref{bif-deg}, it must bifurcate to a periodic orbit of complex eigenvalues.

	\item It is well known that for a fixed energy, any two global minimizers do not 
	cross (see for example, \cite{MF}). We can locally extend the global minimizers for a small interval of energy without them crossing. 

	\item There are at most finitely many families of local minimizers, because 
	they are isolated  and do not accumulate (Lemma~\ref{prop-hyp-min}). 

	\item There are at most finitely many energies on which the global minimizer may be a degenerate periodic orbit. 
\end{itemize}

We have proved the following statement.

\begin{prop}\label{all-local-min}
Given $0 < e < \bar{E}$, there exists an open and dense subset $\cG'$ of $\cG^r$,
such that for each $U\in \cG'$, for the Hamiltonian
$H(\varphi, I) = K(I) - U(\varphi) + \min_\varphi U(\varphi)$, such that the
following hold.

\begin{enumerate}
	\item There are at most finitely many (maybe none) isolated global minimizers 
	$\gamma_j^{E_j,d}$, $j = 1, \cdots, M$.
	\item There are finitely many smooth families of local minimizers
	$$ \gamma_j^E, \quad \bar{a}_j \le E \le \bar{b}_j, \, j = 1, \cdots, N, $$
	with $[e, \bar{E}] \supset \bigcup [\bar{a}_j, \bar{b}_j]$, such that
	$\gamma_j^E$ are hyperbolic for
	$\bar{a}_j < E < \bar{b}_j$. The set $\bar\gamma_j^E$ for $E=\bar{a}_j, \bar{b}_j$ 
	may be hyperbolic or degenerate.
	\item For a fixed energy surface $E$, the sets $\{\gamma_j^{E,d}\}$ and 
	$\bigcup_{\bar{a}_j\le E \le \bar{b}_j}\gamma_j^E$ are pairwise disjoint.
	\item For each $e \le E \le \bar{E}$, the global minimizer is chosen among
	the set of all $\gamma_j^{c_j,d}$  with $E = c_j$, or one of
	the local minimizers $\gamma_j^E$ with $E\in [\bar{a}_j, \bar{b}_j]$.
\end{enumerate}
\end{prop}

\begin{proof}[Proof of Theorem~\ref{local-branch}]
We first show that the set of potentials $U$ satisfying the conclusion of
Theorem~\ref{local-branch} is open. By Lemma~\ref{prop-hyp-min},the family
of local minimizers persists under small perturbation of the potential $U$.
It suffices to show that for sufficiently small perturbation of $U$ satisfying
the conclusion, the global minimizer is still taken at one of the local families.
Assume, by contradiction, that there is a sequence $U_n$ approaching $U$, and
for each $H_n = K-U_n$, there is some global minimizer $\gamma_n^{E_n}$ not from
these families. By picking a subsequence, we can assume that it converges to
a closed curve $\gamma_*$, which belong one of the local families $\gamma_j^E$.
Using local uniqueness from Lemma~\ref{prop-hyp-min}, $\gamma_n^{E_n}$ must
belong to one of the local families as well. This is a contradiction. 

To prove density, it suffices to prove that for a potential $U$ satisfying
the conclusion of Proposition~\ref{all-local-min}, we can make an arbitrarily
small perturbation, such that there are no degenerate global minimizers.

Our strategy is to eliminate the degenerate global minimizers one by one using 
a sequence of perturbations. Suppose $\gamma_E^d$ is an isolated degenerate minimizer. Then by (\cite{CI}, proof of Theorem D, page 40-42), there is a perturbation $\delta U$ with the property that: $\delta U \ge 0$, $U|{\gamma_E^d} = 0$, such that $\gamma_E^d$ is still a global minimizer, but the associated Hamiltonian orbit $\eta_E^d$ becomes hyperbolic (i.e. the perturbation keeps the orbit $\eta_E^d$ intact but changes it's linearization). This perturbation strictly decreases the number of of isolated degenerate minimizers by $1$. 

Suppose $\gamma_E^d$ is the terminal point of one of the local families $\gamma_j^E: E \in [\bar{a}_j, \bar{b}_j]$, let's say $\gamma_E^d = \gamma_j^E(\bar{b}_j)$. Then since  claim $\gamma_j^E$ can not be extended to $E > \bar{b}_j$, the global minimizer for $E \in (\bar{b_j} , \bar{b}_j + \delta)$ must be contained in a different local branch, say $\gamma_i^E: E \in [\bar{a}_i, \bar{b_i}]$. In particular, $\gamma_i^E(\bar{b}_j)$ is another global minimizer. 

Let $V$ be an neighborhood of $\gamma_j^{E_j}$, such that $\bar{V}$ is disjoint from
the set of other global minimizers with the same energy. For $\delta>0$ we define
$U_\delta: \T^2 \to \R$ such that $U_\delta|\gamma_j^{E_j} = \delta$ and
$\supp U_\dt \subset V$. Let $H_\delta = K - U - U_\delta$, and let $l_{E, \delta}$
be the perturbed length function. We have
$$
l_{E_j,\delta}(\gamma_j^{E_j}) = \int_{\gamma_j^{E_j}} \sqrt{2(E_j+U+\delta)K} > l_{E_j}(\gamma_j^{E_j}) = l_{E_j}(h)= l_{E_j,\delta}(h).
$$
As a consequence, $\gamma_j^{E_j}$ is no longer a global minimizer for
the perturbed system. Moreover, for sufficiently small $\delta$, no
new degenerate global minimizer can be created. Hence the perturbation
has decreased the number of degenerate global minimizers strictly.
By repeating this process finitely many times, we can eliminate all degenerate
global minimizers.
\end{proof}

\subsection{Non-degeneracy at high energy}
\label{sec:proof-dr123}

In this section we complete the proof of Theorem~\ref{thm:DR-non-deg-HE}. This amounts to proving
that finite local families of local minimizers, obtained from the previous section,
are ``in general position''.

We assume that the potential $U_0 \in \cG^r$  satisfies the conclusions of 
Theorem~\ref{local-branch}. Let $\gamma_j^{E, U}$ denote the branches of local
minimizers from Theorem~\ref{local-branch}, where we have made the dependence
on $U$ explicit. There exists an neighbourhood $V(U_0)$ of $U_0$, such that
the local branches $\gamma_j^{E,U}$ are defined for
$E \in [a_j - \sigma/2, b_j + \sigma/2]$ and $U \in V(U_0)$.

Define a set of functions
$$ f_j: [a_j - \sigma/2, b_j + \sigma/2] \times V(U_0) \to \R, \quad 
f_j(E,U) = l_E(\gamma_j^{E,U}). $$
Then $\gamma_i^{E,U}$ is a global minimizer if and only if
$$
f_i(E,U) = f_{\min}(E,U) :=\min_j f_j(E,U),
$$
where the minimum is taken over all $j$'s such that $E \in [a_j, b_j]$.

The following proposition implies Theorem~\ref{thm:DR-non-deg-HE}.
\begin{prop}
There exists an open and dense subset $V'$ of $V(U_0)$ such that for every
$U \in V'$, the following hold:
\begin{enumerate}
	\item For each $E \in [E_0, \bar{E}]$, there at at most two $j$'s such that
	$f_j(E,U) = f_{\min}(E,U)$;
	\item There are at most finitely many $E\in [E_0, \bar{E}]$ such that there
	are two $j$'s with
	$f_j(E,U) =    f_{\min}(E,U)$;
	\item For any $E\in [E_0, \bar{E}]$ and $j_1,j_2$ be such that 
	$f_{j_1}(E,U) = f_{j_2}(E,U) = f_{\min}(E,U)$;
	we have
	$$ \frac{\partial}{\partial E}f_{j_1}(E,U) \ne \frac{\partial}{\partial E} f_{j_2}(E,U) .$$
\end{enumerate}
\end{prop}

\begin{proof}
We first note that it suffices to prove the theorem under the additional
assumption that all functions $f_j$'s are defined on the same interval
$(a,b)$ with $f_{\min}(E,U) = \min_j f_j(E,U)$. Indeed,  we may partition
$[E_0, \bar{E}]$ into finitely many intervals, on which the number of local
branches is constant, and prove proposition on each interval.

We define a map
$$
f = (f_1, \cdots, f_N): (a,b)\times V(U_0) \to \R^N,
$$
and  subsets
$$
\Delta_{i_1,i_2,i_3} = \{(x_1, \cdots, x_n); x_{i_1}=x_{i_2}=x_{i_3}\},
\quad 1\le i_1< i_2< i_3\le N,
$$
$$
\Delta_{i_1,i_2} = \{(x_1, \cdots, x_n); x_{i_1}=x_{i_2}\},
\quad 1\le i_1 < i_2 \le N
$$
of $\R^N \times \R^N$. We also write $f^U(E)=f(E,U)$.
The following two claims imply our proposition:
\begin{enumerate}
	\item For an open and dense set of $U\in V(U_0)$, for all
	$1\le i_1< i_2< i_3\le N$, the set $(f^U)^{-1}\Delta_{i_1,i_2,i_3}$ is empty.
	\item For an open and dense set of $U\in V(U_0)$, and all
	$1\le i_1 < i_2 \le N$, the map $f^U:(a,b)\to \R^N$ is transversal to
	the submanifold $\Delta_{i_1,i_2}$.
\end{enumerate}
Indeed, the first claim imply the first statement of our proposition.
It follows from our second claim that there are at most finitely many points
in $(f^U)^{-1}\Delta_{i_1,i_2}$, which implies the second statement. Furthermore,
using the second claim, we have for any $E\in (f^U)^{-1}\Delta_{i_1,i_2}$,
the subspace  $(Df^U(E))\R$ is transversal to $T\Delta_{i_1,i_2}$. This implies
the third statement of our proposition.

For a fixed energy $E$ and $(v_1, \cdots, v_N) \in \R^N$,
let $\delta U : \T^2 \to \R$ be such that $\delta U(\varphi) = v_j$
on an open neighbourhood of $\gamma_j^E$ for each $j =1, \cdots N$.
Let $l_{E,\epsilon}$ and  $\gamma_j^{E,\epsilon}$ denote the arclength and
local minimizer corresponding to the potential $U + \epsilon \delta U$.
For any $\varphi$ in a neighbourhood of $\gamma_j^E$, we have
$$
E + U(\varphi) + \delta U(\varphi) = E + U(\varphi) + \epsilon v_j,
$$
hence for sufficiently small $\epsilon >0$, $\xi_j^{E,\epsilon} =
\xi_j^{E+\epsilon v_j}$.

The directional derivative
$$
D_U f(E,U)\cdot \delta U = \frac{d}{d\epsilon}\Bigr|_{\epsilon =0}
l_{E,\epsilon}(\xi_j^{E,\epsilon}) = \frac{d}{d\epsilon}\Bigr|_{\epsilon =0}
l_{E+ \epsilon v_j}(\xi_j^{E+ \epsilon v_j})=\frac{\partial}{\partial E}f_j(E,U)v_j.
$$
It follows from a direct computation that each $f_j$ is strictly increasing in
$E$ and the derivative in $E$ never vanishes. As a consequence, we can choose
$(v_1, \cdots, v_N)$ in such a way, that $D_U f(E,U)\cdot \dt U$ takes
any given vector in $\R^N$. This implies the map
$$ D_U f:(a,b) \times TV(U_0) \to \R^N $$
has full rank at any $(E,U)$. As a consequence, $f$ is transversal to any
$\Delta_{i_1, i_2, i_3}$ and
$\Delta_{i_1,i_2}$. Using Sard's lemma, we obtain that for a generic $U$,
the image of  $f^U$ is disjoint from
$\Delta_{i_1, i_2, i_3}$ and that $f^U$ is transversal to $\Delta_{i_1, i_2}$.
\end{proof}

\subsection{Unique hyperbolic minimizer at very high energy}
\label{sec:very-high-energy}

In this section we prove Proposition~\ref{prop:DR-very-high-energy}, namely for the mechanical system 
\[
H^s(\varphi, I) = K(I) + U(\varphi), 
\]
the shortest loop in the homology $h \in H_1(\T^2, \Z)$ for  Jacobi metric $g_E(\varphi, v) = 2(E +  U(\varphi)) K^{-1}(v)$, where $K^{-1}(v)$ is the Legendre dual of $K$. Since the metrics
\[
g_E(\varphi, v) = 2(1 + E^{-1} U(\varphi)) K^{-1}(\sqrt{E}v), \quad 
\bar{g}_E(\varphi, v) = 2( 1 + E^{-1} U(\varphi))K^{-1}(v)
\]
differ by a constant multiple, it suffices to study the shortest loops for the metric $\bar{g}_E$. Denote $\delta = E^{-1}$, this is then equivalent to study the mechanical system
\[
 H_\delta(\varphi, I) = K(I) + \delta U(\varphi)
 \] 
 with energy $E= 1$. Without loss of generality, we may assume $h = (1,0)$. Let us write 
 \[
 U(\varphi)= U_1(\varphi_1) + U_2(\varphi_1, \varphi_2), 
 \]
 where $\int_0^1 U_2(\varphi_1, s)ds =0$. We will show an averaging effect that eliminates $U_2$. Let us denote by $\rho_0$ the unique positive number such that $K^{-1}(\rho_0 h) = 1$, and let $c_0 = \partial K^{-1}(\rho_0 h)$.  

 \begin{lem}\label{lem:mech-averaging}
 	 	There is $C>1$ depending only on $\|\partial^2 K\|, \|\partial^2 K^{-1}\|$ such that
 	 	 then the associated Hamiltonian orbit $\eta_h^1 = \{(\varphi, I)(t)\}$ satisfies
  \[
  \|\dot{\varphi} - \rho_0 h\|, \, \|I - c_0\| \le C \delta. 
  \]
 \end{lem}
 \begin{proof}
Let us note the shortest curve $\gamma_h^1$ corresponds to the orbit $\eta_h^1 = (\varphi_h, I_h)$ of the Hamiltonian $H_\delta$, satisfying $\varphi_h(T) = \varphi_h(0) + h$. Note:
\begin{enumerate}
 \item $\ddot\varphi_h(t) = \delta\nabla U(\varphi_h(t)) = O(\delta)$, 
 \item Since the energy is $1$, $K(I_h(t)) = 1 + O(\delta)$, $K^{-1}(\dot{\varphi}_h(t)) = K^{-1}(\partial_I K(I_h(t))) = K(I_h(t)) = 1 + O(\delta)$; 
\end{enumerate}
The first item implies $\dot{\varphi}_h(t) = \frac{1}{T} h + O(\delta)$, when combined with the second item, implies $\dot{\varphi}_h(t) = \rho_0 h + O(\delta)$. apply $\partial K^{-1}$, we get $I_h(t) = c_0 + O(\delta)$. 
 \end{proof}

Note that the condition $\partial K(I) = \rho h$ corresponds precisely to the resonance segment $\Gamma = {I\in \R^2, \, \partial K(I) \cdot (1, 0) =0}$, in other words, we are in precisely the same situation as in single resonance. In fact, since there are fewer degree of freedom, the normal form is much stronger. 
\begin{prop}
For any $M >0$, there is $\delta_0 = \delta_0(K, M)$ and $C = C(K, M) >0$ such that for every $H_\delta$ with $0 < \delta < \delta_0$, there is a $C^\infty$ coordinate change $\Phi$ homotopic to identity such that 
\[
H_\delta \circ \Phi = K(I) + \delta U_1(\varphi_1) + \delta R(\varphi), 
\]
with $\|R\|_{C^2_I} \le  \sqrt{\delta} $, and $\|\pi_\theta(\Phi- Id)\|_{C^2} \le C \delta$ and $\|\pi_p(\Phi - Id)\|_{C^2} \le C \delta$, where the norms are measured on the set $I \in B_{M\sqrt{\delta}}(c_0)$. 

If $U_1(\varphi_1)$ has a unique non-degenerate minimum at $\varphi_1 = \varphi_1^*$, and satisfies
\[
U(\varphi_1) - U(\varphi_1^*) \ge \lambda \|\varphi_1 - \varphi_1^*\|^2. 
 \] 
 Then by choosing $\delta_0 = \delta_0(K, M, \lambda) >0$ smaller, the system $H_\delta \circ \Phi$ has a normally hyperbolic invariant cylinder $\cC$ given by 
\[
(\varphi_1, I_1) = (F(\varphi_2, I_2) + \varphi_1^*, G(\varphi_2, I_2)), \quad \varphi_2 \in \T,  I_2 \in \Gamma \cap B_{M\sqrt{\delta}}(c_0),
\]
satisfying $\|F\|_{C^0} \le C \kappa$, $\|G\|_{C^0} \le C \kappa \sqrt{\delta}$. $\cC$ contains all the invariant sets in 
\[
\|\varphi_1 - \varphi_1^*\| < C^{-1}, \quad \dist(I, \Gamma\cap B_{M\sqrt{\delta}}) < C^{-1}. 
\]

The minimal periodic orbit $\eta_h^1$ is contained in $\cC$ and is a graph over $\varphi_2$. It is therefore unique and hyperbolic. 
\end{prop}
\begin{proof}
The first part is very similar to Theorem~\ref{thm:normal-form-sr}, but in the much simpler case of two-degrees of freedom. In this case, we can give out the coordinate change explicitly. Define 
\[
W(\varphi_1, \varphi_2) =  \rho_0^{-1} \int_0^1 s U_2(\varphi_1, \varphi_2 + s) ds, 
\]
note that $W$ is well defined on $\T^2$ due to $\int_0^1 U_2(\varphi_1, s) ds = 0$. Note that $\partial_{\varphi_2}W = \rho_0^{-1} U_2$.  Consider the coordinate change $I \mapsto I + \delta \nabla W$. Then the new Hamiltonian is 
\[
\begin{aligned}
 &  K(I + \delta \nabla W) - \delta U_1 - \delta U_2  \\
 & = 	K(I) + \delta \partial K(I) \cdot \nabla W - \delta U_1 - \delta U_2 + \delta^2 K(\nabla W) \\
 & = K(I) - \delta U_1 +  \delta(\rho_0 h \cdot \nabla W - U_2) + \delta^2 K(\nabla W) + \delta (\partial K(I) - \rho_0 h) \cdot \nabla W. 
\end{aligned}
 \] 
 The main observation is that the second term in the last equality vanishes, the third term is $O(\delta^2)$ in $C^2$ norm, and the last term is $O(\delta^{\frac32})$ in $C^2_I$ norm. 

 The other parts of the proposition is identical to Theorem~\ref{thm:nhic-sr} and Theorem~\ref{thm:var-local-sr}, and are omitted. 
\end{proof}

Proposition~\ref{prop:DR-very-high-energy} follows immdediately from the above proposition.

\subsection{Proof of Proposition~\ref{prop:DR-non-deg-Crit}}
	\label{sec:non-deg-crit}

We prove Proposition~\ref{prop:DR-non-deg-Crit}, namely non-degenerate conditions $[DR1^c]-[DR4^c]$. This consists of two steps consisting of two localized perturbations of the potential $U$.

First, we perturb $U$ near the origin to achieve properties $[DR1^c] - [DR2^c]$. 
Let $W'$ be a $\rho$-neighborhood of the origin in $\R^2$ for small
enough $\rho>0$ so that it does not intersect sections $\Sigma^s_\pm$
and $\Sigma^u_\pm$. Consider $\xi(\varphi)$ a $C^\infty$-bump function so that
$\xi(\varphi)\equiv 1$ for $|\varphi|<\rho/2$ and $\xi(\varphi)\equiv 0$ for $|\varphi|>\rho$.
Let $Q_\zeta(\varphi)=\sum \zeta_{ij} \varphi_i \varphi_j$ be a symmetric quadratic
form. Consider $U_\zeta(\varphi)=U(\varphi)+\xi(\varphi)(Q_\zeta(\th)+\zeta_0)$.
In $W'\times \R^2$ we can diagonalize both: the quadratic form
$K(p)=\langle Ap,p\rangle$ and the Hessian $\partial^2 U(0)$. Explicit
calculation shows that choosing properly $\zeta$ one can make the minimum
of $U$ at $0$ being unique and eigenvalues to be distinct.

Suppose now that conditions $[DR1^c] - [DR2^c]$ hold. 
Fix a point $\th^* \in \gm^+$ at a distance of order of one
from the origin. In particular, it is away from sections $\Sigma^s_\pm$.
Let $W''$ be its small neighborhood so that intersects only one
homoclinic $\gm^+$. Denote $w^u=W^u\cap \Sigma^u_+$ an unstable curve
on the exit section $\Sigma^u_+$ and $w^s=W^s\cap \Sigma^s_+$ a stable
curve on the enter section $\Sigma^s_+$. Denote on $w^u$ (resp. $w^s$)
the point of intersection $\Sigma^u_+$ (resp. $\Sigma^s_+$) with strong
stable direction $q^{ss}$ (resp. $q^{uu}$). Recall that $q^+$ (resp. $p^+$)
denotes the point of intersection of $\gm^+$ with $\Sigma^u_+$
(resp. $\Sigma^s_+$). We also denote by $T^{uu}(q^+)$ (resp. $T^{ss}(p^+)$)
subspaces the tangent to $w^u$ (resp. $w^s$) at the corresponding points.
The critical energy surface $\{H=K-U=0\}$ is denoted by $S_0$.
In order to satisfy conditions $[DR3^c]-[DR4^c]$ the global map $\Phi_{glob}^+$
has to satisfy
\begin{itemize}
	\item $\Phi_{glob}^+ w^u \cap w^s \neq q^{ss}$ and
	$(\Phi_{glob}^+)^{-1}(w^s) \cap w^u \neq q^{uu}$.

	\item
	$D\Phg^+(q^+)|_{TS_0} T^{uu}(q^+) \pitchfork T^{ss}(p^+),
	\quad  D\Phg^-(q^-)|_{TS_0} T^{uu}(q^-) \pitchfork T^{ss}(p^-).$
\end{itemize}
The first condition can also be viewed as a property of the restriction
of $\Phg^+|_{S_0}$.  Notice that $\Phg^+$, restricted to $S_0$, is
a $2$-dimensional map.

Consider perturbations $\dt U \in \mathcal G^r$ of the potential $U$
localized in $W''$. By Lemma \ref{trans-map} the differential map
$D_U \Phi$ generates a subspace orthogonal to the gradient $\nabla H^U$
and the Hamiltonian vector field $\chi^U(\cdot)$ of $H^U$. Notice that
when we restrict $\Phg^+$ onto $\Sigma^u_+\cap S_0$ we factor
out $\nabla H^U$ and $\chi^U(\cdot)$. Both conditions on $\Phg^+$
(resp. $D\Phg^+|_{S_0}$) are non-equality conditions on images and
preimages for a $2$-dimensional map. Thus, these conditions can be
satisfies by Lemma \ref{trans-map}.


\section{Derivation of the slow mechanical system}
\label{slow-fast-section}

In this section we derive the slow mechanical system. The discussions here applies to arbitrary degrees of freedom, and indeed we will consider 
$$
H_\epsilon(\theta, p, t) = H_0(p) + \epsilon H_1(\theta, p, t), \quad \theta \in \T^n, p \in \R^n, t \in \T, 
$$
and let $p_0 \in \R^n$ be an $n-$resonance, namely, there are linearly independent $k_1, \cdots, k_n \in \Z^{n+1}$ such that
\[
k_i \cdot (\partial H_0(p_0), 1), \quad i = 1, \cdots, n. 
\]
Since the resonance depends only on the hyperplane containing $k_1, \cdots, k_n$, we may choose $k_1, \cdots, k_n$ to be \emph{irreducible}, namely
\[
\Span_\Z \{k_1, \cdots, k_n\} = \Span_\R \{ k_1, \cdots, k_n\} \cap \Z^{n+1}. 
\]
Results in this section applies to the rest of the paper by restricting $n=2$.

First, we will reduce the system near an  $n-$resonance to a normal form. 
After that, we perform a coordinate change on the extended phase space, and 
an energy reduction to reveal the slow system. 

In section~\ref{sec:norm-form-double}, we describe a resonant normal form.

In section~\ref{sec:aff-resc}, we describe the affine coordinate change and 
the rescaling, revealing the slow system.

In section~\ref{sec:rescaled-barriers}  we discuss variational properties of
these coordinate changes.

\subsection{Normal forms near double resonances}
\label{sec:norm-form-double}

Write  $\omega_0 := \partial_p H_0(p_0)$, then the orbit $(\omega_0,1)\,t$
is periodic.  Let
$$ T = \inf_{t>0}\{t(\omega_0,1)\in \Z^3\} $$
be the minimal period, then there exists a constant $T_* = T_*(k_1, \cdots, k_n) >0$, such that 
$T \le T_*(k_1, \cdots, k_n)$.

Given a function $f: \T^n\times \R^n \times \T \to \R$, we define
$$
[f]_{\om_0}(\theta, p, t) = \frac{1}{T} \int_0^T f(\theta + \omega_0 s, p, t + s)ds. $$
$[f]_{\om_0}$ corresponds to the resonant component related 
to the double resonance.

Writing $H_1(\theta, p, t) = \sum_{k \in \Z^n} h_k(p)
e^{2\pi i k \cdot (\theta,t)}$,
and let $\Lambda = \Span_\Z\{k_1, \cdots, k_n\}\subset \Z^n$, then 
$$   
[H_1]_{\om_0}(\theta, p, t) =\sum_{k \in \Lambda}  h_k(p) e^{2\pi i k \cdot (\theta,t)}.
$$

We define a rescaled differential in the action variable by
\be \label{scaled-norm} 
\partial_I f(\theta, p,t) =
\sqrt{\epsilon} \partial_p f(\theta, p,t), 
\ee
and use the notation $C^r_I$ to denote the $C^r$ norm with
$\partial_p$ replaced by $\partial_I$. For a vector or a matrix valued function, 
we take the $C^r$ (or $C_I^r$) norm to be the sum of the norm of all elements. Let  $B_r^n(p)$ denote the $r$-neighborhood of $p_0$ in $\R^n$.

\begin{thm}\label{double-norm-form}
Assume that $r\ge 5$.  
Then for any $M_1 >0$,  there exists $\epsilon_0 = \epsilon_0>0$, 
$C_1>1$ depending only on $\|H_0\|_{C^r}, T_*(k_1, \cdots, k_n), n , M_1$ 
such that for any $0 < \epsilon < \epsilon_0$, 
there exists a $C^\infty$ symplectic coordinate change
$$
\Phi_\eps:  \T^n \times B^n_{M_1\sqrt{\epsilon}}(p_0) \times
\T \to \T^n \times B^n_{2 M_1\sqrt{\epsilon}}(p_0)\times \T,$$
which is the identity in the $t$ component, 
such that
$$
N_\epsilon(\theta,p,t) := H_\epsilon \circ \Phi_\epsilon(\theta,p,t) =
H_0(p) + \epsilon Z(\theta,p, t) + 
\epsilon Z_1(\theta,p, t,\eps)\footnote{formally speaking this term 
contains two terms of the averaging expansion} + \epsilon R(\theta,p,t, \epsilon), $$
where $Z = [H_1]_{\omega_0}$, $[Z_1]_{\omega_0} = Z_1$, and 
\begin{equation}
\label{eq:ZR-bound}
\|Z_1\|_{C^2_I} \le C_1\sqrt{\epsilon}, \quad
\|R\|_{C^2_I} \le C_1 \epsilon^{\frac32},   
\end{equation}
and 
\[
\|\Phi_\eps - Id\|_{C^2_I} \le C_1\epsilon. 
\]
Moreover, $\Phi_\epsilon$ admits the following extensions:
\begin{enumerate}
	\item $\Phi_\epsilon$ can be extended to $\T^n \times \R^n \times \T$ such that $\Phi_\epsilon$ is identity outside of $\T^n \times U_{4M_1\sqrt{\epsilon}}(p_0) \times \T$. 
	\item The extension in item 1 can be further extended to $\widetilde{\Phi}_\epsilon(\theta, p, t, E)$ on $\T^n \times \R^n \times \T \times \R$, such that 
	\[
	\widetilde{\Phi}_\epsilon(\theta, p, t, E) = (\Phi_\epsilon(\theta, p, t), E + \widetilde{E}(\theta, p, t)),
	\]   
	and $\widetilde{\Phi}_\epsilon$ is exact symplectic. 
	\item $\|\widetilde{\Phi}_\epsilon - Id\|_{C^2_I} \le C_1\epsilon$ holds. 
\end{enumerate}
\end{thm}
\begin{rmk}
\begin{itemize}
	\item Our normal form is related to the classical ``partial averaging'', see, for example expansion (6.5) in \cite[Section 6.1.2]{AKN}.
	Our goal here is to obtain precise control of
	the norms with minimal regularity assumptions. In particular, the norm estimate
	of $\Phi_\eps-Id$ is stronger than the usual results, and is needed in the proof
	of Proposition~\ref{var-rel-norm}.
	\item  
	It is possible to lower the regularity assumption to $r \ge 4$, and use the weaker estimate $\|R\|_{C_I^2} \le C_1 \epsilon$. This, however, requires more technical discussions in the next few sections, and we choose to avoid it. 	
	\item Due to existence of the extension, we consider $N_\epsilon = H_\epsilon \circ \Phi_\epsilon$ as defined on all of $\T^n \times \R^n \times \T$, however, the normal form only holds on the local neighborhood. 
\end{itemize}
\end{rmk}

The rest of this section is dedicated to proving Theorem~\ref{double-norm-form}.
Denote $\Pi_\th(\th,p,t)=\th,\ \Pi_p(\th,p,t)=p$ the natural projections. 
For a map $\Phi: \T^n \times U \times \T \to \T^n \times \R^n \times \T$, which is the 
identity in the last component, denote $\Phi=(\Phi_\theta, \Phi_p, Id)$,
where $\Phi_\theta=\Pi_\th \circ \Phi$ and $\Phi_p=\Pi_p \circ \Phi$.

\begin{lem}\label{rescaled-norm} We have the following properties of
the rescaled norm.
\begin{enumerate}
	\item $\|f\|_{C_I^r} \le \|f\|_{C^r}$, $\|f\|_{C^r} \le\epsilon^{-r/2} \|f\|_{C_I^r}$.
	\item $\|\partial_\theta f\|_{C_I^{r-1}} \le \|f\|_{C_I^r}$, 
	$\|\partial_p f\|_{C_I^{r-1}} \le
	\frac{1}{\sqrt{\epsilon}}\|f\|_{C_I^r}$.
	\item There exists $c_{r,n}>1$ such that
	$\|fg\|_{C_I^r} \le c_{r,n}\,\|f\|_{C_I^r} \|g\|_{C_I^r}$.
	\item Let $\Phi$ be as before, and $Id$ denote the identity map. 
	There exists $c_{r,n} >1$ such that if 
	\[
	\max\{\|\Pi_\theta \, (\Phi - \text{Id})\|_{C_I^r},
	\|\Pi_p(\Phi- Id)\|_{C_I^r}/\sqrt{\epsilon} \} < 1
	\]
	we have 
	\[
	\|f \circ \Phi\|_{C^r_I}\le c_{r,n}  \|f\|_{C^r_I}. 
	\]
\end{enumerate}
\end{lem}  
\begin{proof}
The first two conclusions follow directly from definition.
For the third conclusion, we have
$\|\widetilde{f}\|_{C^r} = \|f\|_{C_I^r}$, where
$$ \widetilde{f}(\theta, I) = f(\theta, \sqrt{\epsilon} I).$$
The conclusion then follows from properties of the $C^r-$norm. 

For the fourth conclusion, we note
$$
f \circ \Phi = \widetilde{f} \circ \widetilde{\Phi},
$$
where $\widetilde{f}$ is as before, and $\widetilde{\Phi}(\theta, I) =
(\Phi_\theta(\theta, \sqrt{\epsilon}I),
\Phi_p(\theta, \sqrt{\epsilon}I)/\sqrt{\epsilon})$.
Moreover, let us denote $\Psi = \Phi - Id$, and note that 
$\widetilde{\Psi} = \widetilde{\Phi} - Id$. Then there exists $c>0$ 
depending only on dimension that 
\begin{multline*}
\| D\widetilde{\Phi}\|_{C^{r-1}} \le c + 
\|D(\widetilde{\Phi}-Id)\|_{C^{r-1}} \le c  + 
\max\{ \|\Pi_\theta \widetilde{\Psi}\|_{C^r} , \|\Pi_p \widetilde{\Psi}\|_{C^r} \} \\
\le c + \max\{\|\Pi_\theta \, (\Phi - Id)\|_{C_I^r},
\|\Pi_p(\Phi- Id)\|_{C_I^r}/\sqrt{\epsilon} \} \le c+1,
\end{multline*}
By the Faa-di Bruno formula,  there exists $d_r>0$ such that 
\begin{multline*}
\|f \circ \Phi\|_{C^r_I} 
= \| \widetilde{f} \circ \widetilde{\Phi}\|_{C^r} 
\le d_r \|\widetilde{f}\|_{C^r} \left(  1 + \|D \widetilde{\Phi}\|^r_{C^{r-1}} \right) \\
\le d_r \|f\|_{C_I^r} (1 + (c+1)^r) \le c_{r,n} \|f\|_{C_I^r}
\end{multline*}
where $c_{r,n} = d_r(1 + (c+1)^r)$. 
\end{proof}


For $\rho>0$, denote
$$
D_\rho = \T^n \times U_{\rho}(p_0) \times \T \times \R.
$$

We require another lemma for estimating the norm of a Hamiltonian 
coordinate change. This is an adaptation of Lemma 3.15 from \cite{DH}. 
\begin{lem}\label{bound-Phi-G}

For $0 < \rho' < \rho< \rho_0$, $r\ge 4$, $2 \le \ell \le r-1$, $\eps>0$ 
small enough, and  a $C^r$ Hamiltonian  $G: D_{\rho_0} \to \R$, 
let $\Phi_t^G$ denote the Hamiltonian flow defined by $G$.

Let $c_{r,n}>1$ be the constant from item 4, Lemma~\ref{rescaled-norm}. Assume that 
\[
\frac{1}{\sqrt{\epsilon}}\|G\|_{C_I^{l+1}}  < \frac{1}{c_{r,n}\sqrt{\epsilon}} 
\min\left\{ \sqrt{\epsilon}, \rho- \rho' \right\}, 
\]
then for $0 \le t \le 1$, the flow $\Phi_t^G$ is well defined 
from $D_{\rho'}$ to $D_\rho$. Moreover,
\begin{equation}\label{eq:PhiG}
\| \Pi_\theta (\Phi_t^G - Id) \|_{C_I^l} \le c_{r,n} \| \partial_p G\|_{C_I^l}, 
\quad \| \Pi_p( \Phi_t^G - Id) \|_{C_I^l} \le c_{r,n} \| \partial_\theta G\|_{C_I^l}.
\end{equation}
\end{lem}
\begin{proof}
Define 
\[
A_t = \max_{0 \le \tau \le t} 
\left\{ \| \Pi_\theta (\Phi_\tau^G - Id) \|_{C_I^{l}},  
\frac{1}{\sqrt{\epsilon}} \| \Pi_p( \Phi_\tau^G - Id) \|_{C_I^{l}}  \right\}. 
\]
Let $t_0$ be the largest $t \in [0,1]$ such that the following conditions hold.
\begin{itemize}
	\item[(a)] $A_t \le 1$. 
	\item[(b)] $\Phi_{-t}^G(D_\rho) \supset D_{\rho'}$, in other words, 
	$\Phi_t^G: D_{\rho'} \to D_\rho$ is well defined. 
\end{itemize}
We first show \eqref{eq:PhiG} holds on $0 \le t \le t_0$, then we show $t_0 =1$. 

By Lemma~\ref{rescaled-norm}, item 4, we have 
\[
\left\|\Pi_\theta \left(  \Phi_t^G - Id \right) \right\|_{C_I^l} \le 
\int_0^t \left\| \partial_p G \circ \Phi_\tau^G   \right\|_{C_I^l} d\tau \le c_{r,n} \|\partial_p G\|_{C_I^l}
\]
and 
\[
\left \|\Pi_p \left(  \Phi_t^G - Id \right) \right\|_{C_I^l} \le 
\int_0^t \left\| \partial_\theta G \circ \Phi_\tau^G   \right\|_{C_I^l} d\tau\\
\le c_{r,n}\|\partial_\theta G\|_{C_I^l}. 
\]


		We now show $t_0 =1$. By the estimates obtained, 
		\[
		A_{t_0} \le c_{r,n} \max\{\|\partial_p G\|_{C_I^l},   \frac{1}{\sqrt{\epsilon}} \|\partial_\theta G\|_{C_I^l}\} \le \frac{c_{r,n}}{\sqrt{\epsilon}} \|G\|_{C_I^{l+1}} < 1. 
		\] 
		Moreover, we have 
		\[
		\|\Pi_p \left(  \Phi_{t_0}^G - Id \right) \|_{C_I^l} \le c_{r,n}\|\partial_\theta G\|_{C_I^l} < \sqrt{\epsilon}\frac{c_{r,n}}{\sqrt{\epsilon}}\|G\|_{C_I^{l+1}} < \rho - \rho',
		\] 
		implying $\Phi_{-t_0}D_\rho \supsetneq D_{\rho'}$.  If $t_0 \ne 1$, we can extend conditions (a) and (b) to $t>t_0$, contradicting the maximality 
		of $t_0$. 
		\end{proof}

		We now state our main technical lemma, which  is an adaptation of an inductive
	lemma due to Bounemoura.

	\begin{lem}\label{average-ind}
	There exists a constant $K >2$ depending only on $r$ such that the following hold.   Assume the parameters $r\ge 4$, $\rho>0$,  $\mu>0$ satisfy
	\[
	0<\epsilon\le \mu^2, \quad  KT\mu<
	\min\left\{\frac{1}{4},\ \frac{\rho}{2(r-2)\sqrt{\epsilon}}\right\}.
	\]
	Assume that
	$$ H: \T^n \times U_\rho(p_0) \times \T\times \R \to \R, \quad
	H(\theta, p, t, E) = l + g_0 + f_0,$$
	where $l(p,E)=(\omega_0, 1) \cdot (p,E)$ is linear, $g_0, f_0$ are $C^r$
	and depend only on $(\theta, p, t)$, and 
	\begin{equation}\label{eq:initial-fg}
	\|g_0\|_{C_I^r} \le \sqrt{\epsilon}\mu, \quad  \|f_0\|_{C_I^r(\rho_0)} \le \epsilon , \quad \|\partial_p f_0\|_{C_I^{r-1}} \le \epsilon.
	\end{equation}
Then for 
$j \in \{1, \cdots, r-2\}$ and 
$\rho_j = \rho - j\frac{\rho}{2(r-2)} > \rho/2$,  there exists 
a collection of $C^{r-j}$--symplectic maps
$\Phi_j: D_{\rho_j} \to D_\rho$, of the special form
$$ \Phi_j(\theta, p, t, E) = (\Theta(\theta, p, t), P(\theta, p, t), t, E + \widetilde{E}(\theta, p,t)). $$
The maps $\Phi_j$ have the properties
\begin{equation}\label{eq:phi-id}
\| \Pi_\theta(\Phi_j - Id) \|_{C^{r-j}_I(D_{\rho_j})} 
\le jK^{j+1} (T\mu)^2, \quad 
\|\Pi_p(\Phi_j - Id)\|_{C_I^{r-j}(D_{\rho_j})} 
\le j K^j (T\mu)\sqrt{\epsilon}, 
\end{equation}
and
$$ H \circ \Phi_j = l + g_j + f_j, $$
for each $j \in \{1, \cdots, r-2\}$ satisfying
$g_{j} = g_{j-1} + [f_{j-1}]_{\om_0}$ and
\begin{equation}
\label{eq:bound}
\|g_j\|_{C^{r-j}(D_{\rho_j})} \le   
(2 - 2^{-j}) \sqrt{\epsilon}\mu, \quad \|f_j\|_{C^{r-j}(D_{\rho_j})} \le 
(KT\mu)^j\,\epsilon. 
\end{equation}
\end{lem}
\begin{proof}

The proof is an adaptation of the proof of Proposition 3.2, \cite{Bo}, page 9.

Following \cite{Bo}, we define
\be \label{eq:change} 
\chi_j = \frac1T \int_0^T s(f_j-[f_j]_{\om_0})(\theta + \omega_0s, p, t+s)ds,
\ee
$\Phi_0 = Id$, and 
$$ \Phi_{j+1} = \Phi_j \circ \Phi_1^{\chi_j}, $$
where $\Phi_s^{\chi_j}$ is the time-$s$ map of the Hamiltonian flow of $\chi_j$.

Using the fact that $\chi_j$ is independent of $E$, we have  the map $\Phi^{\chi_j}:=\Phi_1^{\chi_j}$
is independent of $E$ in the $(\theta, p,t)$ components. Furthermore,
$\Phi^{\chi_j}$ is the identity in the $t$ component, and $\Pi_E\Phi^{\chi_j}-E$
is independent of $E$. Hence $\Phi_j$ takes the special format described in
the lemma. The special form of $\Phi_j$ implies that $g_j$ and $f_j$ are
also independent of $E$, allowing the induction to continue.

First of all, assuming the step $j$ of the induction is complete,  we prove the norm estimate \eqref{eq:phi-id}. Using \eqref{eq:bound}, 
\[
\|[f_j]_{\om_0}\|_{C_I^{r-j}} \le \|f_j\|_{C_I^{r-j}} \le 
(KT\mu)^j\epsilon, \quad 
\|\chi_j\|_{C_I^{r-j}} \le 2T
\|f_j\|_{C_I^{r-j}} <  2T(KT\mu)^{j}\epsilon.
\]
We will choose $K$ such that $K \ge 2c_{r,n}$. Then
\begin{equation}\label{eq:chi-j}
\frac{1}{\sqrt{\epsilon}} \|\chi_j\|_{C_I^{r-j}} \le \frac{1}{c_{r,n}} (2Tc_{r,n}\sqrt{\epsilon})
(KT\mu)^{j} \le  \frac{1}{c_{r,n}}(KT\mu)^{j+1} 
< \frac{1}{c_{r,n}} \max\{ \frac{1}{4}, \frac{\rho}{2\sqrt{\epsilon}(r-2)}\}, 
\end{equation}
therefore Lemma~\ref{bound-Phi-G} applies with 
$\rho = \rho_{j}$, $\rho' = \rho_{j+1}$, $G = \chi_j$. 

For $j\ge 1$, using the inductive assumption and $c_{r,n}>1$,
\[
\|\partial_p\chi_j\|_{C_I^{r-j-1}} \le 
\frac{1}{\sqrt{\epsilon}}\|\chi_j\|_{C_I^{r-j}} \le  (KT\mu)^{j+1} ,	
\]
while for $j =0$, the initial assumption on $f_0$ implies
\[
\|\partial_p \chi_0\|_{C_I^{r-1}} \le 2T \|\partial_pf_{0}\|_{C_I^{r-1}} 
\le 2T\epsilon \le 2T \mu^2 < \frac{1}{c_{r,n}} (KT\mu)^2 ,
\]
since $T \ge 1$ and $K \ge 2c_{r,n}$. Apply Lemma~\ref{bound-Phi-G}, we get $\Phi^{\chi_j}$ are well defined
maps from $D_{\rho_{j+1}}$ to $D_{\rho_j}$, and 
\[
\|\Pi_\theta(\Phi^{\chi_j}-Id)\|_{C_I^{r-j-1}} \le c_{r,n} \| \partial_p\chi_j\|_{C_I^{r-j-1}} \le (KT\mu)^2,
\]
\[
\|\Pi_p(\Phi^{\chi_j}-Id)\|_{C_I^{r-j-1}} \le 
c_{r,n}\|\partial_\theta\chi_j\|_{C_I^{r-j-1}} \le c_{r,n} \|\chi_j\|_{C_I^{r-j}} < (KT\mu)^{j+1} \sqrt{\epsilon}
\]
using \eqref{eq:chi-j}. Since $\Phi_1 = \Phi_1^{\chi_0}$, we obtain \eqref{eq:phi-id} for $j=1$. 

For each $j \ge 2$, using the bound on $\Phi^{\chi_j}- Id$ above, Lemma~\ref{rescaled-norm}, item 4 applies for $\Phi = \Phi^{\chi_j}$. Then
\[
\| \Pi_\theta (\Phi_{j-1} - Id) \circ \Phi^{\chi_{j-1}} \|_{C_i^{r-j}} \le c_{r,n} \|\Pi_\theta (\Phi_{j-1} - Id)\| .
\]
We have 
\begin{align*}
&\|\Pi_\theta (\Phi_j - Id)\|_{C_I^{r-j}} \\
& =  \| \Pi_\theta (\Phi_{j-1} \circ \Phi^{\chi_{j-1}}- 
\Phi^{\chi_{j-1}} +\Phi^{\chi_{j-1}}  - Id)\|_{C_I^{r-j}} \\
& \le \| \Pi_\theta (\Phi_{j-1} - Id) \circ \Phi^{\chi_{j-1}} \|_{C_I^{r-j}} + \| \Pi_\theta (\Phi^{\chi_{j-1}} - Id) \|_{C_I^{r-j}} \\
& \le  c_{r,n} \|\Pi_\theta(\Phi_{j-1} - Id)\|_{C_I^{r-j}} + 
\| \Pi_\theta (\Phi^{\chi_{j-1}} - Id) \|_{C_I^{r-j}} \\
& \le  \sum_{i=0}^{j-1} c_{r,n}^{j-i-1}\| \Pi_\theta (\Phi^{\chi_{i}} - Id) \|_{C_I^{r-i-1}} 
< j K^{j-1} (KT\mu)^2,
\end{align*}
noting $c_{r,n} < K$ and $\|\Pi(\Phi^{\chi_j} - Id)\|_{C_r^{r-j-1}} \le (KT\mu)^2$. 

By the same reasoning, we get 
\[
\|\Pi_p (\Phi_j - Id)\|_{C_I^{r-j}} \le  \sum_{i=0}^{j-1} c_{r,n}^{j-i-1}
\| \Pi_p (\Phi^{\chi_{i}} - Id) \|_{C_I^{r-i-1}} \le  jK^{j-1} (KT\mu)\sqrt{\epsilon}. 
\]

We now proceed with the induction in $j$. The inductive assumption \eqref{eq:bound} holds for $j=0$ due to \eqref{eq:initial-fg}, which is in fact stronger. For the inductive step, define
$$
g_{j+1} = g_j + [f_j]_{\om_0}.
$$ 
Since $\|f_j\|_{C_I^{r-j}} \le  (KT\mu)^j \epsilon <  4^{-j} \sqrt{\epsilon}\mu$, we get
\[
\|g_{j+1}\|_{C_I^{r-j-1}} \le \|g_j\|_{C_I^{r-j}} + \|f_j\|_{C_I^{r-j}}  \le ( 2 - 2^{-j} + 4^{-j}) \sqrt{\epsilon}\mu < (2-2^{j+1}) \sqrt{\epsilon}\mu. 
\]

By a standard computation, 
$$
f_{j+1} = \int_0^1 \{f_j^s, \chi_j\} \circ \Phi_s^{\chi_j}ds,
$$
where $f_j^s = g_j + sf_j + (1-s)[f_j]_{\om_0}$. 
Estimate \eqref{eq:bound} implies $\|f_j^s\|_{C_I^{r-j}} \le 3 \|g_j\|_{C_I^{r-j}} \le 6  \sqrt{\epsilon} \mu  $. Using Lemma~\ref{rescaled-norm}, items 2 and 3,  there exists an absolute constant $d>1$ such that for any $f, g \in C^l(D_\rho)$, we have
\[
\|\{f,g\}\|_{C_I^l} = \|\partial_\theta f \cdot \partial_p g + \partial_p f \cdot \partial_\theta g\|_{C_I^l} \le \frac{d}{\sqrt{\epsilon}} \|f\|_{C_I^{l+1}} \|g\|_{C_I^{l+1}}. 
\]
Therefore 
\begin{multline*} \|\{f_j^s, \chi_j\}\|_{C_I^{r-j-1}}  \le \frac{d}{\sqrt{\epsilon}} \|f^s_j\|_{C_I^{r-j}} \|\chi_j\|_{C_I^{r-j}} \\
\le  \frac{d}{\sqrt{\epsilon}} \cdot   3 \sqrt{\epsilon} \mu \cdot 2T \|f_j\|_{C_I^{r-j}}  \le   (6dT\mu) (KT\mu)^{j} \epsilon .
\end{multline*}
Furthermore, by Lemma~\ref{rescaled-norm}, items 4,
\begin{multline*}
\|f_{j+1}\|_{C_I^{r-j-1}} \le  \max_{0 \le s \le 1}\| \{f_j^s, \chi_j\} \circ \Phi_s^{\chi_j} \|_{C_I^{r-j-1}} \le c_{r,n}  \max_{0 \le s \le 1}\|\{f_j^s, \chi_j\}\|_{C_I^{r-j-1}}  \\
\le (6c_{r,n} d) T\mu(KT\mu)^j \epsilon \le (KT\mu)^{j+1} \epsilon
\end{multline*}
if we choose $K \ge 6c_{r,n}d$. 
The induction is complete. 
\end{proof}

\begin{proof}
[Proof of Theorem~\ref{double-norm-form}]
First, we show that 
We write
$$ H(\theta, p, t, E) = H_\epsilon(\theta, p, t) - H_0(p_0) + E = l + g_0 + f_0,$$
where $l(p, E) = (\omega_0, 1)\cdot (p, E)$,
	$$
	g_0(\theta, p, t) = H_0(p) - H_0(p_0) - \omega_0 \cdot p
	$$
	and $f_0 = \epsilon H_1$.

	Define $\rho = 2 M_1 \sqrt{\epsilon}$.  we have the following estimates: $\partial_\theta g_0 =0$, 
	\[
	\|\partial_{p} g_0\|_{C^0} = \|\partial_p H_0(p) - \partial_p H_0(p_0) \|_{C^0}  \le \|H_0\|_{C^2} \cdot \rho = 2\|H_0\|_{C^2}  M_1 \sqrt{\epsilon}
	\]
	and $\|\partial^j_{p^j} g_0\|_{C^0} \le \|H_0\|_{C^{1+j}}$ for all $j \ge 2$. Then for some $\widetilde{C}_1 >1$ 
	depending on $H_0,H_1$,
	\[
	\| \partial_p g_0\|_{C_I^{r-1}(D_\rho)} \le  \max_{j \ge 1}(\sqrt{\epsilon})^{j-1} \|\partial^{j}_{p^{j}} g_0\|_{C^0} \le  {\widetilde C }_1   M_1\sqrt{\epsilon}. 
	\]
	On the other hand, using $\|f\|_{C^r_I} \le \|f\|_{C^r}$, we have
	\[
	\|f_0\|_{C_I^r(D_\rho)} \le \widetilde{C}_1 \epsilon , \quad 
	\|\partial_pf_0\|_{C_I^{r-1}(D_\rho)} \le \widetilde{C}_2 \epsilon .
	\]

	Choose $\mu = \widetilde{C}_1 M_1\sqrt{\epsilon}\, \ge \, \sqrt{\epsilon}$, we have 
	\[
	\|g_0\|_{C_I^r} \le \max\{ \|\partial_\theta g_0\|_{C^{r-1}_I(D_\rho)}, \sqrt{\epsilon} \| \partial_p g_0\|_{C_I^{r-1}(D_\rho)}  \} \le \widetilde{C}_1 M_1\epsilon = \sqrt{\epsilon}\mu,
	\]
	and $\|f_0\|_{C_I^r(D_\rho)}, \|\partial_pf_0\|_{C_I^{r-1}(D_\rho)} \le \sqrt{\epsilon}\mu$. 
	By choosing $\epsilon$ sufficiently small, we ensure
	\[
	KT\mu = KT\widetilde{C}_1 M_1\sqrt{\epsilon} < \min\{\frac14, \frac{\rho}{2(r-2)\sqrt{\epsilon}}\} = \min\{\frac14, \frac{2\bar{E}}{2(r-2)}\}. 
	\]
	The conditions of Lemma~\ref{average-ind} are satisfied with these parameters and we apply
	the lemma with $j = 3$. There exists a map $\Phi_3: D_{\rho/2} \to D_\rho$,
	$$
	\|\Phi_3 - Id\|_{C_I^2(D_{\rho/2})} \le \max\{ 2K^3(T\mu)^2, 2K^2 (T\mu)\sqrt{\epsilon}\} \le \widetilde{C}_2 \epsilon
	$$
	for some $\widetilde{C}_2>1$ depending on $T, \widetilde{C}_1, M_1$. Moreover, 
	\[
	H \circ \Phi_3 = l + g_0 + [f_0]_{\om_0} + 
	[f_1]_{\om_0}  + [f_2]_{\om_0} + f_3,
	\]
	with
	$$
	\|[f_1]_{\om_0}\|_{C_I^2(D_{\rho/2})} \le (KT\mu)\sqrt{\epsilon}\mu \le
	\widetilde{C}_3 \epsilon^{\frac32},
	$$
	$$
	\|[f_2]\|_{C^2_I(D_{\rho/2})} \le (KT\mu)^2 \sqrt{\epsilon}\mu \le 
	\widetilde{C}_3 \epsilon^2,
	$$
	\[
	\|[f_3]\|_{C^2_I(D_{\rho/2})} \le (KT\mu)^3 \sqrt{\epsilon}\mu \le
	\widetilde{C}_3 \epsilon^{\frac52}, 		
	\]
	for some $\widetilde{C}_3$ depending on $T, \widetilde{C}_1, M_1$. 
	Using $l+g_0 = H_0(p) + E - H_0(p_0)$,  
	define $\epsilon Z = [f_0]_{\om_0},\ 
	\epsilon Z_1 = [f_1 + f_2]_{\om_0},$ and 
	$\epsilon R = f_3$, we obtain
	$$
	(H_\epsilon + E) \circ \Phi_2 =
	H_0 + \epsilon Z + \epsilon Z_1 + \epsilon R + E
	$$
	with the desired estimates and constant $C_1 = \max\{\widetilde{C}_2, \widetilde{C}_3\}$. Finally, we define
	$\Phi_\epsilon(\theta, p,t) = \Phi_3(\theta, p, t, E)$.
	This is well defined since $\Phi_3$ is independent of $E$.

	We now use a smooth approximation technique to show the normal form $\Phi$ can be taken to be $C^\infty$. Using standard techniques, for every $\sigma>0$, there is $C = C(r) >0$ such that there is $C^\infty$ functions $H_0'(p)$ and $H_1'(\theta, p, t)$ satisfying:
	\[
	\|H_0 - H_0'\|_{C^2} \le \sigma, \quad \|H_1 - H_1'\|_{C^2} \le \sigma, \quad \|H_0\|_{C^r}, \|H_1\|_{C^r} \le C. 
	\]
	Note that the estimates of the normal form depends on the $C^r$ norm of the smooth approximation. We apply the above procedure to  $H_0' + \epsilon H_1'$, and obtain a $C^\infty$ coordinate change $\Phi'$ such that 
	\[
	(H_0' + \epsilon H_1' + E) \circ \Phi' = H_0 + \epsilon Z' + \epsilon Z_1' + \epsilon R'. 
	\]
	We have $\|\epsilon Z' - \epsilon Z\|_{C^2} \le \epsilon \|H_1 - H_1'\|_{C^2} < \epsilon \sigma$ and $\|Z_1'\|_{C^2_I}\le \tilde{C} \epsilon$, $\|R\|_{C^2_I}\le \tilde{C} \epsilon^{\frac32}$. We now compute
	\[
\begin{aligned}
 & 	(H_0 + \epsilon H_1 + E) \circ \Phi' \\
 & = H_0 + \epsilon Z + \epsilon Z_1' + + \epsilon R' +  (H_0 - H_0') \circ \Phi' + \epsilon(H_1 - H_1') \circ \Phi' + \epsilon (Z' - Z) . 
\end{aligned}
	\]
	The last three terms are bounded in $C^2_I$ norms by $\tilde{C} \sqrt{\epsilon}^{-2}\sigma$. To ensure the remainder is small we can take $\sigma = \epsilon^4$. 
	\end{proof}

	The normal form lemma stated here also applies to double resonances with long period, which then combined with the idea of Lochak, can be used to study single resonance. This is Proposition~\ref{prop:normal-form-dr-gen}, which we prove here. 
	\begin{proof}
	[Proof of Proposition~\ref{prop:normal-form-dr-gen}]
	As in the proof of Theorem~\ref{double-norm-form}, we apply the Lemma~\ref{average-ind} to $l + g_0 + f_0$ where $l + g_0 = H_0 + E$ and $f_0 = \epsilon H_1$. Choose $\rho = K_1 \sqrt{\epsilon}$, then 
	\[
	\|g_0\|_{C^r_I} \le C K_1 \sqrt{\epsilon}, \quad \|f_0\|_{C_I^r} \le \epsilon, \quad \|\partial_p f_0\|_{C_I^{r-1}} \le \epsilon
	\]
	where $C$ depends on $r, n, \|H_0\|_{C^r}$. Choose
	\[
	\mu = C K_1 \sqrt{\epsilon},
	\]
	 by assumption $T \le C_1/(K_1^2 \sqrt{\epsilon})$, then $KT\mu \le (K C C_1) K_1^{-1}$. The condition of Lemma~\ref{average-ind} is satisfied if $K_1$ is sufficiently large. We apply Lemma~\ref{average-ind} with $j =1$, the conclusion is  
	\[
	\left\| \Pi_\theta (\Phi - Id) \right\|_{C^{r-1}_I} \le K^2(KCC_1)^2 K_1^{-2}, \quad 
	\left\| \Pi_p (\Phi - Id) \right\|_{C^{r-1}_I} \le K(KCC_1) K_1^{-1} \sqrt{\epsilon}. 
	\]
	and 
	\[
	\|f_1\|_{C^{r-j}} \le (KCC_1) K_1^{-1} \epsilon. 
	\]
	We can also apply the same smooth approximation to upgrade the coordinate change to $C^\infty$. 
	\end{proof}

	\subsection{Affine coordinate change, rescaling and energy reduction}
	\label{sec:aff-resc}

	\paragraph{Definition of the slow system.}
	Recall that an $n$-tuple of vectors 
	$k_1, \cdots, k_n \in \Z^{n+1}$ defines an irreducible 
	lattice, if the lattice $\Lb=\langle k_1, \cdots, k_n \rangle_\Z$
	is not contained in any lattice of the same rank or, equivalently, 
	$\langle k_1, \cdots, k_n \rangle_\R \cap \Z^{n+1}=\Lb$.
	Let $B_0$ be the $n\times (n+1)$ matrix whose rows 
	are vectors $k_1^T, \cdots, k_n^T$ from $\Z^{n+1}$, and 
	let $k_{n+1}\in \Z^{n+1}$ be such that 
	\begin{equation}
	\label{eq:int-matrix}
	B := \begin{bmatrix}
	k_1^T \\ \vdots \\ k_{n+1}^T
	\end{bmatrix} \in SL(n+1, \Z).   
	\end{equation}
	Such a $k_{n+1}$ exists if and only if $k_1, \cdots, k_n$ is irreducible. 
	By possibly changing the signs of the rows of $B$, we ensure
	\begin{equation}
	\label{eq:beta}
	\beta := k_{n+1} \cdot (\omega_0, 1) >0
	\end{equation}
	In particular, we have 
	\begin{equation}
	\label{eq:w-cancel}
	B \bmat{\omega_0 \\ 1} = \bmat{ 0 \\ \beta}. 
	\end{equation}

	Let $A_0 = \partial^2_{pp}H_0(p_0)$, we define 
	\begin{equation}
	\label{eq:KI}
	K(I) = \frac12 \left(  B_0 A_0 B_0^T \right)I \cdot I, \quad I \in \R^n
	\end{equation}
	and $U(\varphi)$, $\varphi \in \T^n$ by the relation
	\begin{equation}
	\label{eq:Uphi}
	U(k_1 \cdot (\theta, t), \cdots, k_n \cdot (\theta, t)) = - Z(\theta, p_0, t). 
	\end{equation}
	Consider the corresponding autonomous system 
	\[
	G_\epsilon(\theta, p, t, E) = N_\epsilon + E,
	\]
	We show that $G_\epsilon$ can be reduced to 
	\be \label{eq:Hseps} 
	H_\epsilon^s(\varphi, I, \tau) = 
	\frac{1}{\beta}\left(  K(I) - U(\varphi) + \sqrt{\epsilon}P(\varphi, I, \tau,\eps) \right), 
	\quad \varphi \in \T^n, I \in \R^n, \tau \in \sqrt{\epsilon}\T, 
	\ee
	with a coordinate change and an energy reduction. See Proposition \ref{prop:Heps-Hseps}.

	\paragraph{Linear and rescaling coordinate change.}

	Our coordinate change is a combination of a linear and a rescaling coordinate change. Namely, we have 
	\[
	(\theta, p, t, E) = \Phi_1(\varphi, p^s, s, p^f): \quad 
	\bmat{\theta \\ t} = B^{-1} \bmat{\varphi \\ s}, \quad \bmat{p- p_0 \\ E + H_0(p_0)} = B^T \bmat{p^s \\ p^f},
	\]
	\[
	(\varphi, p^s, s, p^f) = \Phi_2(\varphi, I, \tau, F)
	= (\varphi, \sqrt{\epsilon}I, \tau/\sqrt{\epsilon}, \epsilon F). 
	\]
	We then have 
	\be \label{def:Phi-SL}
	\Phi_L = \Phi_1 \circ \Phi_2 : \quad (\varphi, I, \tau, F) \in 
	\T^n \times \R^n \times \sqrt{\epsilon}\T \times \R \mapsto (\theta, p, t, E) 
	\ee
	by the formula
	\begin{equation}
	\label{eq:linear-rescale}
	\bmat{\theta \\ t} = B^{-1} \bmat{\varphi \\ \tau/\sqrt{\epsilon}},
	\qquad
	\bmat{p - p_0 \\ E + H_0(p_0)} = \sqrt{\epsilon} B^T \bmat{  I \\ \sqrt{\epsilon} F}. 
	\end{equation}
	Given any $M_1>0$, there exists $C_B >1$ such that 
	\begin{equation}
	\label{eq:PhiLS-dom}
	\begin{aligned}
	& \{\ \|(I, \sqrt{\epsilon} F)\| < C_B^{-1} M_1\ \} 
	\quad \subset \quad \Phi_L\left( \ \{\|(p-p_0, E + H_0(p_0))\|  < \sqrt{\eps}\, M_1\}\ \right) \\
	& \quad \subset  \{\ \|(I, \sqrt{\epsilon} F)\| < C_B M_1\ \}  .
	\end{aligned}
	\end{equation}

	Define 
	\[
	G_\epsilon^s(\varphi, I, \tau, F) = \frac{1}{\sqrt{\epsilon}} G_\epsilon \circ \Phi_L,
	\]
	then the Hamiltonian flows of $G_\epsilon^s$ and $G_\epsilon$ are conjugate via $\Phi_L$. This can be seen from the equivalence of the following Hamiltonian systems:
	\begin{align*}
	& 	(G_\epsilon,\quad d\theta \wedge dp + dt \wedge dE) \quad \sim \quad (G_\epsilon \circ \Phi_L, \quad \sqrt{\epsilon}(d\varphi \wedge dI + d\tau \wedge dF)) \\
	& \sim \left(  \frac{1}{\sqrt{\epsilon}} G_\epsilon \circ \Phi_L,\quad  d\varphi \wedge dI + d\tau \wedge dF \right)
	\end{align*}

	We have the following lemma. 
	\begin{lem}\label{lem:G-eps-slow}
	Assume that $N_\epsilon$ satisfies \eqref{eq:ZR-bound} on $D_{M_1\sqrt{\epsilon}}$. Then there exists $C_B>1$ 
	and $C_1 = C_1(H_0, \widetilde{C}, M)$ such that  
	\begin{equation}
	\label{eq:G-eps-s}
	G_\epsilon^s = \sqrt{\epsilon}\left(  K(I) - U(\varphi) + (\beta +  \sqrt{\epsilon} \, l(I, F)) \cdot F + \sqrt{\epsilon} P_1(\varphi, I, \tau, F,
	\eps) \right),
	\end{equation}
	where $\beta = k_{n+1} \cdot (\omega(p_0),1)$, and 
	\[
	\|l\|_{C^2}, \|P_1\|_{C^2} \le C_1, 
	\]
	with the norm taken on the set 
	\begin{equation}
	\label{eq:CB}
	\{ \|(I,\sqrt{\epsilon}F)\| < C_B^{-1} M_1\}. 
	\end{equation}
	\end{lem}
	\begin{proof}[Proof of Lemma~\ref{lem:G-eps-slow}]
	Due to \eqref{eq:PhiLS-dom}, $G_\epsilon \circ \Phi_L$ is defined on the $\{ \|(I, \sqrt{\epsilon} F)\|< C_B^{-1}M_1\}$. 

	Consider the expansion of $N_\eps$ from Theorem 
	\ref{double-norm-form}
	$$
	N_\epsilon(\theta,p,t) := H_\epsilon \circ \Phi_\epsilon(\theta,p,t) =
	H_0(p) + \epsilon Z(\theta,p, t) + 
	\epsilon Z_1 + \epsilon R. $$
	Denote 
	$$A = \partial^2(H_0(p) + E)(p_0, - H_0(p_0)) =  \bmat{\partial^2_{pp}H_0(p_0) & 0 \\ 0 & 0},
	$$ 
	where $H_0(p)+E$ is a function of $p$ and $E$.  
	Expanding $H_0(p) + E$ to the third order at $(p_0, -H_0(p_0))$, 
	and $Z(\cdot, p)$ to the first order at $p_0$, we have 
	\[
	G_\epsilon = \bmat{\omega_0 \\ 1} \cdot \bmat{p - p_0 \\ E + H_0(p_0)} +  \frac12 A \bmat{p \\ E} \cdot \bmat{p \\ E} + \epsilon Z(\theta, p_0,t) + \check{H}_0 + \epsilon \check{Z}_1 + \epsilon R,
	\]
	with 
	\[
	\check{H}_0 = H_0(p) - H_0(p_0) - \omega(p_0)\cdot (p-p_0)- \frac12 A \bmat{p \\ E} \cdot \bmat{p \\ E}, \quad \check{Z}_1 = Z(\theta, p, t) - Z(\theta, p_0, t) + Z_1. 
	\]

	On the set $\|p - p_0\| \le M_1 \sqrt{\epsilon}$ and $|E + H_0(p)| \le M_1\sqrt{\epsilon}$, for some $C = C(H_0, \tdC, M_1)$, we have
	\[
	\|\check{H_0}\|_{C^2_I} \le C \epsilon^{\frac32}, \quad \|\check{Z}_1\|_{C^2_I} \le C \sqrt{\epsilon}, \quad \|R\|_{C^2_I} \le C \epsilon^{\frac32}. 
	\]

	From  \eqref{eq:linear-rescale} we have  $U(\varphi_1, \cdots, \varphi_n) = -Z(\theta, p_0,t)$, and from \eqref{eq:w-cancel}, 
	\[
	(\omega_0 , 1) \cdot (p - p_0, E + H_0(p_0)) = (0 , \beta) \cdot (\sqrt{\epsilon} I , \epsilon F) = \epsilon \beta F, 
	\]
	we get 
	\begin{align}
	&G_\epsilon \circ \Phi_L =\qquad \qquad  \nonumber \\
	& \label{eq:N-eps-PhiL}
	\epsilon \left(
	\frac12  BAB^T \bmat{I \\  \sqrt{\epsilon} F} \cdot \bmat{I \\ \sqrt{\epsilon} F} 
	+ \beta F -  U(\varphi) + (\check{H}_0/\epsilon +  \check{Z}_1 +  R)\circ \Phi_L
	\right).
	\end{align}
	Finally, note $K(I) = \frac12 BAB^T \bmat{I \\ 0} \cdot \bmat{I \\ 0}$, and define 
	\[
	l(I, F) \cdot F =  BAB^T \bmat{I \\ 0}\cdot \bmat{0 \\ \sqrt{\epsilon} F} + \frac12 BAB^T \bmat{0 \\ \sqrt{\epsilon} F} \cdot \bmat{0 \\ \sqrt{\epsilon} F},
	\]
	\[
	\sqrt{\epsilon} P = (\check H_0 +  \epsilon \check{Z}_1 + R) \circ \Phi_L,
	\]
	then \eqref{eq:G-eps-s} follows directly from \eqref{eq:N-eps-PhiL} and definition of $\Phi_L$. It is also clear from the definition that $l(I,F)/\sqrt{\epsilon}$ has bounded $C^2$ norm on the set $\{\|I\|, \|\sqrt{\epsilon}F\| \le C_B^{-1} M\}$. 

	For the norm estimates, we note the coordinate change $\Phi_1$ increase $C^2$ norm by a factor $C_B$ depending only on $B$. 
	Note that $\check{Z}_1$ depends only on $k_1 \cdot (\theta, t), \cdots, k_n\cdot (\theta, t), p$, $\check{Z}_1$ depends only on $\varphi, I, F$ (and not in $\tau$). As a result, using the definition of the rescaled norm (\ref{rescaled-norm}), we have 
	\begin{align*}
	&\|  (\check H_0/\epsilon + \check{Z}_1) \circ \Phi_L\|_{C^2}
	\le \|(\check H_0/\epsilon + \check{Z}_1) \circ \Phi_1 (\varphi, \sqrt{\epsilon} I, \epsilon F)\|_{C_I^2} \\
	&\le \|(\check H_0/\epsilon + \check{Z}_1) \circ \Phi_1\|_{C^2} 
	\le C_B C \sqrt{\epsilon}. 
	\end{align*}
	For terms depending on $\tau$, we have 
	\[
	\|R \circ \Phi_L\|_{C^2} = \|R \circ \Phi_1 (\varphi, \tau/\sqrt{\epsilon}, \sqrt{\epsilon}I, \epsilon F)\|_{C^2} \le C_B (\sqrt{\epsilon})^{-2} \|R\|_{C_I^2} \le C_B C \sqrt{\epsilon}, 
	\] 
	where we used \eqref{eq:ZR-bound}. 	The $C^2$-norm estimate of $P$ follows from 
	estimates on $C^2_I$ norms of $\check H_0,\ \check Z_1, R$. 
	\end{proof}

	\paragraph{Reduction on energy surface}

	We perform a standard reduction on the energy surface $G_\epsilon^s =0$, 
	with $\tau$ as the new time, obtaining a time-periodic system.

	\begin{lem}\label{lem:energy-reduction}
	Assume that conclusions of Lemma~\ref{lem:G-eps-slow} hold
	on the set $\{\|(I, \sqrt{\epsilon}F)\|< M_2\}$. Then there exists 
	$\epsilon_0 = \epsilon_0(H_0, C_1, M_2),\ 
	C_3 = C_2(H_0, C_1, M_2, U)>0$, such that for 
	any $0 < \epsilon < \epsilon_0$, there is a function 
	\[
	H_\epsilon^s: \T^n \times B^{n}_{M_2}(0) \times \sqrt{\epsilon}\T \to \R
	\]
	uniquely solving the equation 
	\[
	G_\epsilon^s(\varphi, I, \tau, -  H_\epsilon^s) =0
	\]
	on the set $\{\|(I, \sqrt{\epsilon}F)\| \le M_2\}$. Moreover, $H^s_\eps$ has the form 
	(\ref{eq:Hseps}), i.e. 
	\[
	H_\epsilon^s(\varphi, I, \tau) =  
	\frac{1}{\beta}\left(  K(I) - U(\varphi) + \sqrt{\epsilon}P(\varphi, I, \tau,\eps) \right) 
	\]
	where $\|P\|_{C^2} \le C_3$. 

	In particular, $H_\epsilon^s \to H^s/\beta$ uniformly in $C^2(\T^2 \times \R^2 \times \R)$.
	\end{lem}
	\begin{proof}[Proof of Lemma~\ref{lem:energy-reduction}]

	We can choose $\epsilon_0$ such that for any $0<\eps<\eps_0$
	\[
	\frac{\partial}{\partial F}(G_\epsilon^s) > \frac{\beta}2 >0
	\]
	on $\{\|(I, \sqrt{\epsilon}F)\| \le M_2\}$. Therefore, 
	$H_\epsilon^s$ exists by the implicit function theorem. 
	Moreover, there exists a constant $C'$, depending on $H_0$ 
	and $U$, but independent of $\epsilon$, 
	such that $\|H^s_\epsilon\|_{C^2} \le C'$. 

	Let 
	$Q = H_\epsilon^s - \frac{1}{\beta}(K(I)- U(\varphi)) = H_\epsilon^s - H^s$, then $\|Q\|_{C^2} \le \|H^s\|_{C^2} + C'$. 
	We know
	\begin{align*}
	0 &= G_\epsilon^s(\varphi, I, \tau, -Q - H^s) \\
	&= \sqrt{\epsilon} \left(  - \beta Q + \sqrt{\epsilon}l(I, Q + H^s)(Q + H^s) + \sqrt{\epsilon} P_1(\varphi, I, \tau, Q + H^s,\textr{\eps}) \right) \\
	&=: \sqrt{\epsilon}(-\beta Q + \sqrt{\epsilon}P_2(\varphi, I, \tau, Q + H^s,\eps)).
	\end{align*}
	Therefore,
	\[
	Q = \frac{\sqrt{\epsilon}}{\beta} P_2(\varphi, I, \tau, Q + H^s, \eps).
	\]
	To solve this implicit equation notice that
	there exists $C''>0$ depending only on $C_1$, $H^s$ such that  $\|P_2\|_{C^2} \le C''$. Application of the Faa-di Bruno formula
	show that for some $C>0$ depending only on $n$ we have 
	\[
	\|Q\|_{C^2} = \frac{\sqrt{\epsilon}}{\beta} 
	\|P_2(\varphi, I, \tau, Q + H^s,\eps)\|_{C^2} \le \sqrt{\epsilon} C C''  \|Q\|_{C^2}^2  \le \sqrt{\epsilon} C C''(\|H^s\|_{C^2} + C'),
	\]
	and the lemma follows. 
	\end{proof}

	We have 
	\begin{equation}
	\label{eq:G-eps-H-eps}
	G_\epsilon^s(\varphi, I, \tau, F) = 0,  \quad \Longleftrightarrow \quad   H_\epsilon^s +  F =0.
	\end{equation}

	The following Proposition follows from the standard energy reduction (see for example \cite{Ar4}). 
	\begin{prop}\label{prop:Heps-Hseps}
	Assume the conclusions of Theorem~\ref{double-norm-form} hold 
	on the set 
	$$D_{M_1}:=\{\|(p - p_0, H+ H_0(p_0)\| < M_1 \sqrt{\epsilon}\}
	$$
	for sufficently large $M_1$ depending only on $H_0$. \
	Then for $M_3 = C_B^{-2} M_1$, $0 < \epsilon < \epsilon_0(H_0, C_1, M_1)$, 
	the following are equivalent:
	\begin{enumerate}
		\item The curve $(\theta, t, p, E)(t)$ is an orbit of $N_\epsilon(\theta, p, t) + E$ inside $D_{M_1}$.
		\item The curve $(\varphi, \tau, I,  F)(t) = \Phi_L(\theta, t, p, E)(t)$ is 
		the time change of an orbit of $N_\epsilon^s(\varphi, I, \tau) + F$, with $\tau$ as the new time. 
	\end{enumerate}  
	\end{prop}

	\subsection{Variational properties of the coordinate changes}
	\label{sec:rescaled-barriers}

	We have made two reductions: The normal form
	\[
	H_\epsilon(\theta, p, t) \to N_\epsilon(\theta, p, t) = H_\epsilon \circ \Phi_\epsilon, 
	\]
	and the coordinate change with time change 
	\[
	N_\epsilon(\theta, p, t) + E \to H_\epsilon^s(\varphi, I, \tau) + F.  
	\]

	In this section, we discuss the effect of these reductions on the Lagrangian, barrier function, and the Mather, Aubry, Ma\~ne sets. The main conclusion of this section is the following proposition, which follows directly from Propositions~\ref{var-rel-norm} and \ref{prop:slow-variational} below. 

	\begin{prop}\label{prop:He-He-slow-var}
	For $M_1$ large enough depending only on $H_0, k_1, \cdots, k_n, n$, there exist $C_2>1$ and $\epsilon_0>0$ depending on $H_0, k_1, \cdots, k_n, n, M_1$,  such that the following hold for any 
	$0 < \epsilon < \epsilon_0$. 
	\begin{enumerate}
		\item For $M' = C_B^{-2} M_1/2$ (see \eqref{eq:CB}), and $\|c - p_0\| \le M' \sqrt{\epsilon}$ and $\alpha = \alpha_{H_\epsilon}(c)$, let $\bar{c}$ and $\bar{\alpha}$ satisfy
		\[
		\begin{bmatrix}
		c - p_0 \\ - \alpha + H_0(p_0)
		\end{bmatrix} = 
		B^T
		\begin{bmatrix}
		\sqrt{\epsilon}\,\bar{c} \\ - \epsilon\, \bar{\alpha}
		\end{bmatrix},
		\]
		then $\alpha_{H_\epsilon^s}(\bar{c}) = \bar{\alpha}$. 

		\item Let $c, \alpha_{N_\epsilon}(c), \bar{c}, \alpha_{H_\epsilon^s}(\bar{c})$ satisfy the relation in item 1, and suppose 
		$(\varphi_i, t_i) \in \T^n \times \T$, $(\varphi_i, \tau_i) \in \T^n \times \sqrt{\epsilon}\,\T$, $i =1,2$ satisfies 
		\[
		\begin{bmatrix}
		\theta_i \\ t_i
		\end{bmatrix} = 
		B^{-1} 
		\begin{bmatrix}
		\varphi_i \\ \tau_i/\sqrt{\epsilon}
		\end{bmatrix}
		\mod \Z^n \times \Z, 
		\]
		then 
		\[
		|h_{H_\epsilon, c}(\theta_1, t_1; \theta_2, t_2) -  \sqrt{\epsilon} h_{H_\epsilon^s, \bar{c}}(\varphi_1, \tau_1; \varphi_2, \tau_2)| \le C \epsilon. 
		\]
		\item Let $c, \alpha_{H_\epsilon}(c), \bar{c}, \alpha_{H_\epsilon^s}(\bar{c})$ be as before, then 
	\end{enumerate}
	\[
	\widetilde \cM_{H_\epsilon^s}(\bar{c})  = \Phi_L \circ \Phi_\epsilon (\widetilde \cM_{H_\epsilon}(c)), \ 
	\widetilde \cA_{H_\epsilon^s}(\bar{c})  = \Phi_L \circ \Phi_\epsilon (\widetilde \cA_{H_\epsilon}(c)), \ 
	\widetilde \cN_{H_\epsilon^s}(\bar{c})  = \Phi_L\circ \Phi_\epsilon (\widetilde \cN_{H_\epsilon}(c)).
	\]
	\end{prop}

	\ 

	Recall that $\Phi_\epsilon$ can be extended to the whole phase space, and  $N_\epsilon = H_\epsilon \circ \Phi_\epsilon$ is considered as a function on $\T^n \times \R^n \times \T$. 
	\begin{prop}\label{var-rel-norm} In the setup of Proposition
	\ref{prop:He-He-slow-var} there exists $C_3>1$ depending 
	on $H_0, k_1, \cdots, k_n, n, M_1$ such that  
	we have the following relation. 
	\begin{enumerate}
		\item $\alpha_{H_\epsilon}(c) = \alpha_{N_\epsilon}(c)$,
		$\Phi_\epsilon \, \widetilde{\mM}_{H_\epsilon}(c) = \widetilde{\mM}_{N_\epsilon}(c)$,
		$\Phi_\epsilon \, \widetilde{\cA}_{H_\epsilon}(c) = \widetilde{\cA}_{N_\epsilon}(c)$,
		$\Phi_\epsilon \, \widetilde{\N}_{H_\epsilon}(c) = \widetilde{\N}_{N_\epsilon}(c)$.
		\item $|A_{H_\epsilon, c}(\theta_1, \widetilde{t}_1; \theta_2, \widetilde{t}_2) -
		A_{N_\epsilon, c}(\theta_1, \widetilde{t}_1; \theta_2, \widetilde{t}_2)|
		\le C_3 \epsilon$.
		\item  $|h_{H_\epsilon, c}(\theta_1, {t}_1; \theta_2, {t}_2) -
		h_{N_\epsilon, c}(\theta_1, {t}_1; \theta_2, {t}_2)|  \le C_3\epsilon$.
	\end{enumerate}
	\end{prop}
	\begin{proof}[Proof of Proposition~\ref{var-rel-norm}]
	The symplectic invariance of the alpha function and the Mather, Aubry and
	Ma\~ne sets follows from exactness $\widetilde{\Phi}_\epsilon$
	(see \cite{Be2}). In order to get the quantitative estimate 
	of action, we need more detailed estimates.

	Writing $\widetilde{\Phi}_\epsilon(\theta, p, t, E) = (\Theta, P, t, \Phi_E)$, from Theorem~\ref{double-norm-form} and using the rescaled norm (\ref{rescaled-norm}),
	we have
	$$ \|\widetilde{\Phi}_\epsilon - Id\|_{C^0}\le C_1 \epsilon, \quad \|\Phi_\epsilon - Id\|_{C^1} \le C_1 \sqrt{\epsilon}
	$$
	Denote $\widetilde E=\Phi_E-E$. By exactness of $\widetilde{\Phi}_\epsilon$, we have there exists a function
	$S: \T^n \times \R^n \times \T \times \R \to \R$ such that
	\begin{equation}
	\label{eq:exact-one-form}
	P d \Theta + \Phi_E dt - (p d\theta + E dt) = P d\Theta - p d\theta + \widetilde{E}dt = dS(\theta,p, t, E) = dS(\theta, p, t  )
	\end{equation}
	In particular, given  a curve $(\theta, p, t, E)(t)$, $ t \in [\tdt_1, \tdt_2]$ with $N_\epsilon + E =0$, we have for $(\Theta, P, t, \Phi_E)(t) = \widetilde{\Phi}_\epsilon(\theta, p, t, E)$ (and hence $H_\epsilon(\Theta, P, t) + \Phi_E =0$), we apply \eqref{eq:exact-one-form} to the tangent vector of the curve to get  
	\begin{align*}
	&  \frac{d}{dt} S(\theta, p, t, E) =  P \cdot \dot{\Theta} + \Phi_E - (p \cdot \dot{\theta} + E)  \\
	& =  P \cdot \dot{\Theta} - H_\epsilon(\Theta, P, t) - (p \cdot \dot{\theta} - N_\epsilon(\theta, p, t)) \\
	&= L_{H_\epsilon}(\Theta, \dot{\Theta}, t) - L_{N_\epsilon}(\theta, \dot{\theta}, t)
	\end{align*}
	As a result,
	\begin{equation}
	\label{eq: He-Ne-action}
	\begin{aligned}
	& 	\int_{\tdt_1}^{\tdt_2} \left(  L_{H_\epsilon}(\Theta, \dot{\Theta}, t) - c \cdot \dot{\Theta}\right) dt 
	- \int_{\tdt_1}^{\tdt_2} \left(  L_{N_\epsilon}(\theta, \dot{\theta}, t) - c \cdot \dot{\theta}\right) dt \\
	& = \left(  S(\theta, p, t, E) - c \cdot (\widetilde{\Theta} - \widetilde{\theta}) \right)\Bigr |_{\tdt_1}^{\tdt_2},
	\end{aligned}
	\end{equation}
	where $\widetilde{\Theta}$, $\widetilde{\theta}$ are lifts of $\Theta(t), \theta(t)$ to the universal cover. Moreover,  from 
	$$
	\|\Phi_\epsilon- id\|_{C^0} \le C_1 \epsilon,
	$$ 
	we get  $\Theta - \theta$ is a well defined vector function on 
	$\T^n \times B_{M_1}(p_0)^n \times \T$. In particular, we hve 
	$\widetilde{\Theta} - \widetilde{\theta} = \Theta - \theta$. 

	We now estimate the $C^0-$norm of $S$. Write $S_0 = p \cdot (\Theta -\theta)$, using \eqref{eq:exact-one-form} we have 
	\[
	dS = (P - p)d\Theta + pd(\Theta - \theta) + \widetilde{E}dt = (P - p)d\Theta + dS_0 -  (\Theta - \theta)dp + \widetilde{E}dt. 
	\]
	Since $\|\Theta - \theta\|_{C^0}, \|P - p\|_{C^0}, \|\widetilde{E}\|_{C^0} = \|\Phi_E - E\|_{C^0} \le C_1 \epsilon$, we have 
	$\|d S\|_{C^0} \le C_1 \epsilon$, and $\|S-S(0)\|_{C^0} \le C' \epsilon$ for some $C'$ depending on $C_1$. 

	We now conclude that the integral in \eqref{eq: He-Ne-action} is bounded by $C' \epsilon$.  Apply this estimate to a minimizer, we have 
	$$ |A_{N_\epsilon, c}(\theta(\tdt_1), \widetilde{t}_1; \theta(\tdt_2), \widetilde{t}_2) -
	A_{H_\epsilon,c}(\Theta(\tdt_1), \tdt_1; \Theta(\tdt_2), \tdt_2))|
	\le C'\epsilon. $$
	Since $\|\Phi_\epsilon - Id\|_{C^0} \le C_1\epsilon$, we have
	$\|\theta(\tdt_i) - \Theta(\tdt_i)\|\le C_1 \epsilon$, $i = 1,2$. The estimate follows from the Lipschitz property of $A_{H,c}$.

	Taking limit, we obtain the estimate for the barrier $h_{H,c}$.
	\end{proof}

	We now study the relation between $N_\epsilon$ and $H_\epsilon^s$. We extend the definition of $H_\epsilon^s$
	to $\T^n \times \R^n \times \R$ such that the $\sqrt{\epsilon}P$ term is supported on the set $\{\|(I, \sqrt{\epsilon}F)\| < 2 M_2\}$. 

	\begin{prop}\label{prop:slow-variational}
	Assume the conclusions of Theorem~\ref{double-norm-form} hold 
	on the set 
	$$
	\{\|(p - p_0, H+ H_0(p_0)\| < M_1 \sqrt{\epsilon}\}
	$$ 
	for sufficently large $M_1$ depending only on $H_0$. 
	Then for $M_3 = C_B^{-2} M_1$, we have: 
	\begin{enumerate}
		\item Let $(\theta, p, t, E)$ satisfies $\|(p-p_0, E + H_0(p_0))\|< M_3$, $N_\epsilon(\theta, p, t) + E =0$, 
		and 
		\[
		(\varphi, I, \tau, F) = \Phi_L(\theta, p, t, E). 
		\]
		Let $L_{N_\epsilon}$ and $L_{H_\epsilon^s}$ be the  Lagrangians for $N_\epsilon$ and $H_\epsilon^s$. Then for 
		\[
		v = \partial_p N_\epsilon(\theta, p, t), \quad 
		\bmat{v^s\\ v^f} =  B \bmat{v \\ 1}, 		
		\]
		we have 
		\[
		L_{N_\epsilon}(\theta, v, t) - p_0 \cdot v + H_0(p_0) = \epsilon v^f L_{H_\epsilon^s}\left(  \varphi, \frac{1}{\sqrt{\epsilon}}\frac{v^s}{v^f} , s\sqrt{\epsilon} \right). 
		\]

		\item Suppose $(c, \alpha), (\bar{c}, \bar{\alpha}) \in \R^n \times \R$ satisfies $\|c - p_0\| \le \frac{M_3}2 \sqrt{\epsilon}$ and 
		\[
		\begin{bmatrix}
		c - p_0 \\ - \alpha + H_0(p_0)
		\end{bmatrix} = 
		B^T
		\begin{bmatrix}
		\sqrt{\epsilon}\bar{c} \\ \epsilon\bar{\alpha}
		\end{bmatrix},
		\]
		then $\alpha = \alpha_{N_\epsilon}(c)$ if and only if $\bar{\alpha} = \alpha_{H_\epsilon^s}(\bar{c})$. 
		\item Let $c, \alpha_{N_\epsilon}(c), \bar{c}, \alpha_{H_\epsilon^s}(\bar{c})$ satisfies the relation in item 2, and suppose 
		$(\varphi_i, t_i) \in \T^n \times \T$, $(\varphi_i, \tau_i) \in \T^n \times \sqrt{\epsilon}\,\T$, $i =1,2$ satisfies 
		\[
		\begin{bmatrix}
		\theta_i \\ t_i
		\end{bmatrix} = 
		B^{-1} 
		\begin{bmatrix}
		\varphi_i \\ \tau_i/\sqrt{\epsilon}
		\end{bmatrix}
		\mod \Z^n \times \Z, 
		\]
		then 
		\[
		h_{N_\epsilon, c}(\theta_1, t_1; \theta_2, t_2) = \sqrt{\epsilon} h_{H_\epsilon^s, \bar{c}}(\varphi_1, \tau_1; \varphi_2, \tau_2). 
		\]
		\item Let $c, \alpha_{N_\epsilon}(c), \bar{c}, \alpha_{H_\epsilon^s}(\bar{c})$ be as before, then 
		\[
		\widetilde{\cM}_{H_\epsilon^s}(\bar{c})  = \Phi_L(\widetilde{\cM}_{N_\epsilon}(c)), \quad
		\widetilde{\cA}_{H_\epsilon^s}(\bar{c})  = \Phi_L(\widetilde{\cA}_{N_\epsilon}(c)), \quad
		\widetilde{\cN}_{H_\epsilon^s}(\bar{c})  = \Phi_L(\widetilde{\cN}_{N_\epsilon}(c))
		\]
	\end{enumerate}
	\end{prop}

	\begin{proof}[Proof of Proposition~\ref{prop:slow-variational}]
	The choice of $M_3$ is to ensure that for $M_2 = C_B^{-1}M_1$ as in Lemma~\ref{lem:G-eps-slow},
	\[
	\Phi_L (\{\|(p - p_0, E + H_0(p_0))\| < M_3\}) \subset \{\| (I, \sqrt{\epsilon}F)\| < M_2\}. 
	\]
	\emph{Item 1. } 
	Recall that  $G_\epsilon = H_\epsilon + E$, and 
	\begin{equation}
	\label{eq:LNe}
	\begin{aligned}
	& L_{N_\epsilon} - p_0 \cdot v + H_0(p_0) = \sup_p \left\{  (p-p_0) \cdot v - N_\epsilon(\theta, p, t)  + H_0(p_0) \right\} \\
	& = \sup_{p, E}\left\{  (p-p_0, E + H_0(p_0)) \cdot (v, 1): \quad N_\epsilon(\theta, p, t) + E =0  \right\},
	\end{aligned}
	\end{equation}
	where the supremum is acheived at a unique point since $H_\epsilon + E \le  0$ is a strictly convex set. Let us denote $L_{N_\epsilon}^{p_0}  = L_{N_\epsilon} - p_0 \cdot v + H_0(p_0)$.

	Continuing from  \eqref{eq:LNe} and using \eqref{eq:linear-rescale}, we get 
	\begin{align*}
	L_{N_\epsilon}^{p_0}(\theta, v, t) &= \sup\left\{   (B^T)^{-1} \bmat{p - p_0\\ E + H_0(p_0)} 
	\cdot B \bmat{v \\ 1} : \quad G_\epsilon(\theta, p, t, E)  = 0 \right\} \\
	& =  \sup \left\{ \sqrt{\epsilon} (I, \sqrt{\epsilon} F) \cdot (v^s, v^f) : 
	\quad G_\epsilon^s(\varphi, I, \tau, F) = 0  \right\}  \\
	& =  \sqrt{\epsilon} \sup \left\{ (I,F)\cdot (v^s, v^f\sqrt{\epsilon}) : \quad H_\epsilon^s(\varphi, I, \tau) + F =0     \right\}, \\
	& =  \epsilon v^f \sup \left\{ (I,F)\cdot (v^s/(v^f \sqrt{\epsilon}), 1) : \quad H_\epsilon^s(\varphi, I, \tau) + F =0     \right\} \\
	& = \epsilon v^f L_{H_\epsilon^s}(\varphi, v^s/(v^f \sqrt{\epsilon}), \tau) . 
	\end{align*}
	
	\emph{Item 2.}  We first derive a relation between the integrals of the Lagrangians. Let $(\theta, p, t, E)(t)$, $t\in [t_1, t_2]$  be a solution to $N_\epsilon + E$ with 
	$$
	N_\epsilon + E =0\quad \text{ and }\quad
	\|(p - p_0, E + H_0(p_0))\| < M_1 \sqrt{\epsilon}.
	$$ 
	Let $(\varphi, I, \tau, F)(t) = \Phi_L^{-1} (\theta, p, t, E)(t)$. Then from \eqref{eq:linear-rescale}, 
	\[
	\bmat{ \dot{\varphi}\\ \dot{\tau}/\sqrt{\epsilon} } = B \bmat{ \dot{\theta} \\ 1} := \bmat{v^s \\ v^f}(t) , \quad \frac{d\varphi}{d\tau} = \frac{v^s}{\sqrt{\epsilon} v^f}. 
	\]
	From item 1 we get 
	\begin{align*}
	&\int_{t_1}^{t_2} L_{N_\epsilon}^{p_0}(\theta, \dot{\theta}, t) dt 
	= \int_{t_1}^{t_2} \epsilon v^f(t) L_{H^s_\epsilon}\left(\varphi, v^s/(v^f\sqrt{\epsilon}), \tau  \right) dt \\ 
	& = \int_{\tau_1}^{\tau_2} \epsilon v^f(\tau) L_{H_\epsilon^s}(\varphi(\tau), \frac{d\varphi}{d\tau}(\tau), \tau) \frac{dt}{d\tau}(\tau) d\tau 
	= \sqrt{\epsilon} \int_{\tau_1}^{\tau_2} L_{H_\epsilon^s}(\varphi, \frac{d\varphi}{d\tau}, \tau) d\tau.  
	\end{align*}	

	Let $\widetilde{\theta}(t)$ be a lift of $\theta(t)$ to the universal cover and $\widetilde{\varphi}(t), \widetilde{\tau}(t)$ be a lift of $(\varphi(t), \tau(t))$. Then 
	\[
	\bmat{\widetilde{\varphi}(t) \\ \widetilde{\tau}(t)/\sqrt{\epsilon}} = B \bmat{\widetilde{\theta}(t) \\ t} + const. 
	\]
	Given any $(c, -\alpha) \in \R^n \times \R$, we have 
	\begin{equation}
	\label{eq:N-eps-H-slow-action}
	\begin{aligned}
	&\int_{t_1}^{t_2} \left(  L_{N_\epsilon}(\theta, \dot{\theta}, t) - c\cdot \dot{\theta} + \alpha \right) dt \\
	& = \int_{t_1}^{t_2} L_{N_\epsilon}^{p_0}(\theta, \dot{\theta}, t) dt - (c -p_0, -\alpha + H_0(p_0)) \cdot \left(  (\widetilde{\theta}(t_2)-\widetilde{\theta}(t_1), t_2 - t_1) \right) \\
	& =  \sqrt{\epsilon}\int_{\tau_1}^{\tau_2} L_{H_\epsilon^s} d\tau - (B^T)^{-1} \bmat{ c - p_0\\ - \alpha + H_0(p_0)} \cdot B \bmat{ \widetilde{\theta}(t_2) - \widetilde{\theta}(t_1) \\ t_2 - t_1} \\
	& =  \sqrt{\epsilon}\int_{\tau_1}^{\tau_2} L_{H_\epsilon^s} d\tau  - (\sqrt{\epsilon}\bar{c}, -\epsilon\bar{\alpha}) \cdot (\widetilde{\varphi}(\tau_2) - \widetilde{\varphi}(\tau_1), (\widetilde{\tau}_2 - \widetilde{\tau}_1)/\sqrt{\epsilon}) \\
	& = \sqrt{\epsilon} \int_{\tau_1}^{\tau_2} \left(  L_{H_\epsilon^s} - \bar{c}\cdot \frac{d\varphi}{d\tau} + \bar{\alpha} \right)(\varphi, \frac{d\varphi}{d\tau}, \tau) d\tau. 
	\end{aligned}
	\end{equation}

	We now prove Item 2. Suppose $\alpha = \alpha_{N_\epsilon}(c)$, let  $(\theta, p, t, E)(t)$, $t \in [-\infty, \infty)$ be an orbit in the Mane set $\cN_{N_\epsilon}(c)$. Then there exists $C_*>0$ such that  we have $\|p(t) - c\| \le C^* \sqrt{\epsilon}$, by Proposition~\ref{prop:semi-concave-near-int}. Since $\|c - p_0\| \le M_3\sqrt{\epsilon}/2$, with $M_3$ large enough, we have $\|p(t) - p_0\| \le M_3 \sqrt{\epsilon}$ for all $t$.

	Moreover, using the fact that the orbit is semi-static, we have 
	\[
	-C \le \int_0^T \left(  L_{N_\epsilon} - c \cdot v  + \alpha \right)(\theta(t), \dot{\theta}(t),t) dt \le C.
	\] 
	From \eqref{eq:N-eps-H-slow-action}, we get 
	\[
	-C \le \int_{\tau(0)}^{\tau(T)}   \left(  L_{H_\epsilon^s} - \bar{c} \cdot v  + \bar{\alpha}\right) d\tau \le C
	\]
	for all $T > 0$. This implies $\bar{\alpha}$ is Ma\~ne critical for $\bar{c}$.

	As a result, we obtain that $\alpha$ is Mane critical for $L_{N_\epsilon} - c\cdot v$ implies $\bar{\alpha}$ is Mane critical for $L_{H_\epsilon^s} - \bar{c} \cdot v$. Moreover the converse is true by reversing the above computations. 

	\emph{Item 3}. We apply \eqref{eq:N-eps-H-slow-action} to any one-sided minimizer of $h_{N_\epsilon, c}(\theta_1, t_1; \theta_2, t_2)$, and use the localization of calibrated orbits. 

	\emph{Item 4}. \eqref{eq:N-eps-H-slow-action} implies that an orbit $(\theta, p, t, E)(t)$ is semi-static for $L_{N_\epsilon, c}$ implies $(\varphi, I, \tau, F) = \Phi_L^{-1}(\theta, p, t, E)$ is a reparametrization of a semi-static orbit of $L_{H_\epsilon^s, \bar{c}}$. The converse also holds. This implies the relation between Mane sets. The same applies to Aubry sets. We note that the Mather set is  precisely the support of all invariant measures contained in the Aubry set, and therefore is also invariant. 
	\end{proof}


\section{Variational aspects of the slow mechanical system}
\label{sec:var-slow}

In this section we study the variational properties of the slow mechanical system
\[
H^s(\varphi, I) = K(I) - U(\varphi), 
\]
with $\min U = U(0) = 0$.

The main goal of this section is to derive some properties \ of \ the \
``channel'' \ $\bigcup_{E>0}\LF_\beta(\lambda_h^E h)$, and information about
the Aubry sets for $c\in \LF_\beta(\lambda_h^E h)$. More precisely, we prove Proposition~\ref{slow-localization} and justify the picture Figure~\ref{fig:sim-non-simp-coh}. 

\begin{itemize}
	\item In section~\ref{sec:basic-channel}, we show that each $\LF_\beta(\lambda_h^E h)$ is an segment parallel  to $h^\perp$.
	\item In section~\ref{sec:char-channel}, we provide a characterization of the segment, and provide information about the Aubry sets.
	\item In section~\ref{sec:width}, we provide a condition for the ``width'' of
	the channel to be non-zero.
	\item In section~\ref{sec:critical}, we discuss the limit of the set $\LF_\beta(\lambda_h^E h)$ as $E\to 0$
	which corresponds to the ``bottom'' of the channel.
\end{itemize}

We drop all supscripts ``s'' to simplify the notations.
The results proved in this section are mostly contained in \cite{Ma6} in some form. Here we reformulate some of them for our purpose and also provide some different proofs.

\subsection{Relation between the minimal geodesics and the Aubry sets}
\label{sec:basic-channel}

Assume that $H(\varphi, I)$ satisfies the conditions $[DR1^h]-[DR3^h]$ and $[DR1^c]-[DR4^c]$. Then for $E\ne E_j$, $1\le j \le N-1$, there exists a unique shortest
geodesic $\gamma_h^E$ for the metric $g_E$ in the homology $h$. For the bifurcation
values $E=E_j$, there are two shortest geodesics $\gamma_h^E$ and $\bar{\gamma}_h^E$.

The function $l_E(h)$ denotes the length of the shortest $g_E-$geodesic in homology
$h$. By Lemma \ref{lm:monot-length} the length function $l_E(h)$ is continuous and 
strictly increasing in $E\ge 0$. It is easy to see that 
it is positive homogeneous ($l_E(nh)=nl_E(h)$, $n\in \NN$)
and sub-additive ($l_E(h_1 + h_2) \le l_E(h_1) + l_E(h_2)$) 
in $h$.

Assume that the curves $\gamma_h^E$ are parametrized using the Maupertuis
principle, namely, it is the projection of the associated Hamiltonian orbit. Let $T(\gamma_h^E)$ be the period under this parametrization, and write $\lambda(\gamma_h^E) = 1/(T(\gamma_h^E))$.

We pick another vector $\bar{h}\in H_1(\T^2, \Z)$ such that
$h, \bar{h}$ form a basis of $H_1(\T^2, \Z)$ and for
the dual basis $h^*, \bar{h}^*$ in $H^1(\T^2, \R)$ we have
$\langle h, \bar h^*\rangle =0$.
We denote $\bar h^*$ by $h^\perp$ to emphasise the latter fact.

The main result of this section is:
\begin{thm}
\label{prop-lf}
\begin{enumerate}
	\item For $E=E_j$,
	$$\LF_\beta(\lambda(\gamma_h^E)\cdot h) = \LF_\beta(\lambda(\bar\gamma_h^E)
	\cdot h). $$
	As a consequence, write $\lambda_h^E = \lambda(\gamma_h^E)$, then the set
	$\LF_\beta(\lambda_h^E h)$ is well defined (the definition is independent of
	the choice of $\gamma_h^E$).

	\item
	For each $E>0$, there exists $-\infty \le a_E^-(h) \le a_E^+(h)\le \infty$
	such that
	$$
	\LF_\beta(\lambda_h^E h) = l_E(h) h^* + [a_E^-(h), a_E^+(h)] \ h^\perp.
	$$
	Moreover, the set function $[a_E^-, a_E^+]$ is upper semi-continuous in $E$.

	\item For each $c \in \LF_\beta(\lambda_h^E h)$, $E \ne E_j$, there is a unique $c-$minimal measure supported on $\gamma_h^E$.

	\item For each $c \in \LF_\beta(\lambda_h^{E_j} h)$, there are two $c-$minimal measures supported on $\gamma_h^{E_j}$ and $c$.

	\item For $E>0$, assume that the torus $\T^2$ is not completely foliated by
	shortest closed $g_E-$geodesics in the homology $h$, then
	$a_E^+(h)-a_E^-(h)>0$ and the channel has non-zero width.
\end{enumerate}
\end{thm}

Assume that $\gamma$ is a geodesic parametrized according to the Maupertuis principle. First, we note the following useful relation.
\begin{equation}\label{eq:rel-maupertuis}
L(\gamma, \dot{\gamma}) + E = 2(E+U(\gamma)) = \sqrt{g_E(\gamma, \dot{\gamma})},
\end{equation}
where $L$ denote the associated Lagrangian.

According to  \cite{DC}, the minimal measures for $L$ is
in one-to-one correspondence with the minimal measures of $\frac12 g_E(\varphi, v)$.
On the other hand, any minimal measure $\frac12 g_E$ with a rational rotation number
is supported on closed geodesic. The following lemma characterizes
minimal measures supported on a closed geodesic.

\begin{lem}\label{orthogonal}
\begin{enumerate}
	\item Assume that $c\in H^1(\T^2, \R)$ is such that $\alpha_H(c)=E>0$.
	Then for any $h\in H_1(\T^2, \Z)$,
	$$	l_E(h) - \langle c, h \rangle \ge 0. $$
	\item Let $\gamma$ be a closed geodesic of $g_E$, $E>0$, with
	$[\gamma]=h\in H_1(\T^2, \Z)$. Let $\mu$ be the invariant measure
	supported on the periodic orbit associated to $\gamma$. Then given
	$c\in H^1(\T^2, \R)$  with $\alpha_{H^s}(c)=E$,
	\begin{equation}\label{eq:length-c}
	\mu \text{ is }c-\text{minimal if and only if } \ \
	l_E(h) - \langle c, h \rangle = 0.
	\end{equation}
	\item Let $\gamma$ be a closed geodesic $g_E$, $E\ge 0$, with
	$[\gamma]=h\in H_1(\T^2, \Z)$ and $\alpha_{H^s}(c)=E$.
	Then $\gamma \subset \cA_{H^s}(c)$ if and only if \eqref{eq:length-c} holds.
\end{enumerate}
\end{lem}
\begin{proof}
Let $\gamma$ be a closed geodesic of $g_E$, $E>0$, with $[\gamma]=h$.
Assume that with the Maupertuis
	parametrization, the periodic of $\gamma$ is $T$. Let $\mu$ be
	the associated invariant measure, then $\rho(\mu)=h/T$. Assume that
	$\alpha(c)=E$, by definition, we have
	$$ \int L d\mu + E \ge \beta(h/T) + \alpha(c) \ge \langle c, h/T\rangle. $$
	By \eqref{eq:rel-maupertuis}, we have
	$$ \int L d\mu + E = \frac{1}{T}\int_0^T (L + E)(d\gamma) =
	\frac{1}{T}\int_0^T \sqrt{g_E(d\gamma)} = l_E(\gamma)/T.  $$
	Combine the two expressions, we have $l_E(\gamma)- \langle c, h\rangle\ge 0$.
	By choosing $\gamma$ such that $l_E(\gamma) = l_E(h)$, statement 1 follows.

	To prove statement 2, notice that if $\mu$ is $c-$minimal, then $\alpha(c)=E$ and
	the equality
	$$ \int L d\mu + E = \langle c, h/T\rangle$$
	holds. Equality \eqref{eq:length-c} follows from the same calculation as statement 1.

	For $E>0$, $\gamma \subset \cA_{H^s}(c)$ if an only if $\gamma$ is a minimal measure.
	Hence we only need to prove statement 3 for $E=0$. In this case, $\gamma $ can
	be parametrized as a homoclinic orbit. $\gamma \subset \cA_{H^s}(c)$ if and only if
	$$
	\int_{-\infty}^\infty (L - c \cdot v + \alpha(c))(d\gamma) =0.
	$$
	Since
	$$  \int_{-\infty}^\infty (L - c \cdot v + \alpha(c))(d\gamma) =
	\int_{-\infty}^\infty (L+E)(d\gamma) -
	\langle c, h\rangle = l_E(h) - \langle c, h\rangle, $$
	the statement follows.
	\end{proof}

\begin{proof}[Proof of Theorem~\ref{prop-lf}, item 1-4]
By Lemma~\ref{orthogonal}, if there are two shortest geodesics $\gamma_h^E$
and $\bar{\gamma}_h^E$ for $g_E$, for any $c$, the invariant measure supported
on $\gamma_h^E$ is $c-$ minimal if and only if the measure on $\bar{\gamma}_h^E$
is $c-$minimal. This implies statement 1.

Statement 2 follows from  the fact that $\LF_\beta(\lambda_h^Eh)$ is a closed convex set, and
\eqref{eq:length-c}.

Statement 3 and 4 follows directly from Lemma~\ref{orthogonal}. Item 5 is proved in 
Proposition~\ref{non-zero-width}. 
\end{proof}

We also record the following consequence of Theorem 4.1. 
\blm \label{lm:monot-length}
The period $T_E:=T(\gm^E_h)$ and action 
$l^E_h$ as functions of $E$ are strictly monotone for $E>0$. 
\elm 

\begin{proof}
Let $\gamma_h^E$ be a minimal geodesic of $g_E$, then for any $E' < E$
we have 
\[
l_{E'}(h) \le \int \sqrt{g_{E'}(\dot{\gamma}_h^E)} < \int \sqrt{g_E(\dot{\gamma}_g^E)} = l_E(h), 
\]
therefore $l_E(h)$ is strictly monotone. 

To prove strict monotonicity of $T_E$ application of  
Theorem~\ref{prop-lf} part 2 gives
\[
\cL\mathcal F_\bt (h/T_E)\cdot h= l^E_h h^* \cdot h, 
\]
where $h^*\cdot h>0$. Since $\cL\mathcal F_\bt (h/T_E)$ are distinct for 
different $E$, $T_E \ne T_{E'}$ for $E \ne E'$, proving strict monotonicity. 
\end{proof}

\subsection{Characterization of the channel and the Aubry sets}
\label{sec:char-channel}

In this section we provide a precise characterization of the set
$$ \LF_\beta(\lambda_h^E h) = l_E(h)\,h^* + [a_E^-(h), a_E^+(h)]\,h^{\perp}.  $$
For each $E>0$, define
\[
d_E^\pm(h) = \pm \inf_{n\to \infty} (l_E(nh \pm \bar{h})- l_E(nh)),	
\]
where $h, \bar h$ is a basis in $H_1(T^2,\R)$ and the dual of $\bar h$
satisfies $\langle  \bar h^*, h \rangle =0$ so we denote it $h^\perp$.
Note that the sequence $l_E(nh \pm \bar{h})- l_E(nh)$ is decreasing, so
the infimum coincides with the limit. We will omit
dependence on $h$ when it is not important.

\begin{lem}\label{lem:cont-width}
The function $d_E^\pm(h)$ is continuous in $E>0$. 
\end{lem}
\begin{proof}

From  the sub-linearity of $l_E(h)$, we get,
\[
-l_E(\bar h) \le l_E(nh \pm \bar{h})- l_E(nh) \le l_E(\bar h),
\]
as a result, the family $l_E(nh \pm \bar{h})- l_E(nh)$ as a function of $E \in [E_1, E_2] \subset (0, \infty)$
is equi-continuous and equi-bounded. Then there exists a subsequence such that 
\[
\lim_{n\to \infty} l_E(nh \pm \bar{h})- l_E(nh) = d_E^+(h)
\]
uniformly. This implies $d_E^+(h)$ is continuous. The same argument works for $d_E^-(h)$. 
\end{proof}

\begin{prop}
For each $E>0$, we have
$$ d_E^\pm(h) = a_E^\pm(h). $$
\end{prop}
\begin{proof}
We first show
$$ d_E^-(h) \le a_E^-(h) \le a_E^+(h) \le d_E^+(h). $$
Omit dependence on $h$. Denote $c^+= l_E(h) h^* + a_E^+ h^{\perp}$,
by definition, $l_E(h) - \langle c^+ , h\rangle=0$. By Lemma~\ref{orthogonal},
statement 1, for $n\in \NN$,
\begin{multline*}
0 \le l_E(n h + \bar{h}) - \langle c^+, n h + \bar{h}\rangle =
l_E(n h + \bar{h}) - nl_E(h) - \langle c^+, \bar{h} \rangle \\
= l_E(n h + \bar{h}) - nl_E(h) - a_E^+.
\end{multline*}
Take infimum in $n$, we have $d_E^+ - a_E^+ \ge 0$. Perform the same calculation
with $nh+ \bar{h}$ replaced by $nh - \bar{h}$, we obtain
$0 \le l_E(nh - \bar{h}) - nl_E(h) + a_E^-$, hence $a_E^- - d_E^- \ge 0$.

We now prove the opposite direction. Take any  $c \in l_E(h) h^* + [d_E^-, d_E^+]h^\perp$, we first show that $\alpha(c)=E$.

Take $\rho \in \Q h + \Q \bar{h}$, then any invariant measure $\mu$
with rotation number $\rho$ is supported on some
$[\gamma] =m_1 h + m_2 \bar{h}$ with $m_1, m_2 \in \Z$. Let $T$ denote
the period, by Lemma~\ref{act-lowerbnd} below,
$$
\beta(\rho) + E = l_E(m_1 h + m_2 \bar{h})/T  \ge \langle c, m_1 h + m_2 \bar{h}\rangle /T =
\langle c, \rho\rangle.
$$
Since $\beta$ is continuous, we have $\alpha(c) = \sup \langle c,\
\rho\rangle - \beta(\rho) \le  E$,
where the supremum is taken over all rational $\rho$'s.
Since the equality is achieved at $\rho = h$, we conclude that $\alpha(c)=E$.

By Lemma~\ref{orthogonal}, statement 2, the measure supported on $\gamma_h^E$ is $c-$minimal, and hence $c\in \LF_\beta(\lambda_h^E h)$.
\end{proof}

Recall $h, \bar h$ form a basis in $H_1(T^2,\Z)$ and the dual of $\bar h$ is perpendicular to $h$ and denoted by $h^\perp$.
\begin{lem}\label{act-lowerbnd}
For any $c \in l_E(h)\, h^* + [d_E^-, d_E^+]\, h^\perp$ and
$m_1, m_2 \in \Z$, we have
$$
l_E(m_1 h + m_2 \bar{h}) - \langle c , m_1 h + m_2\bar{h}\rangle \ge 0 .
$$
Moreover, if $c \in l_E(h) h^* + (d_E^-, d_E^+)h^\perp$ and $m_1, m_2 \neq 0$,
there exists $a>0$ such that
$$
l_E(m_1 h + m_2 \bar{h}) - \langle c , m_1 h + m_2\bar{h}\rangle >a> 0.
$$
\end{lem}
\begin{proof}
The inequality for $m_1=0$ or $m_2=0$ follows from positive homogeneity of $l_E$.
We now  assume $m_1, m_0$.

If $m_2 >0$, for a sufficiently large $n\in \NN$,  have
\begin{align*}
&l_E(m_1 h + m_2 \bar{h}) - \langle c, m_1 h + m_2 \bar{h} \rangle \\
&= l_E(m_1 h + m_2 \bar{h})  + l_E((nm_2 - m_1)h) -  \langle c, nm_2 h + m_2 \bar{h} \rangle \\
&\ge l_E(m_2(nh+ \bar{h})) - \langle c, m_2(nh + \bar{h}) = m_2 (l_E(nh + \bar{h}) - \langle c, nh + \bar{h} \rangle) \\
& \ge  l_E(nh + \bar{h}) - \langle c, nh + \bar{h} \rangle .
\end{align*}
Since
$$ l_E(nh + \bar{h}) - \langle c, nh + \bar{h} \rangle = l_E(nh + \bar{h}) - nl_E(h) - \langle c, \bar{h}\rangle, $$
for $c \in l_E(h) h^* + (d_E^-, d_E^+)h^\perp$, then there exists $a>0$ such that for sufficiently large $n$,
$$ \lim_{n\to \infty}l_E(nh + \bar{h}) - nl_E(h) - \langle c, \bar{h}\rangle > a. $$

For $m_2 <0$, we replace the term $(nm_2 - m_1)h$ with $(-nm_2 - m_1)h$ in the above calculation.
\end{proof}

We have the following characterization of the Aubry sets for the cohomologies contained in the channel.

\begin{prop}\label{char-aubry}
For any $E>0$ and $c \in l_E(h) h^* + (d_E^-, d_E^+)h^\perp$, we have
$$ \cA_{H^s}(c) = \gamma_h^E $$
if $E$ is not a bifurcation value and
$$ \cA_{H^s}(c) = \gamma_h^E \cup \bar\gamma_h^E $$
if $E$ is a bifurcation value.
\end{prop}

\begin{proof}
We first consider the case when $E$ is not a bifurcation value.
Since $\gamma_h^E$ is the unique closed shortest geodesic, if $\cA_{H^s}(c)\supseteq \gamma_h^E$,
it must contain an infinite orbit $\gamma^+$. Moreover, as $\gamma_h^E$ supports
the unique minimal measure, the orbit $\gamma^+$ must be biasumptotic to $\gamma_h^E$.
As a consequence, there exists $T_n, T_n'\to \infty$ such that $\gamma^+(-T_n)-\gamma^+(T_n')\to 0$.
By closing this orbit using a geodesic, we obtain a closed piece-wise geodesic curve $\gamma_n$.
Moreover, since $\gamma^+$ has no self-intersection, we can arrange it such that $\gamma_n$ also
have no self-intersection.  We have
$$ \int (L- c\cdot v + \alpha(c))(d\gamma_n) = \int (L+ E)(d\gamma_n) - \langle c, [\gamma_n]\rangle
= l_E(\gamma_n) - \langle c, [\gamma_n]\rangle. $$
By the definition of the Aubry set, and take limit as $n\to \infty$, we have
$$ \lim_{n\to \infty} l_E(\gamma_n) - \langle c, [\gamma_n]\rangle = 0. $$
Since $\gamma_n$ has no self intersection, we have $[\gamma_n]$ is irreducible. 
However, this contradicts to the strict inequality obtained in Lemma~\ref{act-lowerbnd}.

We now consider the case when $E$ is a bifurcation value, and there are two shortest 
geodesics $\gamma_h^E$ and $\bar{\gamma}_h^E$. Assume by contradiction that 
$\cA_{H^s}(c)\supseteq \gamma_h^E \cup \bar{\gamma}_h^E$. For
mechanical systems on $\T^2$, the Aubry set satisfies an ordering property. As a consequence, 
there must exist two infinite orbits $\gamma^+_1$ and $\gamma_2^+$
contained in the Aubry set, where $\gamma_1^+$ is forward asymptotic to
$\gamma_h^E$ and backward asymptotic to $\bar\gamma_h^E$, and $\gamma_2^+$ is
forward asymptotic to $\bar{\gamma}_h^E$ and backward asymptotic to $\gamma_h^E$.
Then there exists $T_n, T_n', S_n, S_n' \to \infty$ such that
$$
\gamma_1^+(T_n')-\gamma_2^+(-S_n), \gamma_2^+(S_n')- \gamma_1^+(-T_n) \to 0
$$
as $n\to \infty$. The curves $\gamma_{1,2}^+$, $\gamma_h^E$, $\bar{\gamma}_h^E$
are all disjoint on $\T^2$. Similar to the previous case, we can construct
a piecewise geodesic, non-self-intersecting closed curve $\gamma_n$ with
$$
\lim_{n\to \infty}\int (L- c\cdot v + \alpha(c))(d\gamma_n) =0.
$$
This, however, lead to a contradiction for the same reason as the first case.
\end{proof}

\subsection{The width of the channel}
\label{sec:width}

We show that under our assumptions, the ``width'' of the channel
$$
d_E^+(h) - d_E^-(h) = \inf_{n\in \NN} (l_E(nh + \bar{h})- l_E(nh)) +
\inf_{n\in \NN} (l_E(nh - \bar{h})- l_E(nh)),
$$
is non-zero.

The following statement is a small modification of a theorem of Mather (see \cite{Ma6}),
we provide a proof using our language.

\begin{prop}\label{non-zero-width}
For $E>0$, assume that the torus $\T^2$ is not completely foliated by shortest
closed $g_E-$geodesics in the homology $h$. Then
$$ d_E^+(h) - d_E^-(h) >0. $$
\end{prop}

\brm
This is the last item of Theorem \ref{prop-lf}.
\erm

\begin{proof}
Let $\M$ denote the union of all shortest closed $g_E-$geodesics in the homology $h$.
We will show that $\M \ne \T^2$ implies $d_E^+(h) - d_E^-(h) >0$. Omit $h$ dependence.
For $n\in \N$, denote
$$
d_n =  (l_E(nh + \bar{h})- l_E(nh)) + (l_E(nh - \bar{h})- l_E(nh)).
$$
Assume by contradiction that $\inf d_n = \lim d_n >0$.

Let $\gamma_0$ be a shortest geodesic in homology $h$. We denote $\widetilde\gamma_0$
its lift to the universal cover, and use ``$\le_{\tilde{\gamma}}$'' to denote the order on $\widetilde\gamma$
defined by the flow. Let $\gamma_1$ and $\gamma_2$ be shortest curves in the homology
$nh + \bar{h}$ and $nh-\bar{h}$ respectively, and let $T_1$ and $T_2$ be their periods.
$\gamma_i$ depends on $n$ but we will not write it down explicitly.

	Let $\widetilde{\gamma}_i$, $i=0,1,2$ denote a lift of $\gamma_i$ to the universal cover.
	Using the standard curve shortening lemma in Riemmanian geometry, it's easy to see
	that $\widetilde{\gamma}_i$ and $\widetilde{\gamma}_j$ may intersect at most once.
	Let $a\in \gamma_0 \cap \gamma_1$ and lift it to the universal cover without changing
	its name. Let $b\in \gamma_0 \cap \gamma_2$, and we choose a lift in $\widetilde\gamma_0$
	by the largest element such that $b\le a$. We now choose the lifts $\widetilde\gamma_i$ of
	$\gamma_i$, $i=1,2$, by the relations $\widetilde\gamma_1(0)=a$ and $\widetilde\gamma_2(T_2)=b$.

	We have for $1\le k \le 2n$, $\widetilde\gamma_2(T_2)+kh > \widetilde\gamma_1(0)$ and
	$$
	\widetilde\gamma_2(0)  + kh  = b - (nh -\bar{h}) + kh \le_{\tilde{\gamma}_0} a + nh + \bar{h} = \widetilde\gamma_1(T_1).
	$$
	As a consequence, $\widetilde\gamma_2 + kh$ and $\widetilde\gamma_1$ has a unique intersection. Let
	$$
	x_k = (\widetilde\gamma_2 + kh)\cap \widetilde\gamma_1, \quad \bar{x}_k =
	(\widetilde\gamma_2+ kh)\cap (\widetilde\gamma_1 - h).
	$$
	We have $x_k$ is in increasing order on $\widetilde\gamma_1$ and $\bar{x}_k$ is in decreasing
	order after projection to $\gamma_2$ (see Figure~\ref{fig:shortes-geod}) . Define
	$$
	\widetilde\gamma_k^* = (\widetilde\gamma_2 + kh)|[\bar{x}_k, x_k] * \widetilde\gamma_1|[x_k, x_{k+1}],
	$$
	and let $\gamma_k^*$ be its projection to $\T^2$. We have $[\gamma_k]=h$ and
	$$ \sum_{k=1}^{2n} l_E(\gamma_k^*) = l_E(\gamma_1) + l_E(\gamma_2).$$
	Using $l_E(\gamma_k)\ge l_E(h)$ and $l_E(\gamma_1) + l_E(\gamma_2) \le 2n l_E(h) + d_n$, we obtain
	$$ l_E(h) \le l_E(\gamma_k) \le l_E(h) + d_n. $$

	\begin{figure}[t]
	\centering 
	\includegraphics[width=5in]{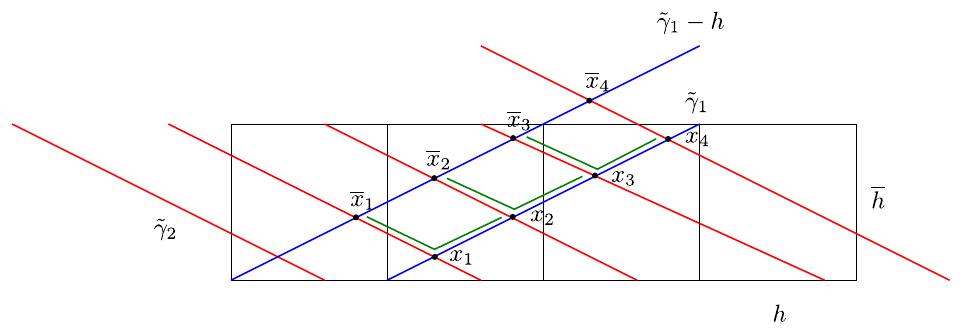} 
	\caption{Proof of Proposition~\ref{non-zero-width}, green curves are $\tilde{\gamma}_k^*$'s. }
	\label{fig:shortes-geod}
	\end{figure}

	Any connected component in the completement of $\M$ is diffeomorphic to an annulus.
	Pick one such annulus, and let $b>0$ denote the distance between its boundaries.
	Since $\gamma_1$ intersects each boundary once, there exists a point $y_n\in \gamma_1$
	such that $d(y_n, \M)= b/2$. Since $\gamma_1 \subset \bigcup_k \gamma_k^*$ there exists
	some $\gamma_k^*$ containing $y_n$. By taking a subsequence if necessary, we may assume $y_n \to y_* \notin \cM$.
	Using the above discussion, we have
	$$ l_E(h) \le \inf_{y_* \in \gamma, [\gamma]=h}l_E(\gamma) \le \inf_n l_E(h) + d_n = l_E(h).$$
	Taking limits, we conclude that there exists
	a rectifiable curve $\gamma_*$ containing $y_*$ with $l_E(\gamma_*) = l_E(h)$, hence $\gamma_*$
	is a shortest curve. But $y_* \notin \M$, leading to a contradiction.
	\end{proof}

	Proposition~\ref{non-zero-width} clearly applies to the slow system as there are either one or two shortest
	geodesics.

	\subsection{The case $E=0$}
	\label{sec:critical}

	We now extend the earlier discussions to the case $E=0$. While the functions $a_E^\pm$
	is not defined at $E=0$, the functions $d_E^\pm$ is well defined at $E=0$.
Recall $h, \bar h$ form a basis in $H_1(T^2,\Z)$ and the dual of $\bar h$ is perpendicular to $h$
and denoted by $h^\perp$.

\begin{prop}\label{prop-E0}
The properties of the channel and the Aubry sets depends on the type of homology $h$.
\begin{enumerate}
	\item Assume $h$ is simple and critical.
	\begin{enumerate}
		\item $d_0^+(h) - d_0^-(h) >0$.
		\item $l_0(h)h^* + [d_0^-(h), d_0^+(h)]\, h^\perp \subset \LF_\beta(0)$.
		\item For $c\in l_0(h)\,h^* + [d_0^-(h), d_0^+(h)]\, h^\perp$, we have $\gamma_h^0 \subset \cA_{H^s}(c)$;

		For $c\in l_0(h)h^* + (d_0^-(h), d_0^+(h))\, h^\perp$, we have $\gamma_h^0 = \cA_{H^s}(c)$.
	\end{enumerate}
	\item Assume $h$ is simple and non-critical.
	\begin{enumerate}
		\item $d_0^+(h) - d_0^-(h) >0$.
		\item $l_0(h)h^* + [d_0^-(h), d_0^+(h)]\,h^\perp \subset \LF_\beta(0)$.
		\item For $c\in l_0(h)h^* + [d_0^-(h), d_0^+(h)]\,h^\perp$, we have $\gamma_h^0 \cup \{0\} \subset
		\cA_{H^s}(c)$;

		For $c\in l_0(h)h^* + (d_0^-(h), d_0^+(h))h^\perp$, we have $\gamma_h^0 \cup \{0\}= \cA_{H^s}(c)$.
		\item The functions $d^\pm_E(h)$ is right-continuous at $E=0$.
	\end{enumerate}
	\item Assume $h$ is non-simple and $h=n_1h_1 + n_2h_2$, with $h_1,h_2$ simple.
	\begin{enumerate}
		\item $d_0^+(h) = d_0^-(h)$. Moreover, let $c^*(h) = l_E(h_1)h^*_1 + l_E(h_2)h_2^*$,
		where $(h_1^*, h_2^*)$ is the dual basis to $(h_1, h_2)$, then
		$$ c^* = l_E(h)h^* + d_0^\pm(h) h^\perp, $$
		where $h^\perp$ is a unit vector perpendicular to $h$.
		\item $\gamma_{h_1}^0 \cup \gamma_{h_2}^0 = \cA_{H^s}(c^*)$.
		\item $d_0^+(h_1)- d_0^-(h_1)>0$ with
		$$ l_E(h_1)h_1^* + d_0^+(h_1)h_2^* = c^*. $$
	\end{enumerate}
\end{enumerate}
\end{prop}

Before proving Proposition~\ref{prop-E0}, we first explain how the proof of
Proposition~\ref{non-zero-width} can be adapted to work even for $E=0$.

\begin{lem}\label{crit-width}
Assume that there is a unique  $g_0-$shortest geodesic in the homology $h$.  Then
$$ d_0^+(h) - d_0^-(h) >0 .$$
\end{lem}
\begin{proof}

We will try to adapt the proof of Proposition~\ref{non-zero-width}.
Let $\gamma_0$, $\gamma_1$ and $\gamma_2$ be shortest geodesics in homologies $h$,
$nh + \bar{h}$ and $nh- \bar{h}$, respectively. We choose an arbitrary
parametrization for $\gamma_i$ on $[0,T]$. Note that the parametrization is only continuous in general.

The proof of Proposition~\ref{non-zero-width} relies only on the property that lifted
shortest geodesics intersects at most once. For $E=0$, we will rely on a weaker property.

Let $\widetilde{\gamma}_i$ be the lifts to the universal cover $\R^2$. The degenerate point $\{0\}$
lifts to the integer lattice $\Z^2$. Since $g_0$ is a Riemannian metric away from the integers,
using the shortening argument, we have: if $\gamma_i$ intersect $\gamma_j$ at more than one point,
then either the intersections occur only at integer points, or the two curve coincide on a segment
with integer end points.

Let $a_0\in \gamma_0 \cap \gamma_1$ and let $\widetilde{\gamma}_0$ and $\gamma_1$ be lifts
with $\widetilde{\gamma}_0(0)=\widetilde{\gamma}_1(0)=a_0$. If $a\notin \Z^2$, then it is
the only intersection between the two curves. If $a_0\in \Z^2$, we define $a_0'$ to be
the largest intersection between $\widetilde{\gamma}_0|[0,T)$ and $\widetilde{\gamma}_1$ according
to the order on $\widetilde{\gamma}_0$. $a_0'$ is necessarily an integer point, and since
$a_0'\in \widetilde{\gamma}_0$, there exists $n_1 < n$ such that $a_0' - a_0 =n_0 h$.
Moreover, using the fact that $\widetilde{\gamma}_0$ is minimizing, we have
$$ l_0(\widetilde\gamma_0|[a_0,a_0']) = l_0(\widetilde{\gamma}_1|[a_0,a_0']). $$
We now apply a similar argument to $\widetilde{\gamma}_0 + \bar{h}$ and $\widetilde{\gamma}_1$.
Let  $a_1 = \widetilde\gamma_0(T) + \bar{h} = \widetilde{\gamma}_1(T)$ and let $a_1'$ be the smallest
intersection between $\widetilde{\gamma}_0|(0,T]$ and $\widetilde{\gamma}_1$. Then there exists
$n_1 \in \mathbb N$, $n_0 + n_1 <n$, such that $a_1 - a_1' = n_1 h$. Moreover,
$$ l_0((\widetilde\gamma_0+\bar{h})|[a_1',a_1]) = l_0(\widetilde{\gamma}_1|[a_1',a_1]). $$
Let $\widetilde{\eta}_1 = \widetilde{\gamma}_1|[a_0', a_1']$ and $\eta_1$ be its projection.
We have $[\eta]=(n-n_0-n_1)h + \bar{h}=:m_1 h + \bar{h}$, and
$$
l_0(\eta_1) - m_1 l_0(h) = l_0(\gamma_1) - nl_0(h) .
$$
The curve $\eta_1$ has the property that it intersects $\gamma_0$ only once.
Apply the same argument to $\gamma_2$, we obtain a curve $\eta_2$ with
$[\eta_2] = m_2 h - \bar{h}$, and
$$ l_0(\eta_2) - m_2 l_0(h) = l_0(\gamma_2) - nl_0(h) .$$

To proceed as in the proof of Proposition~\ref{non-zero-width}, we show that
if $\widetilde{\eta}_1$ and $\widetilde{\eta}_2$ are lifts of $\eta_1$ and $\eta_2$ with
the property that
$$ \widetilde{\eta}_1(0), \widetilde{\eta}_2(T) \in \{\widetilde{\gamma}_0(t)\}, \quad
\widetilde{\eta}_1(T), \widetilde{\eta}_2(0) \in \{\widetilde{\gamma}_0(t) + \bar{h}\}, $$
then $\widetilde{\eta}_1$ and $\widetilde{\eta}_2$ intersects only once. Indeed,
there are no integer points between $\widetilde{\gamma}_0$ and $\widetilde{\gamma}_0 + \bar{h}$.

We have
$$ l_0(\eta_1) - m_1 l_0(h) + l_0(\eta_2) - m_2 l_0(h) = d_n,$$
where $d_n$ is as defined in  Proposition~\ref{non-zero-width}. Assume $\inf d_n =0$,
proceed as in the proof of  Proposition~\ref{non-zero-width}, we obtain curves $[\gamma_k]=h$,
positive distance away from $\gamma_0$, such that
$$
l_0(h) \le l_0(\gamma_k) \le l_0(h) + d_n.
$$
This leads to a contradiction.
\end{proof}

\begin{proof}[Proof of Proposition~\ref{prop-E0}]

\emph{Case 1, $h$ is simple and critical.}

(a) This follows from Lemma~\ref{crit-width}.

(b) We note that Lemma~\ref{act-lowerbnd} depends only on positive homogeinity and
sub-additivity of $l_E(h)$, and hence applies even when $E=0$. We obtain for $c\in l_0(h)h^* + [d_0^-(h),
d_0^+(h)]\bar{h}^* $
$$ l_0(h') - \langle c, h' \rangle\ge 0, \forall h' \in H_1(\T^2, \Z^2). $$
Since $l_E(h)$ is strictly increasing, we obtain $l_E(h') - \langle c, h' \rangle> 0$
for $E>0$. By Lemma~\ref{orthogonal}, there are no $c-$minimal measures with energy $E>0$.
As a consequence, $\alpha(c)=0$. Since $\{0\}$ is a $c-$minimal measure with
rotation number $0$, we conclude $l_0(h)h^* + [d_0^-(h), d_0^+(h)]\bar{h}^* \subset \LF_\beta(0)$.

(c) Since we proved $\alpha(c)=0$, the first conclusion follows from Lemma~\ref{orthogonal}.
For the second conclusion, we verify that the proof of Proposition~\ref{char-aubry} for non-bifurcation val
applies to this case.

(d) The set function $[d_E^-(h), d_E^+(h)]$ is upper semi-continuous at $E=0$ from the right, by definition
We will show that it is continuous. Assume by contradiction that
$$ [\liminf_{E\to 0+}d_E^-(h), \limsup_{E\to 0+}d_E^+(h)] \subsetneq [d_0^-(h), d_0^+(h)]. $$
Then there exists $c\in l_0(h)h^* + (d_0^-(h), d_0^+(h))\bar{h}^*$ and
$$c(E) \notin l_E(h)h^* + [d_E^-(h), d_E^+(h)]h^\perp$$ such that $c(E)\to c$.
By part (c), the Aubry set $\cA_{H^s}(c)$ supports a unique minimal measure.
By Proposition~\ref{cont-barrier}, the Aubry set is upper semi-continuous in $c$.
Hence any limit point of $\cA_{H^s}(c(E))$ as $E\to 0$ is in $\cA_{H^s}(c)$.
This implies that $\widetilde{\cA}_{H^s}(c(E))$ approaches $\gamma_h^E$ as $E\to 0$.
Since $\gamma_h^E$ is the unique closed geodesic in a neighbourhood of itself, we conclude that $\widetilde{\cA}_{H^s}(c(E)) = \gamma_h^E$
for sufficiently small $E$. But this contradicts with $c(E) \notin l_E(h)h^* + [d_E^-(h), d_E^+(h)]h^\perp$

\emph{Case 2, $h$ is simple and non-critical.}

(a) This follows from Lemma~\ref{crit-width}.

(b) The proof is identical to case 1.

(c) For the first conclusion, we can directly verify that
$\gamma_h^0 \subset \cA_{H^s}(c)$ and $\{0\} \subset \cA_{H^s}(c)$.
For the second conclusion, we note that  proof of
Proposition~\ref{char-aubry} for bifurcation values applies to this case.

\emph{Case 3, $h$ is non-simple with $h=n_1h_1 + n_2h_2$.}

(a) Assume that $\bar{h} = m_1h_1 + m_2h_2$ for some $m_1, m_2\in \Z$.
For sufficiently large $n\in \NN$, we have $nh \pm \bar{h} \in \NN h_1 + \NN h_2$.
As a consequence,
\begin{multline*} l_0(nh \pm \bar{h}) - l_0(nh)\\
= (nn_1 \pm m_1)l_0(h_1) + (nn_2 \pm m_2)l_0(h_2) - (nn_1 l_0(h_1) + n n_2 l_0(h_2))\\
=\pm m_1 l_0(h_1) \pm m_2 l_0(h_2).
\end{multline*}
We obtain $d_0^+(h) - d_0^-(h) =0$ by definition.

We check directly that
$$ l_0(h) - \langle c^*, h\rangle =0.$$
Since $l_0(h)h^* + d_0^-(h)h^\perp = l_0(h)h^* + d_0^+(h)h^\perp $ is the unique $c$ with this property. The second claim follows.

(b) We note that any connected component of the complement to $\gamma_{h_1}^0\cup\gamma_{h_2}^0$ is
contractible. If $\cA_{H^s}(c)$ has other components, the only possibility is a contractible orbit bi-asymp
to $\{0\}$.
However, such an orbit can never be minimal, as the fixed point $\{0\}$ has smaller action.

(c) The statement $d_0^+(h_1) - d_0^-(h_1)>0$ follows from part 1(a). for the second claim, we compute
$$ d_0^+(h_1) = \inf_n l_0(nh_1 + h_2) - l_0(nh_2) = l_0(h_2)$$
and the claim follows.
\end{proof}

%
%
%

\subsection*{Acknowledgment}
V.K. has been partially supported of the NSF grant DMS-1402164 and enjoyed the hospitality of 
the ETH Institute for Theoretical Studies. He has been partially supported by the Dr. Max Rssler, 
the Walter Haefner Foundation and the ETH Zurich Foundation. K.Z. is supported by the NSERC Discovery grant, reference number 436169-2013, and thank  the hospitality of the  Institute for Theoretical Studies 
for visiting in 2017. During the spring semester of 2012 the authors visited the Institute
for Advanced Study. They acknowledge hospitality and a highly stimulating research atmosphere. 
The authors are grateful to John Mather for numerous inspiring discussions. The course of lectures Mather
gave in the spring of 2009 was very helpful for the authors.
The authors thank Abed Bounemoura and Chong-Qing Cheng for helpful
comments on preliminary versions of this paper. Bounemoura  pointed out
some errors related to normal forms and Cheng on error related to total disconnectedness of Ma\~ne sets.
The authors are thankful to Giovanni Forni for helpful remarks on exposition.

\section{Notations}
We provide a list of notations for the reader's convenience. 

\subsection{Formulation of the main result}

\begin{longtable}{l p{0.7\textwidth}}

$(\theta, p, t)$ & A point in the phase space $\T^2 \times \R^2 \times \T$.  \\

$\phi_{H}^{s, t}$ & The Hamiltonian flow for an non-autonomous Hamiltonian.  \\

$\Phi_H^t$ & The Hamiltonian flow for an autonomous Hamiltonian. \\ 

$H_\epsilon = H_0 + \epsilon H_1$ & Nearly integrable system.  \\

$D$ & A fixed constant controlling the norm and convexity of $H_0$. \\

$\Z^3_* = \Z^3 \setminus \{(0, 0, 1)\}\Z$ & Set of integer vectors defining a resonance relation.  \\

$k = (k^1, k^2, k^0)$ & Resonance vectors, contained in $\Z^3_*$. \\

$S_k$ & Single resonance surface in the action space given by $k \in \Z^3_*$. \\

$\Gamma_k$ & Singe resonance segment, a closed segment contained in $S_k$.  \\

$S_{k_1, k_2}$ & Double resonance point in the action space given by $k_1, k_2 \in \Z^3_*$. \\

$\cS^r$ & The unit sphere of the $C^r$ functions. \\

$\cP$ & Diffusion path consisting of segments of single resonance segments.  \\

$\cK = \{(k, \Gamma_k)\}$ & Collection of resonances and resonance segments making up a diffusion path. \\

$\cU = \cU(\cP)$ & An open and dense set in $\cS^r$ defining the ``non-degenerate perturbations'' relative to a diffusion path. \\

$\cU_{SR}^\lambda(k_1, \Gamma_{k_1})$ & Set of $H_1 \in \cS^r$ satisfying the quantitative non-resonance conditions $[SR1_\lambda]-[SR3_\lambda]$ relative to the resonant segment $(k_1, \Gamma_{k_1})$. \\

$\cU_{DR}(k_1, k_2)$ & Set of $H_1 \in \cS^r$ satisfying the non-degeneracy conditions $[DR1^h]-[DR3^h]$ and $[DR1^c]-[DR4^c]$. \\

$\cV = \cV(\cU, \epsilon_0)$ & A ``cusp'' set of perturbations, equal to $\{\epsilon H_1:\, H_1 \in \cU, \, 0 < \epsilon < \epsilon_0(H_1)\}$. \\

$\cK^\st(k_1, \Gamma_{k_1}, \lambda)$ & Set of strong additional resonances relative to the resonance segment $\Gamma_{k_1}$, for an perturbation $H_1 \in \cU_{SR}^\lambda(k_1, \Gamma_{k_1})$. \\

$B_\sigma$ & Ball in Euclidean space $\R^n$. \\

$\cV^r_\sigma, \cV_\sigma$ & Ball in the functional space $C^r$. When supscript is not indicated, then stands for $C^r$ where $r$ is from the main theorem. 
\end{longtable}

\subsection{Weak KAM and Mather theory}

\begin{longtable}{l p{0.7\textwidth}}
$\varpi = \varpi(H)$ & The period of a time periodic Hamiltonian, i.e. $H(\theta, p, t + \varpi) = H(\theta, p, t)$. \\
$\T_\varpi = \R /(\varpi \Z)$ & Torus with period $\varpi$. \\
$\bH = \bH(D)$ & A family of Hamiltonians satisfying uniform conditions depending on the parameter $D>1$. \\
$L = L_H$ & The Lagrangian of $H$. \\
$L_{H, c}$ & The ``penalized'' Lagrangian $L_H(\theta, v, t) - c \cdot v$. \\
$A_{H,c}(x, s, y, t)$ & The minimal action for the Lagrangian $L_{H, c}$. \\
$\alpha_H(c), \alpha_L(c)$ & Mather's alpha function. \\
$\beta_H(\rho), \beta_L(\rho)$ & Mather's beta function. \\
$h_{H, c}(x, s, y, t)$ & The time dependent Peierl's barrier function. \\
$h_{H, c}(x, y)$ & The discrete Peierl's barrier function, equal to $h_{H, c}(x, 0, y, 0)$. \\
$d_{H, c}(x,s, y, t)$ & Mather's semi-distance \\
$T^{s, t}_c u(x)$ & Lax-Oleinik semi-group defined on $\T^n$. \\
$\partial^+ u(x)$ & The supergradient of a semi-concave function at $x$. \\
$\cG_{c, w}$, $w = w(\theta, t)$ & The time-dependent pseudograph as a subset of $\T^n \times \R^n \times \R$, or $\T^n \times \R^n \times \T_{\varpi}$. \\
$\tilde{I}(c, w)$ & The maximum invariant set contained in the psudograph $\cG_{c, w}$. \\
$\tcM_H(c), \tcA_H(c), \tcN_H(c)$ & The continuous (Hamiltonian) Mather, Aubry and Ma\~ne set, defined on $\T^n \times \R^n \times \T_{\varpi(H)}$. \\
$\tcM_H^0(c), \tcA_H^0(c), \tcN_H^0(c)$ & The discrete (Hamiltonian) Mather, Aubry and Ma\~ne set, defined on $\T^n \times \R^n$, invariant under the map $\phi_H = \phi_H^{\varpi(H)}$.  \\
$\tcS(H, c)$ & A static class of the Aubry set $\tcA_H(c)$. \\
$c \vdash c'$ & The forcing relation. \\
$c \dashv\vdash  c'$ & The forcing equivalence relation. 
\end{longtable}

\subsection{Single resonance}

\begin{longtable}{l p{0.7\textwidth}}
$[H_1]_{k_1}$ & The average of $H_1$ relative to the resonance $k_1$. \\

$[H_1]_{k_1, k_2}$ &  The average of $H_1$ relative to the double resoannce $k_1, k_2$. \\

$\Phi_\epsilon$ & The averaging coordinate change at single resonance.  \\

$N_\epsilon$ & Normal form under the coordinate change. Same notation is used at double resonance. \\

$Z_{k_1}(\theta^s, p, t)$ & The resonant component of $H_1$ relative to the resonance $k_1$. \\

$R(\theta, p, t)$ & The remainder in the single resonance normal form. \\

$\cK^\st(k_1, \Gamma_{k_1}, K)$ & Set of strong additional resonance $k_2$ intersecting $\Gamma_{k_1}$ with norm at most $K$, plus any resonances in the diffusion path that intersect $\Gamma_{k_1}$. \\

$\Gamma_{k_1}^{SR}$ & The punctured resonance segment after removing $O(\sqrt{\epsilon})$ neighborhoods of strong double resonances.  \\

$\|\cdot\|_{C_I^r}$ & The rescaled $C^r$ norm where the derivatives in the action variable are rescaled by $\sqrt{\epsilon}$. \\

$T_\omega$ & The period of the rational vector $(\omega, 1) \in \R^3$.  \\

$B \in SL(3, \Z)$ & An integer matrix defining the linear coordinate change relative to a double resonance $k_1, k_2$. \\

$(\theta^s, \theta^f, p^s, p^f, t)$ & The coordinate at single resonance after taking a linear coordinate change. The resonance becomes $(1, 0, 0) \cdot (\omega, 1)  =0$ in this coordinate.  \\

$(\Theta^s, P^s)(\theta^f, p^f, t)$ & Normally hyperbolic invariant cylinder for the single resonant normal form, parametrized using the $(\theta^f, p^f, t)$  variables. 

\end{longtable}

\subsection{Double resonance}

\begin{longtable}{l p{0.7\textwidth}}
$\Phi_\epsilon$ & The averaging coordinate change at double resonance. \\
$N_\epsilon$ & The normal form at double resonance. \\
$\Phi_L$ & The linear coordinate change corresponding to the double resonance $k_1, k_2$. \\
$K(I)$ & The kinetic energy of the slow mechanical system.  \\
$U(\varphi)$ & The potential function of the slow mechanical system. \\
$g_E$ & The Jacobi metric at energy $E$ of the slow mechanical system. \\
$h$ & A homology class in $H^1(\T^2, \Z)$. \\
$\gamma_h^E$ & Shortest curves of the Jacobi metric in homology class $h$. \\
$\eta_h^E$ & The Hamiltonian periodic orbit corresponding to the geodesic $\gamma_h^E$. \\
$\Phl^{ij}$, $i,j \in \{+, -\}$ & Local maps near the saddle fixed point of $H^s$. \\
$\Phg$ & Global map along a homoclinic orbit $\eta$ of $H^s$. \\
$\bar{c}_h(E)$ & Curve of cohomologies chosen in the channel of $h$. \\
$\bar{\Gamma}_h$ & Choice of cohomologies along the homology $h$. \\
$\bar{\Gamma}_h^e$ & In the non-simple case, choice of cohomology curve above energy $e$. \\
$\bar{\Gamma}_{h_1}^{e, \mu}$ & In the non-simple case, choice of cohomology along the adjacent simple homology $h_1$. \\
$\Phi^*_L$ & The relation between cohomology class and alpha function after the coordinate change $\Phi_L$. \\
$\Phi^*_{L, H_\epsilon^s}$ & The relation between the cohomology classes for the coordinate change $\Phi_L$. Depends on the alpha function of the system $H_\epsilon^s$. \\
$\Gamma^{DR}_{k_1, k_2}$ & The choice of cohomology classes at a double resonance $k_1, k_2$. \\

\end{longtable}

\bibliographystyle{plain}
\bibliography{strong-diffusion}

\end{document}